\documentclass[11pt,twoside,reqno]{amsart}
\usepackage[margin=1.1in]{geometry}
\usepackage{graphicx} 
\usepackage{amssymb}
\usepackage{mathrsfs}
\usepackage{amsaddr}
\usepackage[hidelinks]{hyperref}
\usepackage{bm}

\newcommand{\ri}{\mathrm{i}}
\newcommand{\sym}{\mathrm{Sym}}

\newtheorem{theorem}{Theorem}[section]
\newtheorem{proposition}[theorem]{Proposition}
\newtheorem{corollary}[theorem]{Corollary}
\newtheorem{lemma}[theorem]{Lemma}
\theoremstyle{remark}
\newtheorem{remark}[theorem]{Remark}
\newtheorem{definition}[theorem]{Definition}

\title[Nonlinear Hydrodynamic Manifolds]{Nonlinear Hydrodynamic Manifolds for a BGK Equation with Density-Dependent Collision Rate}

\author{Florian Kogelbauer}
\email{florian.kogelbauer@kaust.edu.sa}
\address{King Abdullah University of Science and Technology (KAUST),
CEMSE Division\\
Thuwal, 23955-6900, Saudi Arabia}

\date{\today}

\begin{document}

\begin{abstract}
We give a constructive proof of the existence of nonlinear hydrodynamic
manifolds for a Bhatnagar--Gross--Krook (BGK)-type kinetic equation with
a density-dependent collision rate on $\mathbb R^d$, $d\geq3$, near a
positive homogeneous Maxwellian. Our construction starts from the slow
spectral modes of the linearized dynamics and yields a convergent generalized functional expansion whose terms can be computed recursively. Sums of slow eigenvalues can meet the essential spectrum at the combined frequency, creating external resonances that obstruct stationary inversion of the homological equations. To overcome this obstruction, we introduce the new technique of parametrization of time, an intrinsic clock determined by the slow spectral amplitudes that allows us to solve the hierarchy of homological equations by evolution along the reduced trajectories.
The construction applies directly to the full kinetic equation, without finite-dimensional truncation or a global spectral gap separating the slow part of the spectrum from the origin or the remaining spectrum. The parametrization is differentiable at equilibrium and is a local $C^1$ embedding on the punctured open ball around each Maxwellian, while, on a smaller forward-invariant domain, the mass density provides local
coordinates on the hydrodynamic manifold.
These coordinates yield an exact nonlinear hydrodynamic closure in the form of a nonlocal equation for the macroscopic mass density.
For sufficiently small kinetic states, the reconstructed distributions remain positive, satisfy uniform two-sided Maxwellian bounds, and have finite, non-increasing relative entropy with respect to the
background Maxwellian.
\end{abstract}

\maketitle

\tableofcontents

\newpage

\section{Introduction}

\begin{figure}
    \centering
    \includegraphics[width=0.9\linewidth]{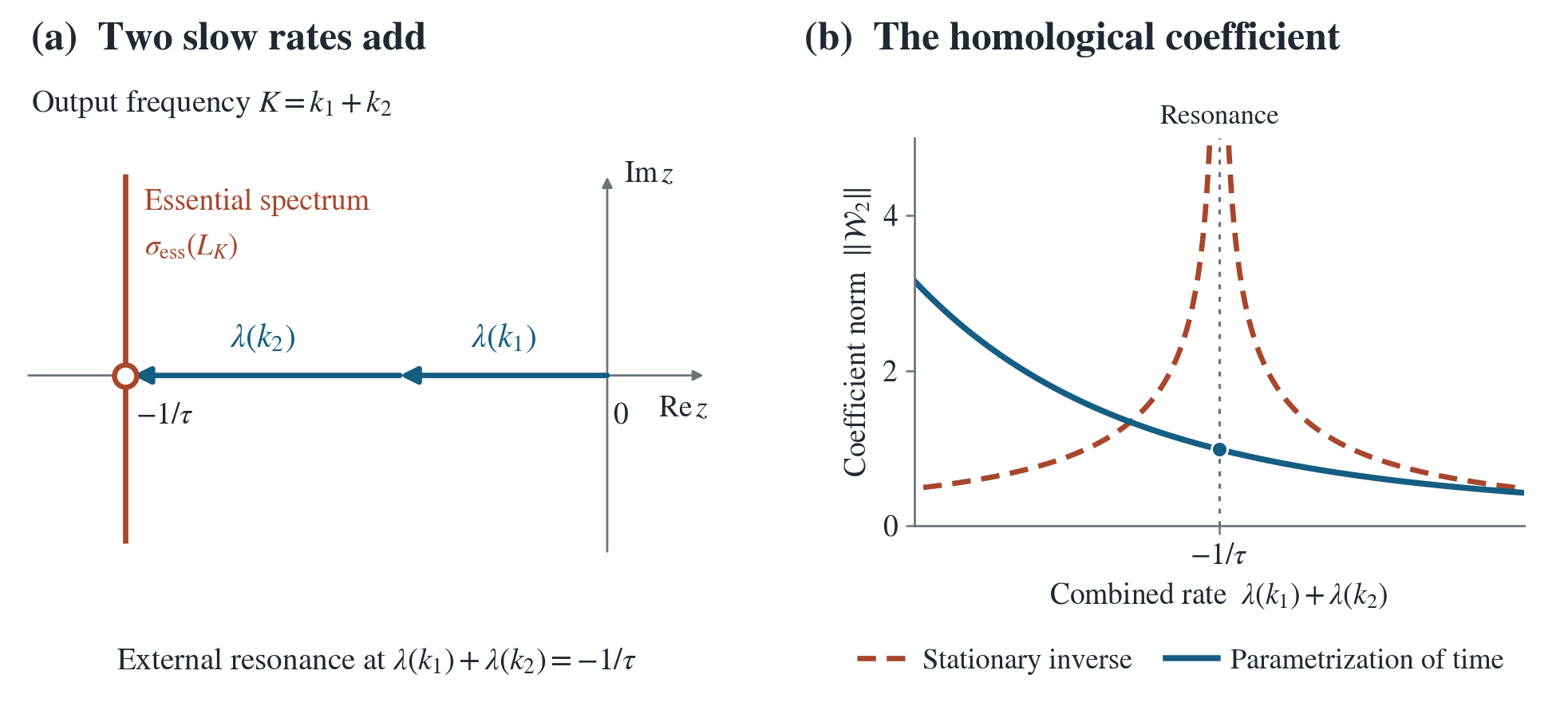}
    \caption{Illustration of external resonance in a quadratic mode interaction.
(a) Two slow modes interact at frequency $K=k_1+k_2$. When their rates
add to a point $\lambda(k_1)+\lambda(k_2)=-1/\tau$ of the essential
spectrum of $L_K$, the stationary homological operator $(-1/\tau)I-L_K$ has
no bounded inverse, see Section \ref{subsection_homological} for technical details.
(b) In the example shown, the norm of the stationary coefficient
(dashed orange) diverges as we approach the external resonance, while the parametrization of time gives a coefficient whose norm remains finite through resonance
(solid blue).}
    \label{fig_ext_res}
\end{figure}

What is the relation between the statistical description of gases and the field equations of continuum mechanics? This question lies at the heart of Hilbert's Sixth Problem and, in its kinetic formulation, asks how the Boltzmann equation and related kinetic models give rise to macroscopic evolution equations \cite{Hilbert1901}. The classical approach to this problem uses moment equations and asymptotic closures, most prominently the Chapman--Enskog expansion, which formally yields the Euler, Navier--Stokes and Fourier constitutive laws. Rigorous hydrodynamic limits instead compare suitably scaled kinetic solutions with fluid solutions as the Knudsen number tends to zero, see \cite{caflisch1979fluid,BardosGolseLevermore1993,GolseSaintRaymond2009,LevermoreMasmoudi2010,SaintRaymond2009HydrodynamicLimits,gallagher2020convergence} and the recent breakthroughs on the passage from many-particle or wave systems to kinetic models \cite{DengHani2023WaveKinetic,DengHaniMaHardSpheres,deng2025hilbert}.\\
In this work, we follow a different paradigm for the passage from kinetic equations to fluid dynamics: we construct a family of special forward-invariant solutions in kinetic space around each global Maxwellian, called the \textit{hydrodynamic manifold}. This manifold is a local nonlinear deformation of the slow spectral subspace of the linearized kinetic operator, and defines a local $C^1$ embedding away from the manifold of equilibria. On a smaller forward-invariant domain, the macroscopic mass density gives local coordinates on the hydrodynamic manifold, which leads to an exact closure relation and hydrodynamic reduction. This viewpoint is related to Bogolyubov's separation of kinetic and hydrodynamic time scales and to geometric interpretations of the Chapman--Enskog procedure as discussed in \cite{mckean1969simple,kuramoto1989reduction,GorbanKarlin2014Hilbert}, see also the exact linear hydrodynamic manifolds constructed in \cite{kogelbauer2020slow,KarlinChikatamarlaKooshkbaghi2014,kogelbauer2021non,KogelbauerKarlinClosure2025}. On such a manifold, the kinetic equation is satisfied exactly, while the macroscopic variables provide coordinates for a closed reduced evolution. The central question following this paradigm is therefore not only how to approximate a kinetic equation by a fluid equation, but whether an invariant family of kinetic states carrying exact fluid dynamics can be constructed. Related constructions include invariant manifolds for steady Boltzmann flows \cite{liu2013invariant}, invariant-manifold theory for PDEs \cite{bates1989invariant,chen1997invariant,wayne1997invariant}, and spectral submanifolds for systems with discrete spectral structure \cite{kogelbauer2018rigorous,buza2024spectral,buza2026spectral}.\\
On the technical side, our construction is motivated by the parametrization method
\cite{cabre2003parameterization}, which aims to continue a selected linear spectral subspace at equilibrium to a nonlinear
invariant manifold. In the classical parametrization method, one determines
the coefficients of the parametrization successively by solving stationary
linear equations, called homological equations. For the whole-space
kinetic equation considered here, no spectral gap separates the selected
hydrodynamic modes from the remaining spectrum of the full operator.
The spectral separation used in classical inertial-manifold constructions
and in spectral-submanifold constructions for dissipative PDEs
\cite{Zelik2014,buza2024spectral} is therefore unavailable.\\
The obstruction to stationary inversion arises when
hydrodynamic modes interact. At an \emph{internal resonance}, the sum of the
interacting eigenvalues equals the hydrodynamic eigenvalue at the combined
frequency. If their spectral values add to a point of
the essential spectrum at the combined frequency, the corresponding
homological operator has no bounded inverse. Such an \emph{external
resonance} can cause the stationary construction with prescribed linear
reduced dynamics to fail already at quadratic order, even after integration
over the interacting frequencies
(Proposition~\ref{prop_stationary_quadratic_obstruction}). To overcome this obstruction we introduce a new technique, the \emph{parametrization of time}, in which we evolve the coefficients through the resonance and recover elapsed time from the slow amplitudes. Along a nontrivial reduced orbit we define a scalar functional $\theta[\xi]$ such that
\begin{equation}\label{intro_clock}
    \theta[\xi(t)]=\theta[\xi(0)]+t.
\end{equation}
Ideas related to \eqref{intro_clock} have been developed in \cite{WilsonMoehlis2016}, where the authors recover the elapsed time from a leading decaying amplitude of a PDE solution. Resonant spectral dynamics have been analyzed in \cite[Thm.~4]{FoiasSaut1987}, where it is shown that resonances can produce polynomial factors multiplying decaying exponentials for the Navier--Stokes equation, as well as in \cite[Sec.~4.2]{gallay2002invariant}, where corrections of polynomial-logarithmic type were analyzed on a vorticity manifold. We stress that our construction \eqref{intro_clock} uses the full frequency domain of slow modes of the kinetic spectrum and resolves resonances for the full nonlinear PDE without scaling arguments or self-similar coordinate changes. Figure \ref{fig_ext_res} illustrates the problem of resonances of the homological equation and the regularization by the parametrization of time schematically.\\
Allowing the coefficient kernels of a generalized functional expansion to depend on this intrinsic clock replaces the singular steady homological  equation by an unsteady equation. Its Duhamel representation passes through external resonances without division by the vanishing spectral denominator. Quantitative estimates for the resulting kernels, together with the compact frequency support of the slow parameter functions, allow us to sum the nonlinear expansion term by term. This is a construction in the full infinite-dimensional kinetic equation, without a spatial Fourier discretization, a global spectral gap, or a small-Knudsen-number limit. Although the spectral clock becomes singular at zero, the decay of the coefficients proves differentiability of the summed parametrization at the origin, while continuous differentiability is proved only away from zero. We stress that the time dependence in the coefficients is not an additional physical time variable since, through \eqref{intro_clock}, it defines an autonomous invariant parametrization.\\
We illustrate this mechanism for the nonlinear BGK-type equation based on the relaxation model of Bhatnagar, Gross and Krook~\cite{BhatnagarGrossKrook1954}, 
\begin{equation}\label{maineq_intro}
    \frac{\partial F}{\partial t}+v\cdot\nabla_xF
    =\frac{\rho[F]^\alpha}{\tau}\bigl(\rho[F]\phi-F\bigr),
\end{equation}
where $F:\mathbb R_x^d\times\mathbb R_v^d\times[0,\infty)\to [0,\infty)$ is a distribution function, $\rho[F]=\int_{\mathbb R^d}F\,dv$ is its mass density, $\phi$ is the normalized Maxwellian, $\tau>0$ is the relaxation parameter, and $\alpha$ describes the density dependence of the collision rate. Details, including connections with granular media and related density-dependent relaxation models, are given in Section~\ref{sec_main_equation}, while we refer to \cite{evans2026density} for a study of \eqref{maineq_intro} on the torus $\mathbb T^d$, including exponential relaxation to equilibrium. Earlier analyses establish nonlinear diffusion limits for related relaxation
operators in finite-velocity models \cite{LionsToscani1997,GolseSalvarani2007}
and in radiative transfer with opacity depending on the angular
average \cite{BardosGolsePerthameSentis1988}.\\
We stress that our construction of the spectral clock \eqref{intro_clock} depends only on the exponential evolution of the slow amplitudes in \eqref{reduced_dyn_linear}, so its construction naturally extends to other collision operators with such reduced dynamics. We use the BGK model \eqref{maineq_intro} to illustrate the method because its rank-one linearized gain yields a scalar dispersion relation and explicit spectral projections, while the density-dependent rate gives an explicit Bell polynomial formula for the nonlinear forcing in \eqref{hat_N_n_explicit}.\\

Our principal existence result can be summarized informally as follows.

\begin{theorem}[Main theorem, informal version]
Assume $d\geq 3$ and $\alpha\geq 0$. Near each positive homogeneous
Maxwellian $\rho_0\phi$, equation~\eqref{maineq_intro} admits a
forward-invariant family of real kinetic states parametrized by the slow spectral modes of the
linearized kinetic operator. We call this family the
\emph{nonlinear hydrodynamic manifold}. Its parametrization is
differentiable at any global Maxwellian,
and is a local $C^1$ embedding into kinetic phase space at every nonzero amplitude in the open
parameter ball around the global Maxwellian.\\
On a smaller forward-invariant domain of parameter functions, defined in \eqref{invariant_density_domain_s9}, the mass density
$\rho$ locally determines the full kinetic state and the resulting closure map
$\rho\mapsto F_{\rm spec}[\rho]$ is $C^1$. The mass flux is therefore
determined by $\rho$ alone, and the density satisfies the exact, generally
nonlocal closed equation
\begin{equation}\label{dynamcis_rho}
    \partial_t\rho
    = -\nabla_x\cdot\int_{\mathbb R^d} v\,F_{\rm spec}[\rho]\,dv .
\end{equation}
\end{theorem}

To formulate the main theorem in technical terms, we introduce three
function space settings. For convergence of the kinetic expansion on a small closed ball of slow
amplitudes, we use the weighted Fourier space $\mathcal Z^{(2)}$ defined
in \eqref{def_Xv2_aux}. In this norm the parametrization is differentiable
at equilibrium, with derivative given by the linear lifting of the slow
modes. To prove local $C^1$ invertibility on the image of the punctured
open ball, we work under the stronger clock assumption below and use the
kinetic space $\mathcal X_{k,v}$ defined in
\eqref{closure_continuous_spaces_s9}, which also requires continuity in
frequency. Finally, to recover the kinetic state locally from its mass
density on the smaller forward-invariant domain
\eqref{invariant_density_domain_s9}, we use the space of Fourier
densities $\mathcal Y_\rho$ defined in \eqref{closure_continuous_spaces_s9}.
This space includes frequencies beyond the cutoff on the slow amplitudes
because nonlinear interactions generate higher harmonics. The following theorem gives a precise formulation of the existence and regularity claims for the hydrodynamic manifold.

\begin{theorem}[Main theorem]\label{mainthm_formal}
Consider the normalized perturbation equation \eqref{maineqf} with
$d\geq 3$, $\alpha\geq 0$ and $\tau>0$. Fix
$0<\varepsilon_B<k_c=\tau^{-1}\sqrt{\pi/2}$  as in \eqref{def_k_crit}, and let $\lambda(k)$
and $\hat e(k,v)$ be the hydrodynamic eigenvalue and eigenvector of
\eqref{deflambda} and \eqref{eigenvector}. Let $1\leq p<\infty$ and
let $w:\mathbb R^d\to(0,\infty)$ be a measurable weight satisfying the conditions of \eqref{ass_phi_w},
$$
 \begin{gathered}
 (1+|v|)^n\phi\in L_v^p(\mathbb R^d,w),\qquad n\geq 0,\\
 w^{-1/(p-1)}\in L^1(\mathbb R^d)\quad(p>1),\qquad
 w^{-1}\in L^\infty(\mathbb R^d)\quad(p=1).
 \end{gathered}
$$
We use the velocity space $\mathcal X_v$ from \eqref{defX} and the
weighted Fourier space $\mathcal Z^{(2)}$ from \eqref{def_Xv2_aux}, and define the underlying parameter space of slow amplitudes in \eqref{defXi_eps_C0} as
$$
 \Xi(\varepsilon_B)=
 \left\{\xi\in C^0(\overline{B_{k_c}(0)};\mathbb C):
 \begin{gathered}
 \xi(-k)=\xi(k)^*,\\ \operatorname{supp}\xi\subseteq\overline{B_{k_c-\varepsilon_B}(0)}
 \end{gathered}\right\},
 \qquad \|\xi\|_{C^0}=\sup_{|k|\leq k_c}|\xi(k)|.
$$
We regard this as a real Banach space and extend its elements by zero
outside $B_{k_c}(0)$. Throughout, we use the spectral clock $\theta[\xi]$ of \eqref{deftheta} with
$\theta_0=0$ and $\nu\geq -2\lambda(k_c-\varepsilon_B)$, and fix an even clock weight
$\mu\in L^\infty(B_{k_c}(0))$ with $\mu\geq \mu_0>0$ almost everywhere
and $\int_{B_{k_c}(0)}\mu=1$ as in \eqref{clock_weight_assumption}.\\
There exists a radius $\mathfrak r_*>0$ such that, for every fixed
$0<\mathfrak r<\mathfrak r_*$, the following assertions hold.\footnote{
All intersections below carry the sum norm, and all derivatives are
real Fr\'echet derivatives.}
\begin{enumerate}
\item \emph{Existence, invariance and the derivative at equilibrium.}
The expansion \eqref{defhydro} converges absolutely and uniformly
in $\mathcal Z^{(2)}$ for $0<\|\xi\|_{C^0}\leq\mathfrak r$.
Setting $\hat{\mathcal W}_{\mathfrak r}[0]=0$ extends its sum continuously
to equilibrium and defines the parametrization
$$
\hat{\mathcal W}_{\mathfrak r}:
\{\xi\in\Xi(\varepsilon_B):\|\xi\|_{C^0}\leq\mathfrak r\}
\longrightarrow
\{\widehat f\in\mathcal Z^{(2)}:
    \widehat f(-k,\cdot)=\widehat f(k,\cdot)^*\},
$$
which is differentiable at equilibrium with derivative
$$
 D\hat{\mathcal W}_{\mathfrak r}[0]h(k,v)
 =h(k)\frac{\phi(v)}{1+\tau\lambda(k)+\ri\tau k\cdot v},
 \qquad |k|<k_c-\varepsilon_B,
$$
and zero outside this ball. Along
each solution $\xi(t,k)=e^{t\lambda(k)}\xi(k)$ of the reduced equation
\eqref{reduced_dyn_linear}, the reconstructed perturbation
$\hat{\mathcal W}_{\mathfrak r}[\xi(t)]$ is a mild kinetic solution
for every $t\geq 0$, as expressed by \eqref{prop82_mild_identity}, implying that the image of $\hat{\mathcal W}_{\mathfrak r}$, i.e., the nonlinear hydrodynamic
manifold, is forward invariant. After inverse spatial Fourier transformation, the series converges in
$L_x^\infty(\mathcal X_v)$ and, when $1\leq p\le2$, also in $L_x^2(\mathcal X_v)$.
In each of these norms the physical parametrization is Fr\'echet differentiable
at zero and satisfies the same Duhamel identity in physical space along
every reduced orbit.

\item \emph{Local embedding away from equilibrium.}
If $\nu\geq 2/\tau$, then, after decreasing $\mathfrak r_*$ if necessary,
$$
 \hat{\mathcal W}_{\mathfrak r}:
 \{\xi\in\Xi(\varepsilon_B):0<\|\xi\|_{C^0}<\mathfrak r\}\longrightarrow\mathcal X_{k,v}
$$
is a local real $C^1$ embedding into the kinetic space $\mathcal X_{k,v}$
defined in \eqref{closure_continuous_spaces_s9}. By
Theorem~\ref{prop_kinetic_embedding}, every nonzero amplitude in this
open ball has a neighborhood $V$ on which the parametrization
is a $C^1$ diffeomorphism onto the submanifold
$\hat{\mathcal W}_{\mathfrak r}(V)$ of $\mathcal X_{k,v}$.
On the same punctured ball, the parametrization is also $C^1$
in the physical-space norms of part~1.

\item \emph{Local density coordinates on a smaller invariant domain.}
Assume $\nu\geq 2/\tau$ and fix $0<\delta<1$ and $G>0$. The domain
$$
 \mathcal U_{\delta,G}^{(\mathfrak r)}=
 \left\{\xi\in\Xi(\varepsilon_B):
 \begin{gathered}
 0<\|\xi\|_{C^0}<\delta\mathfrak r,\\
 \displaystyle\int_{|k|<k_c-\varepsilon_B}
 \frac{\mu(k)}{-\lambda(k)}
 \left(\frac{e^{-\lambda(k)\theta[\xi]}|\xi(k)|}{\mathfrak r}\right)
 ^{-\nu/\lambda(k)-1}\,dk<G
 \end{gathered}\right\},
$$
is open in $\Xi(\varepsilon_B)$ and forward invariant under the
reduced dynamics \eqref{reduced_dyn_linear}
(Lemma~\ref{lem_invariant_clock_domain_s9}). For
$0<\mathfrak r<\mathfrak r_{\rm inv}(\delta,G)\leq\mathfrak r_*$, the full
Fourier density map of \eqref{density_chart_section9}
$$
 \mathcal H_\rho\colon
 \mathcal U_{\delta,G}^{(\mathfrak r)}\longrightarrow\mathcal Y_\rho,
 \qquad \mathcal H_\rho(\xi)=\rho[\hat{\mathcal W}_{\mathfrak r}[\xi]]
$$
is a local $C^1$ embedding
for the supremum norm on the density space $\mathcal Y_\rho$ of
\eqref{closure_continuous_spaces_s9}. Near each amplitude, the inverse
on the corresponding density image extends to a $C^1$ map
on an open neighborhood in $\mathcal Y_\rho$ and, consequently, the density locally
determines the full kinetic state with $C^1$ dependence.
The constitutive law and the exact autonomous density equation
\eqref{density_spectral_closure_section9}--\eqref{density_reduced_dynamics_section9}
therefore continue along every finite reduced orbit segment in
$\mathcal U_{\delta,G}^{(\mathfrak r)}$ through successive local charts. For $G>G_{\rm log}$, with $G_{\rm log}$ defined in
\eqref{explicit_clock_class_constant_s9}, the logarithmic profiles of
Remark~\ref{rem_logarithmic_class_s9} form an explicit forward-invariant
subset of $\mathcal U_{\delta,G}^{(\mathfrak r)}$.

\item \emph{Real states, positivity and entropy.}
The full distribution
$F=\rho_0(\phi+f)$, where $f$ is the inverse Fourier transform of
$\hat{\mathcal W}_{\mathfrak r}[\xi(t)]$ and $\rho_0>0$, satisfies
the uniform two-sided Maxwellian bounds of
Theorem~\ref{thm_selected_maxwellian_bounds}, given explicitly in
\eqref{selected_uniform_maxwellian_bounds}. Moreover,
$f\in C([0,\infty);L_x^2L_v^2(\phi^{-1}))$, and the relative entropy of $F$ defined in
\eqref{def_H_Phi_relative}, with $\Phi(y)=y\log y$, is finite and satisfies the integrated form of
\eqref{relative_convex_entropy_dissipation} in the normalization of
\eqref{maineqf} for all forward times. Corollary~\ref{cor_selected_decay} gives the algebraic decay rates for the reconstructed perturbation and its relative entropy.
\end{enumerate}
\end{theorem}

\begin{figure}
    \centering
    \includegraphics[trim={1.7cm 2cm 0.7cm 2cm},clip, width=0.66\linewidth]{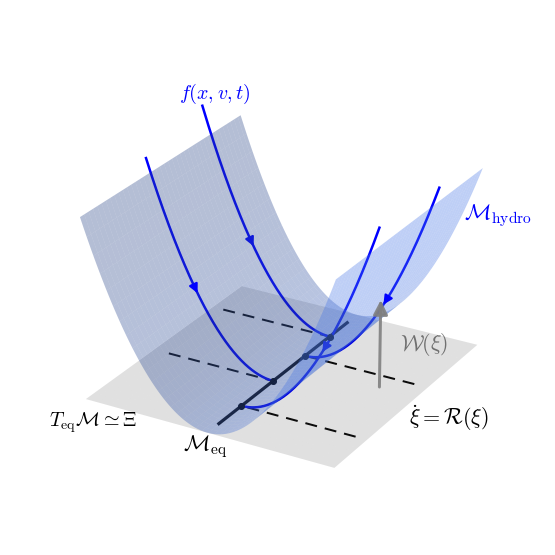}
    \caption{Schematic depiction of the selected hydrodynamic family $\mathcal{M}_{\rm hydro}$ constructed around the straight line of equilibria $\mathcal{M}_{\rm eq}$, with parametrization $\mathcal W$ and reduced dynamics $\dot\xi=\mathcal R(\xi)$. }
\end{figure}

\subsection{Structure of the Paper}
Section~\ref{Notation} fixes notation and introduces the plasma dispersion
function used in the spectral analysis. Section~\ref{sec_main_equation}
presents the density-dependent BGK model, its kinetic entropy and the
behavior of deviations from equilibrium. Section~\ref{sec_CE_series}
derives the moment hierarchy and the first-order Chapman--Enskog closure
and identifies its diffusivity near a positive background.
Section~\ref{sec_inv_equation_general} formulates the invariance equation
for a slow spectral submanifold and explains how the reduced dynamics
yield an exact closure when the macroscopic observables provide local
coordinates on the manifold. Section~\ref{sec_linear_part} identifies
the hydrodynamic eigenvalue and its critical wave number, constructs
the linear hydrodynamic manifold, and derives the resolvent and semigroup
estimates needed for the nonlinear construction.
Section~\ref{sec_cohomology_eq} derives the generalized functional
expansion and the unsteady homological equations, and shows at quadratic
order how an external resonance can obstruct the stationary construction.
Section~\ref{sec_proof_main_theorem} selects the initial conditions
\eqref{Vn_selected_choice} for the unsteady equations and proves that
the resulting expansion converges, that its sum is differentiable at
equilibrium, and that it is $C^1$ on the punctured open ball.
Section~\ref{sec_slow_spectral_closure} constructs local kinetic charts
on this ball, establishes local density coordinates on the invariant
domain \eqref{invariant_density_domain_s9}, and derives the exact density
closure. It also proves uniform Maxwellian bounds and dissipation of
the finite relative entropy for the selected real states.
Finally, Section~\ref{sec_discussion} discusses the limitations of the construction
and possible extensions to other kinetic and dispersive equations.

\section{Basic Notation and General Definitions}\label{Notation}

\subsection{Basic Notation}
Throughout, let $d\in\mathbb N$ denote the spatial and velocity dimension, $\mathbb R^+=(0,\infty)$, and all function spaces be modeled over $\mathbb C$ unless explicitly declared real. For $z=a+\ri b\in\mathbb C$, we write $z^*=a-\ri b$.\\
The normalized global Maxwellian is given by
\begin{equation}\label{defphi}
    \phi(v) = (2\pi)^{-\frac{d}{2}}e^{-\frac{|v|^2}{2}}.
\end{equation}
The Fourier transform of an integrable function $f:\mathbb{R}^d\to\mathbb{C}$ is defined as
\begin{equation}\label{def_Fourier_transform}
    \hat{f}(k) = (2\pi)^{-\frac{d}{2}}\int_{\mathbb{R}^d} f(x) e^{-\ri x\cdot k}\, dx,
\end{equation}
which can be extended to square integrable functions by Plancherel's Theorem. We sometimes use the subscripts $x$, $v$, and $k$ on norms and function spaces to indicate the corresponding integration variable.\\
We write $\kappa=|k|$ for the wave number.
Let $w:\mathbb{R}^d\to(0,\infty)$ be a measurable weight and let $1\leq  p\leq \infty$. We write $L^p(\mathbb{R}^d,w)$, or simply $L^p(w)$ when the integration domain is clear, for the weighted Lebesgue space. For $1\leq p<\infty$ its norm is
$$
    \|f\|_{L_v^p(w)} = \left(\int_{\mathbb{R}^d} |f(v)|^p w(v)\, dv\right)^{1/p},
$$
with the usual essential-supremum convention for $p=\infty$. Analogous notation is used in the $x$- and $k$-variables. We write
$$
   \langle f,g \rangle_w =  \int_{\mathbb{R}^d}  f(v) g^*(v) w(v)\, dv,
$$
for the corresponding sesquilinear form that endows $L^2_v(w)$ with a Hilbert space structure. As mentioned in the introduction, the two natural spaces for the velocity variable are given by
$$
    p=1,\, w(v) = 1,\quad p=2,\, w(v)= \phi(v)^{-1}.
$$
For a Banach space $X$ and $1\leq q\leq\infty$, $L_x^q(X)$ denotes the Bochner space of strongly measurable $X$-valued functions, with its standard norm. In particular, for scalar functions of $(x,v)$ we identify
$$
 L_x^q(L_v^p(w))\equiv L_x^qL_v^p(w),\qquad
 \|f\|_{L_x^qL_v^p(w)}
 =\left\|x\mapsto\|f(x,\cdot)\|_{L_v^p(w)}\right\|_{L_x^q}.
$$
Note that the order of the subscripts specifies the order of integration: $L_v^p(w;L_x^q)$ generally differs from $L_x^q(L_v^p(w))$. In particular, Plancherel's theorem in $x$ identifies $L_v^p(w;L_x^2)$ with $L_v^p(w;L_k^2)$, but does not in general identify $L_x^2(L_v^p(w))$ with $L_k^2(L_v^p(w))$ unless $p=2$, while for $1\leq p\leq  2$, Minkowski's inequality immediately gives $\|f\|_{L_x^2(L_v^p(w))}\leq\|f\|_{L_v^p(w;L_x^2)}$. 
The majorant space used later is defined as 
\begin{equation}\label{def_majorant_space_Y}
 \mathcal Y:=L_k^1(\mathbb R^d)\cap L_k^2(\mathbb R^d),
 \qquad \|h\|_{\mathcal Y}:=\|h\|_{L_k^1}+\|h\|_{L_k^2}.
\end{equation}

For a Hilbert space $X$, we set
$$
 \begin{aligned}
 H_x^s(X)&:=\left\{f\in\mathcal S'(\mathbb R_x^d;X):(1+|k|^2)^{s/2}\hat f(k)\in L_k^2(X)\right\},\\
 \|f\|_{H_x^s(X)}&:=\|(1+|k|^2)^{s/2}\hat f\|_{L_k^2(X)}.
 \end{aligned}
$$
Thus, in particular, $H_x^sL_v^2(\phi^{-1})$ means $H_x^s(L_v^2(\phi^{-1}))$. We use $\Lambda=(1-\Delta_x)^{1/2}$, so that $\|f\|_{H_x^s(X)}=\|\Lambda^sf\|_{L_x^2(X)}$ with the Fourier convention above.\\
Let $K\subseteq \mathbb{R}^d$ be compact. For $r\in\mathbb N_0$, $C^r(K)$ denotes the space of functions whose derivatives up to order $r$ extend continuously to $K$, equipped with the norm
$$
    \|f\|_{C^r(K)} = \max_{|\beta|\leq  r}\sup_{x\in K}|\partial^\beta f(x)|.
$$
Similarly, we denote the space of smooth compactly supported test functions by $C_c^\infty(\mathbb{R}^d)$, and $C_0(\mathbb R^d)$ denotes the continuous functions vanishing at infinity, equipped with the supremum norm.\\
We write $B_r(x_*)$ for the closed ball of radius $r$ centered at $x_*$ and $B_r=B_r(0)$.
For two nonnegative functions $f,g$, we write $f\lesssim g$ if there exists a constant $C>0$, independent of the variables under consideration, such that $f\leq  Cg$. We write $f\asymp g$ if both $f\lesssim g$ and $g\lesssim f$ hold, with constants independent of the variables under consideration.\\
With $\delta$ denoting the Dirac distribution, the repeated convolution can be written as
\begin{equation}\label{convolution_rewrite}
    f_1\ast f_2 \ast \cdots \ast f_n = \int_{\mathbb{R}^{d\times n}}\prod_{j=1}^n f_j(x_j) \delta\left(x-\sum_{j=1}^n x_j\right)\, \prod_{j=1}^n dx_j.
\end{equation}

For a linear operator $T$, we split the spectrum into the discrete spectrum and the essential spectrum,
$$
    \sigma(T) = \sigma_{\rm disc}(T) \cup \sigma_{\rm ess}(T).
$$
The discrete spectrum consists of all isolated eigenvalues of finite algebraic multiplicity, while the essential spectrum is defined in this work as its complement. In particular, embedded eigenvalues or eigenvalues with infinite multiplicity are part of the essential spectrum. This fixes the convention for $\sigma_{\rm ess}$ used below. For a finite set of eigenvalues $\Lambda \subseteq \sigma_{\rm disc}(T)$, denote the sum of their eigenspaces by
$$
    \mathrm{eig}_T(\Lambda) = \bigoplus_{\lambda\in \Lambda}\ker(T-\lambda I),
$$
and define the corresponding Riesz projection as
\begin{equation}\label{def_Riesz_proj}
    P_{\Lambda} = -\frac{1}{2\pi \ri} \oint_{\Gamma(\Lambda)} (T-zI)^{-1}\, dz,
\end{equation}
where $\Gamma(\Lambda)$ is a positively oriented simple closed contour enclosing $\Lambda$ and no other point of $\sigma(T)$. The operator \eqref{def_Riesz_proj} is well-defined by analytic spectral calculus and for semisimple eigenvalues its range equals $\mathrm{eig}_T(\Lambda)$, while in general its range is the corresponding generalized eigenspace, see \cite{hislop2012introduction}.\\
For an operator $A$ generating a $C_0$-semigroup $T(t)=e^{tA}$, let
$$
\omega_0(A)
:=
\inf\left\{\omega\in\mathbb R:
\exists M_\omega\geq1\ \text{such that}\
\|T(t)\|\leq  M_\omega e^{\omega t}\ \forall t\geq0
\right\}
$$
denote its growth bound.\\
For a sufficiently regular map $F:X\to Y$ between Banach spaces, $D^nF(x)$ denotes its $n$th Fr\'echet derivative. When $F$ is analytic near the origin, we write its homogeneous Taylor expansion as
$$
    F(x)=\sum_{n=0}^\infty \mathcal L_F^{(n)}(x),
    \qquad
    \mathcal L_F^{(n)}(x):=\frac1{n!}D^nF(0)[x,\ldots,x].
$$
The Laplace transform of an integrable function $f:[0,\infty) \to \mathbb{C}$ is defined as
$$
    \tilde{f}(s) = \int_{0}^\infty f(t) e^{-st}\, dt.
$$

Let $\mathcal{S}_n$ denote the symmetric group on $n$ elements and write $P\in\mathcal{S}_n$ for a permutation. The symmetric group acts on frequency vectors
$$
    P.\bm{k}_n = \Big(k_{P(1)},...,k_{P(n)}\Big),
$$
and we denote the symmetrization of a function $f:\mathbb{R}^{d\times n}\to \mathbb{C}$ as
\begin{equation}\label{def_symmetrization}
    \sym[f](\bm{k}_n) = \frac{1}{n!}\sum_{P\in\mathcal{S}_n} f(P.\bm{k}_n).
\end{equation}

We denote the falling factorials as
\begin{equation}\label{falling_factorials}
(x)_n = \prod_{j=0}^{n-1} (x-j).
\end{equation}
For a formal scalar power series,
$$
    W(z) = \sum_{n=0}^{\infty} W_n z^n,
$$
we write
\begin{equation}
    [z^n]W(z) = W_n,
\end{equation}
for the $n^{\rm th}$ Taylor coefficient. For a manifold $\mathcal{M}$ and $p\in\mathcal{M}$, let $T_p\mathcal{M}$ denote the tangent space at $p$.

\subsection{The Plasma Dispersion Function}

Let
\begin{equation}\label{defZ}
Z(\zeta) = \frac{1}{\sqrt{2\pi}}\int_{\mathbb{R}} \frac{e^{-\frac{v^2}{2}}}{v-\zeta}\, dv,\quad \Im\zeta\neq 0,
\end{equation}
denote the plasma dispersion function \cite{fried2015plasma}. An explicit formula for $Z$ in terms of more elementary functions is given \begin{equation}\label{Zexplicit}
    Z(\zeta) = \ri \sqrt{\frac{\pi}{2}} e^{-\frac{\zeta^2}{2}}\left[\text{sign}(\Im \zeta)- \text{erf}\left(\frac{-\ri \zeta}{\sqrt{2}}\right)\right],\quad \Im\zeta\neq 0,
\end{equation}
\cite{ReviewAbramowitzStegun1964}. It satisfies the following conjugation and symmetry relations,
\begin{equation}
    Z(-\zeta) = - Z(\zeta),\quad Z(\zeta^*) = Z(\zeta)^*.
\end{equation}
These identities relate the two half-plane restrictions of the Cauchy transform. Each restriction has a distinct entire continuation, denoted by $Z^+$ and $Z^-$, respectively. This also follows directly from the explicit formula \eqref{Zexplicit} for either sign of $\Im\zeta$ such that
\begin{equation}
    Z^{\pm}(\zeta) = \ri \sqrt{\frac{\pi}{2}} e^{-\frac{\zeta^2}{2}}\left[\pm 1- \mathrm{erf}\left(\frac{-\ri \zeta}{\sqrt{2}}\right)\right],\quad \zeta\in\mathbb{C}. 
\end{equation}
On the real line, the principal value differs from either boundary trace:
$$
 Z^\pm(x)=\frac{1}{\sqrt{2\pi}}\operatorname{PV}
 \int_{\mathbb R}\frac{e^{-v^2/2}}{v-x}\,dv
 \ \pm\ri\sqrt{\frac{\pi}{2}}e^{-x^2/2},\qquad x\in\mathbb R.
$$
In particular, $Z^+(x)-Z^-(x)=\ri\sqrt{2\pi}e^{-x^2/2}$ for $x\in\mathbb{R}$. The function $Z$ is a solution to the complex differential equation,
\begin{equation}
    \frac{dZ}{d\zeta} = -1-\zeta Z,
\end{equation}
as follows by integration by parts in \eqref{defZ}, which can be used as a polynomial recurrence relation for the derivatives of $Z$.
The plasma dispersion function satisfies the asymptotics 
\begin{equation}\label{asymptoticsZ}
Z(\ri y)=-\frac{1}{\ri y}-\frac{1}{(\ri y)^3}+\mathcal O(y^{-5}),\qquad y\to+\infty,
\end{equation}
on the positive imaginary axis, which is the asymptotic direction needed below. This follows by expanding the Cauchy kernel and using the Gaussian moments. Note that no branch-independent expansion over the whole complex plane can be asserted.\\
For later calculations, we note that the lower bound 
\begin{equation}\label{lower_Bound_Z}
    |-\ri y Z(\zeta)-1|\geq1-y\sqrt{\frac{\pi}{2}}\,e^{y^2/2}\left[1-\mathrm{erf}\left(\frac{y}{\sqrt{2}}\right)\right]>0,\quad y= \Im\zeta>0,
\end{equation}
holds. Indeed,
$$
 \Re[-\ri yZ(x+\ri y)]
 =\frac{1}{\sqrt{2\pi}}\int_{\mathbb R}
 \frac{y^2e^{-v^2/2}}{(v-x)^2+y^2}\,dv<1.
$$
The convolution of the two nonnegative symmetric decreasing functions in this integral is maximal at $x=0$ (as follows, for example, by writing their level sets as centered intervals). Thus
$|1+\ri yZ(x+\ri y)|\geq1-\Re[-\ri yZ(x+\ri y)]\geq1+\ri yZ(\ri y)>0$ and
evaluating the last expression with \eqref{Zexplicit} proves \eqref{lower_Bound_Z}.

\section{A BGK Kinetic Equation with Density-Dependent Collision Rate}\label{sec_main_equation}

In this section, we introduce the underlying kinetic model of our analysis, specialize to Maxwell molecules, and establish mass conservation, entropy dissipation, and the characterization of equilibria. We then linearize about a global Maxwellian and express the nonlinear perturbation equation in Fourier space as a symmetrized convolution series, preparing the spectral and invariant-manifold analysis below.

\subsection{A BGK-Type Kinetic Equation with a Density-Dependent Collision Rate}

Granular media motivate kinetic models in which the local density
controls the rate of collisional relaxation. In a dilute granular
gas, the collision frequency per particle is proportional to the
number density and a characteristic relative speed: for hard
spheres, this gives a frequency proportional to
$\rho\sqrt{T}$, up to constants determined by the particle
properties. Inelastic encounters additionally dissipate kinetic
energy, coupling collisional relaxation to the evolution of the
granular temperature
\cite{Haff1983,Goldhirsch2003,BrilliantovPoeschel2004,ReviewGarzo2019,Villani2006}.
To isolate the effect of density-dependent relaxation, we study
a mass-conserving BGK equation with collision rate
$\rho[F]^\alpha/\tau$ and prescribed Maxwellian target
$\rho[F]\phi$. At a fixed thermal speed, $\alpha=1$ captures
the dilute-gas density scaling, while other positive exponents
provide a phenomenological family of density-dependent rates.
The prescribed Maxwellian represents relaxation towards a
background with fixed velocity and temperature. Consequently,
momentum and energy are exchanged with that background, and
the model does not describe freely cooling granular gases.
This simplification retains a nonlinear coupling between density
and kinetic relaxation while making the hydrodynamic spectral
branch and the invariant-manifold construction explicit.
It therefore provides a tractable setting in which to study
exact nonlinear hydrodynamic reduction and its relation to the
nonlinear diffusion equation obtained from the first-order
Chapman--Enskog approximation.\\
Thus consider the kinetic equation
\begin{equation}\label{maineq_full}
    \frac{\partial F}{\partial t}+v\cdot\nabla_xF
    =\frac{1}{\tau} Q_\alpha[F],
\end{equation}
for a non-negative distribution function
$F:\mathbb{R}^d\times\mathbb{R}^d\times[0,\infty)\to[0,\infty)$ and a
Knudsen number $\tau>0$.
We take the collision operator to be
\begin{equation}\label{defQ}
   Q_\alpha[F]
    =\rho[F]^\alpha\int_{\mathbb{R}^d}
    \bigl[F(v^*)\phi(v)-F(v)\phi(v^*)\bigr]
    |v-v^*|^\gamma\,dv^*,
\end{equation}
where $\alpha\geq0$ determines the density dependence of the collision rate,
$\gamma\geq0$ describes the velocity dependence of the collision kernel,
$\phi$ is the normalized global Maxwellian \eqref{defphi}, and
\begin{equation}\label{defrho_section3}
    \rho[F](x,t)=\int_{\mathbb{R}^d}F(x,v,t)\,dv
\end{equation}
is the mass density. The macroscopic velocity field is defined as
\begin{equation}\label{defu}
u(x,t) = \frac{1}{\rho(x,t)} \int_{\mathbb{R}^d} v F(x,v,t)\, dv. 
\end{equation}
The density dependence in \eqref{defQ} is motivated by the increased collision frequency in dense media. One-dimensional hard-rod systems provide a classical
setting in which excluded-volume effects lead to density-dependent macroscopic transport, see \cite{boldrighini1983hardrods,boldrighini1997navierstokes,tonks1936complete}. Nonlinear BGK equations with density-dependent relaxation rates, including their
well-posedness, hypocoercive relaxation and nonlinear-diffusion limits, have been studied in \cite{dolbeault2007nonlinear,evans2026density}. A related density-dependent Fokker--Planck model, with a different collision operator, is studied in \cite{anceschi2024spatially}. Related nonlinear relaxation mechanisms also occur in radiative-transfer models in the Rosseland
regime, where the opacity depends on the material state
\cite{clouet1997rosseland,golse2008rosseland}. The model considered here should therefore be viewed as a deliberately simple relaxation model that retains this density dependence, rather than as a direct microscopic derivation from any one of these systems.\\
In the remainder of the paper we consider the analogue of Maxwell molecules,
namely $\gamma=0$. In this case, \eqref{defQ} reduces to
\begin{equation}\label{Qalpha}
    Q_\alpha[F]=\rho[F]^\alpha\bigl(\rho[F]\phi-F\bigr),
\end{equation}
and \eqref{maineq_full} becomes the nonlocal partial differential equation
\begin{equation}\label{maineq}
    \frac{\partial F}{\partial t}+v\cdot\nabla_xF
    =\frac{\rho[F]^\alpha}{\tau}\bigl(\rho[F]\phi-F\bigr).
\end{equation}
Note that the nonlocality in \eqref{maineq} enters only through the velocity integral \eqref{defrho_section3}.\\
While the full Boltzmann equation has $d+2$ collision invariants (five when $d=3$), the collision operator \eqref{defQ} only preserves the zeroth moment in velocity space
\begin{equation}\label{colinv1}
    \int_{\mathbb{R}^d} Q_\alpha[F]\,dv = 0, 
\end{equation}
which follows from exchanging $v$ and $v^*$ in one of the two terms. As a consequence of \eqref{colinv1}, assuming that the integrations below are justified, for example by sufficient spatial decay and regularity, the total mass \begin{equation}\label{rhobar}
    \overline{\rho}
    =\int_{\mathbb{R}^d}\int_{\mathbb{R}^d}F(x,v,t)\,dv\,dx
\end{equation}
is conserved whenever it is finite. We stress, however, that for the reduced collision operator \eqref{Qalpha}, momentum and kinetic energy are not collision invariants: their relaxation, related to the phenomenon of cooling in granular media, is discussed in Section~\ref{sec_moment_closure}. \\
Spatially homogeneous Maxwellians
\begin{equation}
    F^{\rm eq}(v)=\rho_0\phi(v),\qquad \rho_0\geq0,
\end{equation}
are equilibria of \eqref{maineq}. The positive homogeneous Maxwellians form the equilibrium manifold
\begin{equation}\label{defMeq}
    \mathcal{M}_{\rm eq}
    =\{\rho_0\phi(v):\rho_0\in\mathbb{R}^+\}.
\end{equation}
The closure of $\mathcal M_{\rm eq}$ in the density parameter is
$\mathcal M_{\rm eq}\cup\{0\}$, where $0$ denotes the vacuum equilibrium.
For signed densities the formula makes algebraic sense whenever $\rho^\alpha$ is defined.
For $\rho_0>0$, the homogeneous Maxwellian $\rho_0\phi$ has infinite mass on
$\mathbb R_x^d\times\mathbb R_v^d$. We therefore study spatially localized deviations
from a positive homogeneous Maxwellian.

\subsection{Relative Entropy, Lyapunov Functions and Deviations from Equilibrium}

Define the local relative-entropy functional\footnote{We employ the mathematical convention for entropy. For spatially homogeneous solutions on $\mathbb R^d_x$, monotonicity refers to the entropy per unit spatial volume, under the finiteness assumptions below.}
\begin{equation}\label{defentropy}
    \mathcal{S}[F]
    =\int_{\mathbb{R}^d}\int_{\mathbb{R}^d}
      F\log\left(\frac{F}{\rho[F]\phi}\right)\,dv\,dx.
\end{equation}
At every fixed $(x,t)$ with $\rho[F](x,t)>0$, the normalized function $F/\rho[F]$ is a probability density in velocity. Thus $\mathcal S[F]$ is the spatial integral of $\rho[F]$ times the Kullback--Leibler divergence of the velocity probability densities $F/\rho[F]$ and $\phi$. It is nonnegative by Gibbs' inequality \cite{cover1999elements}, with equality precisely for local Maxwellians $F(x,v)=\rho(x)\phi(v)$. On vacuum sets the integrand is defined to be zero. The spatially inhomogeneous transport term, however, prevents \eqref{defentropy} from
being a Lyapunov functional in general and the H-Theorem for equation \eqref{maineq} takes the form of the following balance law.

\begin{proposition}[Local Relative-Entropy Balance: H-Theorem]
Let $F>0$ be a sufficiently regular solution of \eqref{maineq} with finite $\mathcal S[F]$, integrable terms in the balance below, and vanishing spatial boundary fluxes, so that differentiation under the integral and integration by parts are justified, and let $u$ be its velocity field as defined in \eqref{defu}. Then
\begin{equation}\label{entropy_balance}
\begin{split}
    \frac{d}{dt}\mathcal{S}[F]
    &=-\int_{\mathbb{R}^d}(\rho u)\cdot\nabla_x\log\rho[F] \,dx\\
    &\quad-\frac{1}{\tau}\int_{\mathbb{R}^d}\rho[F]^\alpha
       \int_{\mathbb{R}^d}
       \bigl(F-\rho[F]\phi\bigr)
       \log\left(\frac{F}{\rho[F]\phi}\right)\,dv\,dx.
\end{split}
\end{equation}
The collision contribution (the second term) is non-positive. Moreover, $\mathcal S[F]\geq 0$, the functional $\mathcal S$ is convex on its natural positive cone, and $\mathcal S[F]=0$ if and only if $F$ is a local Maxwellian,
    \begin{equation}
        \phi_{\rm loc}(x,v) = \rho(x) \phi(v). 
    \end{equation}
On any time interval on which the transport contribution in \eqref{entropy_balance}
vanishes for almost every time, $t\mapsto\mathcal S[F(t)]$ is non-increasing.
If the solution exists globally, $\mathcal S[F(0)]<\infty$, and this contribution
vanishes for almost every $t\geq 0$, then $\mathcal S[F(t)]$ converges to a limit
in $[0,\mathcal S[F(0)]]$ as $t\to\infty$. For spatially homogeneous solutions on $\mathbb R^d_x$, the same assertion applies to the entropy per unit spatial volume, $\int F\log(F/(\rho\phi))\,dv$, when finite, while its integral over all of space is generally infinite.
\end{proposition}

\begin{proof}
Differentiating \eqref{defentropy} and using $\rho[F]=\int F\,dv$ gives
$$
    \frac{d}{dt}\mathcal{S}[F]
    =\int_{\mathbb{R}^d}\int_{\mathbb{R}^d}
      \partial_tF\log\left(\frac{F}{\rho[F]\phi}\right)\,dv\,dx,
$$
while the terms produced by differentiating $F$ and $\rho[F]$ outside the logarithm cancel. Indeed, $\int \partial_t F\,dv=\partial_t\rho$ and $\int(F/\rho)\,\partial_t\rho\,dv=\partial_t\rho$, so the two terms cancel pointwise in $x$.  For the transport part, integration by parts in $x$ gives
\begin{align*}
    -\int\!\!\int v\cdot\nabla_xF
       \log\left(\frac{F}{\rho\phi}\right)\,dv\,dx
    =\int\!\!\int Fv\cdot\nabla_x
       \log\left(\frac{F}{\rho\phi}\right)\,dv\,dx
    =-\int_{\mathbb{R}^d}(\rho u)\cdot\nabla_x\log\rho\,dx,
\end{align*}
where $\int\!\int v\cdot\nabla_xF\,dv\,dx=0$ under the assumed spatial boundary behavior and $\phi$ is independent of $x$.  For the collision part, \eqref{maineq} and
\begin{equation}\label{a_b_inequality}
    (a-b)\log\frac{a}{b}\geq0,\qquad a,b>0,
\end{equation}
give the second term in \eqref{entropy_balance} and show that it is non-positive.\\
For completeness, convexity follows from the second variation.  If $G$ is an admissible variation and $\rho[G]=\int G\,dv$, then
$$
    D^2\mathcal S[F](G,G)
    =\int_{\mathbb R^d_x}
      \left(
      \int_{\mathbb R^d_v}\frac{G^2}{F}\,dv
      -\frac{\rho[G]^2}{\rho[F]}
      \right)dx\geq 0,
$$
by Cauchy--Schwarz in $v$ or, equivalently, Jensen's inequality for the probability measure $F\,dv/\rho[F]$.  Finally, for each $x$, $\mathcal S$ is $\rho[F](x)$ times the relative entropy of the probability density $F/\rho[F]$ with respect to $\phi$.  Gibbs' inequality therefore gives $\mathcal S[F]\geq 0$, with equality exactly when $F=\rho[F]\phi$ almost everywhere. When the transport contribution vanishes, \eqref{entropy_balance} gives
$\frac{d}{dt}\mathcal S[F(t)]\le0$ almost everywhere on the interval in question.
On $[0,\infty)$, monotonicity and the lower bound give the asserted limit.
\end{proof}

\begin{remark}
Let $\Phi\colon\mathbb R^+\to\mathbb R$ be a continuously differentiable convex function.
For $F>0$, set $g=F/\phi$. The relative entropy with respect to a homogeneous
background $\rho_0\phi$, $\rho_0>0$, is
\begin{equation}\label{def_H_Phi_relative}
 \mathcal H_{\Phi,\rho_0}^{\rm rel}[F]
 :=\iint\phi(v)\left[\Phi(g)-\Phi(\rho_0)
       -\Phi'(\rho_0)(g-\rho_0)\right]\,dv\,dx.
\end{equation}
Assume that this nonnegative integral is finite, its collision production is integrable, and the corresponding renormalized spatial flux vanishes, with differentiation under the integral justified. Since $\int\phi(g-\rho[F])\,dv=0$, subtraction of the affine term does not change the production, and
\begin{equation}\label{relative_convex_entropy_dissipation}
 \frac{d}{dt}\mathcal H_{\Phi,\rho_0}^{\rm rel}[F]
 =-\frac1\tau\int\rho[F]^\alpha
       \int\phi\bigl(\Phi'(g)-\Phi'(\rho[F])\bigr)
                  \bigl(g-\rho[F]\bigr)\,dv\,dx\le0.
\end{equation}
For $\Phi(z)=z\log z$, its relation to the local relative entropy is
\begin{equation}\label{relative_convex_entropy_decomposition}
 \mathcal H_{\Phi,\rho_0}^{\rm rel}[F]
 =\mathcal S[F]+\int\left[\rho[F]\log\frac{\rho[F]}{\rho_0}
                         -\rho[F]+\rho_0\right]dx.
\end{equation}
These renormalized identities avoid subtracting divergent whole-space integrals. Nonnegative states with vacuum are included only when the entropy has a suitable extension at zero and a positive regularization justifies passage to the limit in these identities.
\end{remark}

The entropy \eqref{defentropy} can be used to characterize global equilibria of \eqref{maineq}. Indeed, under the regularity and spatial boundary assumptions stated below, the scalar multiples of the normalized global Maxwellian are the only nonnegative equilibria of equation \eqref{maineq}.
\begin{proposition}[Characterization of Equilibria]
\label{prop_equilibria}
 Let $F \geq0$ be a sufficiently regular stationary solution, with finite local entropy production and sufficient integrability to justify the entropy identity below, of the stationary kinetic equation
\begin{equation}
    v\cdot\nabla_x F
    =
    \frac{\rho[F]^\alpha}{\tau}
    \bigl(\rho[F]\phi-F\bigr). 
    \label{eq:stationary-bgk}
\end{equation}
Assume that the spatial entropy flux vanishes at infinity, in the sense that
\begin{equation}
    \int_{\mathbb{R}^d\times\mathbb{R}^d}
    v\cdot\nabla_x
    \left[
        F\log\left(\frac{F}{\rho_0\phi}\right)-F
    \right]
    \,dv\,dx
    =0
    \label{eq:entropy-flux-assumption}
\end{equation}
for some constant $\rho_0>0$. Then $F$ is a scalar multiple of the global Maxwellian. 
\end{proposition}
\begin{proof}
For the details of this elementary proof, see Appendix \ref{app_equilibria_characterization}. 
\end{proof}

\subsection{Deviations from Equilibrium and Long-Time Behavior of Solutions}

We proceed with the analysis of solutions to \eqref{maineq} close to the manifold of global equilibria. To this end, we write a general solution to \eqref{maineq} as the sum of an equilibrium distribution and a distribution measuring the deviation from equilibrium. Fix $\rho_0>0$ and write
\begin{equation}\label{deviation_ansatz}
    F(x,v,t)=\rho_0\phi(v)+g(x,v,t).
\end{equation}
The Fr\'echet derivative of the collision operator at $\rho_0\phi$ is
$$
DQ_\alpha[\rho_0\phi]\,g
    =\left.\frac{d}{d\varepsilon}Q_\alpha[\rho_0\phi+\varepsilon g]\right|_{\varepsilon=0}
    =\rho_0^\alpha\bigl(\rho[g]\phi-g\bigr),
$$
and we write 
\begin{equation}\label{def_LQ}
    \mathcal{L}_Q[g]=\rho[g]\phi-g,
\end{equation}
for the linearized collision operator. The deviation $g$ satisfies the equation
\begin{equation}
    \frac{\partial g}{\partial t}+v\cdot\nabla_xg
    =\frac{1}{\tau}\bigl(\rho_0+\rho[g]\bigr)^\alpha
      \mathcal{L}_Q[g].
\end{equation}
Setting $g=\rho_0 f$ and introducing the effective Knudsen number
$\tau_{\rm eff}=\tau/\rho_0^\alpha$, then relabelling
$\tau_{\rm eff}$ as $\tau$, gives the normalized equation
\begin{equation}\label{maineqf}
    \frac{\partial f}{\partial t}+v\cdot\nabla_xf
    =\frac{1}{\tau}\bigl(1+\rho[f]\bigr)^\alpha
      \mathcal{L}_Q[f].
\end{equation}
The reference density changes the effective relaxation parameter to $\tau_{\rm eff}=\tau/\rho_0^\alpha$, and hence changes the collision rate relative to transport at fixed spatial and velocity scales. A rescaling of time alone would also rescale the transport term. Splitting \eqref{maineqf} into its linear and nonlinear parts gives
\begin{equation}\label{maineqf_split}
    \frac{\partial f}{\partial t}
    =\mathcal{L}[f]+\mathcal{N}_\alpha(f),
\end{equation}
where
\begin{equation}
    \mathcal{L}=-v\cdot\nabla_x+\frac{1}{\tau}\mathcal{L}_Q
\end{equation}
is the linear part of \eqref{maineq} around a global Maxwellian and
\begin{equation}\label{defNalpha}
    \mathcal{N}_\alpha(f)
    =\frac{1}{\tau}\sum_{n=1}^{\infty}
      \binom{\alpha}{n}\rho[f]^n\mathcal{L}_Q[f],
\end{equation}
is the nonlinear part.

\begin{remark}
For non-integer $\alpha$, the binomial series in \eqref{defNalpha} has radius one in the scalar variable $\rho[f]$ and is therefore pointwise analytic
whenever $\|\rho[f]\|_{C_x^0}<1$. Accordingly, the analytic germ used in the invariant-manifold construction is
naturally restricted to a function-space neighborhood whose spatial topology controls the $C_x^0$-norm and supports the required multiplication estimates, such as the Wiener algebra of scalar functions with integrable Fourier transform and its kinetic module. We stress that the fixed-bandwidth parameter space of Section~\ref{sec_cohomology_eq} is not an algebra: products of functions with Fourier support in $K$ have support in $K+K$. The construction therefore retains all nonlinear summed frequencies and uses convolution estimates in the ambient Fourier spaces, since an
$L_x^2$ assumption alone does not provide the required Banach-algebra property. In the normalized variables, the same condition also guarantees that the total density $1+\rho[f]$ remains bounded away from zero for
non-negative physical distributions. Possible continuation beyond this local analytic chart is discussed in Section~\ref{sec_discussion}. \end{remark}

\begin{proposition}[Convergence to a global Maxwellian]
\label{prop_global_existence}
Let $\alpha\geq0$, $\tau>0$ and assume that
$$
    s>\frac{d}{2}+1. 
$$
Then there exists an $\varepsilon_0>0$ such that the following holds. Let
$$
F_0(x,v)=\rho_0\phi (v)+\rho_0 f_0(x,v),
    \qquad \rho_0>0,
$$
with $F_0\geq0$ such that
$$
    f_0\in
    H^s_xL^2_v(\phi ^{-1})
    \cap
    L^2_v(\phi ^{-1};L^1_x),
    \qquad
    \|f_0\|_{H^s_xL^2_v(\phi ^{-1})}
    \leq \varepsilon_0.
$$
Then the normalized perturbation equation 
\eqref{maineqf} has a unique global small mild solution in
$C([0,\infty);H_x^sL_v^2(\phi^{-1}))$ satisfying
$$
    \sup_{t\geq0}
    \|f(t)\|_{H^s_xL^2_v(\phi ^{-1})}
    \leq 
    C\|f_0\|_{H^s_xL^2_v(\phi ^{-1})}.
$$
Moreover,
\begin{equation}\label{decay_estimate}
    \|f(t)\|_{L^2_xL^2_v(\phi ^{-1})}
    \leq 
    C(1+t)^{-d/4}
    \left(
        \|f_0\|_{H^s_xL^2_v(\phi ^{-1})}
        +
        \|f_0\|_{L^2_v(\phi ^{-1};L^1_x)}
    \right),
\end{equation}
and consequently 
$$
\|F(t)-\rho_0\phi\|_{L_x^2L_v^2(\phi^{-1})}
    =\rho_0\|f(t)\|_{L_x^2L_v^2(\phi^{-1})}
    \longrightarrow0\qquad\text{as }t\to\infty.
$$
\end{proposition}
The proof of Proposition~\ref{prop_global_existence} is given in Appendix~\ref{proof_global_existence}. Estimate \eqref{decay_estimate} gives relaxation to the background Maxwellian for
localized Sobolev data, whose exponent agrees with the heat rate for linear
mass-conserving kinetic equations on the whole space \cite{BGKproofBouin}.

\begin{remark}
The polynomial decay reflects the accumulation of arbitrarily slow spatial Fourier modes near $k=0$, not algebraic decay of each fixed nonzero mode. Indeed, the hydrodynamic eigenvalue computed below satisfies $\lambda(k)=-\tau|k|^2+\mathcal O(|k|^4)$ by \eqref{Taylorlambda}. Thus, for bounded spectral amplitudes near the origin, the corresponding linear contribution has the scaling
$$
 \int_{|k|<\delta} e^{2\lambda(k)t}\,dk\asymp (1+t)^{-d/2},
$$
which yields the $L^2$ decay scale $(1+t)^{-d/4}$ after taking a square root. We emphasize that this back-of-the-envelope calculation is a spectral explanation of the rate, not an assertion that the nonlinear solution is exactly a superposition of linear eigenmodes or that the stated estimate is optimal for every datum. On $\mathbb R^d$, the limiting Maxwellian must be the prescribed background $\rho_0\phi$: a nonzero constant difference between two Maxwellians does not belong to the assumed spatial $L^2$ class.
\end{remark}

\subsection{Set-Up for the Construction of Hydrodynamic Manifolds in Frequency Space}

We first fix the functional analytic setting for the subsequent analysis by introducing several Banach spaces. Let $1\leq  p<\infty$ and set
\begin{equation}\label{defX}
    \mathcal{X}_v=L_v^p(w).
\end{equation}
Fourier inversion maps $L_k^1(\mathcal X_v)$ into
$L_x^\infty(\mathcal X_v)$ for every $1\leq p<\infty$.
A general $L_k^2(\mathcal X_v)$ to $L_x^2(\mathcal X_v)$ Plancherel
bound is available here only when $p=2$. For the particular kernels
constructed below, the assertion for $1\leq p<2$ follows instead from
the auxiliary Hilbert-space argument in
Proposition~\ref{prop_majorant_functional_taylor}.
For $p>2$, we retain only the explicitly proved frequency-space and
$L_x^\infty(\mathcal X_v)$ conclusions.\\
We assume throughout that the weight satisfies
\begin{equation}\label{ass_phi_w}
 \begin{gathered}
 (1+|v|)^n\phi\in L_v^p(w),\qquad n\geq 0,\\
 w^{-1/(p-1)}\in L_v^1\quad(1<p<\infty),\qquad
 w^{-1}\in L_v^\infty\quad(p=1).
 \end{gathered}
\end{equation}
These assumptions include the natural spaces $L_v^2(\phi^{-1})$ and $L_v^1$ for which, in particular,
$$
    \|f\|_{L_v^1}\leq  \|f\|_{L_v^2(\phi^{-1})}.
$$
We use the weighted velocity space and the two Fourier spaces
\begin{equation}\label{def_Xv2_aux}
 \begin{gathered}
 \mathcal X_v^{(2)}:=\{f:(1+|v|)^2f\in\mathcal X_v\},\qquad
 \|f\|_{\mathcal X_v^{(2)}}:=\|(1+|v|)^2f\|_{\mathcal X_v},\\
 \mathcal Z:=L_k^1(\mathcal X_v)\cap L_k^2(\mathcal X_v),\qquad
 \mathcal Z^{(2)}:=L_k^1(\mathcal X_v^{(2)})\cap L_k^2(\mathcal X_v^{(2)}).
 \end{gathered}
\end{equation}
Both intersections carry the sum norm. With $\mathcal Y$ from
\eqref{def_majorant_space_Y}, this gives
$$
 \|f\|_{\mathcal Z^{(2)}}
 =\bigl\|\,\|f(\cdot)\|_{\mathcal X_v^{(2)}}\,\bigr\|_{\mathcal Y},
 \qquad \|f\|_{\mathcal Z}\leq\|f\|_{\mathcal Z^{(2)}}.
$$
The two velocity norms remain distinct: multiplication by $k\cdot v$
maps $\mathcal X_v^{(2)}$ into $\mathcal X_v$ at each fixed $k$,
as quantified in \eqref{C1_generator_weight_loss}.\\
We next reformulate our main equation \eqref{maineq_full} in frequency space by passing to Fourier variables in $x$. With the Fourier convention \eqref{def_Fourier_transform}, equation \eqref{maineqf_split} becomes
\begin{equation}\label{maineqfk}
    \frac{\partial\hat f}{\partial t}
    =\mathcal{L}_k\hat f+\widehat{\mathcal{N}}_\alpha(\hat f),
\end{equation}
with
\begin{equation}\label{def_Lk_section3}
    \mathcal{L}_k[\hat f](v)
    =-\ri(v\cdot k)\hat f(v)-\frac{1}{\tau}\hat f(v)
      +\frac{1}{\tau}\hat\rho\,\phi(v),
\end{equation}
the \emph{fiber operator}, where $\hat\rho(k)=\rho[\hat f(k,\cdot)]$. Note that the multiplier $-\ri k\cdot v$ and the convolutions with real coefficients
in \eqref{def_Nhat} preserve $\hat f(-k,v)=\hat f(k,v)^*$. Thus the Fourier
equation preserves the reality condition, while its individual velocity
fibers can be studied in the complex space $\mathcal X_v$.\\
The nonlinear term takes the form
\begin{equation}\label{def_Nhat}
    \widehat{\mathcal{N}}_\alpha(\hat f)
    =\frac{1}{\tau}\sum_{n=1}^{\infty}
      \binom{\alpha}{n}(2\pi)^{-\frac{nd}{2}}
      \hat\rho^{\ast n}\ast\mathcal{L}_Q[\hat f].
\end{equation}
The above series is well-defined, for example, for $\hat f\in \mathcal Z$ satisfying
$$(2\pi)^{-d/2} C_\rho\|\hat f\|_{L_k^1(\mathcal X_v)}<1,$$ where $C_\rho$ is the
norm bound for the density functional in \eqref{defCrho}.
The convolution algebra estimates then give convergence in the same
space $\mathcal Z$. Note that membership in physical-space $L_x^2$
alone does not ensure that the nonlinear products or the binomial
series are defined.
Using \eqref{convolution_rewrite} and \eqref{def_symmetrization}, we rewrite the frequency-space nonlinearity in symmetric form as

\begin{equation}
\begin{split}
\widehat{\mathcal{N}}_\alpha(\hat f)(k,v)
    &=\frac{1}{\tau}\sum_{n=1}^{\infty}
      \binom{\alpha}{n}\frac{(2\pi)^{-\frac{nd}{2}}}{n+1}
\int_{\mathbb{R}^{d\times(n+1)}}
      \delta\bigl(k-\Sigma\bm{k}_{n+1}\bigr)
      \sum_{j=1}^{n+1}
      \mathcal{L}_Q[\hat f](k_j,v)
      \prod_{\substack{\ell=1\\ \ell\neq j}}^{n+1}
      \hat\rho(k_\ell)\,d\bm{k}_{n+1}.
\end{split}
\end{equation}
The factor $1/(n+1)$ is the normalization in the symmetrization, and the
integral contains $n+1$ Fourier variables because the term
$\rho[f]^n\mathcal{L}_Q[f]$ is homogeneous of degree $n+1$ in $f$.\\
Equation \eqref{maineqfk} is the starting point for our further analysis, in particular the spectral analysis carried out in Section~\ref{sec_linear_part}. We emphasize that the deviation $f$ is taken in a spatial class admitting a Fourier transform, whereas the full distribution $F=\rho_0\phi+\rho_0f$ contains an $x$-independent equilibrium component and need not itself be spatially integrable.

\section{The Moment Closure Problem and the Chapman--Enskog Series for the Mass Density}\label{sec_CE_series}

For completeness, we recall the first-order Chapman--Enskog approximation to the mass flux
and compare the resulting diffusion equation near a positive homogeneous
background with the low-frequency kinetic dynamics. Towards the end, we comment on the finite existence time of solutions to the fast diffusion equation. 

\subsection{The Moment Closure Problem}\label{sec_moment_closure}

Taking velocity moments of a kinetic equation produces an infinite hierarchy of balance laws: the equation for a moment of order $j$ contains a flux of order $j+1$ \cite{Cercignani1988,Grad1949,Struchtrup2005}. A macroscopic model therefore requires a closure relation expressing the unresolved moments in terms of a finite set of retained variables.
Physically, this corresponds to specifying higher-order tensors through a constitutive law, such as Fourier's law of heat conduction, relating the heat flux to the gradient of temperature, or Newton's law of viscosity, relating the deviatoric stress tensor to the gradient of the velocity field. The \emph{Chapman--Enskog series} constructs such a relation perturbatively in the Knudsen number by expanding around local equilibrium \cite{chapman1990mathematical}. For the Boltzmann equation it yields the Euler equations at leading order, the Navier--Stokes--Fourier system at first order, and Burnett-type corrections at higher orders.\\
For the BGK equation \eqref{maineq}, define the pressure tensor
$$
    \Pi=\int_{\mathbb{R}^d}(v-u)\otimes(v-u)F\,dv,
$$
as well as the temperature and the heat flux, 
$$
    T:=\frac{1}{d\rho}\int_{\mathbb{R}^d}|v-u|^2F\,dv,
    \qquad
    q:=\frac12\int_{\mathbb{R}^d}|v-u|^2(v-u)F\,dv.
$$
Assuming sufficient decay in $x$ and $v$, the first three balance laws are then given by
\begin{equation}\label{eq_moment_hierarchy}
\begin{aligned}
    \partial_t\rho+\nabla_x\cdot(\rho u)&=0,\\
    \partial_t(\rho u)+\nabla_x\cdot(\rho u\otimes u+\Pi)
        &=-\frac{\rho^{\alpha+1}}{\tau}u,\\
    \frac d2\partial_t(\rho T)
        +\nabla_x\cdot\left(\frac d2\rho T u+q\right)
        +\Pi:\nabla_xu
        &=\frac{\rho^{\alpha+1}}{\tau}
          \left(\frac{|u|^2}{2}-\frac{dT}{2}+\frac d2\right). 
\end{aligned}
\end{equation}
The last constant term is essential as the collision operator relaxes the distribution toward the fixed Maxwellian $\rho\phi$, whose temperature is one. Equation \eqref{eq_moment_hierarchy} is not closed because $\Pi$ and $q$ require additional evolution equations involving still higher moments. The first-order Chapman--Enskog approximation below supplies the corresponding constitutive law for the mass flux under the assumption of small Knudsen number $\tau$.

\begin{remark}[Velocity Relaxation and Cooling]
While the total mass \eqref{rhobar} of equation \eqref{maineq} is conserved, the relaxation terms in \eqref{eq_moment_hierarchy} show that the collision operator does not conserve
momentum or kinetic energy. In the spatially homogeneous setting, where the density is
constant, the macroscopic velocity satisfies
$$
\frac{d}{dt}u=-\frac{\rho^\alpha}{\tau}u,
$$
and therefore $u(t)=u_0e^{-at}$, where $a=\rho^\alpha/\tau$. The temperature equation in \eqref{eq_moment_hierarchy} instead gives
$$
 T'=a\left(1-T+\frac{|u|^2}{d}\right),\qquad
 T(t)=1+(T_0-1)e^{-at}
       +\frac{|u_0|^2}{d}\left(e^{-at}-e^{-2at}\right).
$$
The mean-velocity contribution cannot be omitted unless $u_0=0$. In particular, $T$ need not relax monotonically. The normalized total second moment
$\Theta=T+|u|^2/d$ satisfies $\Theta'=a(1-\Theta)$ and relaxes exponentially to one. This behavior should be distinguished from the unconditional 
\emph{granular cooling} observed in kinetic models for granular gases \cite{BrilliantovPoeschel2004}, where inelastic collisions
dissipate thermal energy and drive the temperature itself to zero according to Haff's law \cite{Haff1983,Goldhirsch2003}. The present BGK model instead describes collisional relaxation towards a fixed Maxwellian with exponential damping of the mean velocity and relaxation of the excess thermal energy. 
\end{remark}

\subsection{The Chapman--Enskog Series}

For completeness, we perform the leading-order Chapman--Enskog calculation in this section to derive the scalar diffusion closure associated with equation \eqref{maineq}. We seek a formal Taylor expansion of the form 
\begin{equation}
    F=F^{(0)}+\tau F^{(1)}+\tau^2F^{(2)}+\cdots,
    \qquad
\int_{\mathbb{R}^d}F^{(j)}\,dv=0\quad\text{for} \quad j\geq1,
\end{equation}
where the normalization assigns the full mass density to $F^{(0)}$. The order $\tau^{-1}$ equation gives
$$
    F^{(0)}(x,v,t)=\rho(x,t)\phi(v),
$$
implying that the leading-order velocity and temperature are $u^{(0)}=0$ and $T^{(0)}=1$. Since the continuity equation then gives $\partial_t\rho=\mathcal{O}(\tau)$, the order-one kinetic equation is
$$
    v\cdot\nabla_x(\rho\phi)=-\rho^\alpha F^{(1)},
$$
whence
\begin{equation}
    F^{(1)}=-\rho^{-\alpha}(v\cdot\nabla_x\rho)\phi.
\end{equation}
Using $\int v\otimes v\,\phi(v)\,dv=I$, the associated mass flux is
$$
    \rho u=\int_{\mathbb{R}^d}vF\,dv
    =-\tau\rho^{-\alpha}\nabla_x\rho+\mathcal{O}(\tau^2).
$$
Substitution into the continuity equation yields the first-order closure
\begin{equation}\label{prorous_media_eq}
    \partial_t\rho
    =\tau\nabla_x\cdot\bigl(\rho^{-\alpha}\nabla_x\rho\bigr)
    =\tau\Delta\Phi_\alpha(\rho),
\end{equation}
where
\begin{equation}\label{eq_Phi_alpha}
    \Phi_\alpha(r):=
    \begin{cases}
       \dfrac{r^{1-\alpha}-1}{1-\alpha},&\alpha\neq1,\\
       \log r,&\alpha=1.
    \end{cases}
\end{equation}

Since the additive constant in \eqref{eq_Phi_alpha} does not contribute to the dynamics, for $\alpha\neq1$, setting
$$
    m=1-\alpha
$$
we can rewrite \eqref{prorous_media_eq} as
\begin{equation}\label{eq_m_diffusion}
    \partial_t\rho=\frac{\tau}{m}\Delta(\rho^m).
\end{equation}
The diffusivity $\tau\rho^{-\alpha}$ is positive for $\rho>0$,
including when $m<0$. At $\alpha=1$, \eqref{eq_Phi_alpha} gives logarithmic
diffusion.

\subsection{Diffusion near a Positive Background}\label{sec_CE_wellposedness}

For initial densities in $C_b^{2+\beta}(\mathbb R^d)$, $0<\beta<1$, with
$0<c_0\leq \rho(\cdot,0)\leq  C_0$, comparison with the constant solutions
preserves these bounds. Equation~\eqref{prorous_media_eq} thus remains
uniformly parabolic, and standard parabolic estimates give a global classical
solution, smooth for positive times \cite{Lieberman1996}. Linearizing at a
constant density $\rho_0>0$ gives, for the perturbation $r$,
$$
    \partial_t r=\tau\rho_0^{-\alpha}\Delta r.
$$
The diffusivity $\tau\rho_0^{-\alpha}=\tau_{\rm eff}$ is the effective
relaxation parameter in \eqref{maineqf} and agrees with the quadratic
low-frequency term of the hydrodynamic eigenvalue in \eqref{Taylorlambda}.\\
For $d\geq3$ and $0<m<(d-2)/d$, suitable finite-mass solutions of
\eqref{eq_m_diffusion} become extinct in finite time
\cite{BonforteVazquez2010}. Such densities cannot have a uniform positive
lower bound on $\mathbb R^d$. Theorem~\ref{mainthm_formal} concerns
perturbations of a positive homogeneous Maxwellian and gives an exact
local density closure at fixed relaxation time through a convergent
expansion in the slow amplitudes.

\section{The Invariance Equation, Spectral Submanifolds and the Slow Spectral Closure}
\label{sec_inv_equation_general}

In this section, we formulate the geometric framework underlying the construction of the nonlinear hydrodynamic manifold.
We begin by recalling the notion of invariant manifolds for infinite-dimensional dynamical systems and their characterization through the invariance equation, emphasizing that this formulation remains meaningful even when a global flow is not available. We then specify slow spectral manifolds and introduce the concept of \emph{spectral closure}. This perspective recasts the classical moment-closure problem as the construction of an invariant manifold together with its associated reduced vector field, providing a dynamically consistent constitutive law. 

\subsection{Invariant Manifolds}
\begin{figure}
    \centering
    \includegraphics[trim={1.5cm 2cm 2cm 2cm},clip, width=0.66\linewidth]{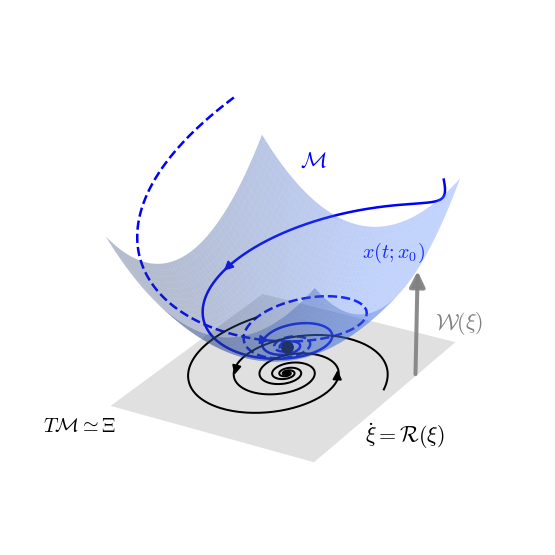}
    \caption[Schematic depiction of an invariant manifold]{Schematic depiction of an invariant manifold $\mathcal{M}$ with parametrization $\mathcal W$ and reduced dynamics $\dot\xi=\mathcal R(\xi)$.}
    \label{plot_slow_manifold}
\end{figure}

Before constructing nonlinear hydrodynamic manifolds for the kinetic equation \eqref{maineqf}, we recall the formulation of invariant-manifold reduction that will be used throughout the paper. Consider the autonomous evolution equation
\begin{equation}\label{eqx}
    \frac{dx}{dt}=\mathfrak{F}(x),\qquad x\in\mathcal{X},
\end{equation}
where $\mathcal{X}$ is a Banach space and $\mathfrak{F}:D(\mathfrak{F})\subseteq\mathcal{X}\to\mathcal{X}$ is a sufficiently regular vector field. Whenever \eqref{eqx} generates a local semiflow, we denote it by $\Phi^t$. For evolution equations governed by unbounded operators, this semiflow may be defined through mild rather than classical solutions and the differential invariance equation introduced below remains meaningful whenever the relevant compositions and derivatives are well-defined.\\
We normalize an equilibrium to the origin, $\mathfrak{F}(0)=0$, and write
\begin{equation}
    \mathfrak{F}(x)=\mathcal{A}x+\mathcal{N}(x),\qquad
    x\in D(\mathcal A),\qquad
    \mathcal{N}(0)=D\mathcal{N}(0)=0.
\end{equation}
Here $\mathcal A:D(\mathcal A)\subseteq\mathcal X\to\mathcal X$ generates a strongly continuous semigroup and $\mathcal N$ is $C^1$ on an open neighborhood of zero in $\mathcal X$. The linearization is $\mathcal A$ as a generator of a $C_0$-semigroup, but an unbounded $\mathcal A$ is not a Fr\'echet derivative $\mathcal X\to\mathcal X$. The differential invariance equation below requires $\mathcal W(\xi)\in D(\mathcal A)$ and the stated derivatives. Otherwise we use the mild flow identity. A set $\mathcal{M}\subseteq\mathcal{X}$ is \emph{positively invariant} if $\Phi^t(\mathcal{M})\subseteq\mathcal{M}$ for all $t\geq0$. Below, ``invariant'' means forward invariant in this sense, since two-sided invariance would additionally require a two-sided flow and equality for all $t\in\mathbb R$.\\
Throughout, an embedded $C^1$ submanifold in a Banach space
$\mathcal X$ is a subset that, near each of its points and in suitable
ambient coordinates, is the graph of a $C^1$ map from an open subset of
a closed subspace $E$ into a closed complement $F$, with $\mathcal X=E\oplus F$.\\
Let $\Xi$ be a Banach space, let $\mathcal{U}\subseteq\Xi$ be open (and contain the origin when a chart through equilibrium is assumed) and let
$$
    \mathcal{W}:\mathcal{U}\to\mathcal{X}
$$
be a \emph{local $C^1$ embedding}. Thus each $\xi_0\in\mathcal U$ has an open
neighborhood $V$ for which $\mathcal W|_V$ is a $C^1$ diffeomorphism onto
a $C^1$ submanifold $\mathcal W(V)$ of $\mathcal X$ with its induced topology.
We write $\mathcal M=\mathcal W(\mathcal U)$ for the full parametrized family.
On each local image $\mathcal W(V)$, the parameters $\xi\in V$ are intrinsic
coordinates, called \emph{reduced variables}, which satisfy the \emph{reduced dynamics},
\begin{equation}\label{reduced_dynamics_general}
    \frac{d\xi}{dt}=\mathcal{R}(\xi),\qquad \mathcal{R}:\mathcal{U}\to\Xi.
\end{equation}
Let $\Phi_{\mathcal{R}}$ denote the flow map associated to \eqref{reduced_dynamics_general}. For the reduced dynamics to be well-defined for all times, the open neighborhood $\mathcal{U}$ has to be forward-invariant under the dynamics of $\xi$,
$$
    \Phi^t_{\mathcal{R}}(\mathcal{U})\subseteq \mathcal{U},\quad t\geq0.
$$
See Figure \ref{plot_slow_manifold} for a schematic depiction of the construction of an invariant manifold in phase space. The parametrization $\mathcal{W}$ and the reduced vector field $\mathcal{R}$ then satisfy the \emph{invariance equation}
\begin{equation}\label{invariance_equation}
    D\mathcal{W}(\xi)\,\mathcal{R}(\xi)
    =\mathfrak{F}(\mathcal{W}(\xi)).
\end{equation}
Indeed, if $\xi(t)$ solves \eqref{reduced_dynamics_general}, the chain rule implies that $x(t)=\mathcal{W}(\xi(t))$ solves \eqref{eqx}. Conversely, whenever $\mathcal{W}$ is locally invertible on its image, the full vector field on $\mathcal W(V)$ determines $\mathcal R$ there. In flow-map form, \eqref{invariance_equation} becomes the semiconjugacy relation, i.e., a conjugacy wherever $\mathcal W$ is invertible on its image, 
\begin{equation}\label{flow_invariance}
    \Phi^t\!\left(\mathcal{W}(\xi)\right)
    =\mathcal{W}\!\left(\Phi^t_{\mathcal{R}}(\xi)\right),
\end{equation}
for all times for which both sides are defined. The parametrization method treats \eqref{invariance_equation} as the primary equation and solves simultaneously for the embedding and the reduced vector field \cite{cabre2003parameterization,cabre2005parameterization}.\\
For a parametrization extending through equilibrium we impose the normalizations
\begin{equation}
    \mathcal{W}(0)=0,\qquad \mathcal{R}(0)=0.
\end{equation}

If $\mathcal W$ is differentiable into $D(\mathcal A)$ with its graph norm at zero and the preceding invariance identity holds, linearization at the origin gives
\begin{equation}
    \mathcal{A}\,D\mathcal{W}(0)
    =D\mathcal{W}(0)\,D\mathcal{R}(0).
\end{equation}
Thus the tangent space $T_0\mathcal{M}=\operatorname{ran}D\mathcal{W}(0)$ is invariant under $\mathcal{A}$, and the leading-order reduced dynamics are conjugate to the restriction of $\mathcal{A}$ to this tangent space. After identifying $\Xi$ with a chosen invariant subspace $E\subseteq\mathcal{X}$, one may use the graph-style normalization
\begin{equation}
    D\mathcal{W}(0)=\iota_E,\qquad
    D\mathcal{R}(0)=\mathcal{A}|_E,
\end{equation}
where $\iota_E:E\hookrightarrow\mathcal{X}$ is the canonical inclusion map.

\begin{remark}
The invariance equation is equivariant under the action of the diffeomorphism group of the parameter space. Indeed, the pair $(\mathcal{W},\mathcal{R})$ is not unique as a parametrized object. If $\Theta$ is a local diffeomorphism of $\Xi$ fixing the origin, then
$$
\widetilde{\mathcal{W}}=\mathcal{W}\circ\Theta,
    \qquad
    \widetilde{\mathcal{R}}(\eta)
    =D\Theta(\eta)^{-1}\mathcal{R}(\Theta(\eta))
$$
parametrizes the same geometric manifold. Normal-form choices for $\mathcal{R}$ therefore act as gauge conditions. Inner resonances determine which nonlinear terms cannot be removed by such reparametrizations \cite{cabre2003parameterization}. 
\end{remark}

\begin{remark}
For the invariant family constructed for the kinetic model in this work, the normalization $\hat{\mathcal W}[0]=0$ and derivative $D\hat{\mathcal W}[0]h=\widehat e_h$ are proved in Proposition~\ref{prop_selected_tangency}. The flow identity holds on the parameter ball, including zero after extension. We stress, however, that graph-norm differentiability and a local $C^1$ manifold chart at zero are not guaranteed by the results in the subsequent sections.
\end{remark}

The invariance equation alone has many solutions: individual trajectories, the full phase space, stable and unstable manifolds, center manifolds and slow manifolds all satisfy it in an appropriate sense. A useful reduction is selected by additional geometric and dynamical requirements, such as prescribed tangency, attraction, smoothness, spectral localization, thermodynamic consistency or computability.
In kinetic theory, the central requirement is that the manifold provide an autonomous constitutive relation between microscopic distributions and a prescribed set of macroscopic observables \cite{gorban2005invariant}.

\subsection{Slow Spectral Submanifolds}

A spectral submanifold is a nonlinear invariant continuation of a selected invariant spectral subspace of the linearization. In finite-dimensional dissipative systems, slow spectral submanifolds provide reduced dynamics associated with the least rapidly decaying modes. Under suitable spectral separation, nonresonance and regularity assumptions, a distinguished manifold can be characterized by maximal smoothness among invariant manifolds with the same tangent space \cite{cabre2003parameterization,HallerPonsioen2016,SzalaiEhrhardtHaller2017}. This characterization is conditional: tangency alone does not imply uniqueness, and neither smoothness nor exponential attraction follows merely from the designation ``slow''. Related invariant-manifold constructions for PDEs and semigroups are discussed in  \cite{wayne1997invariant,buza2024spectral}, as well as in \cite{FoiasSaut1984}, where nonlinear spectral manifolds are used to classify Navier--Stokes data by the leading decay rate of the solution. \\
Let $E\subseteq\mathcal X$ be an invariant spectral subspace of $\mathcal A$, with a spectral projection when one is available. We call a sufficiently regular parametrized invariant manifold $\mathcal M_E$ \emph{associated with} $E$ if
\begin{equation}
    T_0\mathcal M_E=E.
\end{equation}
Writing
$$
    \mathcal W(\xi)=\iota_E\xi+\widetilde{\mathcal W}(\xi),\qquad
    \mathcal R(\xi)=\mathcal A|_E\xi+\widetilde{\mathcal R}(\xi),
$$
where $\widetilde{\mathcal W}(0)=D\widetilde{\mathcal W}(0)=0$ and $\widetilde{\mathcal R}(0)=D\widetilde{\mathcal R}(0)=0$, the invariance equation reads
\begin{equation}\label{nonlinear_invariance_expanded}
    \bigl(\iota_E+D\widetilde{\mathcal W}(\xi)\bigr)
    \bigl(\mathcal A|_E\xi+\widetilde{\mathcal R}(\xi)\bigr)
    =\mathcal A\bigl(\iota_E\xi+\widetilde{\mathcal W}(\xi)\bigr)
    +\mathcal N\bigl(\iota_E\xi+\widetilde{\mathcal W}(\xi)\bigr).
\end{equation}
The slow spectral manifold is associated with those parts of the spectrum with decay rates closest to zero and does not, in general, give an assertion of a uniform transverse decay rate. 
\begin{remark}
When complementary invariant spectral subspaces with ordered real parts exist, one may express a uniform slow-fast separation by
\begin{equation}
 \inf\{\Re z:z\in\sigma(\mathcal A|_{E_{\rm slow}})\}
 >\sup\{\Re z:z\in\sigma(\mathcal A|_{E_{\rm fast}})\}.
\end{equation}
We emphasize that this additional separation is not assumed for the whole-space kinetic construction below since its cutoff on the individual frequencies does not supply such a global splitting.
\end{remark}
To formally solve the invariance equation \eqref{nonlinear_invariance_expanded}, consider first a homogeneous expansion of the form 
$$
\mathcal W(\xi)=\iota_E\xi+\sum_{n\geq2}W_n(\xi,\ldots,\xi),
$$
where $W_n$ are symmetric $n$-linear maps with factorials absorbed into their definition. The degree-$n$ part of \eqref{invariance_equation} gives
\begin{equation}\label{ssm_cohomology_general}
\sum_{j=1}^n W_n(\xi,\ldots,\mathcal{A}|_E\xi,\ldots,\xi)
-\mathcal A W_n(\xi,\ldots,\xi)
=F_n(\xi),
\end{equation}
where $F_n$ depends only on lower-order coefficients and the Taylor coefficients of $\mathcal N$. If nonlinear terms $R_n$ are retained in the reduced dynamics, the right-hand side is instead $F_n-\iota_E R_n$ with the corresponding lower-order contributions included in $F_n$. Equation \eqref{ssm_cohomology_general} defines a hierarchy requiring inversion of shifted linear operators.
\emph{Resonances} occur when sums of selected spectral values meet the spectrum of the full operator, obstructing this inversion \cite{cabre2003parameterization}. In that case the equation may fail to be solvable in the desired function space. We refer to these linear equations for successive homogeneous coefficients as \emph{homological equations}.\\
In the classical nonresonant manifold theory, spectral quotients and smoothness requirements can single out a unique sufficiently regular manifold \cite{cabre2003parameterization}. In the resonant setting, however, such a selection need not be available: the reduced dynamics may have to retain resonant terms, or the parametrization may require weaker regularity or a different functional ansatz. Center-manifold theory likewise illustrates that prescribed tangency does not by itself select a unique invariant manifold \cite{Carr1981}.\\
The kinetic setting considered here differs from finite-dimensional SSM theory in two essential ways. The slow space is a continuously parametrized family of hydrodynamic eigendirections, hence infinite-dimensional, and the hydrodynamic spectral values of the full spatial operator fill the whole interval $[-1/\tau,0]$. The branch has no positive decay rate uniform at low frequencies, and the retained portion has no spectral separation from the omitted portion at the cutoff. Within each selected fiber, however, the distance from $\lambda(k)$ to the essential spectrum is $\tau^{-1}+\lambda(k)\geq \tau^{-1}+\lambda(k_c-\varepsilon_B)>0$, while nonlinear convolution produces summed frequencies beyond this band. The resulting resonances thus cannot be excluded by the linear spectral restriction. The construction we outline in the following therefore selects an invariant hydrodynamic family by a specified normalization of the generalized functional coefficients and solves the resonant hierarchy using a spectral clock, the parametrization of time.

\subsection{Contrast to Inertial Manifolds: Spectral Gap Condition}

We contrast spectral submanifolds with inertial manifolds. An inertial manifold is a global, finite-dimensional, positively invariant and exponentially attracting graph containing the global attractor of a dissipative evolution equation \cite{FoiasSellTemam1988,ConstantinFoiasNicolaenkoTemam1989}.
Inertial manifolds were introduced as a rigorous finite-dimensional reduction theory for PDEs and have been applied most successfully to dissipative parabolic equations such as reaction-diffusion systems, the Kuramoto--Sivashinsky equation, Cahn--Hilliard-type equations, nonlocal Burgers equations, and some modified or hyperviscous fluid models \cite{ConstantinFoiasNicolaenkoTemam1989,Zelik2014,Guo2024,JollyKevrekidisTiti1990}.\\
Classical constructions of inertial manifolds typically require global well-posedness on an absorbing set and a quantitative separation, or squeezing property, between resolved and unresolved modes. For semilinear parabolic equations this is commonly expressed as a \emph{spectral gap} condition between consecutive eigenvalues large enough relative to the Lipschitz constant of the cut-off nonlinearity. For example, for $z_t+A_0z=F(z)$ on a Hilbert space, with $A_0$ positive self-adjoint with compact inverse and $F$ globally Lipschitz with constant $L$ (after a suitable cutoff), a standard sufficient condition is
\begin{equation}\label{spectral_gap}
 \mu_{N+1}-\mu_N>2L,
\end{equation}
where $\mu_j$ are the eigenvalues of $A_0$, counted with multiplicity \cite[Theorem~2.1]{Zelik2014}. This separates the retained and discarded modes. The compact-resolvent hypotheses of this criterion do not hold for the whole-space kinetic operator considered here. Its continuous low-frequency spectrum reaches the origin, and the selected hydrodynamic family is infinite-dimensional, as shown in Section~\ref{sec_linear_part}. The separation used below is the fiberwise one on the fixed subcritical band. It yields the estimates in Section~\ref{sec_proof_main_theorem} even though the full kinetic operator has no gap separating the selected slow modes from their complement.
In particular, the existence of an inertial manifold for the classical Navier--Stokes equations remains a long-standing open problem, motivating approximate inertial manifolds, exponential attractors, spatial averaging methods, and modified models \cite{Zelik2014,Guo2024}. The hydrodynamic manifolds constructed in this work are instead local in amplitude, infinite-dimensional over physical space, spectrally selected at each wave number and constructed without assuming a global attractor or a uniform gap to the essential spectrum.

\subsection{The Spectral Closure}
\label{sec_spectral_closure_general}

\begin{figure}[p]
    \centering
    \rotatebox{90}{%
        \includegraphics[
            width=0.85\textheight,
            height=\textwidth,
            keepaspectratio
        ]{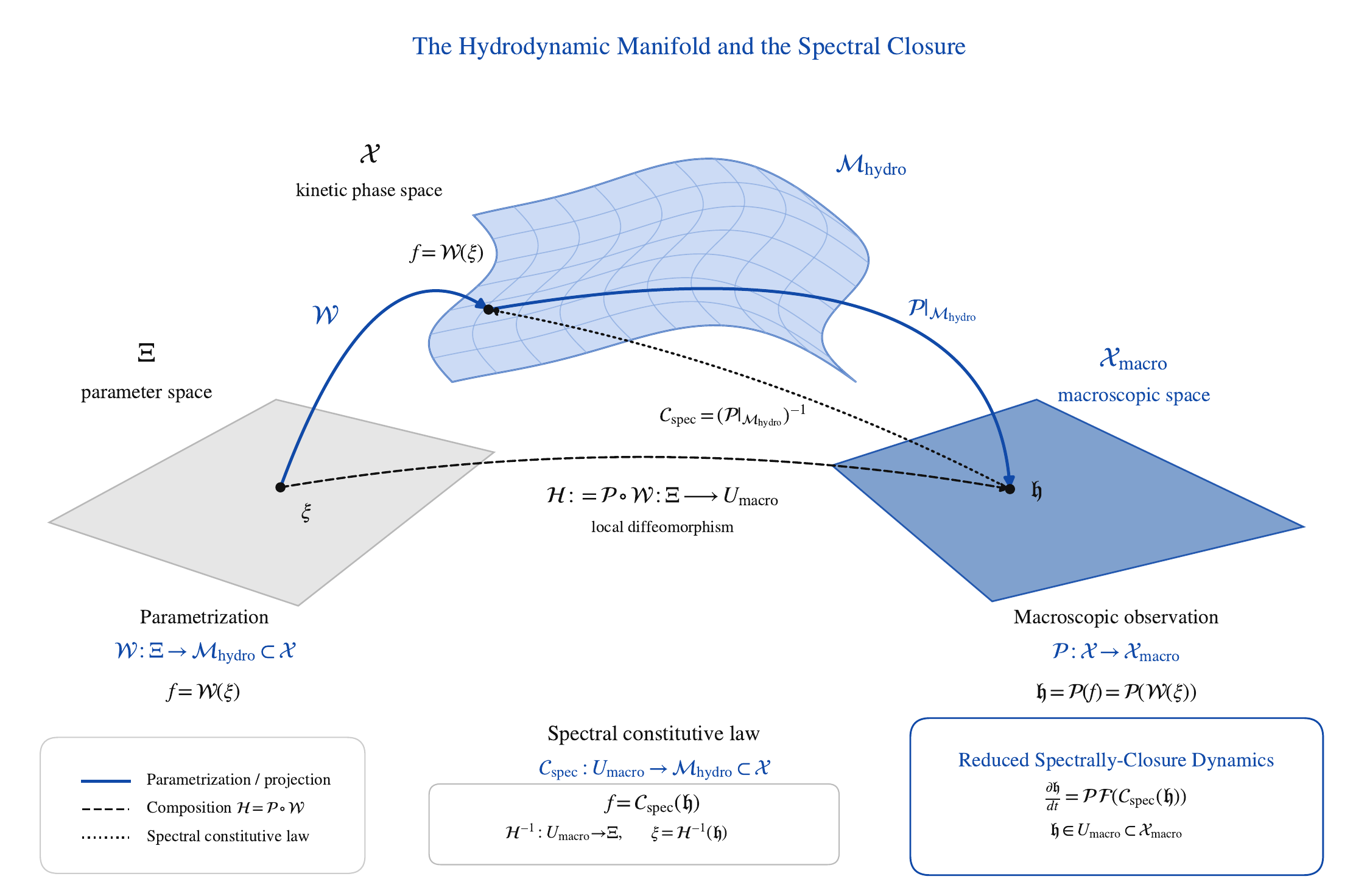}%
    }
    \caption[Local spectral closure]{Schematic relation between the slow spectral parameters $\xi$, the kinetic parametrization $\mathcal W$, and the macroscopic observation $\mathcal P$. The constitutive law $\mathcal C_{\mathrm{spec}}$ and reduced field $\mathcal R_{\mathrm{spec}}$ are defined on each kinetic chart $\mathcal W(V)$ where $\mathcal P|_{\mathcal W(V)}$ is invertible.}
    \label{m_hydro_spectral_closure}
\end{figure}

How is the physical constitutive law, corresponding to a closure, derived from the slow manifold?\\
As discussed in Section \ref{sec_moment_closure}, a moment closure prescribes unresolved kinetic information as a functional of selected macroscopic variables. These macroscopic variables are, abstractly, defined through a $C^1$ \emph{observation map}
$$
\mathcal{P}:\mathcal{X}\to\mathcal{X}_{\rm macro},\quad \mathfrak{h} = \mathcal{P}[f],
$$
where $\mathcal{X}_{\rm macro}$ is the Banach space of physically observable variables or the space of \emph{macroscopic variables}. For the density closure studied here, e.g., $\mathcal{P}f=\rho[f]$ is simply the macroscopic mass density. More generally, the macroscopic variables are the mass density, the velocity field and temperature, see \cite{kogelbauer2024rigorous}. When the realizable macroscopic image is open in an ambient Banach space, its model space must be isomorphic to $\Xi$. More generally the image may be a $C^1$ submanifold modeled on $\Xi$. In the present construction, each local density image has this form because the nonlinear harmonics can extend beyond the compact cutoff of the parameters.\\
A constitutive law is a lifting $\mathcal{C}:\mathcal{X}_{\rm macro}\to\mathcal{X}$ satisfying
\begin{equation}
\mathcal{P}\mathcal{C}=\operatorname{Id}_{\mathcal{X}_{\rm macro}}.
\end{equation}
The closed macroscopic dynamics induced by $\mathcal{C}$ are
\begin{equation}
\frac{d\mathfrak{h}}{dt}
    =D\mathcal{P}(\mathcal{C}(\mathfrak{h}))\,
      \mathfrak{F}(\mathcal{C}(\mathfrak{h})).
\end{equation}
For a generic constitutive law, however, the lifted trajectory $\mathcal{C}(\mathfrak h(t))$ does not solve the kinetic equation: the vector field has a nonzero component normal to the graph of $\mathcal{C}$. This prompts the following definition:

\begin{definition}[Nonlinear Slow Spectral Closure]
Let $\mathcal{M}_{\mathrm{hydro}}=\mathcal W(\mathcal U)$ be a forward-invariant slow spectral family parametrized by a local $C^1$ embedding. Fix a kinetic chart $\mathcal W(V)\subseteq\mathcal M_{\rm hydro}$ on which the macroscopic observation $\mathcal P$ is a diffeomorphism onto its image $\mathcal U_{\rm macro}$, which is either open in $\mathcal{X}_{\rm macro}$ or a $C^1$ submanifold of that space with its induced topology. The \emph{nonlinear spectral closure} is the constitutive law
\begin{equation}
    \mathcal{C}_{\mathrm{spec}}
    :=\bigl(\mathcal P|_{\mathcal W(V)}\bigr)^{-1}:
    \mathcal{U}_{\rm macro}\to\mathcal W(V)\subseteq\mathcal{X},
\end{equation}
and the associated closed hydrodynamic vector field is
\begin{equation}
\mathcal{R}_{\mathrm{spec}}(\mathfrak{h})
    :=D\mathcal{P}(\mathcal{C}_{\mathrm{spec}}(\mathfrak{h}))\,
      \mathfrak{F}(\mathcal{C}_{\mathrm{spec}}(\mathfrak{h})).
\end{equation}
\end{definition}
Thus the spectral closure is the exact constitutive law obtained by restricting the kinetic vector field to the nonlinear hydrodynamic manifold and expressing that restriction in macroscopic coordinates.
In the linear BGK setting, the same principle yields an exact nonlocal hydrodynamic system valid beyond a small-Knudsen expansion and has been termed spectral closure \cite{KogelbauerKarlinSpectral2024,KogelbauerKarlinClosure2025}.\\ 
For a parametrization $f=\mathcal W(\xi)$ with reduced vector field
$\mathcal R$, the construction amounts to locally inverting
$\mathcal P\circ\mathcal W$ on $V$, with both inverses below taken on this chart:
\begin{equation}\label{nonlinear_constitutive_law_density}
 \xi=(\mathcal P\circ\mathcal W)^{-1}(\mathfrak h),\qquad
 \mathcal C_{\mathrm{spec}}(\mathfrak h)
 =\mathcal W\!\left((\mathcal P\circ\mathcal W)^{-1}(\mathfrak h)\right).
\end{equation}
Every velocity moment defined on the kinetic space is then a functional
of $\mathfrak h$. Differentiating $\mathfrak h=\mathcal P(\mathcal W(\xi))$
and using \eqref{invariance_equation} gives
\begin{equation}
 \frac{d\mathfrak h}{dt}
 =D(\mathcal P\circ\mathcal W)(\xi)\,\mathcal R(\xi)
 =\mathcal R_{\mathrm{spec}}(\mathfrak h).
\end{equation}
Thus the constitutive law and the closed dynamics follow from invariance
and local invertibility of the observation map on the manifold.
For mild solutions, the corresponding statement follows from
\eqref{flow_invariance} in these macroscopic coordinates. Figure \ref{m_hydro_spectral_closure} gives a schematic overview of the hydrodynamic manifold, the spectral closure and the relation between the parameter space and the space of macroscopic observables.

\section{Analysis of the Linear Part}
\label{sec_linear_part}
In this section, we analyze the Fourier fiber operator \eqref{def_Lk_section3} in detail, calculate its spectrum, identify its hydrodynamic eigenvalue and derive the estimates used in the nonlinear construction. Throughout this section, $\mathcal X_v=L_v^p(w)$, $1\leq p<\infty$, satisfies \eqref{ass_phi_w}. The spectral calculations are independent of the weight and the operator norms used here need not be rotationally invariant. We distinguish estimates at a fixed frequency from estimates uniform up to the critical wave number, and keep the order of the mixed spatial and velocity norms explicit.

\subsection{Spectral Analysis and the Hydrodynamic Eigenvalue Function}

The spectral theory of linear kinetic operators with finite-rank collision terms has been developed for the one-dimensional linear BGK model in \cite{carty2017elementary,kogelbauer2021non}, for the three-dimensional linear BGK model in \cite{KogelbauerKarlinSpectral2024}, and for the linear Shakhov model in \cite{kogelbauer2024spectral}, relying on the Weinstein--Aronszajn determinant
\cite{aronszajn1948rayleigh,aronszajn1948rayleighII,weinstein1974nonselfadjoint},
see also
\cite{kato1966perturbation}. The spectral analysis presented here is an instance of the general theory of non-normal operators arising as perturbations to multiplication operators as derived in \cite{ljance1970completely}.  The analyses in \cite{carty2017elementary,kogelbauer2021non} use weighted Hilbert spaces, while  a further spectral-decomposition and decay analysis is given in \cite{sukhtayev2024spectral}.\\
We set 
\begin{equation}\label{Lk}
 \mathcal L_k f=-\ri(k\cdot v)f-\tau^{-1}(f-\rho[f]\phi),
 \qquad \rho[f]=\int_{\mathbb R^d}f(v)\,dv,
\end{equation}
with
\begin{equation}\label{domain_Lk}
 D(\mathcal L_k)=\{f\in\mathcal X_v:(k\cdot v)f\in\mathcal X_v\}.
\end{equation}
In particular, $D(\mathcal L_0)=\mathcal X_v$. Note that the domain of $\mathcal{L}_k$ is dense: for any $f\in\mathcal X_v$, the truncations $\mathbf1_{\{|v|\leq N\}}f$ belong to $ D(\mathcal L_k)$ and converge to $f$. This argument does not require local integrability of the weight or density of smooth test functions.\\
The density functional and the Maxwellian projection satisfy
\begin{equation}\label{defCrho}
 |\rho[f]|\leq\|f\|_{L_v^1}\leq C_\rho\|f\|_{\mathcal X_v},\qquad
 C_\rho=
 \begin{cases}
 (\int w^{-1/(p-1)}\,dv)^{(p-1)/p},&1<p<\infty,\\
 \|w^{-1}\|_\infty,&p=1,
 \end{cases}
\end{equation}
by H\"older's inequality. Consequently, the Maxwellian projection associated with the density functional $\rho$,
\begin{equation}\label{def_P_rho}
 P_\rho f:=\rho[f]\phi,\qquad P_\rho^2=P_\rho,
\end{equation}
satisfies the operator norm bound
\begin{equation}\label{Prho_bound}
 \|P_\rho\|_{\rm op}\leq C_\rho\|\phi\|_{\mathcal X_v}.
\end{equation}
Thus $\mathcal L_k$ is a bounded rank-one perturbation of the closed multiplication operator 
$$M_kf=-(\tau^{-1}+\ri k\cdot v)f.$$
For the characterization of the spectrum of $\mathcal{L}_k$ and the following calculations, we define the \emph{critical wave number} as 
\begin{equation}\label{def_k_crit}
    k_c=\frac{1}{\tau}\sqrt{\frac{\pi}{2}}. 
\end{equation}

\begin{figure}
 \centering
 \includegraphics[width=0.45\linewidth]{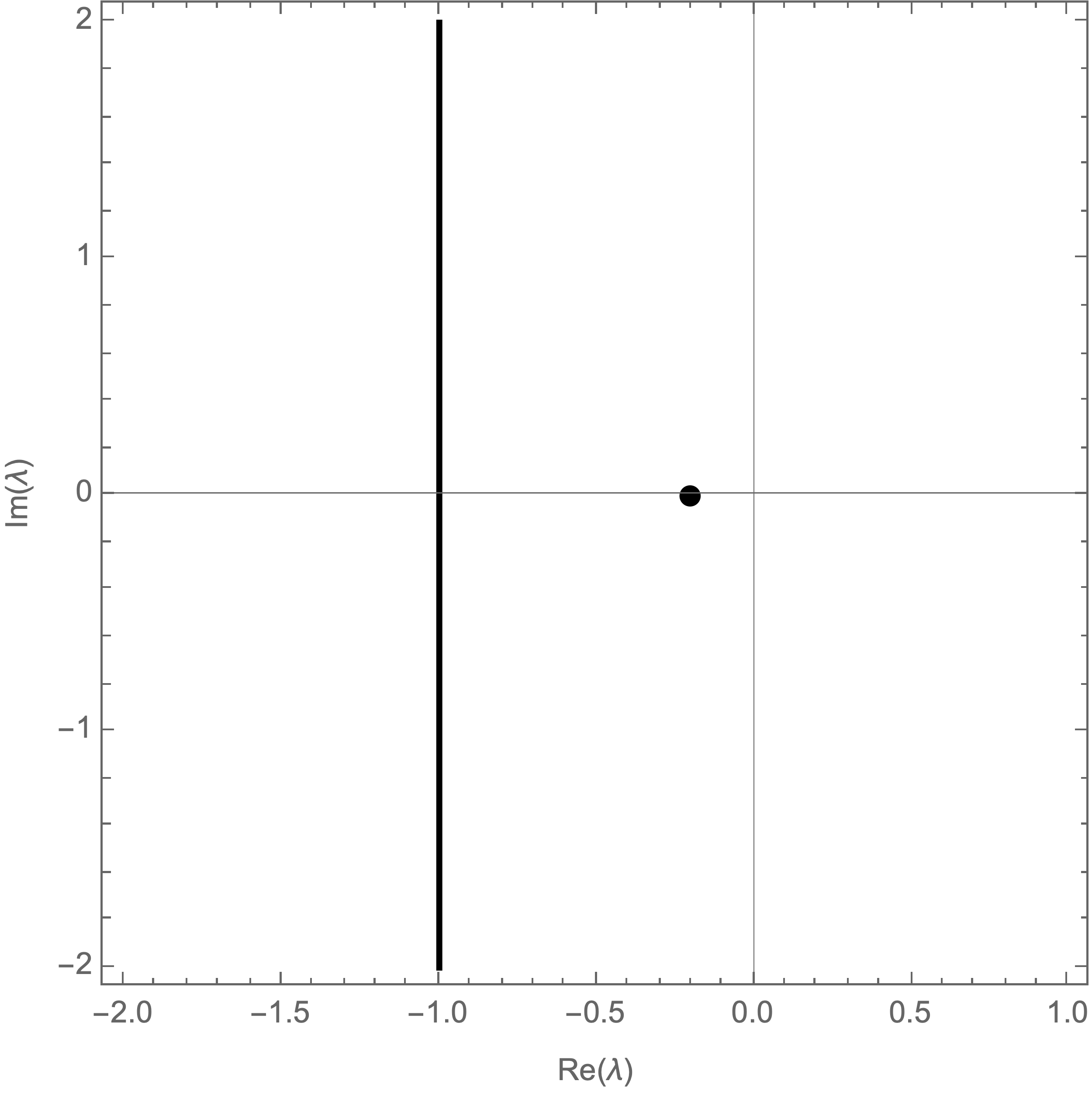}
 \includegraphics[width=0.45\linewidth]{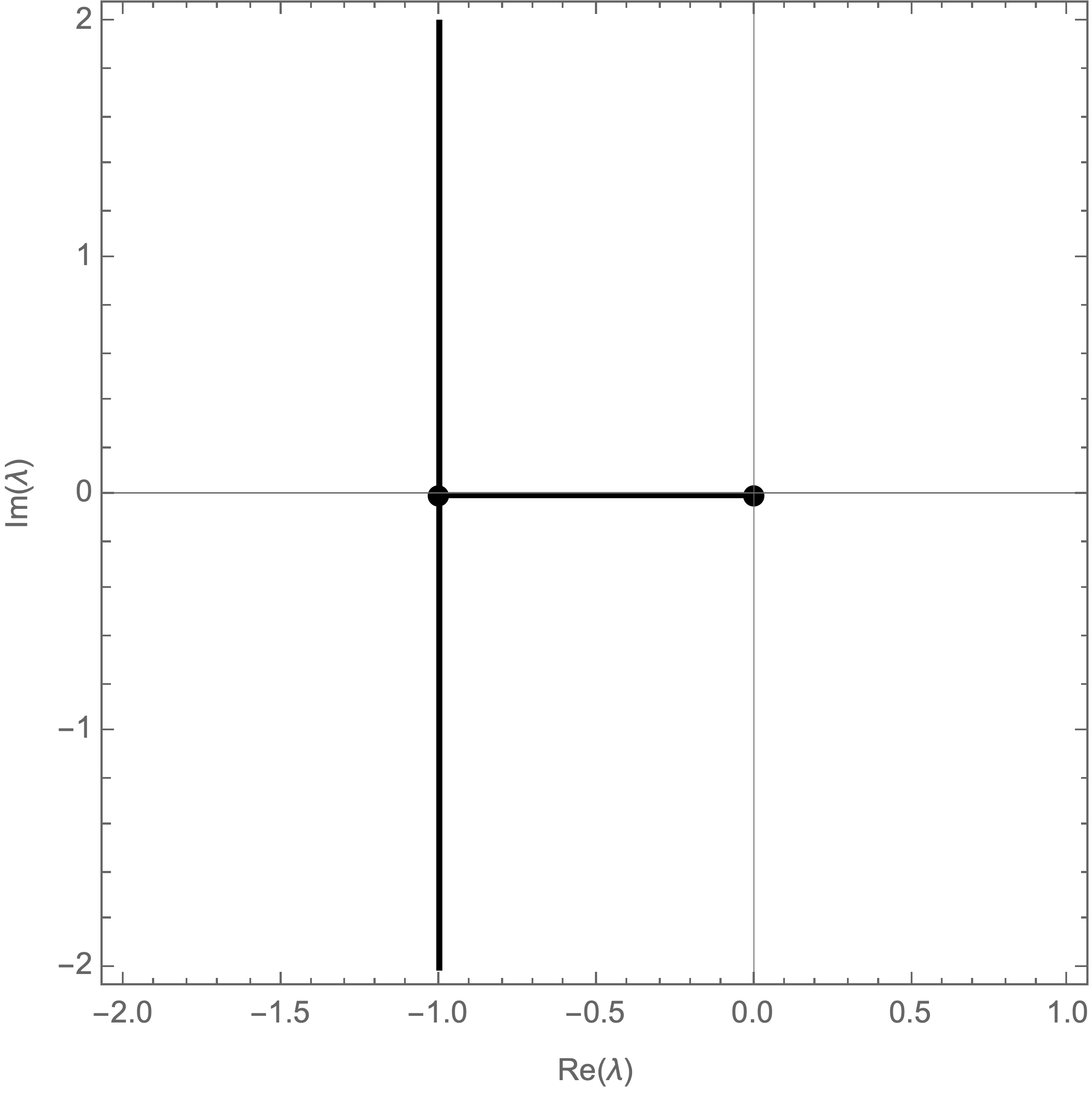}
 \caption{Spectral sets exhibited for a fiber and for the full spatial operator at $\tau=1$. For $0<|k|<k_c$ the hydrodynamic fiber eigenvalue is separated from the essential line. The full spatial operator has arbitrarily slow modes as $k\to0$; a spectral point of an individual fiber need not be a point eigenvalue of the full operator.}
 \label{plot_spectra}
\end{figure}

The exact spectral properties of the operator \eqref{Lk} are summarized in the following proposition, see also Figure \ref{plot_spectra}.

\begin{proposition}\label{prop_spectrum}
Consider the operator $\mathcal L_k:D(\mathcal L_k)\subset\mathcal X_v\to\mathcal X_v$ with domain \eqref{domain_Lk}. For $k=0$, we have that
\begin{equation}
    \sigma(\mathcal{L}_0) = \left\{-\frac{1}{\tau},0\right\},
\end{equation}
where zero defines an isolated eigenvalue of multiplicity one, while $-1/\tau$ defines an eigenvalue of infinite multiplicity, thus constituting the essential spectrum at this wave number.\\
For each wave vector $k\in\mathbb{R}^d$, $k\neq 0$,  the essential spectrum of the operator $\mathcal{L}_k$ is given by 
$$
    \sigma_{\rm ess}(\mathcal{L}_k) = -\frac{1}{\tau } + \ri\mathbb{R}.
$$
For $0<|k|<k_c$, the only spectral point off the essential line is a simple real eigenvalue $\lambda(k)\in(-\tau^{-1},0)$ to its right, determined by
\begin{equation}\label{deflambda}
    Z\left(\ri\frac{\tau  \lambda(k)+1}{\tau  |k|}\right) = \ri \tau  |k|,
\end{equation}
for the plasma dispersion function \eqref{defZ}. The corresponding eigenspace is spanned by the function \begin{equation}\label{eigenvector}
    \hat{e}(k,v) = \frac{\phi(v)}{1+\tau  \lambda(k) +\ri \tau  (k\cdot v)},\quad |k|< k_c,
\end{equation}
which is normalized to 
\begin{equation}
    \rho[\hat{e}](k) = 1,\quad |k|< k_c,\quad \hat{e}(0,v) = \phi(v),
        \end{equation}
and satisfies the bound
\begin{equation}\label{bound_e_hat}
    \|\hat e(k,\cdot)\|_{\mathcal{X}_v}
    \leq  \frac{\|\phi\|_{\mathcal{X}_v}}{1+\tau\lambda(k)},\quad |k|<k_c.
\end{equation}
For $|k|\geq k_c$, we have that $\sigma(\mathcal{L}_k) = \sigma_{\rm ess}(\mathcal{L}_k)$. 
\end{proposition}

\begin{proof}
At $k=0$, $\mathcal L_0=-\tau^{-1}(I-P_\rho)$ gives the assertion.
For $k\ne0$, Kato's compact-perturbation theorem
\cite[Theorem~IV.5.35]{kato1966perturbation} uses the complement of the
semi-Fredholm domain as essential spectrum. To verify the spectral
line directly under our convention, fix $b\in\mathbb R$ and take a unit vector of zero mass supported in $\{|k\cdot v+b|<1/n\}$. Such functions are obtained by subtracting two suitably weighted restrictions of $\phi$ to disjoint bounded subsets of positive measure in the slab. The gain vanishes and the residual for $-\tau^{-1}+\ri b$ is at most $1/n$. Every point of the line is consequently a non-isolated spectral point.\\
For $z=-\tau^{-1}+\ri b$ the operator $\mathcal L_k-zI$ is injective.
Indeed, $(\mathcal L_k-zI)f=0$ implies
$\ri(k\cdot v+b)f=\tau^{-1}\rho[f]\phi$.
If $\rho[f]\ne0$, then $|f|$ has a nonintegrable singularity across the
hyperplane $k\cdot v=-b$, contradicting $\mathcal X_v\hookrightarrow L_v^1$.
Thus $\rho[f]=0$, and the hyperplane has measure zero, so $f=0$.
If the range were closed, the bounded inverse theorem on the graph domain
would give $\|f\|_{\mathcal X_v}\leq C\|(\mathcal L_k-zI)f\|_{\mathcal X_v}$.
The unit vectors above contradict this estimate. Hence the range is not
closed and $\mathcal L_k-zI$ is non-Fredholm.
Off this line,
\eqref{resolvent_formula} and \eqref{defZ} reduce the spectral condition to
$$
 Z\!\left(\ri\frac{1+\tau z}{\tau|k|}\right)=\ri\tau|k|.
$$
For $y\ne0$, \eqref{defZ} shows that $\Im Z(x+\ri y)$ has the sign of $y$,
whereas
$$
 \Re Z(x+\ri y)=\frac1{\sqrt{2\pi}}\int_0^\infty
 \frac{s\bigl(e^{-(s+x)^2/2}-e^{-(s-x)^2/2}\bigr)}{s^2+y^2}\,ds
$$
has the strict sign of $-x$ when $x\ne0$. The spectral condition therefore
forces its argument onto the positive imaginary axis, so $z$ is real
and $z>-\tau^{-1}$. By \eqref{Zexplicit},
$$
 -\ri Z(\ri y)=\int_0^\infty e^{-ys-s^2/2}\,ds,\qquad y>0,
$$
has strictly negative derivative and decreases from $\sqrt{\pi/2}$ to zero. Hence there is exactly
one root for $0<|k|<k_c$, and none for $|k|\geq  k_c$. Since the derivative of $y\mapsto-\ri Z(\ri y)$ at
$y=(1+\tau\lambda(k))/(\tau|k|)$ is strictly negative, the scalar
dispersion relation has a simple zero at $z=\lambda(k)$. The resolvent
formula \eqref{resolvent_formula} then has a simple pole with rank-one
residue, so the eigenvalue is algebraically simple. Rescaling gives
\begin{equation}\label{kernel_characteristic}
 \frac1\tau\int_0^\infty
 e^{-(\lambda(|k|)+\tau^{-1})s-|k|^2s^2/2}\,ds=1.
\end{equation}
At $\lambda=0$ the integral is less than one, so $\lambda(k)<0$.
The respective endpoint limits, from \eqref{asymptoticsZ} and \eqref{Zexplicit},
are $\lambda(0)=0$ and $\lambda(k_c)=-\tau^{-1}$, and
\eqref{kernel_characteristic} remains valid at both endpoints.
Finally, rotating only the Gaussian integral gives
\begin{equation}
 \rho[\hat e]=\frac1{\ri\tau|k|}
 Z\!\left(\ri\frac{1+\tau\lambda(k)}{\tau|k|}\right)=1,
 \qquad 0<|k|<k_c.
\end{equation}
The pointwise bound $|\hat e|\leq\phi/(1+\tau\lambda)$ proves
\eqref{bound_e_hat}, and
$\ri\tau(k\cdot v)\hat e=\phi-(1+\tau\lambda)\hat e$
puts $\hat e$ in $D(\mathcal L_k)$. 
\end{proof}

Note that the above arguments use only the assumptions \eqref{ass_phi_w} on the velocity space. Because the rotation group acts invariantly on the defining relation for $\lambda$, as also used in the above proof, see \cite{kogelbauer2021non} and the subsequent analysis, we call
$$
    \lambda(k)= \lambda(|k|),
$$
the \emph{eigenvalue function}. For later calculations and to make some notation more compact, we also define the \emph{extended eigenvalue function} as 
\begin{equation}\label{lambda_extended}
    \hat{\lambda}(k) = \begin{cases}
        \lambda(k) & \quad \text{if} \quad |k|<k_c,\\
        -\frac{1}{\tau} & \quad \text{if} \quad |k|\geq k_c. 
    \end{cases}
\end{equation}
\begin{remark}
For $|k|<k_c$, we call $\hat e(k,\cdot)$ the hydrodynamic eigenvector
of the fiber operator $\mathcal L_k$. As $k$ varies, $\lambda(k)$ traces
a continuous spectral branch of the full spatial operator, for which
$\hat e(k,v)e^{\ri k\cdot x}$ is a distributional generalized eigenfunction.
\end{remark}

\begin{proposition}[Properties of the Eigenvalue Function]\label{lemma_lambda}
The radial function $\lambda$ extends smoothly and evenly through $\kappa=0$ and continuously to $k_c$. It is strictly decreasing and strictly concave on $(0,k_c)$, with
\begin{equation}\label{Taylorlambda}
 \lambda(\kappa)=-\tau\kappa^2+\tau^3\kappa^4+\mathcal O(\tau^5\kappa^6),
 \qquad \tau\kappa\to0.
\end{equation}
For $0\leq\kappa\leq k_c$,
\begin{equation}\label{lambda_upper_Bound}
 -\sqrt{\frac2\pi}\kappa\leq\lambda(\kappa)
 \leq\frac{\pi-2}{2\tau}-\sqrt{\frac\pi2}\kappa,
\end{equation}
\begin{equation}\label{lambda_quadratic_bounds}
 -\tau\kappa^2\leq\lambda(\kappa)\le-\frac{2\tau}{\pi}\kappa^2.
\end{equation}
\end{proposition}
\begin{proof}
Introduce the Gaussian Mills ratio
\begin{equation}\label{defm}
 m(a)=\int_0^\infty e^{-ay-y^2/2}\,dy,\qquad a\geq 0.
\end{equation}
With $a=(1+\tau\lambda)/(\tau\kappa)$, changing variables in
\eqref{kernel_characteristic} gives
\begin{equation}\label{mkappa}
 \tau\kappa=m(a),\qquad \tau\lambda=am(a)-1=m'(a).
\end{equation}
Integration by parts gives $m'=am-1$, while differentiation under the integral gives
\begin{equation}
 (-1)^jm^{(j)}(a)=\int_0^\infty y^je^{-ay-y^2/2}\,dy>0.
\end{equation}
Since $m$ decreases from $\sqrt{\pi/2}$ to zero, it parametrizes the entire branch. Differentiating \eqref{mkappa} yields
\begin{equation}
 \lambda'(\kappa)=\frac{m''(a)}{m'(a)}<0,
\end{equation}
\begin{equation}
 \lambda''(\kappa)
 =\tau\frac{m'''(a)m'(a)-m''(a)^2}{m'(a)^3}<0.
\end{equation}
The strict inequality follows from strict Cauchy--Schwarz applied to
$y^{1/2}$ and $y^{3/2}$ in the measure $e^{-ay-y^2/2}dy$.\\
For real $\kappa$ and $|\lambda|<(2\tau)^{-1}$, the left-hand side of
\eqref{kernel_characteristic} is $C^\infty$ in $(\lambda,\kappa)$ by
dominated differentiation. Its $\lambda$-derivative at $(0,0)$ is
$-\tau\ne0$, so the smooth implicit-function theorem gives a unique
smooth even branch $\lambda(\kappa)$ near zero. Expanding the integral
to fourth order gives
$$
 1=\frac1{1+\tau\lambda}
   -\frac{\tau^2\kappa^2}{(1+\tau\lambda)^3}
   +\frac{3\tau^4\kappa^4}{(1+\tau\lambda)^5}
   +\mathcal O(\tau^6\kappa^6).
$$
Comparison of coefficients yields \eqref{Taylorlambda}, without
assuming convergence of the Taylor series.\\
The radial function $k\mapsto\lambda(|k|)$ is consequently smooth at the origin. At the other endpoint, $a=0$ gives
$\lambda(k_c)=-\tau^{-1}$ and $\lambda'(k_c)=-\sqrt{\pi/2}$.
Concavity bounds the graph between the chord from $0$ to $k_c$ and the tangent at $k_c$, proving \eqref{lambda_upper_Bound}.

For the quadratic estimates set $J(a)=-m'(a)$. The substitutions $t=x+y$ and $(x,y)=(r\cos\vartheta,r\sin\vartheta)$, respectively, give
$$
 \begin{aligned}
 m(a)^2
 &=\int_0^\infty e^{-at}\int_0^t
    e^{-[x^2+(t-x)^2]/2}\,dx\,dt\geq  J(a),\\
 m(a)^2
 &=\int_0^\infty r e^{-r^2/2}\int_0^{\pi/2}
    e^{-ar(\cos\vartheta+\sin\vartheta)}\,d\vartheta\,dr
 \leq\frac\pi2 J(a).
 \end{aligned}
$$
Here $x^2+(t-x)^2\leq t^2$ and $\cos\vartheta+\sin\vartheta\geq 1$.
Thus
\begin{equation}
 -m(a)^2\leq m'(a)\le-\frac2\pi m(a)^2,
\end{equation}
and \eqref{mkappa} proves \eqref{lambda_quadratic_bounds}. Equivalently, the upper quadratic estimate implies
\begin{equation}
 m(a)\leq\frac{2}{a+\sqrt{a^2+8/\pi}}.
\end{equation}
\end{proof}

\begin{figure}
 \centering
 \includegraphics[width=0.45\linewidth]{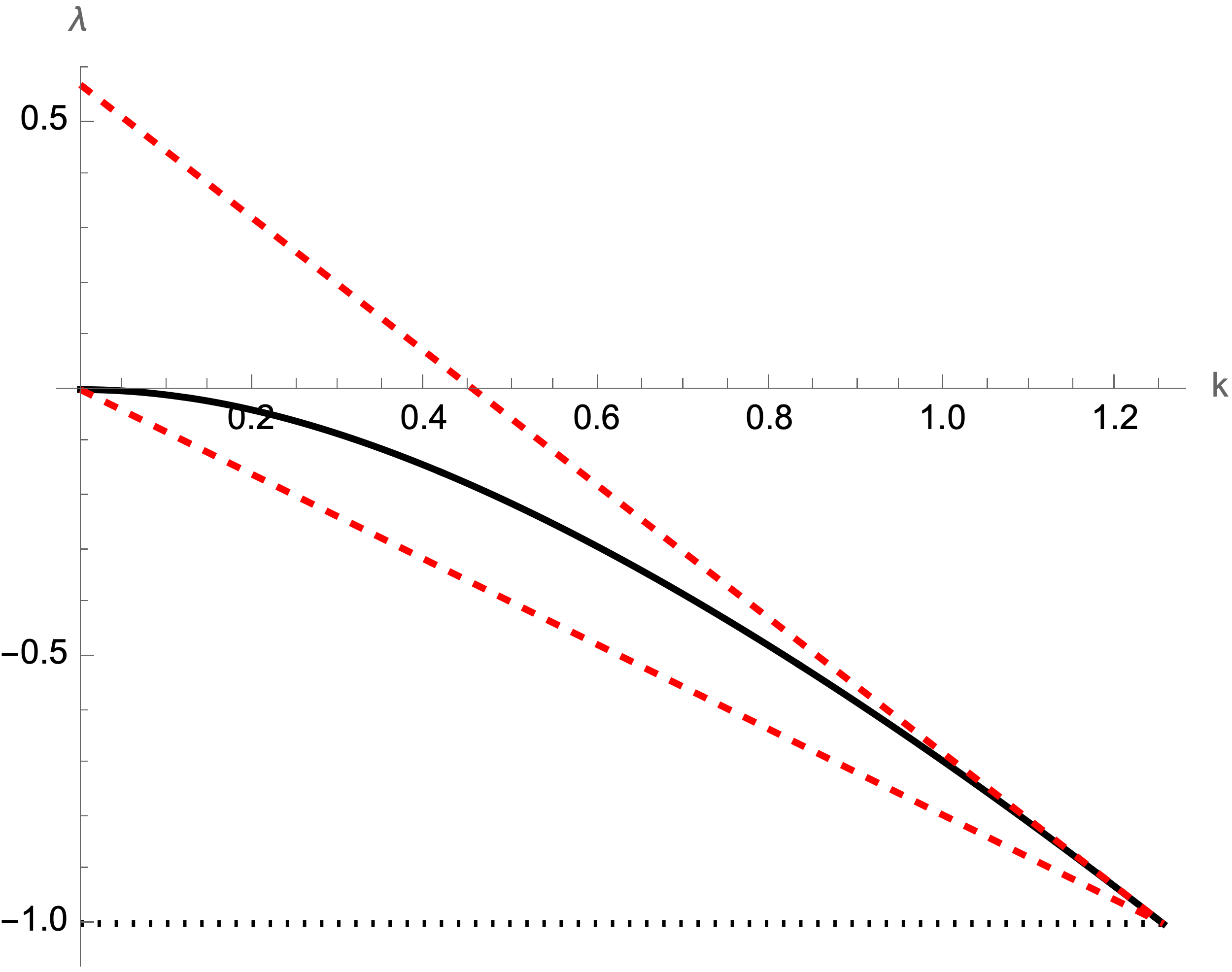}
 \includegraphics[width=0.45\linewidth]{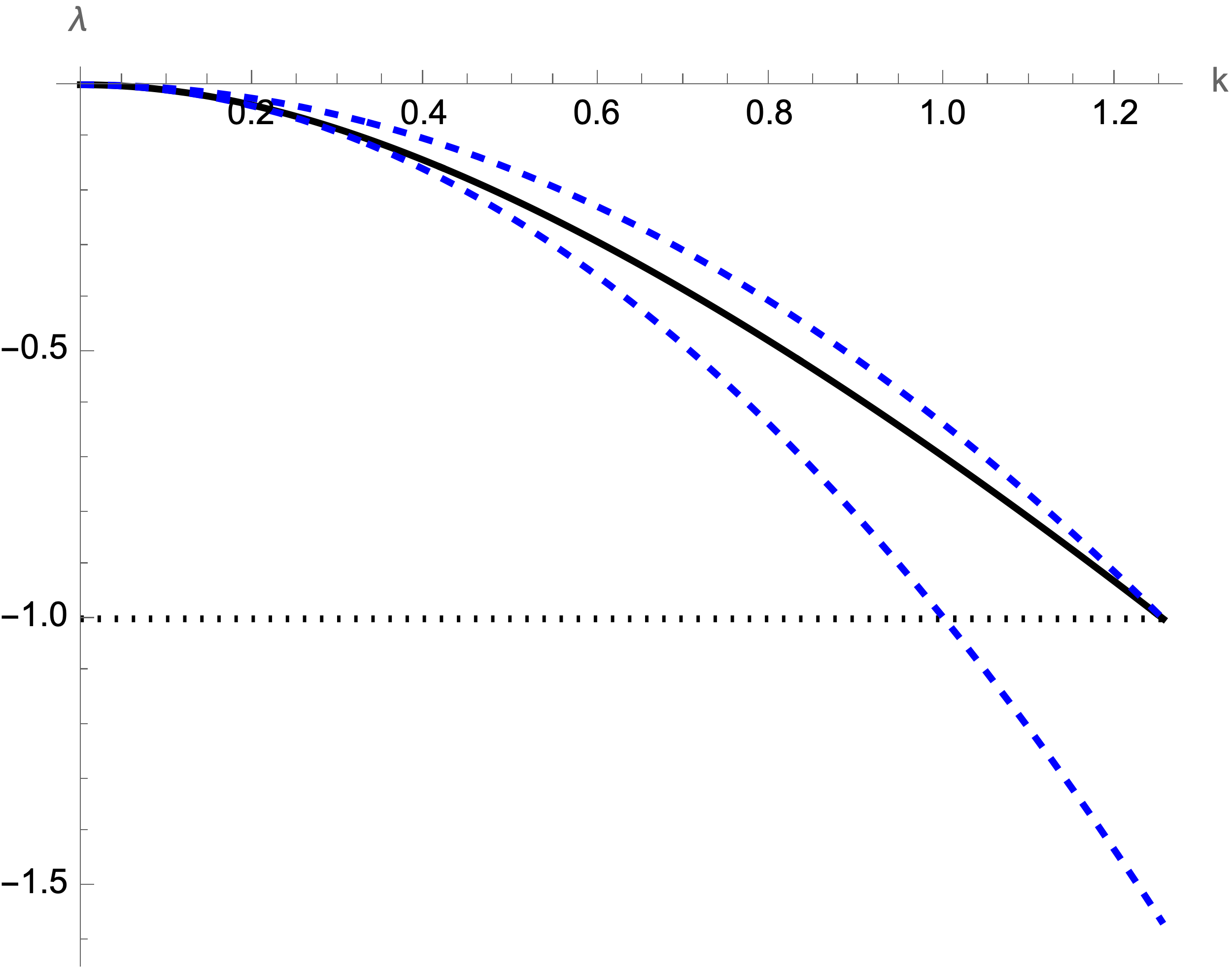}
 \caption{The hydrodynamic eigenvalue for $\tau=1$, with the linear bounds \eqref{lambda_upper_Bound} (left) and quadratic bounds \eqref{lambda_quadratic_bounds} (right). The value at $k_c$ is a continuous endpoint value.}
\end{figure}

\begin{remark}
Note that expansion \eqref{Taylorlambda} is a finite-order Taylor asymptotic, not a convergent power-series representation. Indeed, smoothness at the origin is compatible with the zero radius of convergence derived in \cite{kogelbauer2025relation}. 
\end{remark}

\subsection{Resolvent and Semigroup Bounds}

For $\Re z\ne-\tau^{-1}$ define, including at $k=0$,
\begin{equation}\label{def_Delta_k}
 \Delta_k(z)=1-\int\frac{\phi(v)}{1+\tau z+\ri\tau k\cdot v}\,dv.
\end{equation}
For $k\ne0$, this equals
$1-(\ri\tau|k|)^{-1}Z(\ri(1+\tau z)/(\tau|k|))$, where
the branch of $Z$ is the one specified by \eqref{defZ}.
When $\Delta_k(z)\ne0$, set
\begin{equation}
 \mathcal I_k[g,z]=-
 \frac1{\Delta_k(z)}\int\frac{g(v)}{1+\tau z+\ri\tau k\cdot v}\,dv.
\end{equation}

\begin{lemma}[Resolvent formula and bounds]
For $z\in\rho(\mathcal L_k)$ and $k\ne0$,
\begin{equation}\label{resolvent_formula}
 (\mathcal L_k-z)^{-1}g
 =\frac{\mathcal I_k[g,z]\phi-g}{\tau^{-1}+z+\ri k\cdot v}.
\end{equation}
Whenever $|1+\tau\Re z|>0$ and $\Delta_k(z)\ne0$, one has
\begin{equation}
 \|(\mathcal L_k-z)^{-1}\|_{\rm op}
 \leq\frac{\tau}{|1+\tau\Re z|}
 \left(1+\frac{C_\rho\|\phi\|_{\mathcal X_v}}{|1+\tau\Re z|}\frac1{|\Delta_k(z)|}\right).
\end{equation}
In particular, this bound is uniform on compact sets of pairs $(k,z)$ in the resolvent set satisfying
\begin{equation}\label{resolvent_distance_multiplier}
 \inf |1+\tau\Re z|>0.
\end{equation}
For $k=0$ the resolvent simplifies to 
$$
 (\mathcal L_0-z)^{-1}=-z^{-1}P_\rho-(z+\tau^{-1})^{-1}(I-P_\rho).
$$
\end{lemma}

\begin{proof}
Solving $(\mathcal L_k-z)f=g$ gives
\begin{equation}\label{solf}
 f=\frac{\rho[f]\phi-\tau g}{1+\tau z+\ri\tau k\cdot v},
 \qquad \rho[f]=\tau\mathcal I_k[g,z].
\end{equation}
The estimate follows from \eqref{defCrho} and
\begin{equation}
 |1+\tau z+\ri\tau k\cdot v|^{-1}\le|1+\tau\Re z|^{-1}.
\end{equation}
The constructed $f$ is in $D(\mathcal L_k)$ by rearranging \eqref{solf}. Dominated convergence makes $\Delta_k(z)$ jointly continuous wherever $|1+\tau\Re z|>0$, including $k=0$. Its modulus therefore has a positive minimum on each compact set in the assertion.
\end{proof}

\begin{remark}
Resolvent continuity for a fixed operator does not imply operator-norm continuity in $k$, since the domains of the unbounded multiplication operators vary with $k$. In particular, at $k=0$ points with $\Re z=-\tau^{-1}$ and $\Im z\ne0$ are in the resolvent, whereas the same points lie in the spectrum for every $k\ne0$. We thus use the explicit uniform estimate only under \eqref{resolvent_distance_multiplier} in the following.
\end{remark}

\begin{proposition}\label{semi_group_Lk}
Let $\mathcal{X}_v$ be as in \eqref{defX}, under the assumptions \eqref{ass_phi_w}. For every $k\in\mathbb{R}^d$, the operator
$\mathcal{L}_k$ with domain \eqref{domain_Lk} generates a strongly continuous semigroup $\{e^{t\mathcal{L}_k}\}_{t\geq0}$ on $\mathcal{X}_v$. For every fixed wave number $k$ with $0<|k|<k_c$, there exists a finite constant $C_{\rm sg}(k)>0$ such that
\begin{equation}\label{boundnormsemigroup}
		\|e^{t\mathcal{L}_k}\|_{\rm op}
		\leq 
		C_{\rm sg}(k)e^{\lambda(k)t},
		\qquad t\geq0.
	\end{equation}
The function $C_{\rm sg}$ is locally bounded on
$\{0<|k|<k_c\}$, and every choice of constants in \eqref{boundnormsemigroup} diverges as $|k|\uparrow k_c$.\\
In addition, there exists a constant
	$C_{\mathcal L}=C_{\mathcal L}(\tau,\mathcal{X}_v)>0$,
	independent of $k$, such that
	\begin{equation}\label{boundnormsemigroup_uniform}
		\|e^{t\mathcal{L}_k}\|_{\rm op}
		\leq 
		C_{\mathcal L}
		\left(1+\frac{t}{\tau}\right)^2
		e^{\widehat{\lambda}(k)t},
		\qquad
		k\in\mathbb{R}^d,\quad t\geq0,
	\end{equation}
	where $\widehat{\lambda}$ denotes the extended eigenvalue function
	defined in \eqref{lambda_extended}. In particular, for every
	$k\neq0$ the semigroup is exponentially stable up to a polynomial
	prefactor of degree two.
\end{proposition}

\begin{proof}
We divide the proof into the generation statement, the frequency-wise exponential estimate up to the critical wave number, and the uniform polynomial-exponential estimate.\\

\noindent 
\underline{Generation of the semigroup}:\\
For $k\neq0$, write
$$
\mathcal L_k=M_k+\tau^{-1}P_\rho,
\qquad
M_kg:= -\left(\frac1\tau+\ri k\cdot v\right)g,
$$
with
$$
D(M_k)=D(\mathcal L_k) =
		\left\{
		g\in L_v^p(w):(k\cdot v)g\in L_v^p(w)
		\right\}.
	$$
The operator $M_k$ generates the multiplication semigroup
$$
		(S_k(t)g)(v)
		=
		e^{-t/\tau}e^{-\ri t k\cdot v}g(v),
		\qquad t\geq0,
	$$
which satisfies
\begin{equation}
		\|S_k(t)\|_{\rm op}=e^{-t/\tau}.
\end{equation}
By \eqref{Prho_bound}, the rank-one operator $P_\rho$ is bounded on
$\mathcal X_v$. Hence the bounded perturbation theorem \cite[Theorem~III.1.3]{engel2000one} implies that $\mathcal L_k=M_k+\tau^{-1}P_\rho$ generates a strongly continuous semigroup on $\mathcal X_v$. For $k=0$, the operator is bounded and the semigroup is given explicitly by
\begin{equation}\label{semigroup_k_zero}
		e^{t\mathcal L_0}
		=
		P_\rho+e^{-t/\tau}(I-P_\rho).
	\end{equation}
This proves the generation statement for every $k\in\mathbb R^d$.\\

\noindent 
\underline{Proof of the singular exponential bound
		\eqref{boundnormsemigroup}}:\\
Fix $k$ with $0<|k|<k_c$. By Proposition~\ref{prop_spectrum}, $\lambda(k)>-1/\tau$ is the unique spectral point of $\mathcal L_k$ in the half-plane $\left\{z\in\mathbb C:\Re z>-1/\tau\right\},$ and $\lambda(k)$ is algebraically simple. Let $P_k$ denote the corresponding Riesz projection as defined in \eqref{def_Riesz_proj}, which is bounded for every fixed $0<|k|<k_c$. Since $P_k$ commutes with the semigroup, \begin{equation}\label{semigroup_spectral_decomposition}
		e^{t\mathcal L_k}
		=
		e^{\lambda(k)t}P_k
		+
		e^{t\mathcal L_k}(1-P_k).
	\end{equation}
We next estimate the complementary semigroup. Choose any number $\beta(k)$ satisfying
	\begin{equation}
		-\frac1\tau<\beta(k)<\lambda(k).
	\end{equation}
Since $P_k$ is a spectral projection, $(1-P_k)\mathcal X_v$ is invariant
under $e^{t\mathcal{L}_k}$, and Proposition~\ref{prop_spectrum} implies
that
\begin{equation}
	\sigma\left(
		\mathcal L_k\big|_{(1-P_k)\mathcal X_v}
	\right)
	\subset
	\left\{
		z\in\mathbb C:\Re z\leq -\frac1\tau
	\right\}.
\end{equation}
The free phase $e^{t/\tau}S_k(t)$ is a surjective isometry on the
infinite-dimensional velocity space, so its essential norm is one, while
the rank-one perturbation $\tau^{-1}P_\rho$ preserves the essential
growth bound \cite[Proposition~IV.2.12]{engel2000one}. Thus
\begin{equation}
	\omega_{\rm ess}(\mathcal L_k)
	=
	-\frac1\tau.
\end{equation}
Moreover, since $P_k$ has finite rank and
$$
	e^{t\mathcal L_k}
	-
	e^{t\mathcal L_k}(1-P_k)
	=
	e^{\lambda(k)t}P_k,
$$
the restriction to $(1-P_k)\mathcal X_v$ has the same essential growth
bound,
\begin{equation}
	\omega_{\rm ess}\left(
		\mathcal L_k\big|_{(1-P_k)\mathcal X_v}
	\right)
	=
	-\frac1\tau.
\end{equation}
The characterization of the growth bound in terms of the spectral
bound and the essential growth bound
\cite[Corollary~IV.2.11]{engel2000one} therefore yields
$$
	\omega_0\left(
		\mathcal L_k\big|_{(1-P_k)\mathcal X_v}
	\right)
	=
	-\frac{1}{\tau}.
$$
Hence, for the fixed $k$ under consideration and every
$\beta(k)>-1/\tau$, there exists a finite constant
$C_{\rm rem}(k)$ such that
\begin{equation}
	\|e^{t\mathcal L_k}(1-P_k)\|_{\rm op}
	\leq 
	C_{\rm rem}(k)e^{\beta(k)t},
	\qquad t\geq0.
\end{equation}
Since $\beta(k)<\lambda(k)$, it follows from
\eqref{semigroup_spectral_decomposition} that
$$
	\begin{aligned}
		\|e^{t\mathcal L_k}\|_{\rm op}
		&\leq 
		e^{\lambda(k)t}\|P_k\|_{\rm op}
		+
		C_{\rm rem}(k)e^{\beta(k)t}\\
		&\leq 
		\left(
		\|P_k\|_{\rm op}
		+
		C_{\rm rem}(k)
		\right)e^{\lambda(k)t}.
	\end{aligned}
$$

Thus the optimal constant in \eqref{boundnormsemigroup} satisfies
\begin{equation}\label{definition_gamma_k}
 C_{\rm sg}(k):=\sup_{t\geq 0}e^{-\lambda(k)t}\|e^{t\mathcal L_k}\|_{\rm op}
 \leq\|P_k\|_{\rm op}+C_{\rm rem}(k).
\end{equation}
Its local boundedness is proved after Proposition~\ref{lem_weighted_semigroup_aux}.\\
\noindent
\underline{Proof of the uniform polynomial-exponential bound \eqref{boundnormsemigroup_uniform}}:\\
We now derive an estimate which remains uniform as the hydrodynamic eigenvalue approaches the essential spectrum. To this end, we first derive a representation formula that reduces the estimate of the full semigroup to the mass density evaluated at the semigroup. For $g\in\mathcal{X}_v$, Duhamel's formula gives
\begin{equation}
\label{Duhamel_Tk}
\begin{split}
    e^{t\mathcal L_k}g =
    e^{-\frac{t}{\tau}-\ri tk\cdot v}g+ \frac{1}{\tau}
    \int_0^t
    e^{-\frac{(t-s)}{\tau}-\ri(t-s)k\cdot v}
    \rho[e^{s\mathcal L_k}g] \phi\,ds .
\end{split}
\end{equation}
We define
$$
r_k(t):=e^{t/\tau}\rho[e^{t\mathcal L_k}g], 
$$
multiply \eqref{Duhamel_Tk} by $e^{t/\tau}$ and take the velocity integral to obtain the following linear Volterra equation of the second kind for $r_k$:
\begin{equation}
\label{Volterra_mk}
   r _k(t)
    =
    \hat{g}_k(t)
    +
    \int_0^t \phi_k(t-s)r_k(s)\,ds,\quad t\geq0,
\end{equation}
where 
$$
    \hat{g}_k(t) =
    \int_{\mathbb R^d}
    e^{-itk\cdot v}g(v)\,dv, 
$$
and, using the Fourier transform of the normalized Maxwellian,
\begin{equation}
    \phi_k(t)
    :=
    \frac1\tau
    \int_{\mathbb R^d}
    e^{-itk\cdot v}\phi(v)\,dv
    =
    \frac1\tau e^{-t^2\frac{|k|^2}{2}}.
\end{equation}
By the boundedness of the density functional, we have that
\begin{equation}
\label{bound_hk}
    |\hat{g}_k(t)|
    \leq  C_\rho\|g\|_{\mathcal X_v},
    \qquad t\geq0.
\end{equation}
We will now derive two bounds on the solution $r_k$ to equation \eqref{Volterra_mk}. To ease notation, we set 
$$
    \mu_k = \lambda(k)+\frac{1}{\tau}\geq0 ,\quad 0<|k|<k_c. 
$$
First, assume that $|k|<k_c$. The defining equation for the hydrodynamic eigenvalue \eqref{deflambda} is equivalent
to
\begin{equation}
    \int_0^\infty
    e^{-\mu_kt}\phi_k(t)\,dt
    =
    1,
    \qquad 0\leq  |k| <k_c,
\end{equation}
which follows simply from 
$$
    \begin{aligned}
1 & = \frac{1}{\tau}
\int_{\mathbb R^d}
\frac{\phi(v)}
{\mu_k+\ri k\cdot v}\,dv = \frac{1}{\tau} \int_{\mathbb R^d}
\phi(v) \int_0^\infty e^{-\mu_kt}e^{-\ri tk\cdot v}\,dt\,dv\\
& = \frac{1}{\tau}
\int_0^\infty
e^{-\mu_kt}
\left(
\int_{\mathbb R^d}
\phi(v)e^{-\ri tk\cdot v}\,dv
\right)dt =  \int_0^\infty
    e^{-\mu_kt}\phi_k(t)\,dt,\\
    \end{aligned}
$$
by the definition of $\phi_k$. The function $\Phi_k(t) = e^{-\mu_kt}\phi_k(t)$ is thus a probability density on $[0,\infty)$ for $0<|k|<k_c$. Setting $\tilde{r}_k(t)= e^{-\mu_kt}r_k(t)$ and $\tilde{g}_k(t) = e^{-\mu_kt}\hat{g}_k(t)$, the Volterra equation \eqref{Volterra_mk} can be rewritten for the kernel $\Phi_k$ as 
\begin{equation}
    \tilde{r}_k(t) = \tilde{g}_k(t) + \int_0^t \Phi_k(t-s) \tilde{r}_k(s)\, ds. 
\end{equation}
The classical renewal expansion \cite[Chapter~2]{Gripenberg1990} gives
\begin{equation}\label{sol_rk}
    \tilde{r}_k = \sum_{j=0}^\infty \Phi_k^{*j}\ast \tilde{g}_k,
\end{equation}
where $\Phi_k^{*0}=\delta_0$. For $\sigma>0$, positivity gives
$\int_{[0,t]}\Phi_k^{*j}(ds)\leq e^{\sigma t}\widehat\Phi_k(\sigma)^j$,
where $\widehat\Phi_k$ is the Laplace transform. By \eqref{bound_hk} and
the geometric series,
$$
 |\tilde r_k(t)|\leq
 \frac{e^{\sigma t}C_\rho\|g\|_{\mathcal X_v}}{1-\widehat\Phi_k(\sigma)}.
$$
The probability density $\Phi_k$ is bounded by $1/\tau$, so at least half
its mass lies in $[\tau/2,\infty)$ and we may take $\sigma=(t+\tau)^{-1}$.
Since $0<\sigma\tau/2\leq1/2$, concavity gives
$$
 1-\widehat\Phi_k(\sigma)\geq \tfrac12(1-e^{-\sigma\tau/2})
 \geq \tfrac12(1-e^{-1/2})\sigma\tau, 
$$
while using $e^{\sigma t}\leq e$ and undoing the tilt proves the following
bound with $C_{\rm ren}=2e/(1-e^{-1/2})$.
\begin{equation}\label{bound_rk_subcritical}
		|r_k(t)|
		\leq 
		C_{\rm ren}C_\rho
		\left(1+\frac t\tau\right)
		e^{\mu_kt}
		\|g\|_{\mathcal X_v},
		\qquad
		0<|k|<k_c,\quad t\geq0. 
	\end{equation}
At $|k|=k_c$, $\phi_{k_c}$ is itself a probability density, so the same
estimate holds with $\mu_k=0$. For $|k|\geq k_c$, positivity and
$\phi_k\leq\phi_{k_c}$ give $\phi_k^{*j}\leq\phi_{k_c}^{*j}$ for
every $j\geq1$. The renewal expansion and \eqref{bound_hk} therefore give
\begin{equation}\label{bound_rk_supercritical}
		|r_k(t)|
		\leq 
		C_{\rm ren}C_\rho
		\left(1+\frac t\tau\right)
		\|g\|_{\mathcal X_v},
		\qquad |k|\geq k_c.
	\end{equation}
	Combining \eqref{bound_rk_subcritical} and
	\eqref{bound_rk_supercritical}, and using the extended eigenvalue
	function \eqref{lambda_extended}, yields
	\begin{equation}\label{bound_rk_all_k}
		|r_k(t)|
		\leq 
		C_{\rm ren}C_\rho
		\left(1+\frac t\tau\right)
		e^{(\widehat\lambda(k)+1/\tau)t}
		\|g\|_{\mathcal X_v},
		\qquad
		k\neq0,\quad t\geq0.
	\end{equation}
We now return to the full semigroup. From \eqref{Duhamel_Tk},
	$$
		\begin{aligned}
			\|e^{t\mathcal L_k}g\|_{\mathcal X_v}
			&\leq 
			e^{-t/\tau}\|g\|_{\mathcal X_v}
			+
			\frac{\|\phi\|_{\mathcal X_v}}{\tau}
			e^{-t/\tau}
			\int_0^t|r_k(s)|\,ds.
		\end{aligned}
	$$
	Set
	$$
		\mu_*(k)
		:=
		\widehat\lambda(k)+\frac1\tau
		\geq0.
	$$
	Using \eqref{bound_rk_all_k} gives
	\begin{equation}
		\begin{split}
			\|e^{t\mathcal L_k}g\|_{\mathcal X_v}
			&\leq 
			e^{-t/\tau}\|g\|_{\mathcal X_v} +
			C_{\rm ren}C_\rho\|\phi\|_{\mathcal X_v}\frac{e^{-t/\tau}}{\tau}
			\int_0^t
			\left(1+\frac s\tau\right)
			e^{\mu_*(k)s}\,ds
			\|g\|_{\mathcal X_v},
		\end{split}
	\end{equation}
	and, since $\mu_*(k)\geq0$, we may estimate
	$$
		\begin{aligned}
			\int_0^t
			\left(1+\frac s\tau\right)
			e^{\mu_*(k)s}\,ds
			&\leq 
			e^{\mu_*(k)t}
			\int_0^t
			\left(1+\frac s\tau\right)\,ds =
			e^{\mu_*(k)t}
			\left(
			t+\frac{t^2}{2\tau}
			\right) \leq 
			\tau
			\left(1+\frac t\tau\right)^2
			e^{\mu_*(k)t}.
		\end{aligned}
	$$
	Consequently,
	$$
		\begin{aligned}
			\|e^{t\mathcal L_k}g\|_{\mathcal X_v}
			&\leq 
			e^{-t/\tau}\|g\|_{\mathcal X_v}+
			C_{\rm ren}C_\rho\|\phi\|_{\mathcal X_v}
			\left(1+\frac t\tau\right)^2
			e^{(-1/\tau+\mu_*(k))t}
			\|g\|_{\mathcal X_v}.
		\end{aligned}
	$$
Since $\mu_*\geq0$, the two terms give the following bound with
$C_{\mathcal L}:=1+C_{\rm ren}C_\rho\|\phi\|_{\mathcal X_v}$:
	$$
		\|e^{t\mathcal L_k}\|_{\rm op}
		\leq 
		C_{\mathcal L}
		\left(1+\frac t\tau\right)^2
		e^{\widehat\lambda(k)t},
		\qquad
		k\neq0,\quad t\geq0.
	$$
For $k=0$, the same estimate follows immediately from
	\eqref{semigroup_k_zero} and $\widehat\lambda(0)=0$, with this same
	$C_{\mathcal L}$. This proves
	\eqref{boundnormsemigroup_uniform}.
\end{proof}

\begin{remark}[Critical loss and non-uniformity of the exact rate]
	For every fixed $0<|k|<k_c$, the isolated hydrodynamic eigenvalue
	gives the sharper frequency-wise estimate
	$$
		\|e^{t\mathcal L_k}\|_{\rm op}
		\leq 
		C_{\rm sg}(k)e^{\lambda(k)t}.
	$$
	At the critical wave number the fiber gap closes:
$$
\lambda(k)+\frac{1}{\tau}\to 0 ,\quad |k|\uparrow k_c.
$$
The spectral decomposition
	implies
	$$
		e^{-\lambda(k)t}e^{t\mathcal L_k}
		\longrightarrow
		P_k
	$$
	in operator norm as $t\to\infty$, and hence
	$$
		C_{\rm sg}(k)\geq\|P_k\|_{\rm op}.
	$$
	
	To see the divergence, write $a=1+\tau\lambda(k)\downarrow0$. Integrating
	the modulus of \eqref{eigenvector} on a fixed velocity ball, in coordinates
	with first axis parallel to $k$, gives
	$$
	 \|\hat e(k)\|_{L^1_v}\geq  c\int_0^1
	        (a^2+\tau^2 k_c^2s^2)^{-1/2}\,ds\longrightarrow\infty .
	$$
	By \eqref{defCrho} the $\mathcal X_v$ norm diverges as well. The functions
	$g(v)=(a+\ri\tau k\cdot v)\phi(v)$ are uniformly bounded in $\mathcal X_v$
	by \eqref{ass_phi_w}. In \eqref{Pk_rank_one_corrected} their scalar numerator is one
	and the denominator is at most one by \eqref{Bk_bounds_s6} and hence
	$\|P_k\|_{\rm op}\to\infty$.
	The uniform estimate \eqref{boundnormsemigroup_uniform} replaces this
	singular dependence on the frequency by a polynomial loss in time.
\end{remark}

\begin{remark}[On the polynomial loss]
	The quadratic polynomial factor in
	\eqref{boundnormsemigroup_uniform} is a consequence of estimating the
	mass density first and subsequently estimating the Duhamel integral
	in absolute value. In particular, at the critical wave number the
	renewal estimate gives
	$$
		|r_{k_c}(t)|
		\lesssim
		\left(1+\frac t\tau\right)\|g\|_{\mathcal X_v},
	$$
	and a direct estimate of the Duhamel integral introduces one further
	integration in time. The estimate
	$$
		\int_0^t
		\left(1+\frac s\tau\right)e^{\mu s}\,ds
		\lesssim
		\tau
		\left(1+\frac t\tau\right)e^{\mu t},
		\qquad \mu\geq0,
	$$
	with a constant uniform in $\mu$, fails at $\mu=0$.\\
The present absolute-value argument does not yield a polynomial prefactor of degree one. Retaining the oscillatory factors $e^{-\ri(t-s)k\cdot v}$ may permit a sharper estimate, but no such refinement is proved here.
\end{remark}

\begin{remark}
    The Volterra equation \eqref{Volterra_mk} may alternatively be solved
    by the Laplace transform. Note that the renewal representation \eqref{sol_rk} remains valid at the critical
wave number, where the tilt $\mu_k=\lambda(k)+1/\tau$ vanishes, see also \cite{Cox1962RenewalTheory} for the probabilistic interpretation. We stress that the classical general estimates for Volterra
    equations \cite{Gripenberg1990} do not directly retain the precise
    dependence on the critical parameter required here.
\end{remark}

\begin{remark}
The explicit estimates above are needed to control the frequency dependence of the semigroup bounds. In the Hilbert-space setting, the Gearhart--Pr\"{u}ss theorem does not yield uniform exponential decay over all nonzero $k$, since $\lambda(k)\to0$ as $k\to0$. Its Hilbert-space formulation does not apply to general $L^p$ spaces.
\end{remark}

\begin{proposition}[Auxiliary weighted semigroup bounds]\label{lem_weighted_semigroup_aux}
Let $\mathcal X_v^{(2)}$ be as in \eqref{def_Xv2_aux}. The semigroup $e^{t\mathcal L_k}$ preserves $\mathcal X_v^{(2)}$ and is strongly continuous there. The constant $C_{\mathcal L}^{\rm wgt}:=1+C_{\rm ren}C_\rho\|\phi\|_{\mathcal X_v^{(2)}}$, independent of $k$, satisfies
\begin{equation}\label{weighted_semigroup_general}
 \|e^{t\mathcal L_k}g\|_{\mathcal X_v^{(2)}}
 \leq C_{\mathcal L}^{\rm wgt}\left(1+\frac t\tau\right)^2\|g\|_{\mathcal X_v^{(2)}},
 \qquad k\in\mathbb R^d,\quad t\geq 0.
\end{equation}
Fix $0<\kappa<k_c$. For each $0<\omega<1/\tau$, there exist $C_{\rm far}(\kappa),C_{\rm rem}(\kappa,\omega)>0$ such that, for all $g\in\mathcal X_v^{(2)}$ and $t\geq 0$,
\begin{align}
 \|e^{t\mathcal L_k}g\|_{\mathcal X_v^{(2)}}
 &\leq C_{\rm far}(\kappa)e^{\lambda(\kappa)t/2}\|g\|_{\mathcal X_v^{(2)}},
 &&|k|\geq  \kappa,\label{weighted_semigroup_far}\\
 \|e^{t\mathcal L_k}(I-P_k)g\|_{\mathcal X_v^{(2)}}
 &\leq C_{\rm rem}(\kappa,\omega)e^{-\omega t}\|g\|_{\mathcal X_v^{(2)}},
 &&|k|\leq \kappa.\label{weighted_semigroup_rem}
\end{align}
Here $P_k$ is the hydrodynamic Riesz projection and $P_0=P_\rho$.
\end{proposition}

\begin{proof}
We first apply Proposition~\ref{semi_group_Lk} on the weighted space, and then prove uniformity of the complementary estimate on the closed subcritical ball. All operator norms in this proof are taken on $\mathcal X_v^{(2)}$ to simplify notation.\\

\noindent
\underline{Generation and the full semigroup estimates}:\\
The weight of $\mathcal X_v^{(2)}$ is $(1+|v|)^{2p}w(v)$, which again satisfies \eqref{ass_phi_w}. Indeed, its Maxwellian moment condition follows by increasing the polynomial exponent, while its reciprocal weight is bounded by the original one. In particular,
$$
 \int (1+|v|)^2|g(v)|\,dv
 \leq C_\rho\|g\|_{\mathcal X_v^{(2)}},
 \qquad \phi\in\mathcal X_v^{(2)}.
$$
The multiplication semigroup $S_k(t)$ from Proposition~\ref{semi_group_Lk} has norm $e^{-t/\tau}$ on this space, and $P_\rho$ is bounded there. The bounded perturbation theorem therefore gives a strongly continuous semigroup on $\mathcal X_v^{(2)}$. Its Duhamel formula is again \eqref{Duhamel_Tk} and uniqueness in $\mathcal X_v$ identifies it with the restriction of $e^{t\mathcal L_k}$. Applying \eqref{boundnormsemigroup_uniform} with the new weight gives
$$
 \|e^{t\mathcal L_k}\|_{\rm op}
 \leq C_{\mathcal L}^{\rm wgt}\left(1+\frac t\tau\right)^2
 e^{\widehat\lambda(k)t}.
$$
Since $\widehat\lambda(k)\le0$, this proves \eqref{weighted_semigroup_general}. For $|k|\geq  \kappa$, monotonicity gives $\widehat\lambda(k)\leq\lambda(\kappa)<0$, and hence
$$
 \left(1+\frac t\tau\right)^2e^{\widehat\lambda(k)t}
 \leq\left(\sup_{s\geq 0}\left(1+\frac s\tau\right)^2
 e^{\lambda(\kappa)s/2}\right)e^{\lambda(\kappa)t/2}.
$$
Thus we may fix
$$
 C_{\rm far}(\kappa):=C_{\mathcal L}^{\rm wgt}
 \sup_{s\geq 0}(1+s/\tau)^2e^{\lambda(\kappa)s/2}<\infty.
$$
It is independent of $\omega$. This proves \eqref{weighted_semigroup_far}.\\

\noindent
\underline{The complementary semigroup estimate}:\\
For each fixed $|k|\leq \kappa$, the spectral argument in the proof of Proposition~\ref{semi_group_Lk} applies on $\mathcal X_v^{(2)}$. Indeed, the resolvent has the same scalar dispersion denominator, and its hydrodynamic eigenvector belongs to the weighted space $\mathcal X_v^{(2)}$. Removing the eigenspace therefore leaves the spectrum in $\{\Re z\le-1/\tau\}$. The collision operator has rank one, so the essential growth bound is at most $-1/\tau$. The growth-bound characterization used in that proposition consequently gives
$$
 \lim_{t\to\infty}e^{\omega t}
 \|e^{t\mathcal L_k}(I-P_k)\|_{\rm op}=0,
 \qquad 0<\omega<1/\tau.
$$
At $k=0$, this also follows directly from \eqref{semigroup_k_zero}.\\
To pass from fixed $k$ to a uniform estimate, we use the continuity of the collision contribution in Duhamel's formula. First, \eqref{Pk_rank_one_corrected} and \eqref{Bk_bounds_s6}, proved below directly from the resolvent, show that $P_k$ is norm-continuous and uniformly bounded for $|k|\leq \kappa$. Note that the reciprocal multipliers satisfy
$$
 \left|\frac1{1+\tau\lambda(j)+\ri\tau j\cdot v}
       -\frac1{1+\tau\lambda(k)+\ri\tau k\cdot v}\right|
 \leq\frac{\tau\bigl(|\lambda(j)-\lambda(k)|+|j-k||v|\bigr)}
 {(1+\tau\lambda(\kappa))^2}.
$$
Together with the preceding moment estimate and the Maxwellian moment assumptions, this proves continuity of both factors in the rank-one formula, including at $k=0$. Its scalar denominator is bounded below by $2/\pi$. We use the constant $C_P$ in \eqref{Pk_weighted_bound} with its cutoff taken to be $\kappa$ for the resulting uniform projection bound.\\
Although $k\mapsto S_k(t)$ is not norm-continuous, its composition with the density functional satisfies
$$
 \|\rho S_j(t)-\rho S_k(t)\|_{(\mathcal X_v^{(2)})^*}
 \leq C_\rho t|j-k|e^{-t/\tau}.
$$
Indeed, $|e^{-\ri t j\cdot v}-e^{-\ri t k\cdot v}|\leq t|j-k||v|$. Taking the density in \eqref{Duhamel_Tk} gives the explicit Volterra equation
$$
 \rho[e^{t\mathcal L_k}g]
 =\rho[S_k(t)g]+\frac1\tau\int_0^t
 e^{-(t-s)/\tau-|k|^2(t-s)^2/2}\rho[e^{s\mathcal L_k}g]\,ds.
$$
Set $d(t)=\|\rho e^{t\mathcal L_j}-\rho e^{t\mathcal L_k}\|_{(\mathcal X_v^{(2)})^*}$.
On $[0,T]$ the Gaussian kernel is Lipschitz in frequency on $|k|\leq\kappa$,
and \eqref{weighted_semigroup_general} bounds the density functionals.
Subtracting the two Volterra equations gives
$$
 d(t)\leq C_{T,\kappa}|j-k|+\tau^{-1}\int_0^t d(s)\,ds,
$$
which, after applying the integral version of Gr\"onwall's inequality \cite[Lemma~2.7]{BGKproofTeschl}, yields
$$
 \sup_{0\leq t\leq T}
 \|\rho e^{t\mathcal L_j}-\rho e^{t\mathcal L_k}\|_{(\mathcal X_v^{(2)})^*}
 \leq C_{\rho,\rm Lip}(T,\kappa)|j-k|.
$$
Here $C_{\rho,\rm Lip}(T,\kappa)$ depends on the finite horizon $T$, the cutoff $\kappa$, and the fixed velocity data.
Moreover, $k\mapsto S_k(t)\phi$ is continuous in $\mathcal X_v^{(2)}$, uniformly on bounded time intervals. Returning to Duhamel's formula,
$$
 (e^{t\mathcal L_k}-S_k(t))g
 =\frac1\tau\int_0^t S_k(t-s)\phi\,\rho[e^{s\mathcal L_k}g]\,ds,
$$
shows that $k\mapsto e^{t\mathcal L_k}-S_k(t)$ is norm-continuous for each fixed $t$. The moment estimate used above remains valid for $p=1$, so this argument covers the full range of exponents in \eqref{ass_phi_w}.\\
Fix now $|k|\leq \kappa$. The fixed-frequency decay allows us to choose $t_k>0$ so large that
$$
 \|e^{t_k\mathcal L_k}(I-P_k)\|_{\rm op}
 +2(1+C_P)e^{-t_k/\tau}<\tfrac12e^{-\omega t_k}.
$$
The free parts satisfy
$$
 \|S_j(t_k)(I-P_j)-S_k(t_k)(I-P_k)\|_{\rm op}
 \le2(1+C_P)e^{-t_k/\tau},
$$
while $(e^{t_k\mathcal L_j}-S_j(t_k))(I-P_j)$ is norm-continuous in $j$. Thus, for $j$ in a neighborhood of $k$ within the closed ball,
$$
 \|e^{t_k\mathcal L_j}(I-P_j)\|_{\rm op}<e^{-\omega t_k}.
$$
A finite number of these neighborhoods covers $|j|\leq \kappa$. Denote the corresponding times by $t_1,\ldots,t_N$. For a given $j$, choose one of its neighborhoods and write $t=mt_i+r$, where $m\in\mathbb N_0$ and $0\leq r<t_i$. Since $P_j$ commutes with the semigroup,
$$
 e^{t\mathcal L_j}(I-P_j)
 =e^{r\mathcal L_j}(I-P_j)
 \bigl(e^{t_i\mathcal L_j}(I-P_j)\bigr)^m.
$$
Consequently, \eqref{weighted_semigroup_general} gives
$$
 \|e^{t\mathcal L_j}(I-P_j)\|_{\rm op}
 \leq C_{\mathcal L}^{\rm wgt}\left(1+\frac{t_i}\tau\right)^2
 (1+C_P)e^{\omega t_i}e^{-\omega t}.
$$
Taking the maximum of the finitely many prefactors proves \eqref{weighted_semigroup_rem}.
\end{proof}

\noindent To complete the local bound in \eqref{definition_gamma_k}, fix
$0<\kappa<k_c$ and choose $-\lambda(\kappa)<\omega<1/\tau$.
The slow projection is uniformly bounded on $|k|\leq\kappa$ by
\eqref{Pk_rank_one_corrected}--\eqref{Bk_bounds_s6}; thus
\eqref{weighted_semigroup_rem} and the spectral decomposition give
$\|e^{s\mathcal L_k}\phi\|_{\mathcal X_v}\leq C e^{\lambda(k)s}$
there, since $\phi\in\mathcal X_v^{(2)}$. The alternate order of
Duhamel's formula for the bounded perturbation $\tau^{-1}P_\rho$ gives
$$
 e^{t\mathcal L_k}g=S_k(t)g+
 \tau^{-1}\int_0^t e^{(t-s)\mathcal L_k}\phi\,\rho[S_k(s)g]\,ds .
$$
Now $|\rho[S_k(s)g]|\leq C_\rho e^{-s/\tau}\|g\|_{\mathcal X_v}$.
After factoring out $e^{\lambda(k)t}$, the time integral is bounded
by a constant times $(\tau^{-1}+\lambda(\kappa))^{-1}\|g\|_{\mathcal X_v}$.
This proves local uniformity on the original velocity space, including $p=1$.\par

\subsection{The Linear Hydrodynamic Manifold}

Let $B_c=B_{k_c}(0)$, and disregard its boundary in integrals involving $\hat e$. With the Fourier convention \eqref{def_Fourier_transform}, the linear parametrization is
\begin{equation}\label{defexv}
 e_\xi(x,v)=(2\pi)^{-d/2} \int_{|k|<k_c}\xi(k)\hat e(k,v)e^{\ri k\cdot x}\,dk. 
\end{equation}
Thus $\widehat e_\xi=\xi\hat e$ and $\rho[\widehat e_\xi]=\xi$, extended by zero. The linear flow acts by $\xi(k)\mapsto e^{\lambda(k)t}\xi(k)$.
In $L_k^2(L_v^2(\phi^{-1}))$, normalized lifts supported in shrinking
pairs of reflected balls about $\pm k_*$, $|k_*|<k_c$, are approximate eigenvectors at
$\lambda(k_*)$: their residual norm is bounded by the oscillation of
$\lambda$ on these balls, which tends to zero by evenness and continuity.

\begin{proposition}[Linear parametrization]\label{prop_linear_hydro_mf}
Let $\xi\in C^1(\overline{B_c})$ vanish on $\partial B_c$ and satisfy
$\xi(-k)=\xi(k)^*$. Then
$$
 e_\xi\in L_v^p(w;L_x^2)\cap L_x^\infty(\mathcal X_v)
             \cap L_x^\infty(L_v^\infty),
$$
and
\begin{equation}\label{relative_e_bound}
 |e_\xi(x,v)|\leq C\|\xi\|_{C^1}\phi(v).
\end{equation}
For $1\leq p\le2$ one also has $e_\xi\in L_x^2(\mathcal X_v)$, with norm bounded by $C\|\xi\|_{C^1}$. If, in addition,
\begin{equation}\label{ass_compact_support_xi}
 \xi\in C_c^{2M}(\operatorname{int}B_c),\qquad M\in\mathbb N,\quad 2M>d,
\end{equation}
then, for every $1\leq p<\infty$,
\begin{equation}\label{L1_bound_linear_manifold}
 \|e_\xi\|_{L_x^1(\mathcal X_v)}
 \leq C_{M,\operatorname{supp}\xi}\|\xi\|_{C^{2M}}.
\end{equation}
The physical profile $e_\xi$ is real-valued.
\end{proposition}

\begin{proof}
\noindent
\underline{The $L^2$ and $L^\infty$ estimates}:\\
The lower linear bound in \eqref{lambda_upper_Bound} gives
\begin{equation}
 1+\tau\lambda(k)\geq \frac{k_c-|k|}{k_c}.
\end{equation}
Since $\xi$ vanishes on $\partial B_c$, integration along a radial segment gives, for $k\ne0$,
$$
 \xi(k)=-\int_{|k|}^{k_c}
 \nabla\xi\left(r\frac{k}{|k|}\right)\cdot\frac{k}{|k|}\,dr.
$$
Thus, also at $k=0$ by continuity,
\begin{equation}
 |\xi(k)|\leq\sqrt d\,\|\xi\|_{C^1}(k_c-|k|).
\end{equation}
The factor $\sqrt d$ comes from the maximum-of-partial-derivatives norm fixed in Section~\ref{Notation}. Combining the two estimates shows explicitly how the boundary condition compensates for the eigenvector denominator
$$
|\xi(k)\hat e(k,v)|
 \leq\frac{|\xi(k)|}{1+\tau\lambda(k)}\phi(v)
 \leq\sqrt d\,k_c\|\xi\|_{C^1}\phi(v).
$$
In particular,
\begin{equation}
 \int_{B_c}\frac{|\xi(k)|^2}{(1+\tau\lambda(k))^2}\,dk
 \leq d k_c^2|B_c|\|\xi\|_{C^1}^2.
\end{equation}
Applying scalar Plancherel in $x$ for each fixed $v$, and only then taking the velocity norm, gives
$$
 \|e_\xi\|_{L_v^p(w;L_x^2)}
 =\left\|\left(\int_{B_c}|\xi(k)\hat e(k,\cdot)|^2\,dk\right)^{1/2}\right\|_{\mathcal X_v} \leq\|\phi\|_{\mathcal X_v}
 \left(\int_{B_c}\frac{|\xi(k)|^2}{(1+\tau\lambda(k))^2}\,dk\right)^{1/2}.
$$
For $1\leq p\le2$, Minkowski's inequality gives
$$
 \|e_\xi\|_{L_x^2(\mathcal X_v)}
 \leq\|e_\xi\|_{L_v^p(w;L_x^2)}.
$$
Absolute integration in \eqref{defexv} also yields
$$
 |e_\xi(x,v)|
 \le(2\pi)^{-d/2}\phi(v)
 \int_{B_c}\frac{|\xi(k)|}{1+\tau\lambda(k)}\,dk
 \le(2\pi)^{-d/2}\sqrt d\,k_c|B_c|\|\xi\|_{C^1}\phi(v).
$$
This proves \eqref{relative_e_bound}. Since $\phi\in\mathcal X_v\cap L_v^\infty$, it also proves both $L_x^\infty$ assertions.\\
 
\noindent
\underline{Spatial integrability and reality}:\\
Under \eqref{ass_compact_support_xi}, the denominator $1+\tau\lambda(k)$ has a positive minimum on $\operatorname{supp}\xi$, and all derivatives of $\lambda$ are bounded there, for example,
$$
 \partial_{k_j}\hat e(k,v)
 =-\frac{\tau\bigl(\partial_{k_j}\lambda(k)+\ri v_j\bigr)\phi(v)}
 {(1+\tau\lambda(k)+\ri\tau k\cdot v)^2}.
$$
Repeated differentiation and the product rule therefore give
\begin{equation}\label{derivative_ehat_bound}
 \|\partial_k^\beta[\xi(k)\hat e(k,v)]\|_{L_k^1}
 \leq C_{\beta,\operatorname{supp}\xi}\|\xi\|_{C^{|\beta|}}
 \phi(v)(1+|v|)^{|\beta|},\qquad |\beta|\le2M.
\end{equation}
The zero extension of $\xi\hat e$ is $C^{2M}$ in $k$ and compactly supported in the ball. Integration by parts thus has no boundary terms and gives
$$
 (1+|x|^2)^M e_\xi(x,v)
 =(2\pi)^{-d/2}\int_{B_c}
 (1-\Delta_k)^M[\xi(k)\hat e(k,v)]e^{\ri k\cdot x}\,dk.
$$
Taking the velocity norm and using \eqref{derivative_ehat_bound}, we obtain
$$
 \|e_\xi(x,\cdot)\|_{\mathcal X_v}
 \leq C_{M,\operatorname{supp}\xi}(1+|x|^2)^{-M}
 \|\xi\|_{C^{2M}}\|(1+|v|)^{2M}\phi\|_{\mathcal X_v}.
$$
The last velocity norm is finite by \eqref{ass_phi_w}, and $(1+|x|^2)^{-M}$ is integrable since $2M>d$. This proves \eqref{L1_bound_linear_manifold}.\\
Finally, $\lambda(-k)=\lambda(k)$ implies $\hat e(-k,v)=\hat e(k,v)^*$. Conjugating \eqref{defexv}, changing $k$ to $-k$, and using $\xi(-k)=\xi(k)^*$ gives $e_\xi(x,v)^*=e_\xi(x,v)$.
\end{proof}

\begin{remark}
For the reversed mixed $L^2$ norm, and hence for $L_x^2(\mathcal X_v)$ when $p\le2$, one may replace the $C^1$ boundary condition by $\xi\in H_0^1(\operatorname{int}B_c)$: Hardy's inequality gives
$$
 \int_{B_c}\frac{|\xi(k)|^2}{(k_c-|k|)^2}\,dk
 \leq C\int_{B_c}|\nabla_k\xi(k)|^2\,dk
$$
\cite{maz2013sobolev}. No general mixed-norm Plancherel identity for $p\ne2$ is used. The nonlinear construction instead assumes $C^0$ parameters with a fixed support cutoff strictly inside $B_c$.
\end{remark}

\begin{corollary}[Algebraic decay]
Let $0<\delta<k_c$ and $\xi\in L^\infty(B_c)$ with
$\xi(-k)=\xi(k)^*$ almost everywhere and
$\operatorname{supp}\xi\subset B_{k_c-\delta}(0)$. The linear evolution of \eqref{defexv} is given by 
\begin{equation}
 e_{e^{t\lambda}\xi}(x,v)
 =(2\pi)^{-d/2}\int_{|k|<k_c}
 e^{\lambda(k)t}\xi(k)\hat e(k,v)e^{\ri k\cdot x}\,dk,
\end{equation}
and satisfies
\begin{equation}\label{linear_decay_estimate}
 \|e_{e^{t\lambda}\xi}\|_{L_v^p(w;L_x^2)}
 \leq C_\delta(1+t)^{-d/4}\|\xi\|_\infty,\qquad t\geq 0.
\end{equation}
For $1\leq p\le2$, the same estimate holds in $L_x^2(\mathcal X_v)$. If $|k|\geq \kappa_0>0$ on $\operatorname{supp}\xi$, these estimates hold with exponential decay $e^{-2\tau\kappa_0^2t/\pi}$ in place of $(1+t)^{-d/4}$.
\end{corollary}

\begin{proof}
On the support of $\xi$, the pointwise eigenvector bound gives
$$
 |\hat e(k,v)|\leq\frac{\phi(v)}{1+\tau\lambda(k_c-\delta)},
 \qquad 1+\tau\lambda(k_c-\delta)>0.
$$
Applying scalar Plancherel as in Proposition~\ref{prop_linear_hydro_mf}, followed by \eqref{lambda_quadratic_bounds}, therefore yields
$$ \|e_{e^{t\lambda}\xi}\|_{L_v^p(w;L_x^2)}^2
 \leq\frac{\|\phi\|_{\mathcal X_v}^2\|\xi\|_\infty^2}
 {(1+\tau\lambda(k_c-\delta))^2}
 \int_{|k|<k_c-\delta}e^{-4\tau|k|^2t/\pi}\,dk.
$$
For $0\leq t\le1$, the integral is bounded by $|B_{k_c-\delta}(0)|$. For $t\geq 1$, extending the integral to $\mathbb R^d$ gives
$$
 \int_{|k|<k_c-\delta}e^{-4\tau|k|^2t/\pi}\,dk
 \leq\int_{\mathbb R^d}e^{-4\tau|k|^2t/\pi}\,dk
 =\left(\frac{\pi^2}{4\tau t}\right)^{d/2}.
$$
Combining the two ranges of $t$ and taking square roots proves \eqref{linear_decay_estimate}. Minkowski's inequality gives the assertion in $L_x^2(\mathcal X_v)$ when $p\le2$. If the support is separated from zero, we retain $\operatorname{supp}\xi$ as the integration domain in the first estimate. Then $e^{-4\tau|k|^2t/\pi}\leq e^{-4\tau\kappa_0^2t/\pi}$ throughout that domain, and taking square roots gives the stated exponential rate.
\end{proof}

For $p=2$ and $w=\phi^{-1}$, \eqref{linear_decay_estimate}
recovers the same $L_x^2L_v^2(\phi^{-1})$ decay rate $(1+t)^{-d/4}$
as \eqref{decay_estimate} in Proposition~\ref{prop_global_existence}.
Thus the linear hydrodynamic evolution is consistent with the decay
towards the background Maxwellian established there. In both cases,
the rate reflects the quadratic behavior of $\lambda(k)$ near zero.

\begin{remark}[Positive linear states]
We note that non-negativity of the distribution functions constructed through the slow part of the spectrum can already be guaranteed on the linear level. Assume the hypotheses of Proposition~\ref{prop_linear_hydro_mf}. Then
$$
 \phi(v)+\mathfrak r e_\xi(x,v)
 \geq \phi(v)(1-\mathfrak r C\|\xi\|_{C^1})>0
$$
whenever $\mathfrak r C\|\xi\|_{C^1}<1$, thus showing that, for sufficiently small parametrizations, the linear kinetic subspace defines physical distribution functions. 
\end{remark}

\begin{lemma}[Low-frequency Riesz cancellation]
Let $P_k$ be the Riesz projection for $\lambda(k)$, with $P_0=P_\rho$, and fix $0<k_0<k_c$. The fixed choice
$$
 C_P:=\max\left\{
 \frac{\pi C_\rho\|\phi\|_{\mathcal X_v^{(2)}}}{2(1+\tau\lambda(k_0))^2},
 \frac{\pi\tau C_\rho\|\phi\|_{\mathcal X_v^{(2)}}}{(1+\tau\lambda(k_0))^3}
 \right\}
$$
works in both estimates below:
\begin{equation}\label{Pk_weighted_bound}
 \|P_kf\|_{\mathcal X_v^{(2)}}\leq C_P\|f\|_{\mathcal X_v^{(2)}},\qquad |k|\leq k_0.
\end{equation}
If $\rho[f]=0$, then
\begin{equation}\label{Pk_zero_mass_gain}
 \|P_kf\|_{\mathcal X_v^{(2)}}
 \leq C_P|k|\|f\|_{\mathcal X_v^{(2)}},\qquad |k|\leq k_0.
\end{equation}
The same bounds hold with $\mathcal X_v$ on the left, since $\|g\|_{\mathcal X_v}\leq\|g\|_{\mathcal X_v^{(2)}}$.
As $k\to0$, uniformly for $f\in\mathcal X_v^{(2)}$,
\begin{equation}\label{Pk_low_frequency_expansion}
 \begin{split}
 P_kf={}&\rho[f]\phi
 -\ri\tau\phi(v)\left((k\cdot v)\rho[f]
          +k\cdot\int v f(v)\,dv\right) +\mathcal O_{\mathcal X_v}(|k|^2\|f\|_{\mathcal X_v^{(2)}}).
 \end{split}
\end{equation}
\end{lemma}

\begin{proof}
\noindent
\underline{The Riesz projection and its uniform bound}:\\
Taking the residue of \eqref{resolvent_formula}, with the sign convention \eqref{def_Riesz_proj}, gives
\begin{equation}\label{Pk_rank_one_corrected}
 P_kf(v)=
 \frac{\displaystyle\int\frac{f(v')}{1+\tau\lambda(k)+\ri\tau k\cdot v'}\,dv'}
 {\displaystyle\int\frac{\phi(v')}{(1+\tau\lambda(k)+\ri\tau k\cdot v')^2}\,dv'}
 \,\hat e(k,v).
\end{equation}
Indeed, differentiating \eqref{def_Delta_k} yields
\begin{equation}
 \Delta_k'(\lambda(k))
 =\tau\int\frac{\phi(v)}{(1+\tau\lambda(k)+\ri\tau k\cdot v)^2}\,dv,
\end{equation}
and the factors of $\tau$ cancel in the residue. At $k=0$, the formula reduces to $P_0f=\rho[f]\phi$.\\
We next bound the scalar denominator explicitly. Since $\nabla_v\phi=-v\phi$, Gaussian integration by parts gives, for $k\ne0$,
$$
 \int\frac{(k\cdot v)\phi(v)}{1+\tau\lambda(k)+\ri\tau k\cdot v}\,dv
 =-\ri\tau|k|^2\int
 \frac{\phi(v)}{(1+\tau\lambda(k)+\ri\tau k\cdot v)^2}\,dv.
$$
Multiplying out the denominator in $\int\phi(v)\,dv=1$ and using $\rho[\hat e](k)=1$, we obtain
$$
 1=1+\tau\lambda(k)+\tau^2|k|^2
 \int\frac{\phi(v)}{(1+\tau\lambda(k)+\ri\tau k\cdot v)^2}\,dv.
$$
Consequently, \eqref{lambda_quadratic_bounds} implies
\begin{equation}\label{Bk_bounds_s6}
 \frac2\pi\le
 \int\frac{\phi(v)}{(1+\tau\lambda(k)+\ri\tau k\cdot v)^2}\,dv
 =-\frac{\lambda(k)}{\tau|k|^2}\le1,
 \qquad 0<|k|<k_c.
\end{equation}
At $k=0$ the integral equals one. For $|k|\leq k_0$, we also have
$$
 |1+\tau\lambda(k)+\ri\tau k\cdot v|
 \geq 1+\tau\lambda(k_0)>0.
$$
Using this estimate in both numerator factors of \eqref{Pk_rank_one_corrected} gives
$$
 \|P_kf\|_{\mathcal X_v^{(2)}}
 \leq\frac{\pi C_\rho\|\phi\|_{\mathcal X_v^{(2)}}}
 {2(1+\tau\lambda(k_0))^2}\|f\|_{\mathcal X_v^{(2)}}.
$$
This proves \eqref{Pk_weighted_bound}.\\

\noindent
\underline{The zero-density gain}:\\
If $\rho[f]=0$, we may subtract the constant $(1+\tau\lambda(k))^{-1}$ inside the numerator integral. Thus
$$
 \int\frac{f(v)}{1+\tau\lambda(k)+\ri\tau k\cdot v}\,dv
 =-\frac{\ri\tau}{1+\tau\lambda(k)}
 \int\frac{(k\cdot v)f(v)}{1+\tau\lambda(k)+\ri\tau k\cdot v}\,dv.
$$
In particular,
$$
 \left|\int\frac{f(v)}{1+\tau\lambda(k)+\ri\tau k\cdot v}\,dv\right|
 \leq\frac{\tau C_\rho|k|}{(1+\tau\lambda(k_0))^2}
 \|f\|_{\mathcal X_v^{(2)}},
$$
where we used $\int |v||f(v)|\,dv\leq C_\rho\|f\|_{\mathcal X_v^{(2)}}$. Substitution into \eqref{Pk_rank_one_corrected}, together with \eqref{Bk_bounds_s6}, proves \eqref{Pk_zero_mass_gain} with the stated $C_P$, since $\|P_kf\|_{\mathcal X_v}\leq\|P_kf\|_{\mathcal X_v^{(2)}}$.\\
 
\noindent
\underline{The expansion at zero frequency}:\\
To control the remainder without a power series in the unbounded variable $k\cdot v$, we use the exact identity
$$
 \frac1{1+\tau\lambda(k)+\ri\tau k\cdot v}
 =1-\tau\lambda(k)-\ri\tau k\cdot v
 +\frac{(\tau\lambda(k)+\ri\tau k\cdot v)^2}
 {1+\tau\lambda(k)+\ri\tau k\cdot v}.
$$
The denominator is bounded away from zero on $|k|\leq k_0$, and \eqref{Taylorlambda} gives $\lambda(k)=\mathcal O(|k|^2)$. Hence, for $k$ sufficiently small,
$$
 \left|\frac1{1+\tau\lambda(k)+\ri\tau k\cdot v}
       -1+\ri\tau k\cdot v\right|
 \leq C|k|^2(1+|v|)^2.
$$
Since $\int(1+|v|)^2|f(v)|\,dv\leq C_\rho\|f\|_{\mathcal X_v^{(2)}}$, integration gives
$$
 \int\frac{f(v)}{1+\tau\lambda(k)+\ri\tau k\cdot v}\,dv
 =\rho[f]-\ri\tau k\cdot\int vf(v)\,dv
 +\mathcal O(|k|^2\|f\|_{\mathcal X_v^{(2)}}).
$$
Multiplication of the same pointwise estimate by $\phi$ gives
$$
 \hat e(k,v)=\phi(v)-\ri\tau(k\cdot v)\phi(v)
 +\mathcal O_{\mathcal X_v}(|k|^2),
$$
by the Maxwellian moment assumption. Finally, \eqref{Bk_bounds_s6} and \eqref{Taylorlambda} imply
$$
 \int\frac{\phi(v)}{(1+\tau\lambda(k)+\ri\tau k\cdot v)^2}\,dv
 =1+\mathcal O(|k|^2).
$$
Inserting these three expansions into \eqref{Pk_rank_one_corrected} proves \eqref{Pk_low_frequency_expansion}, with the remainder uniform in $f\in\mathcal X_v^{(2)}$.
\end{proof}

\section{The Hierarchy of Homological Equations Derived from a Generalized Functional Expansion}\label{sec_cohomology_eq}

In this section, we fix the reduced dynamics to be linear and, as a crucial step, introduce the \emph{reparametrization of time} as an extended-valued functional of the coordinate functions of the reduced dynamics, see Section \ref{sec_reparametrize_time}. We then allow the parametrization functional $\mathcal{W}$ to depend on $\theta$, which provides additional flexibility for solving the homological  equation in the presence of resonances. Next, we expand the parametrization of the hydrodynamic manifold into a generalized functional expansion and derive the resulting hierarchy of \emph{unsteady homological  equations}, see Section \ref{sec_func_Taylor_cohomology}. Finally, we examine resonant hypersurfaces in frequency space, study the second-order homological  equation in more detail and distinguish pointwise resonance obstructions from integrability of the reconstructed stationary coefficient, including the kinetic $L^2$ obstruction, see Section~\ref{what_does_not_work} and Proposition~\ref{remark_stationary_quadratic_integrability}. \\
Throughout this section, let $\xi: \mathbb{R}^d\times [0,\infty)\to \mathbb{C}$ denote the reduced variable as discussed in Section \ref{sec_linear_part}. We write 
$$
    \bm{k}_n = (k_1,...,k_n),\quad d\bm{k}_n = \prod_{j=1}^n dk_j, \quad \Sigma\bm{k}_n = \sum_{j=1}^n k_j,
$$
for frequency vectors of wave vectors, the corresponding infinitesimal integration element and the row-sum vector of the matrix $\bm{k}_n$, respectively. Furthermore, we simply write
$$
    |\bm{k}_n| = \sqrt{\sum_{j=1}^n |k_j|^2},
$$
for the norm in  $\mathbb{R}^{d\times n}\cong \mathbb{R}^{dn}$ as vector spaces.\\
For any function $f:\mathbb{R}^d\times[0,\infty)\to\mathbb{C}$ and any set of vectors $\bm{k}_n=(k_1,...,k_n)\in\mathbb{R}^{d\times n}$, we write 
\begin{equation}\label{def_prod_xi}
    f(\bm{k}_n,t) = (f(k_1,t),...,f (k_n,t)),\quad \Pi f(\bm{k}_n,t) = \prod_{j=1}^n f(k_j,t),
\end{equation}
for the vector of the function $f$ evaluated at several frequencies and their product. 
For the subsequent calculations, we denote 
$$
    \Lambda_n(\bm{k}_n) = -\sum_{j=1}^n \lambda(k_j)\geq0,\quad \bm{k}_n\in B_c^n. 
$$
For the nonlinear construction the natural parameter topology is the uniform topology. We therefore choose the real Banach space
\begin{equation}
 \Xi:=\{\xi\in C^0(\overline{B_c};\mathbb C):\xi(-k)=\xi(k)^*\},
 \qquad \|\xi\|_{C^0}:=\|\xi\|_{C^0(\overline{B_c})}.
\end{equation}
No frequency derivatives of $\xi$ are used in the nonlinear homological  hierarchy or in its convergence argument. The stronger $C^r$ assumptions in Proposition~\ref{prop_linear_hydro_mf} serve only to obtain additional physical-space localization of the linear hydrodynamic state and are not part of the nonlinear manifold construction.\\
Fix $\varepsilon_B\in(0,k_c)$ and restrict the individual input frequencies to
$|k_j|\leq k_c-\varepsilon_B$. Their sum $\Sigma\bm k_n=\sum_{j=1}^n k_j$
may reach or cross the critical sphere $|k|=k_c$.
The separation at each input frequency is needed because
$$
 \|\hat e(k,\cdot)\|_{\mathcal X_v}\to\infty,\qquad |k|\uparrow k_c.
$$
We therefore work on the closed Banach subspace
\begin{equation}\label{defXi_eps_C0}
 \Xi(\varepsilon_B):=\bigl\{\xi\in\Xi:
\operatorname{supp}\xi\subseteq\overline{B_{k_c-\varepsilon_B}(0)}\bigr\},
 \qquad \|\xi\|_{\Xi(\varepsilon_B)}=\|\xi\|_{C^0}.
\end{equation}
Thus every admissible leaf vanishes on a fixed collar of the spectral threshold. To ease notation, we denote
$$
 B_c(\varepsilon_B)=B_{k_c-\varepsilon_B}(0).
$$
Since the support has finite measure, $\Xi(\varepsilon_B)\hookrightarrow \mathcal Y$ continuously, or, explicitly, in terms of norms
$$
   \|\xi\|_{L^p}\leq |B_c(\varepsilon_B)|^{1/p}\|\xi\|_{C^0},
$$
for $p=1,2$. 
As outlined before, we proceed in several steps. First, we prescribe linear reduced dynamics and seek a parametrization intertwining them with the kinetic dynamics. Local conjugacy then requires local inversion of the parametrization:
\begin{equation}\label{reduced_dyn_linear}
    \frac{\partial \xi }{\partial t}(k,t) = \lambda(k) \xi(k,t),\quad  \xi(\cdot,0)=\xi_0\in\Xi(\varepsilon_B).
\end{equation}
This is the simplest possible reduced dynamics and all the complexity of the construction of a spectral manifold is transferred to the parametrization functional $\mathcal{W}$. The invariance equation for \eqref{maineqfk} with respect to the reduced dynamics \eqref{reduced_dyn_linear} thus takes on the simple form
\begin{equation}\label{inv_equation}
    \frac{D \hat{\mathcal{W}}}{D \xi}[\lambda\xi](k) = \mathcal{L}_k[\hat{\mathcal{W}}](k) +  \hat{\mathcal{N}}_\alpha[\hat{\mathcal{W}}](k).
\end{equation}
We note that formulation \eqref{inv_equation} of the invariance equation is only well defined for those parameter values where $\mathcal{W}$ is continuously differentiable, see also the discussion in Section \ref{sec_proof_main_thm_proof}. \\
The central difficulty in solving \eqref{inv_equation} is that, although the reduced dynamics \eqref{reduced_dyn_linear} are diagonal in frequency space, the
nonlinearity couples different frequencies through convolution. At higher orders, the resulting homological  equations therefore involve the difference between the spectral
value of the linearized kinetic operator at a summed frequency and sums of hydrodynamic eigenvalues at the constituent frequencies. These spectral differences may vanish, or may bring the essential spectrum to the imaginary axis, so that the
corresponding steady homological  operator ceases to be invertible. Rather than excluding such resonant frequency configurations, we retain the evolution along the reduced flow
as an additional degree of freedom in the coefficient functions of the parametrization. For this purpose, we first reconstruct elapsed time directly from the reduced coordinates
$\xi$. The resulting functional will provide an intrinsic clock on every nontrivial reduced orbit and will subsequently allow us to replace the steady homological  equations
by unsteady ones, in which resonant contributions can be accommodated by controlled time dependence.

\subsection{Parametrization of Time}\label{sec_reparametrize_time}

The linear reduced dynamics \eqref{reduced_dyn_linear} have a special property that will be essential for the treatment of resonances: time along a nontrivial reduced trajectory can be recovered from the instantaneous state $\xi$ itself. Indeed, each
Fourier coordinate evolves exponentially with its own decay rate $\lambda(k)$, and hence
its amplitude carries information about the elapsed time. A single coordinate, however, is not suitable for this purpose globally, since it may vanish for a given parameter function $\xi$. We therefore combine the information contained in all frequencies by means of a
positive weighted average to define an expression of time in terms of parameter functions. We thus seek a scalar functional $\theta[\xi]$ with the cocycle property
\begin{equation}\label{def_spec_time_section}
    \theta[\xi(t)]=\theta[\xi(0)]+t,
\end{equation}
so that differentiation along the reduced flow gives
$d\theta[\xi(t)]/dt=1$. This identity will turn the dependence of the
coefficient kernels on $\theta$ into a time derivative in the homological  equations. The particular choice below is adapted to the exponential decay rates $\lambda(k)$ and,
at the same time, yields the quantitative compensation estimates needed later to control
the growth of resonant coefficient functions as $\theta[\xi]\to\infty$ when $\xi\to 0$. The finite-dimensional counterpart of \eqref{def_spec_time_section} is the logarithm of a weighted sum of powers as discussed in~\cite[Sec.~4.2]{KogelbauerDelaLlave2026Resonant}.

\begin{definition}[Parametrization of Time]\label{def_global_clock}
Let $\nu>0$, $\theta_0\in\mathbb{R}$, and let $\mu\in L^\infty(B_c)$ satisfy
\begin{equation}\label{clock_weight_assumption}
    \mu(-k)=\mu(k)\geq\mu_0>0\quad\text{a.e. on }B_c,\qquad 
    \int_{B_c}\mu(k)\,dk=1. 
\end{equation}
For $\xi\in\Xi$ with $\|\xi\|_{C^0}\leq \mathfrak r$, define
\begin{equation}\label{deftheta}
 S_\mu[\xi]
 :=\int_{B_c}\mu(k)
 \left(\frac{|\xi(k)|}{\mathfrak r}\right)^{-\frac{\nu}{\lambda(k)}}\,dk,
 \qquad
 \theta[\xi]:=\theta_0-\frac1\nu\log S_\mu[\xi],
\end{equation}
whenever $S_\mu[\xi]>0$, and set $\theta[0]=+\infty$ otherwise. The value of the integrand at $k=0$, where $\lambda(0)=0$, is immaterial and may be chosen arbitrarily. We call $\theta$ the \emph{parametrization of time}.
\end{definition}

On physical amplitudes the factor multiplying $\mu(k)$ in
\eqref{deftheta} is even, so replacing a positive normalized
weight by its even part leaves the clock unchanged. With the even choice
in \eqref{clock_weight_assumption}, $\theta[\xi(-\cdot)^*]=\theta[\xi]$
also on unrestricted Fourier amplitudes wherever the clock is finite.\\
Substituting $\xi(t)=e^{\lambda t}\xi(0)$ into \eqref{deftheta} gives
\begin{equation}\label{clock_cocycle}
\begin{split}
S_\mu[\xi(t)]
&=\int_{B_c}\mu(k)
 \left(e^{\lambda(k)t}\frac{|\xi(k,0)|}{\mathfrak r}\right)^{-\frac{\nu}{\lambda(k)}}dk\\
&=e^{-\nu t}S_\mu[\xi(0)],
\end{split}
\end{equation}
so $\theta[\xi(t)]=\theta[\xi(0)]+t$ on each nonzero orbit. We call
\eqref{clock_cocycle} the \emph{cocycle identity}. In the terminology of \cite{NakaoMezic2020}, $S_\mu$ is a Koopman eigenfunctional of the reduced flow with exponent $-\nu$.

\begin{remark}
Fix $\theta_0=0$. For nonzero $\xi\in\Xi(\varepsilon_B)$ the compensated
parameter function $\zeta[\xi]=e^{-\lambda\theta[\xi]}\xi$ stays fixed along the
orbit, while the clock records elapsed time. Substitution in \eqref{deftheta} thus gives
\begin{equation}\label{compensated_clock_normalization_s9}
 \int_{B_c(\varepsilon_B)}\mu(k)
 \left(\frac{|\zeta[\xi](k)|}{\mathfrak r}\right)^{-\nu/\lambda(k)}\,dk=1 .
\end{equation}
Since under \eqref{assumption_nu} the exponent is at least two on $B_c(\varepsilon_B)$, the normalization thus controls the $L^1$ and $L^2$ norms of $\zeta[\xi]$ even when its supremum exceeds $\mathfrak r$.
\end{remark}
The constant $\theta_0$ fixes the origin of the auxiliary time
coordinate and can be regarded as a gauge freedom. Indeed, changing
$\theta_0$ translates $\theta$ by a constant. The coefficient functions
and their initial data in the unsteady homological equations have then
to be translated accordingly, implying that $\theta_0$ is a
reference-time freedom rather than an additional dynamical parameter.
We henceforth fix
$$
    \theta_0=0,
$$
for which the regularity and clock-compensation estimates are collected
in the following lemma.

\begin{lemma}[Estimates and Regularity of the Parametrization of Time]\label{lem_theta}
Let $\xi\in\Xi$ with $\|\xi\|_{C^0}\leq \mathfrak r$. Then
$$
 0\leq  S_\mu[\xi]\leq 1,
$$
and, if $\xi\neq0$,
\begin{equation}\label{global_clock_lower_Bound}
 \theta[\xi]\geq\tau\log\!\left(\frac{\mathfrak r}{\|\xi\|_{C^0}}\right)\geq0.
\end{equation}
The functional $S_\mu$ is continuous on the closed ball
$\{\xi\in\Xi:\|\xi\|_{C^0}\leq\mathfrak r\}$ in its relative topology. If
\begin{equation}
\nu\geq2/\tau,
\end{equation}
then $S_\mu$ is real $C^1$ on the open ball
$\{\xi\in\Xi:\|\xi\|_{C^0}<\mathfrak r\}$, and the real-valued clock
$\theta$ is real $C^1$ on the punctured open ball
$\{\xi\in\Xi:0<\|\xi\|_{C^0}<\mathfrak r\}$.\\
For the compensation estimate, fix $0<\varepsilon_B<k_c$.
On $\Xi(\varepsilon_B)$ it suffices to assume
\begin{equation}\label{assumption_nu}
\nu\geq-2\lambda(k_c-\varepsilon_B).
\end{equation}
Then for $0<\|\xi\|_{C^0}\leq \mathfrak r$,
$0\leq \sigma\leq 1$, and $p=1,2$,
\begin{equation}\label{bound_compensated_xi}
 \|e^{-\sigma \lambda(k)\theta[\xi]}|\xi(k)|\|_{L^p_k} \leq  C_{\rm comp}^{L^p}\,\mathfrak r^\sigma
 \|\xi\|_{C^0}^{1-\sigma}.    
\end{equation}
At $\xi=0$ take the compensated profile to be zero.
\end{lemma}

\begin{proof}
For $k\ne0$, set 
$$
    q(k): = -\frac{\nu}{\lambda(k)},
$$
to ease notation. We have that $0\leq  |\xi(k)|/\mathfrak r\leq  1$ and $ q(k)\geq\nu\tau$. Hence
$$
 0\leq  S_\mu[\xi]\leq 1,
 \qquad
 S_\mu[\xi]
 \leq \left(\frac{\|\xi\|_{C^0}}{\mathfrak r}\right)^{\nu\tau},
$$
where we used that $\int_{B_c}\mu=1$. For $\xi\neq0$, continuity of $\xi$ and $\mu\geq\mu_0>0$ imply that $S_\mu[\xi]>0$, and the last estimate, together with $\theta=-(1/\nu)\log S_\mu$ gives \eqref{global_clock_lower_Bound}.\\
If $\xi_j\to\xi$ uniformly in the closed ball, then the integrands
in \eqref{deftheta} defining $S_\mu[\xi_j]$ converge pointwise for $k\ne0$
and are bounded by $\mu(k)$. Dominated convergence gives continuity on the closed ball. If $\nu\geq2/\tau$, then $-\nu/\lambda(k)\geq2$ and, viewing $C^0(B_c;\mathbb C)$ as a real Banach space,
$$
 DS_\mu[\xi]h
 =\int_{B_c}\mu(k)\,
 \frac{q(k)}{\mathfrak r^{q(k)}}
|\xi(k)|^{q(k)-2}\Re\!\bigl(\xi^*(k)h(k)\bigr)\,dk.
$$
Fix $0<\rho<1$. On $\|\xi\|_{C^0}\leq \rho\mathfrak r$ the integrand is bounded by
$$
 \frac{\mu(k)}{\mathfrak r}
 q(k)\rho^{q(k)-1}\|h\|_{C^0},
$$
and $\sup_{p\geq2}p\rho^{p-1}<\infty$. The analogous bound with $p^2\rho^{p-2}$ controls the variation of the derivative, so differentiation under the integral sign is justified and $S_\mu$ is real $C^1$ on the open ball.
Since $S_\mu[\xi]>0$ for $\xi\ne0$, the chain rule gives the asserted
regularity of $\theta$ on the punctured open ball.\\
For the estimate, assume first $\xi\neq0$.
Since $\mu\geq\mu_0$ by assumption, we can bound 
\begin{equation}\label{clock_unweighted_g_bound}
\frac{1}{S_\mu}\int_{B_{c}(\varepsilon_B)}
 \left(\frac{|\xi(k)|}{\mathfrak r}\right)^{-\frac{\nu}{\lambda(k)}}\,dk
 \leq \frac1{\mu_0}.
\end{equation}
Because $\theta=-(1/\nu)\log S_\mu$, for $p=1,2$, we have that 
$$
 \left(\frac{e^{-\sigma \lambda(k)\theta[\xi]}|\xi(k)|}{\mathfrak r}\right)^p
 =\left(\frac{|\xi(k)|}{\mathfrak{r}}\right)^{p(1-\sigma)}
 \left(
\frac{\left(|\xi(k)|/\mathfrak r\right)^{q(k)}}{S_\mu[\xi]}
\right)^{p\sigma/q(k)}.
$$
On $B_c(\varepsilon_B)$, \eqref{assumption_nu} gives
$0\leq  p\sigma/q(k)\leq 1$, also when equality holds in
\eqref{assumption_nu}. The identity above applies equally when
$\|\xi\|_{C^0}=\mathfrak r$, since $S_\mu[\xi]>0$. Thus the elementary inequality $y^\gamma\leq 1+y$ for $y\geq0$ and $0\leq \gamma\leq 1$, together with \eqref{clock_unweighted_g_bound}, yields
$$
\begin{aligned}
\left\|
e^{-\sigma\lambda(k)\theta[\xi]}|\xi(k)|
\right\|_{L_k^p(B_c(\varepsilon_B))}^p & \leq 
\mathfrak{r}^{p\sigma}\|\xi\|_{C^0}^{p(1-\sigma)}
\left[
|B_c(\varepsilon_B)|
+
\frac{1}{S_\mu[\xi]}
\int_{B_c(\varepsilon_B)}
\left(\frac{|\xi(k)|}{\mathfrak{r}}\right)^{q(k)}
\,dk
\right]
\nonumber\\
&\quad\leq 
\mathfrak{r}^{p\sigma}\|\xi\|_{C^0}^{p(1-\sigma)}
\left(
|B_c(\varepsilon_B)|+\frac{1}{\mu_0}
\right).
\end{aligned}
$$
Taking the $p$th root gives \eqref{bound_compensated_xi} with
$$C_{\rm comp}^{L^p}:=(|B_c(\varepsilon_B)|+\mu_0^{-1})^{1/p},\qquad p=1,2.$$ The case $\xi=0$ then  follows from the zero extension.
\end{proof}

An even weight satisfying \eqref{clock_weight_assumption} is
\begin{equation}
    \mu(k) = |B_c|^{-1}. 
\end{equation}

\subsection{Generalized Functional Expansions and the Homological Equation}\label{sec_func_Taylor_cohomology}
We now introduce the hydrodynamic manifold through a functional expansion over the space of parameter functions. To this end, the Fourier products \eqref{def_prod_xi} play the role of monomials in a Taylor expansion. Under the linear reduced flow \eqref{reduced_dyn_linear}, their time-derivative satisfies 
	\begin{equation}
		\begin{split}
			\frac{d}{dt} \Pi \xi(\bm{k}_n,t)
			&=
			\sum_{j=1}^n
			\left(
			\prod_{\substack{\ell=1\\ \ell\neq j}}^n
			\xi(k_\ell,t)
			\right)
			\frac{\partial \xi}{\partial t}(k_j,t)
			=
			\left(\sum_{j=1}^n \lambda(k_j)\right)
			\Pi \xi(\bm{k}_n,t).
		\end{split}
	\end{equation}
We seek a hydrodynamic parametrization in the form
	\begin{equation}\label{defhydro}
		\hat{\mathcal{W}}[\xi](k,v)
		=
		\sum_{n=1}^\infty \frac{1}{n!}
		\int_{\mathbb{R}^{d\times n}}
		\hat{\mathcal{W}}_n(k,v;\bm{k}_n,\theta[\xi])
		\Pi \xi(\bm{k}_n,t)
		\, d\bm{k}_n,
	\end{equation}
	where, following the notational convention \eqref{convolution_rewrite}, the coefficient functions can be written as
	\begin{equation}
		\hat{\mathcal{W}}_n(k,v;\bm{k}_n,\theta)
		=
		\delta(k-\Sigma\bm{k}_n)
		\hat{\mathfrak{W}}_n(\bm{k}_n,v,\theta)
		\Pi \mathbf 1_{B_c}(\bm{k}_n),
	\end{equation}
	for some functions
	$\hat{\mathfrak{W}}_n:
	\mathbb{R}^{d\times n}\times\mathbb{R}^d\times [0,\infty)
	\to\mathbb{C}$.
	We assume the normalization
	\begin{equation}\label{W1}
		\hat{\mathcal{W}}_1(k,v;k_1)
		=
		\mathbf 1_{B_c}(k)
		\hat{e}(k,v)
		\delta(k-k_1),
	\end{equation}
which fixes the linear term. To compute a higher coefficient, we vary an
auxiliary scalar $\varepsilon$ with the clock held fixed and define
\begin{equation}\label{defhydro_epsilon}
		\hat{\mathcal{W}}_\varepsilon[\xi;\theta](k,v)
		:=
		\sum_{n=1}^\infty \frac{\varepsilon^n}{n!}
		\int_{\mathbb{R}^{d\times n}}
		\hat{\mathcal{W}}_n(k,v;\bm{k}_n,\theta)
		\Pi\xi(\bm{k}_n,t)
		\, d\bm{k}_n.
	\end{equation}

\begin{remark}
At fixed $\theta$, \eqref{defhydro_epsilon} has the algebraic form of a
Volterra expansion \cite{Volterra1887,LeeSchetzen1965,boyd1984analytical},
while we stress again that the actual map \eqref{defhydro} is not an
ordinary analytic Volterra series or a Fr\'echet Taylor series at
equilibrium, since its kernels depend on $\xi$ through $\theta[\xi]$.
This distinction is also reflected in the regularity asserted by the
main theorem.
\end{remark}

\begin{remark}
Convergence of an expansion of the form \eqref{defhydro}, as an infinite
sum of integral terms, is slightly more subtle than convergence of a
finite-dimensional Taylor series. If the coefficient functionals were
uniformly bounded in $\theta$, a sufficient condition for a Banach target
space $\mathcal{X}$ would be the usual analytic estimate
\begin{equation}
    \sum_{n=1}^{\infty}
    \frac{1}{n!}
    \|\mathcal{W}_n(\theta)\|_{\mathcal{L}^n(\Xi(\varepsilon_B),\mathcal{X})}
    \|\xi\|_{C^0}^n
    <\infty,
\end{equation}
uniformly in the relevant range of $\theta$. In the present construction,
however, the coefficient functions may grow with $\theta$ close to
resonant frequency vectors. Consequently, convergence cannot be deduced
from a uniform bound on the kernels alone. Instead, one has to estimate
the products
$\hat{\mathfrak{W}}_n(\bm{k}_n,v,\theta[\xi])\Pi\xi(\bm{k}_n)$
directly and show that the growth generated by the unsteady homological
equations is compensated by the dependence of $\theta[\xi]$ on the
amplitudes of the reduced coordinates. At this stage, \eqref{defhydro}
is therefore understood as a formal order-by-order expansion. Absolute
convergence in a sufficiently small neighborhood of the origin in the
parameter space will be established in
Section~\ref{sec_proof_main_theorem}.
\end{remark}
    
Because the slow-fast splitting is defined frequency by frequency, the expansion is naturally formulated in Fourier space. For $p=2$, Hilbert-space Plancherel transfers $L_k^2(\mathcal X_v)$ estimates to $L_x^2(\mathcal X_v)$. The conclusion for $1\leq p<2$ will instead use the auxiliary Hilbert space in Proposition~\ref{prop_majorant_functional_taylor}. Moreover, the compact frequency support of the coordinate functions ensures that the convolution integrals at every finite order are taken over bounded subsets of frequency space. The stronger estimates required for convergence of the full infinite expansion are established in Section~\ref{sec_proof_main_theorem}. For the formal coefficient calculation it is convenient to invert the spatial Fourier transform. Defining the spatial coordinate functions
	$$
		\eta(x,k,t)
		=
		\xi(k,t)e^{\ri k\cdot x},
	$$
	we obtain, formally,
\begin{equation}\label{W_real_parametrization}
		\begin{split}
			\mathcal{W}[\eta](x,v)
			&=
			\sum_{n=1}^\infty
			\frac{\mathcal{W}_n[\eta;\theta[\eta]](x,v)}{n!},
			\\
			\mathcal{W}_n[\eta;\theta](x,v)
			&=
			(2\pi)^{-d/2} \int_{\mathbb{R}^{d\times n}}
			\hat{\mathfrak{W}}_n(\bm{k}_n,v,\theta)
			\Pi\eta(x,\bm{k}_n,t)
			\, d\bm{k}_n,
		\end{split}
	\end{equation}
In analogy with \eqref{defhydro_epsilon}, for fixed $\theta$ we introduce the bookkeeping expansion
	$$
		\mathcal{W}_\varepsilon[\eta;\theta]
		:=
		\sum_{n=1}^\infty
		\frac{\varepsilon^n}{n!}
		\mathcal{W}_n[\eta;\theta],
	$$
	so that
\begin{equation}\label{Wn_bookkeeping_derivative}
		\left.
		\frac{d^n}{d\varepsilon^n}
		\mathcal{W}_\varepsilon[\eta;\theta]
		\right|_{\varepsilon=0}
		=
		\mathcal{W}_n[\eta;\theta].
	\end{equation}
We stress that the derivative in \eqref{Wn_bookkeeping_derivative} is taken at fixed $\theta$ and is not the $n$-th derivative of the final map $\eta\mapsto\mathcal{W}[\eta]$ after the substitution $\theta=\theta[\eta]$. Moreover, $\theta[\xi]=\theta[\eta]$, since the dependence of $\theta$ on the coordinate functions is only through their modulus and
	$$
		|\eta(x,k,t)|=|\xi(k,t)|.
	$$
	The reduced dynamics translate to
	$$
		\frac{d\eta}{dt}(x,k,t)
		=
		\lambda(k)\eta(x,k,t).
	$$
Finally, the representation \eqref{W_real_parametrization} follows directly from
	$$
		\begin{aligned}
			 \int_{\mathbb{R}^d}
			\hat{\mathcal{W}}_n(k,v;\bm{k}_n,\theta)
			\Pi\xi(\bm{k}_n,t)
			e^{\ri k\cdot x}\, dk & =
			 \int_{\mathbb{R}^d}
			\delta(k-\Sigma\bm{k}_n)
			\hat{\mathfrak{W}}_n(\bm{k}_n,v,\theta)
			\Pi\xi(\bm{k}_n,t)
			e^{\ri k\cdot x}\, dk
			\\
			&=
			 \hat{\mathfrak{W}}_n(\bm{k}_n,v,\theta)
			\Pi\eta(x,\bm{k}_n,t).
		\end{aligned}
	$$
Assume now that the parametrization can be expanded in a generalized functional expansion of the form \eqref{defhydro} and \eqref{W_real_parametrization}. We can perform the formal calculations outlined in Section \ref{sec_inv_equation_general} more explicitly. First, we calculate the time-derivative as
$$
 \begin{aligned}
 \frac{d}{dt}\mathcal W[\eta(t)] 
 &\quad=\frac{d}{dt}(2\pi)^{-d/2}
  \sum_{n=1}^\infty\frac1{n!}\int_{\mathbb R^{d\times n}}
  \hat{\mathfrak W}_n(\bm k_n,v,\theta[\eta(t)])\,\Pi\eta(\bm k_n,t)\,d\bm k_n\\
 &\quad=(2\pi)^{-d/2}\sum_{n=1}^\infty\frac1{n!}\int_{\mathbb R^{d\times n}}
  \left[\frac{\partial\hat{\mathfrak W}_n}{\partial\theta}(\theta)
  \frac{d\theta[\eta]}{dt}\,\Pi\eta(\bm k_n,t)
  \mathbin{+}\hat{\mathfrak W}_n(\theta)
  \frac{d}{dt}\Pi\eta(\bm k_n,t)\right]d\bm k_n\\
 &\quad=(2\pi)^{-d/2}\sum_{n=1}^\infty\frac1{n!}\int_{\mathbb R^{d\times n}}
  \left[\frac{\partial\hat{\mathfrak W}_n}{\partial\theta}(\theta)
  +\left(\sum_{j=1}^n\lambda(k_j)\right)\hat{\mathfrak W}_n(\theta)\right]
  \Pi\eta(\bm k_n,t)\,d\bm k_n .
 \end{aligned}
$$
where the clock identity \eqref{clock_cocycle} gives
$d\theta[\eta(t)]/dt=1$. This is the term absent from the steady
parametrization method and lies at the heart of the parametrization
of time approach to the homological equations.
Since the collision nonlinearity is pointwise in $x$, products of
spatial Fourier modes produce convolution in frequency. At fixed
$\theta$, the Leibniz and Fa\`a-di-Bruno formulas give
\begin{equation}\label{Faa-di-Bruno}
  \frac{d^n}{d\varepsilon^n} f(\varepsilon)^\alpha =  \sum_{l=1}^n (\alpha)_l f(\varepsilon)^{\alpha-l} B_{n,l}(\partial_\varepsilon f,...,\partial_\varepsilon^{n-l+1}f),
\end{equation} 
where $(\alpha)_l$ are the falling factorials \eqref{falling_factorials} and $B_{n,l}$ are the partial Bell polynomials, 
\begin{equation}
    B_{n,l}(x_1,x_2,\dots,x_{n-l+1})  =\sum \frac{n!}{j_1!j_2!\cdots j_{n-l+1}!}
    \prod_{m\geq 1}\Big(\frac{x_m}{m!}\Big)^{j_m},
\end{equation}
with the sum over nonnegative integers satisfying $\sum_m j_m=l$ and $\sum_m m j_m=n$. The combinatorial coefficients are nonnegative, although the factors $(\alpha)_l$ need not be.
Using the spatial representation  \eqref{W_real_parametrization} together with the standard Leibniz rule and the Fa\`{a}-di-Bruno formula \eqref{Faa-di-Bruno}, we can write for the bookkeeping parameter $\varepsilon$: 
\begin{equation}
\begin{aligned}
 &\left.\frac{d^n}{d\varepsilon^n}\mathcal N_\alpha(\mathcal W_\varepsilon[\eta;\theta])\right|_{\varepsilon=0}\\
 &\quad=\left.\frac1\tau\frac{d^n}{d\varepsilon^n}
   \left(\bigl[(1+\rho[\mathcal W_\varepsilon[\eta;\theta]])^\alpha-1\bigr]
   \mathcal L_Q[\mathcal W_\varepsilon[\eta;\theta]]\right)\right|_{\varepsilon=0}\\
 &\quad=\frac1\tau\sum_{j=1}^{n-1}\binom nj
   \left.\frac{d^j}{d\varepsilon^j}
   \bigl[(1+\rho[\mathcal W_\varepsilon[\eta;\theta]])^\alpha-1\bigr]\right|_{\varepsilon=0}
   \mathcal L_Q[\mathcal W_{n-j}][\eta;\theta]\\
 &\quad=\frac1\tau\sum_{j=1}^{n-1}\binom nj\mathcal L_Q[\mathcal W_{n-j}][\eta;\theta]
   \sum_{l=1}^j(\alpha)_l
   B_{j,l}\bigl(\rho[\mathcal W_1][\eta;\theta],\ldots,
               \rho[\mathcal W_{j-l+1}][\eta;\theta]\bigr).
\end{aligned}
\end{equation}
Both factors vanish at $\varepsilon=0$, so the endpoint terms in the Leibniz sum vanish. In particular, the linear coefficient is zero. We write
$$
 \begin{aligned}
 \mathcal N_\alpha(\mathcal W_\varepsilon[\eta;\theta])
 &=\sum_{n=2}^\infty\frac{\varepsilon^n}{n!}\mathfrak N_n[\eta;\theta],\\
 \mathfrak N_n[\eta;\theta](x,v)
 &=(2\pi)^{-d/2}\int_{\mathbb R^{d\times n}}
   \widehat{\mathfrak N}_n(k_1,\ldots,k_n,v,\theta)
   \Pi\eta(x,\bm k_n,t)\,d\bm k_n,
 \end{aligned}
$$
and denote the density coefficients as 
$$
    \widehat{\rho}_m(\bm{k}_m,\theta)
    =
    \int_{\mathbb R^d}
    \widehat{\mathfrak W}_m(\bm{k}_m,v,\theta)\,dv. 
$$
For $1\leq  l\leq  j$, we define the corresponding Bell polynomial in frequency space as the symmetrized convolution kernels
\begin{equation}\label{freq_Bell_poly}
    \widehat{B}_{j,l}
    \bigl[\widehat{\rho}_1,\ldots,\widehat{\rho}_{j-l+1}\bigr]
    (\bm{k}_j,\theta)
    :=
    \sum_{\pi\in\Pi_{j,l}}
    \prod_{B\in\pi}
    \widehat{\rho}_{|B|}
    \bigl((k_i)_{i\in B},\theta\bigr),
\end{equation}
where $\Pi_{j,l}$ denotes the set of all partitions of $\{1,\ldots,j\}$ into $l$ nonempty blocks. Equivalently, \eqref{freq_Bell_poly} is the Fourier-space version of the Bell polynomial $B_{j,l}$, with the product in physical space replaced by convolution, or, after fixing the ordered variables $\bm{k}_j$, by summation over all block partitions. The partition formula makes the symmetrization explicit and preserves the nonnegative combinatorial multiplicities. Thus the coefficient functions of the nonlinearity evaluated on the parametrization are given by
\begin{equation}\label{hat_N_n_explicit}
\begin{split}
    \widehat{\mathfrak N}_n(\bm{k}_n,v,\theta)
    =
    \frac{1}{\tau}
    \sum_{j=1}^{n-1}
    \binom{n}{j}
    \sum_{l=1}^{j}
    (\alpha)_l\,(2\pi)^{-dl/2}\,
    \operatorname{Sym}_{\bm{k}_n}
    \Big[
        &
        \widehat{B}_{j,l}
        \bigl[
            \widehat{\rho}_1,\ldots,
            \widehat{\rho}_{j-l+1}
        \bigr]
        (k_1,\ldots,k_j,\theta)
        \\
        &\times
        \mathcal{L}_Q[\hat{\mathfrak W}_{n-j}]
        (k_{j+1},\ldots,k_n,v,\theta)
    \Big],
\end{split}
\end{equation}
for every $n\geq  2$, while $\widehat{\mathfrak N}_1=0$.  Here $\operatorname{Sym}_{\bm{k}_n}$ denotes symmetrization over the $n$ frequency variables, see \eqref{def_symmetrization}. Consequently,
$$
\begin{aligned}
    \widehat{\mathcal N_\alpha(\mathcal W[\eta])}(k,v)
    =
    \sum_{n=2}^{\infty}
    \frac{1}{n!}
    \int_{\mathbb R^{d\times n}}
    \delta\!\left(k-\sum_{i=1}^n k_i\right)
    \widehat{\mathfrak N}_n(\bm{k}_n,v,\theta)
    \Pi \xi(\bm{k}_n,t)
    \,d\bm{k}_n.
\end{aligned}
$$
Equivalently, suppressing the delta distribution imposing the frequency sum, the nonlinear coefficient functional is
$$
    \widehat{\mathfrak N}_n(\bm{k}_n,v,\theta)
    =
    \widehat{\mathcal C}_n
    [\hat{\mathfrak W}_1,\ldots,\hat{\mathfrak W}_{n-1}]
    (\bm{k}_n,v,\theta),
$$
where $\widehat{\mathcal C}_n$ is the polynomial expression defined explicitly by \eqref{hat_N_n_explicit}.\\
Combining the above two calculations leads to the \emph{unsteady homological  equation} at the index vector $\bm{k}_n = (k_1,...,k_n)$:
\begin{equation}\label{unsteady_cohomology_equ}
\frac{\partial \hat{\mathfrak{W}}_n}{\partial \theta} = \left(\mathcal{L}_{\Sigma \bm{k}_n}-\sum_{j=1}^n \lambda(k_j)\right)\hat{\mathfrak{W}}_n + \hat{\mathcal{C}}_n[\hat{\mathfrak{W}}_1,...,\hat{\mathfrak{W}}_{n-1}].
\end{equation}
Typically, the coefficient functions $\hat{\mathfrak{W}}_n$ are assumed to be independent of time, for which \eqref{unsteady_cohomology_equ} reduces to the classical steady homological equation
\begin{equation}\label{steady_cohomolgy_equ}
 \left(\sum_{j=1}^n \lambda(k_j)-\mathcal{L}_{\Sigma \bm{k}_n}\right)\hat{\mathfrak{W}}_n    =  \hat{\mathcal{C}}_n[\hat{\mathfrak{W}}_1,...,\hat{\mathfrak{W}}_{n-1}],
\end{equation}
for $\theta$-independent $\hat{\mathcal{W}}_n$. At this point the hierarchy is formal. The following sections identify its resonant obstruction, solve the unsteady equations by semigroup methods, and establish convergence of the resulting invariant parametrization.

\begin{remark}\label{remark_rho_N_n}
From the representation \eqref{hat_N_n_explicit}, we see that the nonlinear coefficient functional satisfies
$$
    \rho[\hat{\mathfrak{N}}_n] = 0,
$$
since $\rho\circ \mathcal{L}_Q =0$ and every summand in  \eqref{hat_N_n_explicit} carries such a factor. This property will become important in the convergence proof in Section \ref{sec_proof_main_theorem}. 
\end{remark}

\subsection{Homological  Operator and Resonant Frequencies}\label{subsection_homological}

In the following, we refer to $k_1,\ldots,k_n$ as the individual frequencies and to their sum $\Sigma\bm k_n$ as the summed frequency. \\
We call the linear operator
\begin{equation}\label{def_cohomol_operator}
    \mathfrak{L}_{\bm{k}_n} =  \left(\mathcal{L}_{\Sigma \bm{k}_n}-\sum_{j=1}^n \lambda(k_j)\right),
\end{equation}
the \emph{homological  operator}. For $\bm{k}_n\in B_c^n$, we denote its spectral bound, and hence its exponential decay or growth rate in the estimates below, by
\begin{equation}
   \Omega(\bm{k}_n) = \hat{\lambda}(\Sigma \bm{k}_n) - \sum_{j=1}^n\lambda(k_j),
\end{equation}
i.e., $\Omega(\bm{k}_n)=s(\mathfrak{L}_{\bm{k}_n}):=\sup\{\Re z:z\in\sigma(\mathfrak{L}_{\bm{k}_n})\}$.\\

Whenever $\mathfrak{L}_{\bm{k}_n}$ is invertible, we can solve \eqref{steady_cohomolgy_equ} directly for any right-hand side. For certain combinations of frequencies, however, the homological  operator $\mathfrak{L}_{\bm{k}_n}$ becomes singular, thus motivating the following definition. 

\begin{definition}
    An \emph{internal resonance} at order $n$ occurs whenever the sum of $n$ hydrodynamic eigenvalues equals the hydrodynamic eigenvalue at their summed frequency,
\begin{equation}\label{int_resonance}
        I^{(n)}_{\rm int} = \left\{\bm{k}_n \in B_c^n: |\Sigma\bm{k}_n|<k_c,\quad \lambda(\Sigma\bm{k}_n) = \sum_{j=1}^n \lambda(k_j) \right\},
\end{equation}
while an \emph{external resonance} occurs when the sum of $n$ hydrodynamic eigenvalues equals the real part of the essential spectrum,
\begin{equation}\label{ext_resonance}
     I^{(n)}_{\rm ext} = \left\{\bm{k}_n \in B_c^n: \sum_{j=1}^n \lambda(k_j) = -\frac{1}{\tau} \right\}.
\end{equation}
Since external frequencies occur when the shifted essential spectrum passes the imaginary axis, we also call them \emph{essential resonances}.
 We denote the set of all resonant frequency vectors as 
$$
    I^{(n)} = I^{(n)}_{\rm int} \cup I^{(n)}_{\rm ext},
$$
and call them \emph{resonance sets}. At a point where the defining level-set function is smooth and has nonzero gradient, the implicit-function theorem gives a local hypersurface - the resonance hypersurface.   
\end{definition}

\begin{figure}
    \centering
    \includegraphics[width=0.3\linewidth]{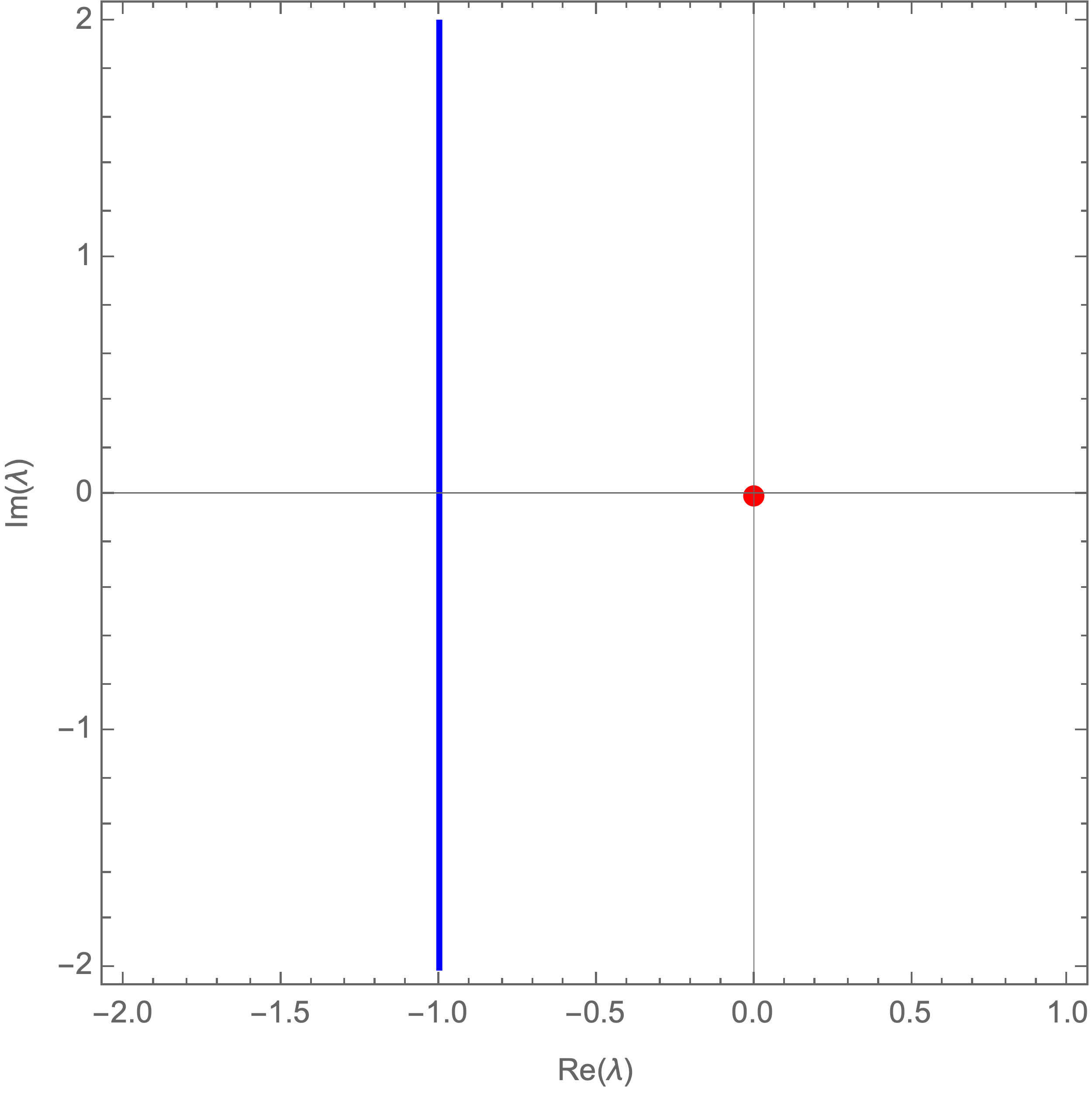}
    \includegraphics[width=0.3\linewidth]{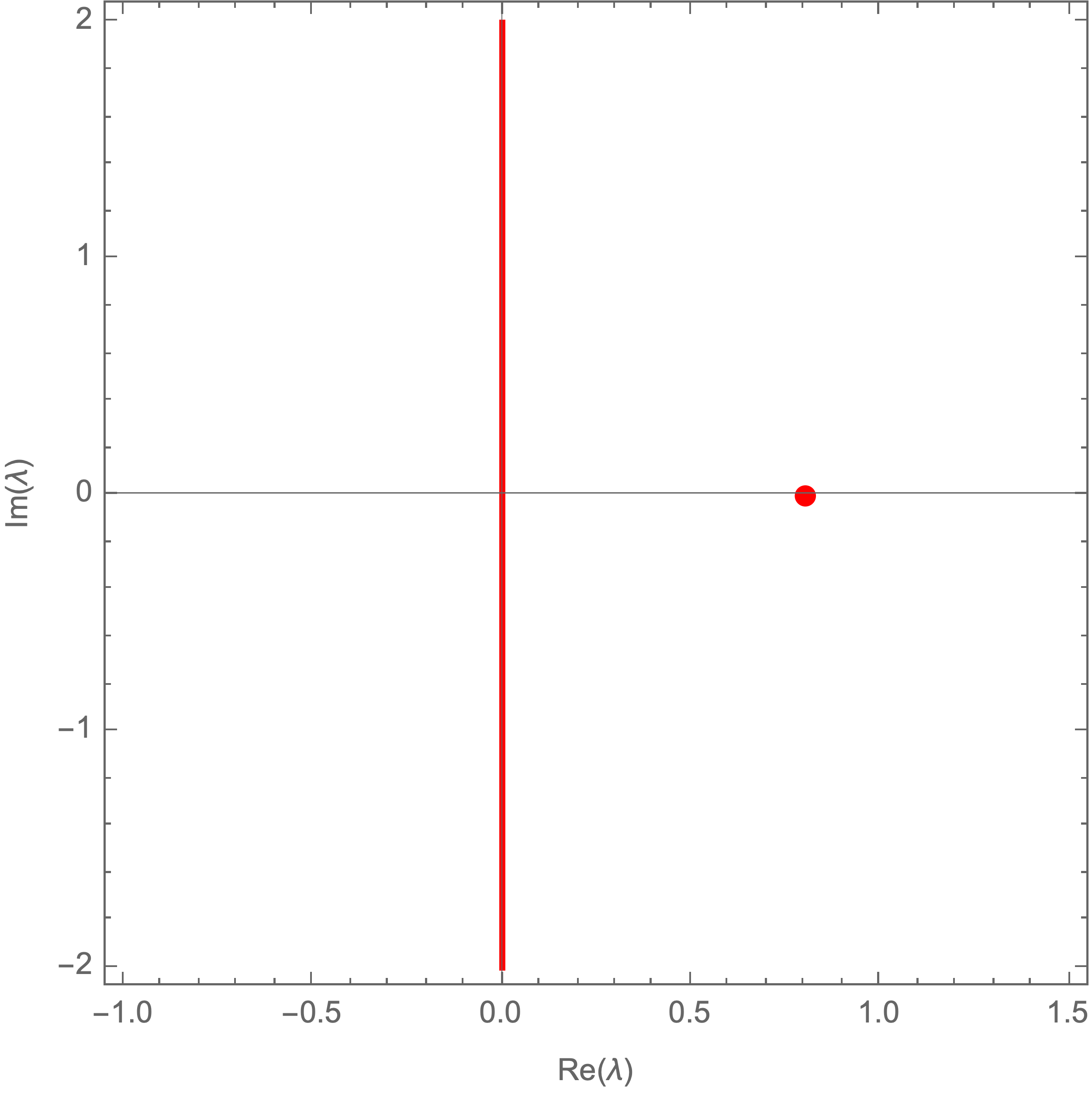}
    \includegraphics[width= 0.3\linewidth]{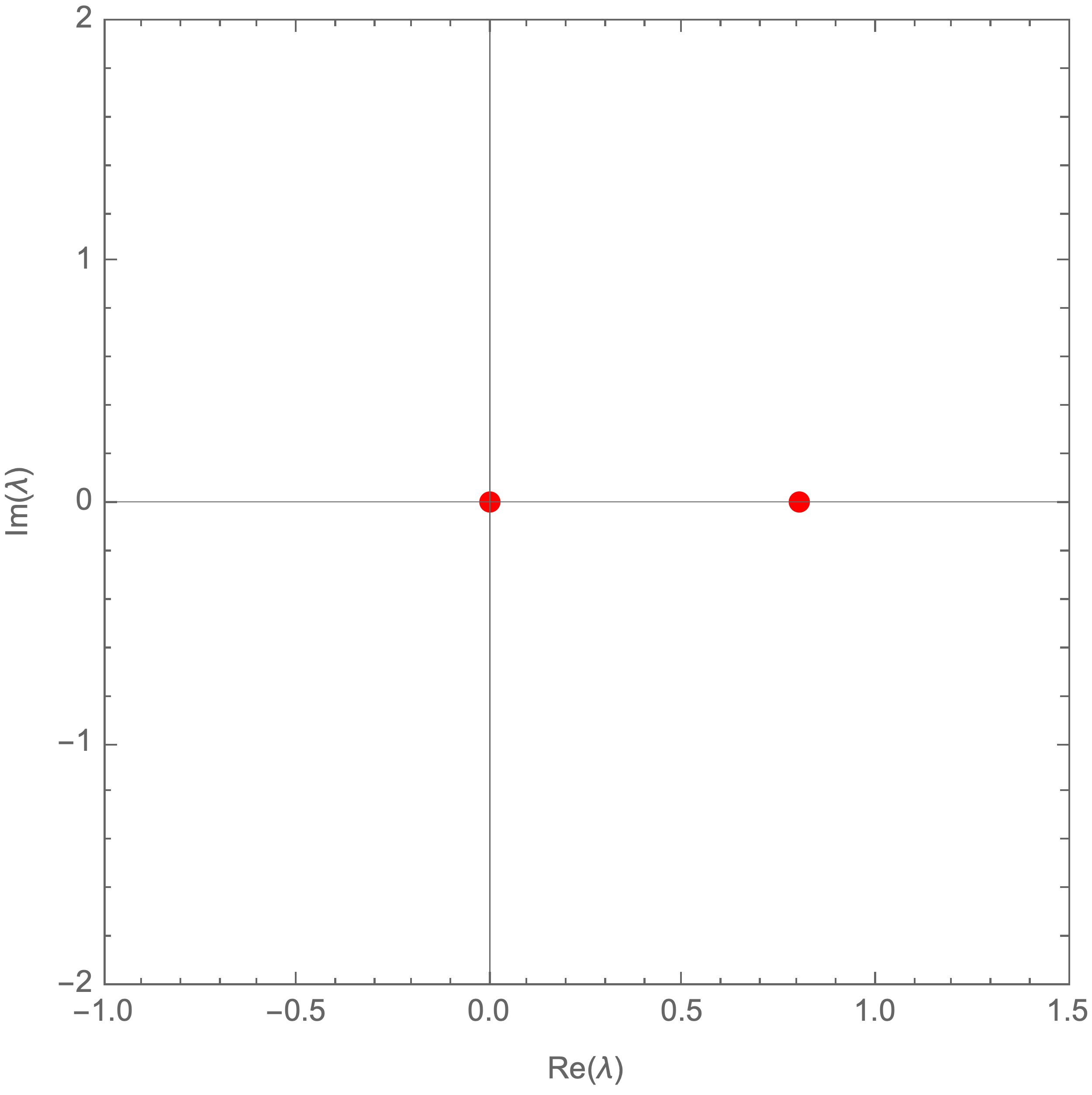}
    \caption{At an internally resonant index (left spectral plot), the isolated eigenvalue is
shifted to the origin and the homological operator $\mathfrak{L}_{\bm{k}_n}$
acquires a one-dimensional kernel.
The middle spectral plot shows an external resonance with
$0<|\Sigma\bm{k}_n|<k_c$. At any external resonance with
$\Sigma\bm{k}_n\ne 0$, the shifted essential spectrum is $\ri\mathbb R$,
and $\mathfrak{L}_{\bm{k}_n}-\ri\omega I$ is non-Fredholm for every
$\omega\in\mathbb R$. The right spectral plot shows an external resonance with
$\Sigma\bm{k}_n=0$, plotted for $\tau=1$. The red point at the origin
is an eigenvalue of infinite multiplicity and constitutes the entire
shifted essential spectrum. The second red point is the simple eigenvalue
$1/\tau$. Thus $\mathfrak{L}_{\bm{k}_n}-zI$ is non-Fredholm only at $z=0$.
Figure~\ref{res_surface} gives the phase diagram for the spectrum of
$\mathfrak{L}_{\bm{k}_2}$ as a function of the frequency pair $(k_1,k_2)$
in one spatial dimension.}
    \label{resonance_spctrum}
\end{figure}

\begin{remark}
Resonance sets are mirror- and rotationally symmetric, i.e.,
$$
\begin{aligned}
    \bm{k}_n \in I^{(n)} &\implies -\bm{k}_n \in I^{(n)},\\
    (k_1,...,k_n)\in I^{(n)} &\implies (Qk_1,...,Qk_n)\in I^{(n)},
    \end{aligned}
$$
for any rotation matrices $Q$, $Q^TQ=\text{Id}$, since $\lambda$ depends on $k$ only through $|k|$. 
\end{remark}

Resonances reflect the stability properties of the homological  operator \eqref{def_cohomol_operator}, where the resonance sets \eqref{int_resonance} and \eqref{ext_resonance} mark changes in the spectral signature of $\mathfrak{L}_{\bm{k}_n}$, see Figure \ref{resonance_spctrum}. We thus define various frequency regions depending on the spectral signature of $\mathfrak{L}_{\bm{k}_n}$.

\begin{definition}
    We call the sets of frequencies
    $$
    K_{\rm stable}^{(n)} := \left\{\bm{k}_n\in B_c^n: \sup_{z\in\sigma(\mathfrak{L}_{\bm{k}_n})}\Re z<0 \right\},\quad  K_{\rm unstable}^{(n)} := \left\{\bm{k}_n\in B_c^n: \inf_{z\in\sigma(\mathfrak{L}_{\bm{k}_n})}\Re z>0 \right\},
    $$
    \emph{stable region} and \emph{unstable region}, respectively, while the set
    $$
      K_{\rm mixed}^{(n)} :=   (B_c)^n\setminus ( K_{\rm stable}^{(n)} \cup  K_{\rm unstable}^{(n)}),
    $$
    is called \emph{mixed region}.
\end{definition}

\begin{figure}
    \centering
\includegraphics[width=1\linewidth]{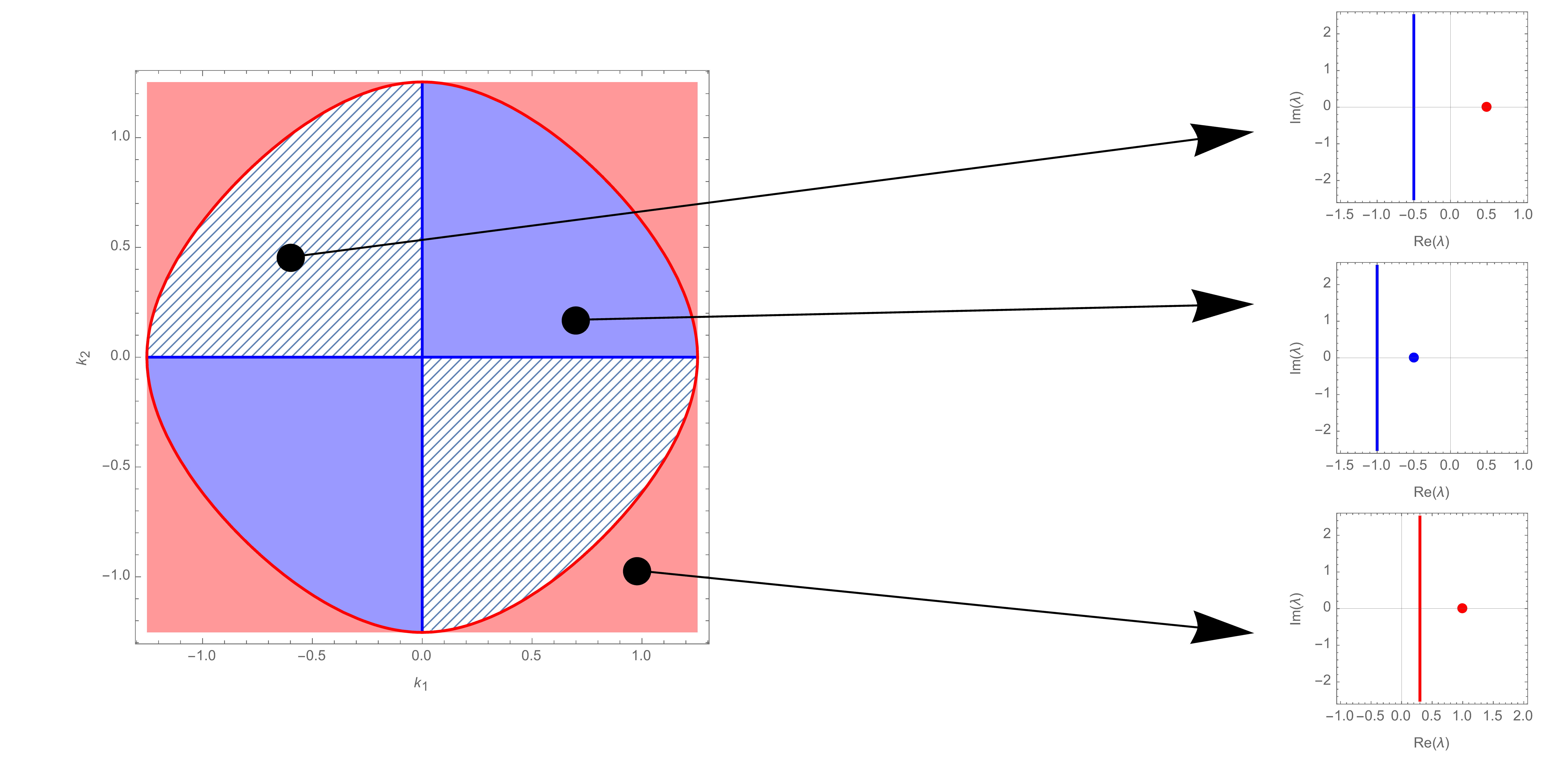}
    \caption{Stable, unstable and mixed regions for the operator $\mathfrak{L}_{\bm{k}_2}$ for $d=1$. Internal resonances (blue lines) from the level sets \eqref{int_resonance} correspond to the boundaries of the stable regions with the mixed regions, while external resonances (red curves) from the level sets \eqref{ext_resonance} correspond to the boundaries of the unstable regions.
    The arrows indicate spectral types within the regions schematically. A vertical essential line corresponds to $k_1+k_2\ne0$, while on the antidiagonal $k_1+k_2=0$, the shifted essential spectrum is the
    single point $-\tau^{-1}-\lambda(k_1)-\lambda(k_2)$.}
  \label{res_surface}
\end{figure}

The two spectral quantities relevant to the subsequent analysis of the resonant homological equations are the spectral bound, equal to the shifted hydrodynamic eigenvalue when the summed frequency is subcritical, $\hat{\lambda}(\Sigma\bm{k}_n)-\sum_j\lambda(k_j)$ and the shifted essential spectrum $-1/\tau-\sum_j\lambda(k_j)$. For subcritical summed frequencies the internal resonance is the zero of the shifted hydrodynamic eigenvalue, while the external resonance is the zero of the shifted essential spectral bound. For a frequency in the stable region we have that 
$$
    \hat{\lambda}\left(\sum_{j=1}^n k_j\right) < \sum_{j=1}^n\lambda(k_j),\quad \bm{k}_n \in  K_{\rm stable}^{(n)},
$$
while for a frequency in the unstable region we have that
$$
    \sum_{j=1}^n\lambda(k_j)<-\frac{1}{\tau} ,\quad \bm{k}_n \in  K_{\rm unstable}^{(n)}.
$$
For a frequency in the mixed region, the estimates
$$
\hat{\lambda}\left(\sum_{j=1}^n k_j\right) \geq  \sum_{j=1}^n\lambda(k_j),\quad \sum_{j=1}^n\lambda(k_j)\geq  -\frac{1}{\tau} ,\quad \bm{k}_n \in K_{\rm mixed}^{(n)},
$$
hold. At regular separating portions, internal resonances mark the stable--mixed boundary and external resonances mark the unstable boundary, while at degenerate points or critical summed frequencies the boundary geometry must be interpreted through the corresponding closures of the sets. The mixed region includes its resonant boundary, hence the non-strict inequalities above. The internal resonance set, however, is disjoint from the external one and their closures can meet at $|\Sigma\bm{k}_n|=k_c$, where the hydrodynamic eigenvalue reaches $-1/\tau$. The cutoff by $\varepsilon_B$ separates the eigenfunctions at the individual frequencies from this threshold, but does not constrain their sums: nonlinear summed frequencies can reach and cross the critical sphere. At order two this behavior can obstruct stationary reconstruction in the kinetic $L^2$-space, see Proposition~\ref{remark_stationary_quadratic_integrability}. Figure \ref{res_surface} shows the phase diagram of different spectral configurations for
$\mathfrak{L}_{\bm{k}_2}$ and $d=1$.\\
At an external resonance, zero belongs to the essential spectrum of
$\mathfrak L_{\bm k_n}$, so the steady homological equation cannot be
solved by a bounded inverse. This does not exclude a solution for a
particular forcing, nor does it force every unsteady solution to grow.
The latter also depends on the forcing and the homological initial
condition. For the quadratic term, the full resolvent and the
convolution against the parameter function must therefore be retained.
As shown in Proposition~\ref{remark_stationary_quadratic_integrability},
at a fixed summed frequency $|k_1+k_2|>k_c$ the singularity across a
regular external resonance is locally integrable in the individual
frequencies. This estimate ceases to be uniform as
$|k_1+k_2|\downarrow k_c$, where the closure of the internal resonance
set meets the external resonance set. The reconstructed stationary
quadratic term can then fail to belong to kinetic $L^2$, even though
the finite-time semigroup formula for the unsteady equation remains
well-defined.\\

\subsection{Explicit Analysis of the Quadratic Term}

Before we proceed with the general analysis, we take a closer look at the quadratic contribution of \eqref{defhydro} by solving the invariance equation explicitly at second order in $\xi$. This shows how the nonlinearity acts on the linear term and gives the quadratic coefficient used to compute the spectral closure in Section~\ref{sec_slow_spectral_closure}. The second-order term also distinguishes locally integrable external singularities from the obstruction at the critical spectral frequency, see Proposition~\ref{remark_stationary_quadratic_integrability}. The calculation also provides the base case for the later inductive bounds on the coefficient functions.\\
The symmetric form of the quadratic term of $\mathcal{N}_\alpha$ reads,
$$
    \mathcal{N}_\alpha^{(2)}[f,g] = \frac{\alpha}{2\tau}\Big(\rho[f]\mathcal{L}_Q[g] + \rho[g]\mathcal{L}_Q[f] \Big),
$$
and the unsteady homological  equation at order two takes the following form:
\begin{equation}
 \begin{split}
 \partial_\theta\hat{\mathcal W}_2
 +[\lambda(k_1)+\lambda(k_2)]\hat{\mathcal W}_2
 =\mathcal L_k\hat{\mathcal W}_2 +2\widehat{\mathcal N}_\alpha^{(2)}
   (\hat{\mathcal W}_1(\cdot;k_1),
    \hat{\mathcal W}_1(\cdot;k_2)).
 \end{split}
\end{equation}
Using \eqref{W1}, the Fourier product rule and $\rho[\hat e]=1$ give
$$
 \begin{aligned}
 \widehat{\mathcal N}_\alpha^{(2)}
 (\hat{\mathcal W}_1(\cdot;k_1),\hat{\mathcal W}_1(\cdot;k_2))
 ={}&\frac{\alpha (2\pi)^{-d/2} }{2\tau}
 \delta(k-k_1-k_2)\mathbf1_{B_c}(k_1)\mathbf1_{B_c}(k_2)\\
 &\times\bigl(2\phi-\hat e(k_1)-\hat e(k_2)\bigr).
 \end{aligned}
$$
After removing the common delta and cutoffs, this yields
\begin{equation}\label{cohomol2_explicit}
 \begin{split}
 \partial_\theta\hat{\mathfrak W}_2
 +[\lambda(k_1)+\lambda(k_2)-\mathcal L_{k_1+k_2}]
     \hat{\mathfrak W}_2
 =\frac{\alpha (2\pi)^{-d/2} }{\tau}
  \bigl(2\phi-\hat e(k_1)-\hat e(k_2)\bigr).
 \end{split}
\end{equation}
We stress that, at the fiber level, $\lambda(k)$ is an isolated eigenvalue of $\mathcal{L}_k$ and $\hat e(k,\cdot)$ its eigenvector for $|k|<k_c$. Only after reconstructing the full spatial operator do the Fourier modes $\hat e(k,v)e^{\ri k\cdot x}$ appear as generalized eigenfunctions.\\
The resonance curves for equation \eqref{cohomol2_explicit} are depicted in Figure \ref{res_surface} for $d=1$. At a non-resonant index we may choose a $\theta$-independent particular solution of \eqref{cohomol2_explicit}, obtaining
\begin{equation}\label{solve_W2_nonres}
    \hat{\mathfrak{W}}_2(k_1,k_2,v) = \frac{\alpha (2\pi)^{-d/2} }{\tau}[\lambda(k_1)+\lambda(k_2)-\mathcal{L}_{k_1+k_2}]^{-1} \Big(2\phi(v) - \hat{e}(k_2,v) - \hat{e}(k_1,v)\Big)
\end{equation}

\begin{remark}[Trivial center direction]
For $k_1=k_2=0$, equation \eqref{cohomol2_explicit} reduces to 
$$
    \frac{\partial\hat{\mathfrak{W}}_2}{\partial \theta}(0,0)  +\frac{1}{\tau}\left(\hat{\mathfrak{W}}_2(0,0)-\rho[\hat{\mathfrak{W}}_2](0,0)\phi\right)   =  0,
$$
which is solved trivially by $\hat{\mathfrak{W}}_2(0,0) = 0$. With the normalization $\hat{\mathfrak V}_n(0,\ldots,0)=0$ for $n\geq2$, the selected higher-order coefficients in the pure center direction vanish as well. We stress that this is a statement about the coefficient kernels at zero frequency, not an embedded center line in $\Xi(\varepsilon_B)$. A continuous function supported only at $k=0$ vanishes identically, and multiplication by $\lambda(k)$ therefore has trivial kernel on that parameter space. A nonzero spatially homogeneous perturbation instead has Fourier transform $(2\pi)^{d/2}\phi(v)\delta_0(k)$, which, being a tempered distribution, is not an admissible continuous parameter function. The equilibrium family \eqref{defMeq} is obtained by varying the independent background density $\rho_0$, while adjoining it as a center coordinate would require an enlarged phase space and a separate construction.
\end{remark}

\subsection{Alternative Approaches and Their Limitations}\label{what_does_not_work}

The steady homological equation \eqref{steady_cohomolgy_equ} requires
inversion of $\mathfrak L_{\bm k_n}$, whereas the unsteady equation
\eqref{unsteady_cohomology_equ} requires its semigroup only over finite
time intervals. We discuss why restricting the individual frequencies
or changing the reduced dynamics does not, in general, justify steady
inversion. In doing so, we distinguish failure of inversion at a
resonant frequency vector from failure of the reconstructed convolution
in the kinetic function space. Proposition~\ref{remark_stationary_quadratic_integrability}
gives an explicit quadratic example of the latter.

\begin{remark}
An annular restriction of the individual frequencies can exclude
external resonances. Indeed, if
$$
 \lambda^{-1}\!\left(-\frac1{2\tau}\right)<|k_j|<k_c,
$$
then monotonicity of $\lambda$ gives
$$
 \sum_{j=1}^n\lambda(k_j)<-\frac{n}{2\tau}\le-\frac1\tau,
 \qquad n\geq 2.
$$
Thus all such frequency vectors belong to $K_{\rm unstable}^{(n)}$.
This restricts the selected parameter space and excludes the slow
frequencies near zero, while it does not supply the inertial-manifold gap
condition \eqref{spectral_gap}.\\
Conversely, no fixed ball of individual frequencies about zero with
positive radius excludes external resonances at every order. For any
$0<k_{\rm max}<k_c$, choose $n$ so large that
$n\lambda(k_{\rm max})<-1/\tau$. Since $\lambda(0)=0$ and $\lambda$
is continuous, there is $s_*\in(0,1)$ such that
$n\lambda(s_*k_{\rm max})=-1/\tau$. Taking
$k_1=\cdots=k_n=s_*k_{\rm max}\omega$, with $|\omega|=1$, gives
an external resonance strictly inside the chosen ball.
\end{remark}

\begin{remark}
Internal resonances can be treated by retaining nonlinear terms in
the reduced dynamics, as in the parametrization method
\cite{cabre2005parameterization}. One then seeks
\begin{equation}
 \mathcal R[\xi](k)
 =\sum_{n=1}^\infty\int_{\mathbb R^{d\times n}}
 \mathcal R_n(k;\bm k_n)\Pi\xi(\bm k_n,t)\,d\bm k_n,
 \qquad
 \mathcal R_1(k;k_1)=\lambda(k)\delta(k-k_1).
\end{equation}
At an internal resonance, projection of the steady homological
equation onto the hydrodynamic eigenspace at $\Sigma\bm k_n$ shows
that the forcing must have zero projection under $P_{\Sigma\bm k_n}$.
When this hydrodynamic direction is included in the reduced space,
a nonlinear reduced term can cancel that projection. Note that the Riesz projection is given by
\eqref{Pk_rank_one_corrected}, and differs from $P_\rho$ at nonzero
frequency. If the projected forcing is nonzero at a regular internal
resonance, steady inversion without this cancellation has a simple pole,
which is not locally absolutely integrable in a transverse frequency variable.\\
Such a correction does not change the essential spectrum of
$\mathcal L_{\Sigma\bm k_n}$ and therefore does not ensure solvability
at an external resonance. Nevertheless, noninvertibility alone does
not exclude a solution for a particular forcing. We next examine
the actual quadratic forcing and its reconstructed stationary term,
with the prescribed linear reduced dynamics \eqref{reduced_dyn_linear}.
\end{remark}

\begin{proposition}[Obstruction to stationary quadratic reconstruction]
\label{prop_stationary_quadratic_obstruction}
\label{remark_stationary_quadratic_integrability}
Assume $\alpha>0$, $d\geq 2$, and
$$
 0<\varepsilon_B<k_c-\kappa_*,\qquad
 \lambda(\kappa_*)=-\frac1{2\tau}.
$$
For the linear reduced dynamics \eqref{reduced_dyn_linear}, there are
nonnegative even $\xi\in C_c^\infty(\operatorname{int}B_c(\varepsilon_B))$
of arbitrarily small amplitude for which the ordinary stationary quadratic
reconstruction \eqref{solve_W2_nonres} is defined as a Bochner integral
for almost every $k$ with $|k|>k_c$ and satisfies
\begin{equation}\label{stationary_quadratic_reconstruction}
 \begin{gathered}
 \hat f_2(k,v)=\frac12\int_{\mathbb R^d}
 \hat{\mathfrak W}_2(k_1,k-k_1,v)\xi(k_1)\xi(k-k_1)\,dk_1,\\
 \int_{|k|>k_c}\|\hat f_2(k,\cdot)\|_{L_v^2(\phi^{-1})}^2\,dk=\infty.
 \end{gathered}
\end{equation}
\end{proposition}

\begin{proof}[Proof]
For fixed $|k|>k_c$ the convolution is locally integrable across a
regular resonance. We show that as $\eta=\tau(|k|-k_c)\downarrow0$ the
density contribution is bounded below by $c/\eta$ on a bounded velocity slab of
thickness $\eta$. The other terms grow at most logarithmically there.
Squaring and integrating over the slab will leave a lower bound $c'/\eta$
for the radial integral.\\
Strict monotonicity and concavity give $k_c/2<\kappa_*<k_c$.
Choose orthonormal vectors $e_1,e_2$ and set
\begin{equation}\label{k10k20}
 k_1^0=\frac{k_c}{2}e_1+\sqrt{\kappa_*^2-k_c^2/4}\,e_2,
 \qquad
 k_2^0=\frac{k_c}{2}e_1-\sqrt{\kappa_*^2-k_c^2/4}\,e_2.
\end{equation}
Both $k_1^0,k_2^0\in B_c(\varepsilon_B)$ and satisfy
$$
 |k_1^0+k_2^0|=k_c,\qquad
 \lambda(k_1^0)+\lambda(k_2^0)=-\frac1\tau.
$$
This pair belongs to $I_{\rm ext}^{(2)}$ and to the closure of
$I_{\rm int}^{(2)}$. Indeed, for $s<\kappa_*$ close to
$\kappa_*$, $\lambda^{-1}(2\lambda(s))$ is strictly between $0$ and
$2s$ and tends to $k_c$. Two leaves of length $s$, chosen symmetrically
about $e_1$ in the $(e_1,e_2)$-plane with this sum length, give internal
resonances tending to \eqref{k10k20}. At the limiting sum there is no
isolated hydrodynamic eigenvalue. For $|k_1+k_2|>k_c$, the essential resonance
$\Lambda_2=1/\tau$ separates $K_{\rm stable}^{(2)}$, where
$\Lambda_2<1/\tau$, from $K_{\rm unstable}^{(2)}$, where
$\Lambda_2>1/\tau$. We approach the displayed pair through these
supercritical summed frequencies.\\
The full stationary resolvent, including its density contribution,
follows from \eqref{cohomol2_explicit} and \eqref{resolvent_formula}, and is given by 
\begin{equation}\label{stationary_quadratic_full_resolvent}
 \begin{aligned}
 \hat{\mathfrak W}_2(k_1,k_2,v)
 &=\frac{\alpha(2\pi)^{-d/2}
       \bigl(2\phi(v)-\hat e(k_1,v)-\hat e(k_2,v)\bigr)
       +\phi(v)\rho[\hat{\mathfrak W}_2](k_1,k_2)}
       {1+\tau\lambda(k_1)+\tau\lambda(k_2)
                  +\ri\tau(k_1+k_2)\cdot v},\\
 \rho[\hat{\mathfrak W}_2](k_1,k_2)
 &=\frac{\alpha(2\pi)^{-d/2}}
       {\Delta_{k_1+k_2}(\lambda(k_1)+\lambda(k_2))}
   \int_{\mathbb R^d}
   \frac{2\phi(w)-\hat e(k_1,w)-\hat e(k_2,w)}
        {1+\tau\lambda(k_1)+\tau\lambda(k_2)
                  +\ri\tau(k_1+k_2)\cdot w}\,dw.
 \end{aligned}
\end{equation}
The scalar integral in the second line is real by conjugation under
$w\mapsto-w$. It has bounded continuous limits from each side of
$\Lambda_2=1/\tau$ near the frequency pair $(k_1^0,k_2^0)$.

To see this, put $k=k_1+k_2$ and $\sigma=\operatorname{sgn}(1-\tau\Lambda_2)$.
The Laplace representation of each of the two denominators and the
Gaussian Fourier transform give
$$
 \int_{\mathbb R^d}\frac{\hat e(k_j,v)}
 {1-\tau\Lambda_2+\ri\tau k\cdot v}\,dv
 =\sigma\int_0^\infty\!\!\int_0^\infty
 e^{-|1-\tau\Lambda_2|s-(1+\tau\lambda(k_j))t
      -\tau^2|\sigma s k+t k_j|^2/2}\,ds\,dt .
$$
For $\phi$ in place of $\hat e(k_j)$ one omits $t$ and the terms
containing it. The $s$-integral of the Gaussian is at most
$\sqrt{2\pi}/(\tau|k|)$ by projection onto $k$. At fixed $t$ it is
continuous up to either side by Gaussian domination in $s$, and the
remaining factor is bounded by $e^{-ct}$ for some $c>0$ on the leaf
supports. This proves the asserted one-sided limits by dominated
convergence. At \eqref{k10k20} one has $k\cdot k_j^0=k_c^2/2$.
On the side $\sigma=1$, integration first in $s$ at the boundary gives
$(\tau k_c)^{-1}e^{-\tau^2\kappa_*^2t^2/2}m(\tau k_ct/2)$ for the mixed
Gaussian and $m(0)/(\tau k_c)=1$ for the single Gaussian. Hence
\par
\begin{equation}\label{stationary_merger_numerator}
 \begin{aligned}
 &\lim_{\substack{k_j\to k_j^0\ (j=1,2)\\\Lambda_2(k_1,k_2)<1/\tau}}
 \int_{\mathbb R^d}
 \frac{2\phi(v)-\hat e(k_1,v)-\hat e(k_2,v)}
      {1+\tau\lambda(k_1)+\tau\lambda(k_2)
                 +\ri\tau(k_1+k_2)\cdot v}\,dv\\
 &\qquad=2-\frac{2}{\tau k_c}\int_0^\infty
 e^{-t/2-\tau^2\kappa_*^2t^2/2}
 m\!\left(\frac{\tau k_ct}{2}\right)\,dt>0,
 \end{aligned}
\end{equation}

where $m$ is again the Mills ratio \eqref{defm}. The strict inequality follows
from $m(\tau k_ct/2)<m(0)=\tau k_c$ for $t>0$ and the dispersion
relation \eqref{kernel_characteristic}, which gives
$\int_0^\infty e^{-t/2-\tau^2\kappa_*^2t^2/2}\,dt=1$.
Consequently, the scalar numerator of the density contribution in
\eqref{stationary_quadratic_full_resolvent} is bounded below by a
positive constant on the nearby stable region.\\
The same Gaussian integration gives, with the already defined
$\Lambda_2=-\lambda(k_1)-\lambda(k_2)$,
\begin{equation}
 \Delta_{k_1+k_2}(-\Lambda_2)
 =1-\frac{\operatorname{sgn}(1-\tau\Lambda_2)}
              {\tau|k_1+k_2|}
 m\!\left(\frac{|1-\tau\Lambda_2|}{\tau|k_1+k_2|}\right),
 \qquad \Lambda_2\ne\frac1\tau.
\end{equation}
For subcritical summed frequencies, the zeros of this denominator
are exactly at the internal resonances. At the chosen pair its
limit from $\Lambda_2<1/\tau$ vanishes at an essential resonance.
We write the summed frequency in the reconstruction as
$k=(k_c+\eta/\tau)\omega$, where $|\omega|=1$ and $\eta>0$.
Since $m'(0)=-1$, the following expansions hold:
\begin{equation}\label{stationary_merger_denominator}
 \begin{aligned}
 \Delta_k(-\Lambda_2)
 &=\frac{\eta}{\tau|k|}
   +\frac{1-\tau\Lambda_2}{\tau^2|k|^2}
   +O\bigl((1-\tau\Lambda_2)^2\bigr)
   \asymp\eta+1-\tau\Lambda_2,
 &&\Lambda_2<\frac1\tau,\\
 \Delta_k(-\Lambda_2)&>1,
 &&\Lambda_2>\frac1\tau.
 \end{aligned}
\end{equation}
For each fixed $|k|>k_c$, the denominator on the stable side stays
positive at the essential resonance. In that case Gaussian
integration gives the local bound
$\|\hat{\mathfrak W}_2\|_{L_v^2(\phi^{-1})}
\leq C_k|1-\tau\Lambda_2|^{-1/2}$, which is integrable across a
regular resonance hypersurface. The difficulty is that this bound
is not uniform as $|k|\downarrow k_c$: the scalar denominator in
\eqref{stationary_merger_denominator} also tends to zero, while its
numerator remains positive by \eqref{stationary_merger_numerator}.\\
Choose nonnegative smooth bumps supported in small balls
about $k_1^0,k_2^0$, positive near their centers, and add their reflections
about the origin to obtain a real even $\xi$. Among sums of two centers,
only $k_1^0+k_2^0$ and its permutation give $k_c e_1$: a reflected center
lowers the first coordinate, and $2k_j^0$ has nonzero second coordinate.
Shrink the supports so only the two original cross terms contribute near
$k_c e_1$ and the stable-side scalar numerator is uniformly positive on both. At fixed summed
frequency $k$, the essential resonance is regular, since
$$
 \left.\nabla_{k_1}\Lambda_2(k_1,k-k_1)
       \right|_{k=k_ce_1,\ k_1=k_1^0}
 =-\frac{\lambda'(\kappa_*)}{\kappa_*}(k_1^0-k_2^0)\ne0.
$$
Replace the $e_2$ component of $k_1$ by $\Lambda_2$. The nonzero
transverse derivative gives a Jacobian bounded above and below on
both cross terms, uniformly for $k$ near $k_c e_1$. Shrinking the
chart, fix a tangential cube of positive measure and an interval about
$1/\tau$ on whose product $\xi(k_1)\xi(k-k_1)\geq c>0$ near the
first pair. These choices work for $\omega$ in a fixed open patch
$\Gamma$ about $e_1$ and $0<\eta<\eta_0$.
Set $v_{\parallel}=\omega\cdot v$ and restrict to $\eta\leq v_{\parallel}\le2\eta$.
After division by $\phi(v)$, the density contribution from the
stable region has nonnegative real part. Restricting its frequency
integral further to $\eta\le1-\tau\Lambda_2\le2\eta$ gives
$$
 \begin{aligned}
 \Re\frac{\hat f_2(k,v)}{\phi(v)}
 &\geq  c\int_{(1-2\eta)/\tau}^{(1-\eta)/\tau}
 \frac{(1-\tau\Lambda_2)\,d\Lambda_2}
 {(\eta+1-\tau\Lambda_2)
       \bigl((1-\tau\Lambda_2)^2+\tau^2|k|^2v_{\parallel}^2\bigr)}
       -C(1+|\log\eta|)\\
 &\geq \frac{c'}{\eta}-C(1+|\log\eta|).
 \end{aligned}
$$
The error bounds the direct forcing contribution on both sides and
the density contribution from the unstable region. Indeed,
$|2\phi-\hat e(k_1)-\hat e(k_2)|\lesssim\phi$ by the cutoff,
the scalar integral in \eqref{stationary_quadratic_full_resolvent}
is bounded, and $\Delta_k(-\Lambda_2)>1$ on the unstable side.
Each of these terms is bounded, after division by $\phi$, by
$$
 C\int_{1/\tau-\epsilon}^{1/\tau+\epsilon}
 \frac{d\Lambda_2}
 {\sqrt{(1-\tau\Lambda_2)^2+\tau^2|k|^2v_{\parallel}^2}}
 \leq C(1+|\log\eta|).
$$
They therefore cannot cancel the positive contribution of size
$\eta^{-1}$.\\
Restricting also to $|v-v_{\parallel}\omega|\le1$ gives a velocity set of
volume proportional to $\eta$, on which $\phi$ is bounded below.
For all sufficiently small $\eta>0$, we obtain
\begin{equation}
 \|\hat f_2((k_c+\eta/\tau)\omega,\cdot)\|_{L_v^2(\phi^{-1})}^2
 \gtrsim\eta^{-1},\qquad \omega\in\Gamma.
\end{equation}
The radial Jacobian $(k_c+\eta/\tau)^{d-1}/\tau$ is bounded below,
so integration in $\eta$ proves
\eqref{stationary_quadratic_reconstruction}.
On every compact subannulus of $|k|>k_c$ the stable-side scalar denominator is bounded away from zero. In the two independent leaf variables the gradient of $\Lambda_2$ is nonzero on the bump supports. Thus $|1-\tau\Lambda_2|^{-1/2}$ is locally integrable there. Fubini's theorem in the variables $(k_1,k_1+k_2)$ gives the ordinary Bochner convolution for almost every such $k$. For $\omega\in\Gamma$ and each fixed $0<\eta<\eta_0$, this convolution exists by the transverse change of variables above. The divergence comes from $\int_0^{\eta_0}d\eta/\eta$. Multiplication
of $\xi$ by any positive scalar preserves it and gives the claimed
arbitrarily small nonnegative even amplitudes.
\end{proof}

The unsteady homological equation, however, avoids precisely this inversion. For example,
at order two with zero homological initial condition, its solution is
$$
 \hat{\mathfrak W}_2(k_1,k_2,\cdot,\theta)
 =\frac{\alpha(2\pi)^{-d/2}}{\tau}
 \int_0^\theta e^{s\mathfrak L_{\bm k_2}}
       \bigl(2\phi-\hat e(k_1)-\hat e(k_2)\bigr)\,ds.
$$
In the stable region, letting $\theta\to\infty$ recovers the
stationary solution. Keeping $\theta$ finite instead gives a
well-defined solution on the stable, mixed and unstable regions,
including both resonance sets, by \eqref{boundnormsemigroup_uniform}.
There is no division by the scalar denominator in
\eqref{stationary_merger_denominator}. At an internal resonance
the projected forcing can contribute a term linear in $\theta$, while
at an essential resonance the finite-time semigroup integral remains
defined even though the stationary inverse does not. The task is
then to control the time-dependent coefficients together with the
parameter functions. Section~\ref{sec_proof_main_theorem} achieves
this through the selected homological initial conditions and the
clock estimates. Since $\theta[\xi(t)]=\theta[\xi(0)]+t$, this time
dependence is expressed by the reduced state itself and yields an
autonomous invariant parametrization.

\begin{remark}[Why Parabolic Similarity Variables do not Remove the Continuous Spectrum]
For parabolic equations on unbounded domains, self-similar variables
can produce an autonomous operator with a confining drift and a
spectrum suitable for invariant-manifold arguments in weighted
spaces \cite{wayne1997invariant,gallay2002invariant}, see also \cite{GallayRaugel1998} where the authors also use scaling variables to identify the parabolic terms in the asymptotic expansion of a damped wave equation. The same change
of variables does not preserve the kinetic equation considered here.
Indeed, substituting
$$
 f(x,v,t)=(1+t)^{-\beta}g(y,v,s),\qquad
 y=(1+t)^{-1/2}x,\qquad s=\log(1+t),
$$
into \eqref{maineqf} gives
$$
 \partial_s g-\frac12y\cdot\nabla_y g-\beta g
 +e^{s/2}v\cdot\nabla_y g
 =\frac{e^s}{\tau}
   \bigl(1+e^{-\beta s}\rho[g]\bigr)^\alpha\mathcal L_Q[g].
$$
No choice of $\beta$ removes the different factors $e^{s/2}$ and
$e^s$ multiplying transport and collisions. Scaling velocity as well
would change the prescribed Maxwellian, the collision projection
and the velocity spaces used in the spectral analysis. Thus this
parabolic rescaling does not turn the kinetic equation into an
autonomous confining problem, and it does not remove the need to
treat the resonances of the kinetic homological operator directly.
\end{remark}

\section{The General Solution to the Unsteady Homological Equation and Convergence of the Invariant Parametrization}\label{sec_proof_main_theorem}

In this section we solve the unsteady homological  equation \eqref{unsteady_cohomology_equ} order by order, choosing the initial data in \eqref{Vn_selected_choice} to cancel the hydrodynamic transient wherever it would decay too slowly, leaving a tail integral over future clock times. We then establish a uniform weighted homological  estimate, sum the resulting majorant, and prove differentiability at equilibrium separately from $C^1$ regularity away from it. Throughout, \eqref{ass_phi_w} is assumed, including its $p=1$ condition. Each input frequency satisfies $|k_j|\leq k_c-\varepsilon_B$.\\
For fixed $\bm k_n$, the general mild solution at order $n$, once the lower orders have been specified, is given by 
\begin{equation}\label{sol_unsteady_cohomol_W}
 \hat{\mathfrak W}_n(\theta)
 =e^{\theta\mathfrak L_{\bm k_n}}\hat{\mathfrak V}_n
 +\int_0^\theta e^{(\theta-s)\mathfrak L_{\bm k_n}}
       \widehat{\mathfrak N}_n(s)\,ds.
\end{equation}
Here $\widehat{\mathfrak N}_n$ is as defined in \eqref{hat_N_n_explicit}, and $\hat{\mathfrak V}_n\in\mathcal X_v$ is a selected homological initial condition to be specified at each order. The integral is a Bochner integral for locally integrable forcing, see \cite[Section~VI.7]{engel2000one}.\\
The natural exponential growth rate, derived from the algebraic structure of the homological  equation, is given by $\Lambda_n=-\sum_j\lambda(k_j)$. We thus define the growth-rate corrected coefficient functions and forcings as 
\begin{equation}\label{def_Un_renormalized}
 \hat{\mathfrak U}_n(\bm k_n,\theta)
 :=e^{-\Lambda_n\theta}\hat{\mathfrak W}_n(\bm k_n,\theta),
 \qquad
 \widehat{\mathfrak n}_n:=e^{-\Lambda_n\theta}\widehat{\mathfrak N}_n.
\end{equation}
Every summand of \eqref{hat_N_n_explicit} partitions the same leaf set $\{k_1,\ldots,k_n\}$ into disjoint blocks. If $B_1,\ldots,B_\ell$ are the Bell blocks and $B_0$ is the block entering the factor $\mathcal L_Q\hat{\mathfrak U}_{|B_0|}$, then
$$
 \Lambda_{|B_0|}((k_i)_{i\in B_0})+
 \sum_{r=1}^{\ell}\Lambda_{|B_r|}((k_i)_{i\in B_r})
 =-\sum_{j=1}^n\lambda(k_j)=\Lambda_n(\bm k_n).
$$ Since the corresponding $\Lambda$'s add to $\Lambda_n$, one obtains 
\begin{equation}\label{eq:renormalized_nonlinearity}
 \widehat{\mathfrak n}_n
 =\hat{\mathcal C}_n[\hat{\mathfrak U}_1,\ldots,\hat{\mathfrak U}_{n-1}],
 \qquad \rho[\widehat{\mathfrak n}_n]=0.
\end{equation}
Thus the original and renormalized forcings are distinguished in the following and the factors $(2\pi)^{-dl/2}$ in \eqref{hat_N_n_explicit} are unchanged by this renormalization. 
The equations and their mild solutions read in the growth-rate corrected variables:
\begin{align}
 \partial_\theta\hat{\mathfrak U}_n
 &=\mathcal L_{\Sigma\bm k_n}\hat{\mathfrak U}_n+\widehat{\mathfrak n}_n,
 \label{equ_U}\\
 \hat{\mathfrak U}_n(\theta)
 &={e^{\theta\mathcal L_{\Sigma\bm k_n}}}\hat{\mathfrak V}_n
   +\int_0^\theta {e^{(\theta-s)\mathcal L_{\Sigma\bm k_n}}}\widehat{\mathfrak n}_n(s)\,ds.
 \label{sol_unsteady_cohomol}
\end{align}

\subsection{Selected Homological Initial Conditions in the Duhamel Formula}

\begin{figure}[tb]
\centering
\includegraphics[width=\textwidth]{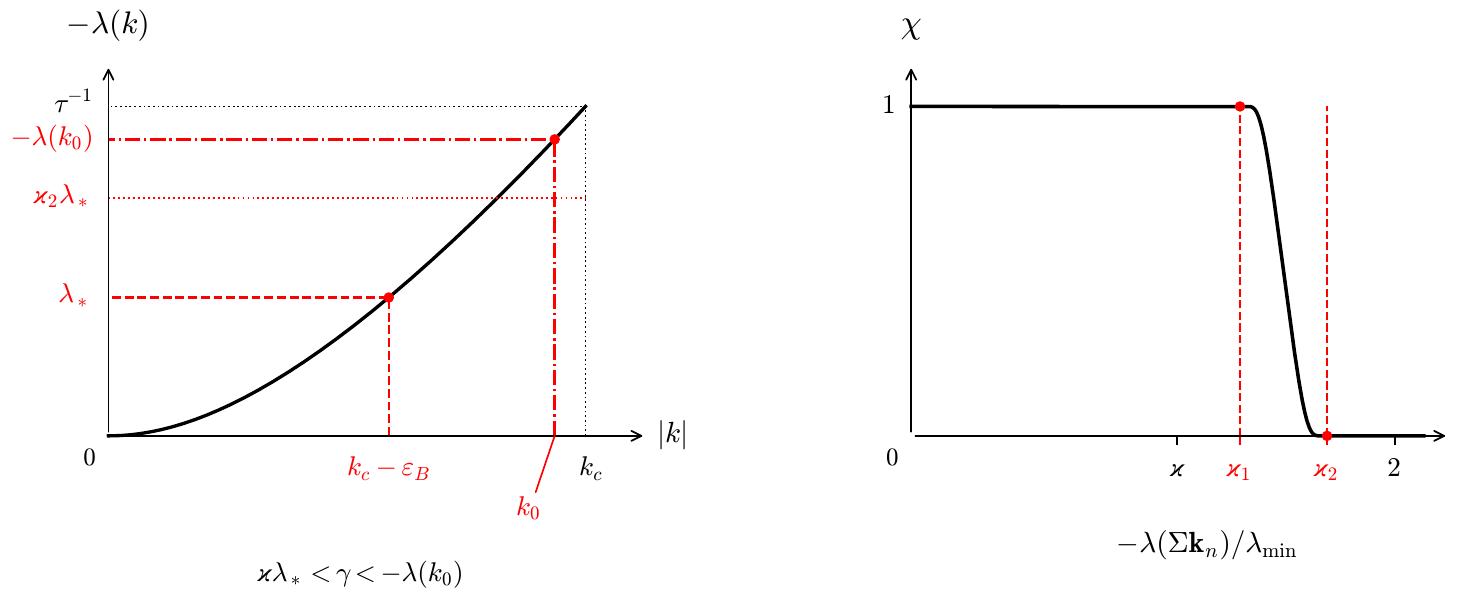}
\caption{Illustration of an admissible choice of the cutoff parameters on the individual and summed frequencies (left) and the smooth selection cutoff for $\lambda_{\min}>0$ (right). The rate $\gamma$ is chosen between
$\varkappa\lambda_*$ and $-\lambda(k_0)$.
\label{fig_selection_cutoffs}}
\end{figure}

To select the initial data in Duhamel's formula \eqref{sol_unsteady_cohomol_W}, we distinguish the following frequency regimes. We set 
\begin{equation}\label{def_lambda_min}
 \lambda_*:=-\lambda(k_c-\varepsilon_B)>0,
 \qquad \lambda_{\min}(\bm{k}_n) :=\min_{1\leq j\leq n}\{-\lambda(k_j)\}\geq 0,
\end{equation}
with $\lambda_*<1/\tau$ and fix constants
\begin{equation}\label{selection_rates}
1<\varkappa<\varkappa_1<\varkappa_2<\min\{2,(\tau\lambda_*)^{-1}\}.
\end{equation}
By Proposition~\ref{lemma_lambda}, $\tau\lambda_*\uparrow1$
as $\varepsilon_B\downarrow 0$ implying that, for each fixed $\varepsilon_B>0$, the
rates in \eqref{selection_rates} can be chosen freely, but must approach one in the limit. In particular, a fixed $\varkappa_2>1$ requires
$\tau\lambda_*<1/\varkappa_2$ and a strict gap cannot be avoided. We then choose a frequency $k_0<k_c$ so that $-\lambda(k_0)>\varkappa_2\lambda_*$,
and fix $\gamma$ with
$\varkappa\lambda_*<\gamma<-\lambda(k_0)$.
Let then $\chi\in C^\infty([0,\infty);[0,1])$ be a smooth cutoff equal to one on
$[0,\varkappa_1]$ and zero on
$[\varkappa_2,\infty)$, see Figure~\ref{fig_selection_cutoffs} for an illustration of the different parameters.\\
Let us now construct the initial condition as indicated.  For $0<|\Sigma\bm k_n|<k_0$, once the proper lower orders have been
constructed, prescribe the \emph{ homological initial condition}
\begin{equation}\label{Vn_selected_choice}
 \hat{\mathfrak V}_n(\bm k_n)
 =-\chi\!\left(\frac{-\lambda(\Sigma\bm k_n)}{\lambda_{\min}}\right)
       \int_0^\infty e^{-\lambda(\Sigma\bm k_n)s}
       P_{\Sigma\bm k_n}\widehat{\mathfrak n}_n(\bm k_n,s)\,ds,
 \qquad n\geq 2,
\end{equation}
whenever $\lambda_{\min}>0$ and the cutoff is nonzero and set the selected homological initial condition to zero otherwise. Where the cutoff is nonzero, the \emph{tail} selected by \eqref{Vn_selected_choice} is the contribution of the slow-mode forcing from clock time $\theta$ onward: 
\begin{equation}
\underbrace{-\chi\!\left(\frac{-\lambda(\Sigma\bm k_n)}{\lambda_{\min}}\right)
  \int_\theta^\infty e^{\lambda(\Sigma\bm k_n)(\theta-s)}
  P_K\widehat{\mathfrak n}_n(\bm k_n,s)\,ds}_{\text{tail contribution at }\theta},
\qquad
2\lambda_{\min}+\lambda(\Sigma\bm k_n)>(2-\varkappa_2)\lambda_{\min}>0. 
\end{equation}

\begin{remark}
On the individual frequencies, the eigenvectors are uniformly bounded and each individual decay rate is at most $\lambda_*$. The cutoff
$k_c-\varepsilon_B<k_0<k_c$ selects where to use the slow projection.
For $n\geq2$ our target in \eqref{prop_unsteady_cohomology_bound_eq} is
a decay at rate $\varkappa \lambda_{\min}(\bm k_n)$. The strict gain $\varkappa>1$ is used for
tangency at equilibrium. For $0<|\Sigma\bm k_n|<k_0$ and $\lambda_{\min}(\bm k_n)>0$, we use the tail contribution of the initial condition when $-\lambda(\Sigma\bm k_n)\leq\varkappa_1 \lambda_{\min}(\bm k_n)$ and forward evolution when
$-\lambda(\Sigma\bm k_n)\geq\varkappa_2 \lambda_{\min}(\bm k_n)$. In the transition region both estimates hold, since
$-\lambda(\Sigma\bm k_n)>\varkappa \lambda_{\min}(\bm k_n)$ for the forward term and $-\lambda(\Sigma\bm k_n)<2\lambda_{\min}(\bm k_n)$ for the tail.
The complementary modes and all summed frequencies for which $|\Sigma\bm k_n|\geq k_0$
decay forward at rate $\gamma>\varkappa\lambda_*$ by
Lemma~\ref{lem_selected_semigroup}. At $\lambda_{\min}(\bm k_n)=0$, we just aim for boundedness and thus prescribe zero initial data. 
\end{remark}

\begin{remark}
Where $\chi=1$, the tail in \eqref{Vn_selected_choice} imposes the Lyapunov--Perron terminal condition $e^{-\lambda(K)\theta}P_K\hat{\mathfrak U}_n(\theta)\to0$ as $\theta\to\infty$, with $K=\Sigma\bm k_n$ (cf.~\cite[Lemma~2.2, Eq.~(2.2)]{Gallay1993CenterStable}). A related choice of initial data appears in \cite[Lemma~23 and Eq.~(114)]{Schlag2009} for the
cubic nonlinear Schr\"odinger equation, where integration of the projected forcing over future times cancels an exponentially growing mode.
\end{remark}

The choice \eqref{Vn_selected_choice} uniquely determines the slow component, while the cutoff interpolates with forward evolution where both choices give the required decay. We refer to the conditions $\hat{\mathfrak V}_n$ prescribed by
\eqref{Vn_selected_choice} at successive orders as the
\emph{homological initial conditions} for the unsteady hierarchy \eqref{unsteady_cohomology_equ}.
Note that \eqref{Vn_selected_choice} is always zero at $\Sigma\bm k_n=0$, since
$P_0\widehat{\mathfrak n}_n=\rho[\widehat{\mathfrak n}_n]\phi=0$.
At order one, we retain
$\hat{\mathfrak U}_1(k,\theta)=e^{\lambda(k)\theta}\hat e(k)$.
The selected homological initial condition is symmetric in the individual frequencies and is independent
of the parameter $\xi$ and $\mathfrak r$. It cancels
sufficiently slow hydrodynamic transients, but does not invert the
stationary kinetic operator. The integral in
\eqref{Vn_selected_choice} will converge by the induction below.

\begin{lemma}[Semigroup bounds at the selected rate]\label{lem_selected_semigroup}
There exists a constant $C_{\rm sel}$ such that, for all $k\in\mathbb R^d$ and $t\geq 0$,
$$
 \begin{aligned}
 \|e^{t\mathcal L_k}(I-P_k)\|_{\mathcal L(\mathcal X_v^{(2)})}
 &\leq C_{\rm sel}e^{-\gamma t},&&|k|\leq k_0,\\
 \|e^{t\mathcal L_k}\|_{\mathcal L(\mathcal X_v^{(2)})}
 &\leq C_{\rm sel}e^{-\gamma t},&&|k|\geq  k_0.
 \end{aligned}
$$
\end{lemma}
\begin{proof}
All operator norms are taken on $\mathcal X_v^{(2)}$.
The first estimate follows from \eqref{weighted_semigroup_rem} in
Proposition~\ref{lem_weighted_semigroup_aux}, with $\kappa=k_0$ and
$\omega=\gamma$. These choices are admissible since
$0<\gamma<-\lambda(k_0)<\tau^{-1}$.\\
For $|k|\geq  k_0$, monotonicity and \eqref{lambda_extended} give
$\widehat\lambda(k)\leq\lambda(k_0)$. Applying
\eqref{boundnormsemigroup_uniform} on the weighted space, as in the
proof of Proposition~\ref{lem_weighted_semigroup_aux}, therefore yields
$$
 \|e^{t\mathcal L_k}\|_{\rm op}
 \leq C_{\mathcal L}^{\rm wgt}(1+t/\tau)^2e^{\lambda(k_0)t}.
$$
Set
$$
 \delta=-\lambda(k_0)-\gamma>0,\qquad
 C_{\rm poly}(\delta)=\sup_{s\geq 0}(1+s/\tau)^2e^{-\delta s}<\infty.
$$
The preceding bound is at most
$C_{\mathcal L}^{\rm wgt}C_{\rm poly}(\delta)e^{-\gamma t}$.
Thus both assertions hold with
$$
 C_{\rm sel}=\max\{C_{\rm rem}(k_0,\gamma),
                    C_{\mathcal L}^{\rm wgt}C_{\rm poly}(\delta)\}.
$$
This estimate remains uniform through $|k|=k_c$ because
\eqref{boundnormsemigroup_uniform} holds for all frequencies.
\end{proof}

To estimate the contributions of the nonlinear functional expansion, we define
\begin{equation}\label{weighted_majorant_def}
 M_n(\bm k_n)=\sup_{\theta\geq 0}e^{2\lambda_{\min}\theta}
       \|\widehat{\mathfrak n}_n(\bm k_n,\theta)\|_{\mathcal X_v^{(2)}},\qquad n\geq 2,
\end{equation}
and the scalar function $E:\mathbb{R}^d\to\mathbb{R}^+$
\begin{equation}\label{def_Q}
 E(x)
 =1+\frac{\mathbf1_{\{0<|x|<k_0\}}}{|x|}
 \quad(x\ne0),\qquad E(0)=1.
\end{equation}
The factor $1/|x|$ records the low-frequency cost of the slow projection.
On zero-density forcing the projection gains $|x|$, while integration
at the slow rate can cost $1/(-\lambda(x))=O(|x|^{-2})$. The resulting
singularity is locally square integrable precisely for $d\geq3$. We use this integrability below to estimate the convolution of the input profiles. 

\begin{proposition}[Weighted regularity and decay of the selected hierarchy]
\label{prop_weighted_hierarchy}
The recursive selection of homological initial conditions of the form \eqref{Vn_selected_choice}
is well-defined. For each fixed vector of individual frequencies, the coefficient kernels
and the forcing are continuous in $\theta$ with values in
$\mathcal X_v^{(2)}$, the forcing has zero density,
\begin{equation}\label{rho_N_zero_sec8}
 \rho[\widehat{\mathfrak n}_n(\bm k_n,\theta)]=0,
\end{equation}
and $M_n(\bm k_n)<\infty$. There is a constant $C_{\rm maj}$,
independent of the order $n$ and the vectors of individual frequencies, such that
\begin{equation}\label{prop_unsteady_cohomology_bound_eq}
 \sup_{\theta\geq 0}e^{\varkappa \lambda_{\min}\theta}
       \|\hat{\mathfrak U}_n(\bm k_n,\theta)\|_{\mathcal X_v^{(2)}}
 \leq C_{\rm maj}E(\Sigma\bm k_n)M_n(\bm k_n),\qquad n\geq 2.
\end{equation}
The same estimate holds with the $\mathcal X_v$ norm on the left,
with the same constant and the same weighted forcing majorant
$M_n$ from \eqref{weighted_majorant_def}.  The constant depends only on the fixed cutoffs on the individual and summed frequencies, velocity space, $\tau$, and selection parameters in
\eqref{selection_rates}, together with the chosen rate $\gamma$, while
it is independent of $\xi$, $\nu$, and $\mathfrak r$.
\end{proposition}

\begin{proof}

We argue by induction on the order, first estimating the forcing
and then the hydrodynamic and complementary components of the solution.
All individual frequencies lie in $B_c(\varepsilon_B)$, so
$0\leq\lambda_{\min}\leq\lambda_*$ by \eqref{def_lambda_min}.
Frequency arguments are suppressed only for the current fixed vector.\\

\noindent
\underline{The forcing estimate}:\\
At order one,
$\hat{\mathfrak U}_1(k,s)=e^{\lambda(k)s}\hat e(k)$ is continuous in
$\mathcal X_v^{(2)}$, and \eqref{eigenvector} gives
$$
 \sup_{s\geq 0}e^{-\lambda(k)s}
       \|\hat{\mathfrak U}_1(k,s)\|_{\mathcal X_v^{(2)}}
 =\|\hat e(k)\|_{\mathcal X_v^{(2)}}
 \leq\frac{\|\phi\|_{\mathcal X_v^{(2)}}}{1-\tau\lambda_*}<\infty.
$$
For each constructed coefficient, define its scalar majorant by
\begin{equation}\label{block_amplitude_constant}
 \begin{aligned}
 U_1(k)&:=\|\hat e(k)\|_{\mathcal X_v^{(2)}},\\
 U_n(\bm k_n)&:=\sup_{s\geq 0}
 e^{\varkappa\lambda_{\min}(\bm k_n)s}
 \|\hat{\mathfrak U}_n(\bm k_n,s)\|_{\mathcal X_v^{(2)}},
 \qquad n\geq 2.
 \end{aligned}
\end{equation}
Suppose the assertions hold at every order below $n\geq 2$, so that
the corresponding $U_{|I|}$ are finite for every proper nonempty block
$I$ of frequencies. Expanding
\eqref{hat_N_n_explicit} by the partition formula
\eqref{freq_Bell_poly} and using
\eqref{eq:renormalized_nonlinearity}, each term of the forcing is a
scalar multiple of
$$
 \left(\prod_{r=1}^{l}
   \rho\!\left[\hat{\mathfrak U}_{|I_r|}
      ((k_i)_{i\in I_r},s)\right]\right)
 \mathcal L_Q\!\left[\hat{\mathfrak U}_{|I_0|}
      ((k_i)_{i\in I_0},s)\right].
$$
Here $I_0,I_1,\ldots,I_l$ partition $\{1,\ldots,n\}$ into
nonempty blocks, with $l\geq 1$. The block $I_0$ enters the
collision factor and the others enter density factors. There are
therefore at least two blocks, all of size less than $n$, so every
factor has already been constructed.
For each such block,
$$
 \lambda_{\min}((k_i)_{i\in I})
 =\min_{i\in I}\{-\lambda(k_i)\}
 \geq \lambda_{\min}(\bm k_n),
$$
and \eqref{block_amplitude_constant}, together with $\varkappa>1$,
gives
$$
 \|\hat{\mathfrak U}_{|I|}((k_i)_{i\in I},s)\|_{\mathcal X_v^{(2)}}
 \leq U_{|I|}((k_i)_{i\in I})e^{-\lambda_{\min}((k_i)_{i\in I})s}
 \leq U_{|I|}((k_i)_{i\in I})e^{-\lambda_{\min}s}.
$$
For $|I|\geq 2$, the inductive estimate is
$$
 U_{|I|}((k_i)_{i\in I})
 \leq C_{\rm maj}E\!\left(\sum_{i\in I}k_i\right)
             M_{|I|}((k_i)_{i\in I}).
$$
Use $C_\rho$ from \eqref{defCrho}. The projection bound
\eqref{Prho_bound}, applied on $\mathcal X_v^{(2)}$, and \eqref{def_LQ}
give the collision bound with $C_Q=1+C_\rho\|\phi\|_{\mathcal X_v^{(2)}}$.
Thus the norm of the displayed block product is at most
$$
 C_QC_\rho^l\left(\prod_{r=0}^{l}U_{|I_r|}((k_i)_{i\in I_r})\right)
 e^{-(l+1)\lambda_{\min}s}
 \leq C_QC_\rho^l\left(\prod_{r=0}^{l}U_{|I_r|}((k_i)_{i\in I_r})\right)
       e^{-2\lambda_{\min}s}.
$$
Summing with the absolute values of the scalar coefficients in
\eqref{hat_N_n_explicit} proves $M_n(\bm k_n)<\infty$ and hence
$$
 \|\widehat{\mathfrak n}_n(\bm k_n,s)\|_{\mathcal X_v^{(2)}}
 \leq M_n(\bm k_n)e^{-2\lambda_{\min}s}.
$$
All order and frequency dependence remains in $M_n$.
The finite product formula also proves continuity in time, and
$\rho\mathcal L_Q=0$, as in Remark~\ref{remark_rho_N_n}, gives
\eqref{rho_N_zero_sec8} term by term.\\

\noindent
\underline{The hydrodynamic component}:\\
Let $0<|\Sigma\bm k_n|<k_0$. For $\lambda_{\min}>0$, write
$$
 \chi_n=\chi\!\left(
          \frac{-\lambda(\Sigma\bm k_n)}{\lambda_{\min}}\right),
$$
and set $\chi_n=0$ when $\lambda_{\min}=0$, without evaluating the ratio.
The cutoff conditions preceding \eqref{Vn_selected_choice} give
$$
 \begin{gathered}
 0\leq\chi_n\le1,\\
 \chi_n\ne0\ \Longrightarrow\
       -\lambda(\Sigma\bm k_n)<\varkappa_2\lambda_{\min},\\
 1-\chi_n\ne0\ \Longrightarrow\
       -\lambda(\Sigma\bm k_n)\geq \varkappa_1\lambda_{\min}.
 \end{gathered}
$$
On the support of $\chi_n$, the integral defining
$\hat{\mathfrak V}_n$ converges absolutely, since, by the zero-density
estimate \eqref{Pk_zero_mass_gain}, we have that 
$$
 \begin{aligned}
 \int_0^\infty e^{-\lambda(\Sigma\bm k_n)s}
       \|P_{\Sigma\bm k_n}\widehat{\mathfrak n}_n(s)
          \|_{\mathcal X_v^{(2)}}\,ds
 &\leq C_P|\Sigma\bm k_n|M_n
       \int_0^\infty
       e^{-[2\lambda_{\min}+\lambda(\Sigma\bm k_n)]s}\,ds\\
 &\quad=\frac{C_P|\Sigma\bm k_n|M_n}
                   {2\lambda_{\min}+\lambda(\Sigma\bm k_n)}<\infty,
 \end{aligned}
$$
because
$$
 2\lambda_{\min}+\lambda(\Sigma\bm k_n)
 >(2-\varkappa_2)\lambda_{\min}>0.
$$
Since the integrand lies in the closed range of $P_{\Sigma\bm k_n}$,
so does $\hat{\mathfrak V}_n$ and on this range the semigroup acts by
$$
 e^{t\mathcal L_{\Sigma\bm k_n}}P_{\Sigma\bm k_n}
 =e^{\lambda(\Sigma\bm k_n)t}P_{\Sigma\bm k_n}.
$$
Substituting the homological initial condition \eqref{Vn_selected_choice} into \eqref{sol_unsteady_cohomol}
and splitting its integral at $\theta$ therefore gives
\begin{equation}\label{selected_green_formula}
 \begin{aligned}
 P_{\Sigma\bm k_n}\hat{\mathfrak U}_n(\theta)
 &=(1-\chi_n)\int_0^\theta
       e^{\lambda(\Sigma\bm k_n)(\theta-s)}
       P_{\Sigma\bm k_n}\widehat{\mathfrak n}_n(s)\,ds-\chi_n\int_\theta^\infty
       e^{\lambda(\Sigma\bm k_n)(\theta-s)}
       P_{\Sigma\bm k_n}\widehat{\mathfrak n}_n(s)\,ds.
 \end{aligned}
\end{equation}
The second term in the above splitting is defined to be zero when $\chi_n=0$ and its
unmultiplied integral is used only where it converges. In particular,
the formula never requires multiplying a divergent integral by zero.
Where $\chi_n=1$ only the second integral remains. Where $1-\chi_n\ne0$, the cutoff gives
$-\lambda(\Sigma\bm k_n)\geq\varkappa_1\lambda_{\min}$, so \eqref{lambda_quadratic_bounds} yields
$$
 -\lambda(\Sigma\bm k_n)-\varkappa\lambda_{\min}
 \geq \left(1-\frac{\varkappa}{\varkappa_1}\right)[-\lambda(\Sigma\bm k_n)]
 \geq \frac{2\tau}{\pi}\left(1-\frac{\varkappa}{\varkappa_1}\right)|\Sigma\bm k_n|^2.
$$
By the zero-density estimate \eqref{Pk_zero_mass_gain} and the definition
of $M_n$ in \eqref{weighted_majorant_def}, the forward integral in
\eqref{selected_green_formula} therefore satisfies, for every $\theta\geq0$,
$$
 \begin{aligned}
 e^{\varkappa\lambda_{\min}\theta}
 \left\|\int_0^\theta e^{\lambda(\Sigma\bm k_n)(\theta-s)}
       P_K\widehat{\mathfrak n}_n(s)\,ds\right\|_{\mathcal X_v^{(2)}} & \leq C_P|\Sigma\bm k_n|M_n\int_0^\theta
       e^{[\lambda(\Sigma\bm k_n)+\varkappa\lambda_{\min}](\theta-s)}
       e^{-(2-\varkappa)\lambda_{\min}s}\,ds\\
 &\leq C_P|\Sigma\bm k_n|M_n
       \frac{1-e^{[\lambda(\Sigma\bm k_n)+\varkappa\lambda_{\min}]\theta}}
            {-\lambda(\Sigma\bm k_n)-\varkappa\lambda_{\min}}
 \leq \frac{\pi C_P}{2\tau(1-\varkappa/\varkappa_1)}
       \frac{M_n}{|\Sigma\bm k_n|},
 \end{aligned}
$$
where we used $2-\varkappa>0$ in the second inequality.
Where $\chi_n\ne0$, we have $-\lambda(\Sigma\bm k_n)<\varkappa_2\lambda_{\min}$,
and hence
$$
 2\lambda_{\min}+\lambda(\Sigma\bm k_n)
 \geq \frac{2-\varkappa_2}{\varkappa_2}[-\lambda(\Sigma\bm k_n)]
 \geq \frac{2\tau(2-\varkappa_2)}{\pi\varkappa_2}|\Sigma\bm k_n|^2.
$$
For the tail integral, the same projection estimate and the substitution
$u=s-\theta$ give
$$
 \begin{aligned}
 e^{\varkappa\lambda_{\min}\theta}
 \left\|\int_\theta^\infty e^{\lambda(\Sigma\bm k_n)(\theta-s)}
       P_{\Sigma\bm k_n}\widehat{\mathfrak n}_n(s)\,ds\right\|_{\mathcal X_v^{(2)}}
 &\leq C_P|\Sigma\bm k_n|M_n e^{-(2-\varkappa)\lambda_{\min}\theta}
       \int_0^\infty e^{-[2\lambda_{\min}+\lambda(\Sigma\bm k_n)]u}\,du\\
 &= C_P|\Sigma\bm k_n|M_n
       \frac{e^{-(2-\varkappa)\lambda_{\min}\theta}}
            {2\lambda_{\min}+\lambda(\Sigma\bm k_n)}
 \leq \frac{\pi\varkappa_2 C_P}{2\tau(2-\varkappa_2)}
       \frac{M_n}{|\Sigma\bm k_n|}.
 \end{aligned}
$$
The forward estimate also covers $\lambda_{\min}=0$, since
then $\chi_n=0$ and its denominator is $-\lambda(\Sigma\bm k_n)>0$.
The above estimates identify the source of the singular part of
$E(\Sigma\bm k_n)$: integration in time costs a factor bounded by
a constant times $|\Sigma\bm k_n|^{-2}$, whereas the zero-density
projection supplies one factor $|\Sigma\bm k_n|$, which, in their product, gives
$|\Sigma\bm k_n|^{-1}$.\\
 
\noindent
\underline{The complementary component}:\\
For $0<|\Sigma\bm k_n|<k_0$, the preceding construction gives
$(I-P_{\Sigma\bm k_n})\hat{\mathfrak V}_n=0$, including when
$\lambda_{\min}=0$ or $\chi_n=0$. At $\Sigma\bm k_n=0$, the
prescription also gives $\hat{\mathfrak V}_n=0$.
Since the Riesz projection commutes with the semigroup, the
complementary component is given by
$$
 \begin{aligned}
 &(I-P_{\Sigma\bm k_n})\hat{\mathfrak U}_n(\theta) =\int_0^\theta e^{(\theta-s)\mathcal L_{\Sigma\bm k_n}}
          (I-P_{\Sigma\bm k_n})\widehat{\mathfrak n}_n(s)\,ds,
 \qquad |\Sigma\bm k_n|<k_0.
 \end{aligned}
$$
For $|\Sigma\bm k_n|\geq  k_0$, the homological initial condition \eqref{Vn_selected_choice} is zero and hence 
$$
 \hat{\mathfrak U}_n(\theta)
 =\int_0^\theta e^{(\theta-s)\mathcal L_{\Sigma\bm k_n}}
         \widehat{\mathfrak n}_n(s)\,ds,
 \qquad |\Sigma\bm k_n|\geq  k_0.
$$
We next apply, respectively, the complementary and full semigroup bounds
of Lemma~\ref{lem_selected_semigroup}, with common constant $C_{\rm sel}$. In both cases the weighted norm is bounded by
$$
 \begin{aligned}
 C_{\rm sel}M_n e^{\varkappa\lambda_{\min}\theta}
       \int_0^\theta e^{-\gamma(\theta-s)}
                       e^{-2\lambda_{\min}s}\,ds & = C_{\rm sel}M_n\int_0^\theta
       e^{-(\gamma-\varkappa\lambda_{\min})(\theta-s)}
       e^{-(2-\varkappa)\lambda_{\min}s}\,ds\\
 &\leq\frac{C_{\rm sel}M_n}{\gamma-\varkappa\lambda_{\min}}
       \leq\frac{C_{\rm sel}M_n}{\gamma-\varkappa\lambda_*},
 \end{aligned}
$$
This bound includes $\lambda_{\min}=0$ and the denominator is
positive by the choice of $\gamma$ following \eqref{selection_rates}.
At $\Sigma\bm k_n=0$, both $P_0\hat{\mathfrak V}_n$ and
$P_0\widehat{\mathfrak n}_n(s)=\rho[\widehat{\mathfrak n}_n(s)]\phi$
vanish. Thus the complementary estimate controls the entire solution
there, consistently with \eqref{semigroup_k_zero}, and $E(0)=1$ is valid.\\
For $0<|\Sigma\bm k_n|<k_0$, we add the bounds on the hydrodynamic and
complementary components. At $\Sigma\bm k_n=0$ and for
$|\Sigma\bm k_n|\geq  k_0$, only the complementary or full semigroup bound is needed. Thus we may take
$$
 C_{\rm maj}:=\max\left\{
 \frac{C_{\rm sel}}{\gamma-\varkappa\lambda_*},\quad
 \frac{\pi C_P}{2\tau}
 \left[\frac{1}{1-\varkappa/\varkappa_1}
             +\frac{\varkappa_2}{2-\varkappa_2}\right]\right\}.
$$
The second entry also covers $\lambda_{\min}=0$. The definition
\eqref{def_Q} now gives \eqref{prop_unsteady_cohomology_bound_eq},
with the same $C_{\rm maj}$ at every induction step.
The convergent integral in \eqref{Vn_selected_choice} shows that
$\hat{\mathfrak V}_n\in\mathcal X_v^{(2)}$ wherever the selection is
nonzero, while elsewhere, we can just set $\hat{\mathfrak V}_n=0$.
Note that the forcing is continuous in that space by the first part of the proof. Strong continuity
of $e^{\theta\mathcal L_{\Sigma\bm k_n}}$ and the finite-time Bochner
formula \eqref{sol_unsteady_cohomol} therefore give continuity of
$\hat{\mathfrak U}_n$ in $\theta$. Multiplication by
$e^{\Lambda_n\theta}$ also gives continuity of
$\hat{\mathfrak W}_n$ and $\widehat{\mathfrak N}_n$ at every finite
time. This completes the induction. The embedding following
\eqref{def_Xv2_aux} gives the unweighted estimate with the same
$C_{\rm maj}$ and $M_n$.
\end{proof}

\begin{remark}\label{rem_weighted_majorant_role}
The factor $|{\Sigma\bm k_n}|$ in the slow projection uses a velocity moment of the zero-density forcing, which cannot be obtained from its unweighted norm alone. The weighted hierarchy, however, supplies this moment automatically. Furthermore, for each fixed frequency vector, it also gives classical differentiability into $\mathcal X_v$: ${\mathcal X_v^{(2)}}\hookrightarrow D(\mathcal L_{\Sigma\bm k_n})$ continuously in the graph norm, since
\begin{equation}\label{C1_generator_weight_loss}
 \|\mathcal L_{\Sigma\bm k_n} f\|_{\mathcal X_v}
 \leq C_{\rm gen}(1+|{\Sigma\bm k_n}|)\|f\|_{\mathcal X_v^{(2)}}.
\end{equation}
Here one may take $C_{\rm gen}:=\max\{1,C_Q/\tau\}$, by \eqref{def_LQ} and the collision bound in the proof of Proposition~\ref{prop_weighted_hierarchy}.
The selected homological initial conditions \eqref{Vn_selected_choice} constructed in Proposition~\ref{prop_weighted_hierarchy} belong to $\mathcal X_v^{(2)}$, by the convergent integral in \eqref{Vn_selected_choice}, and hence to $D(\mathcal L_{\Sigma\bm k_n})$. Differentiating the Duhamel integral using this graph-norm continuity therefore proves \eqref{equ_U} in $\mathcal X_v$, with a continuous right-hand side. This domain argument is needed before differentiating with respect to the clock.
\end{remark}

\subsection{Convergence and Regularity of the Invariant Parametrization}\label{sec_proof_main_thm_proof}

Assume $d\geq 3$, $\theta_0=0$, and $\nu\geq 2\lambda_*$, with $\lambda_*$ as in \eqref{def_lambda_min}. For $\xi\in\Xi(\varepsilon_B)$, defined in \eqref{defXi_eps_C0}, with $0<\|\xi\|_{C^0}\leq\mathfrak r$, define the conjugate-symmetric compensated profile, extended by zero outside ${B_c(\varepsilon_B)}$, by
$$
 \zeta[\xi](k):=e^{-\lambda(k)\theta[\xi]}\xi(k)\mathbf1_{{B_c(\varepsilon_B)}}(k).
$$
We use $\mathcal Y$ with the sum norm in \eqref{def_majorant_space_Y}. The compensation estimate \eqref{bound_compensated_xi} of Lemma~\ref{lem_theta}, with $\sigma=1$, gives
\begin{equation}\label{clock_Y_bound_direct}
\begin{gathered}
 \|\zeta[\xi]\|_{\mathcal Y}\leq C_{\rm comp}^{\mathcal Y}\mathfrak r,
 \\
 C_{\rm comp}^{\mathcal Y}
 :=|B_c(\varepsilon_B)|+\mu_0^{-1}
   +\bigl(|B_c(\varepsilon_B)|+\mu_0^{-1}\bigr)^{1/2}.
\end{gathered}
\end{equation}
Indeed, the proof of that lemma, using \eqref{clock_unweighted_g_bound} and $\nu/(-\lambda(k))\geq 2$ for nonzero $k\in B_c(\varepsilon_B)$, gives
$\|\zeta[\xi]\|_{L_k^r}^r\leq\mathfrak r^r(|{B_c(\varepsilon_B)}|+\mu_0^{-1})$ for $r=1,2$, also when $\|\xi\|_{C^0}=\mathfrak r$.
Here $\mu_0$ is from \eqref{clock_weight_assumption}. In particular, $C_{\rm comp}^{\mathcal Y}=C_{\rm comp}^{L^1}+C_{\rm comp}^{L^2}$.\\
The majorants \eqref{block_amplitude_constant} give the pointwise bounds
\begin{equation}\label{selected_kernel_majorant}
 \begin{aligned}
 \|\hat{\mathfrak U}_1(k,\theta)\|_{\mathcal X_v^{(2)}}
 &=U_1(k)e^{\lambda(k)\theta},\\
 \|\hat{\mathfrak U}_n(\bm k_n,\theta)\|_{\mathcal X_v^{(2)}}
 &\leq U_n(\bm k_n)e^{-\varkappa\lambda_{\min}(\bm k_n)\theta},
 \qquad n\geq 2,\quad\theta\geq 0.
 \end{aligned}
\end{equation}
The order-one weight retains the exact linear decay rate $-\lambda(k)$, while the higher-order weight uses $\varkappa$ from \eqref{selection_rates}. Furthermore, for every nonnegative $z\in \mathcal Y$ supported in ${B_c(\varepsilon_B)}$, define the functional
\begin{equation}\label{def_integrated_majorant}
    \mathfrak u_n(q;z):=\frac1{n!}\int_{{B_c(\varepsilon_B)}^n}\delta(q-\Sigma\bm k_n)
 U_n(\bm k_n)
 \prod_{j=1}^nz(k_j)\,d\bm k_n.
\end{equation}
The estimates below retain the frequency dependence of $U_n$ inside
the convolution integrals.
The small-divisor estimate \eqref{small_divisor_convolution_bound}
controls the factors $E$ after integration, and the scalar majorant
\eqref{integrated_scalar_majorant_equation} controls the order
dependence of the profiles normalized by $1/n!$.

\begin{proposition}[Integrated majorant and convergence]\label{prop_majorant_functional_taylor}
Under the preceding hypotheses, there are $A_n\geq 0$ such that $A(s)=\sum_{n\geq 1}A_ns^n$ is analytic near zero and, for every nonnegative $z\in \mathcal Y$ supported in ${B_c(\varepsilon_B)}$,
\begin{equation}\label{weighted_profile_majorant_final}
\begin{aligned}
 \|\mathfrak u_n(\cdot;z)\|_{\mathcal Y}&\leq A_n\|z\|_{\mathcal Y}^n.
\end{aligned}
\end{equation}
The coefficients are independent of $\xi$, $\mathfrak r$ and $\nu$,
but depend on the fixed selection parameters.\\
There exists $\mathfrak r_*>0$ such that, for $0<\mathfrak r<\mathfrak r_*$,
\eqref{defhydro} converges absolutely and uniformly on
$0<\|\xi\|_{C^0}\leq\mathfrak r$ in $\mathcal Z^{(2)}$ from
\eqref{def_Xv2_aux}. The sum satisfies $\hat{\mathcal W}[\xi](-k,\cdot)=\hat{\mathcal W}[\xi](k,\cdot)^*$ for almost every $k$. Its inverse Fourier transform converges absolutely
in $L_x^\infty(\mathcal X_v)$, and also in $L_x^2(\mathcal X_v)$ if
$1\leq p\le2$.\\
The sum satisfies the mild invariance identity \eqref{prop82_mild_identity}
along every reduced orbit. The radius depends only on the fixed $d$,
$\varepsilon_B$, $\tau$, $k_0$, $\mathcal X_v$, $\alpha$, $\mu$, and
selection parameters in \eqref{selection_rates}.
\end{proposition}

\begin{proof}
We first establish the scalar majorant, then use it to prove
convergence and mild invariance. The physical-space conclusions follow by
Fourier inversion and an auxiliary Hilbert-space estimate.\\

\noindent
\underline{The convolution estimate}:\\
Recall $E$ from \eqref{def_Q}, with $E(0)=1$. Since $d\geq 3$,
\begin{equation}
 \|E-1\|_2^2=\frac{|\mathbb S^{d-1}|}{d-2}k_0^{d-2}<\infty.
\end{equation}
For $r\geq 2$ and $f_1,\ldots,f_r\in \mathcal Y$,
Young's inequality and Cauchy--Schwarz give
\begin{equation}\label{small_divisor_convolution_bound}
\begin{gathered}
 \|E(f_1*\cdots*f_r)\|_{\mathcal Y}
 \leq C_E\prod_{i=1}^r\|f_i\|_{\mathcal Y},\quad C_E=1+\|E-1\|_2.
\end{gathered}
\end{equation}
Indeed, with $g=f_1*\cdots*f_r$ and all norms taken over
frequency space,
$$
\begin{aligned}
 \|g\|_1+\|g\|_2
 &\le(\|f_1\|_1+\|f_1\|_2)\prod_{i=2}^r\|f_i\|_1,\\
 \|g\|_2+\|g\|_\infty
 &\leq\|f_1\|_2(\|f_2\|_1+\|f_2\|_2)
       \prod_{i=3}^r\|f_i\|_1,\\
 \|(E-1)g\|_1+\|(E-1)g\|_2
 &\leq\|E-1\|_2(\|g\|_2+\|g\|_\infty),
\end{aligned}
$$
where we used the convention that empty products equal one. Applying the triangle inequality to
$Eg=g+(E-1)g$ and bounding each individual norm by its sum norm
proves \eqref{small_divisor_convolution_bound}.
Note that $E-1\in \mathcal Y$,
whereas $E\notin \mathcal Y$.\\

\noindent
\underline{The integrated majorant}:\\
Use $C_\rho$ from \eqref{defCrho} and $C_Q$ from the forcing
estimate in the proof of Proposition~\ref{prop_weighted_hierarchy}, which bound
$\rho$ and $\mathcal L_Q$ on $\mathcal X_v^{(2)}$. The order-one estimate in that proof gives
$$
 \|\mathfrak u_1(\cdot;z)\|_{\mathcal Y}\leq A_1\|z\|_{\mathcal Y},\qquad
 A_1\geq \frac{\|\phi\|_{\mathcal X_v^{(2)}}}{1-\tau\lambda_*}.
$$

The coefficient kernels are strongly measurable in frequencies and continuous in time, by the multiplication semigroup, the forward recursion, and the absolutely convergent integrals in \eqref{Vn_selected_choice}. Measurability of the latter follows by truncating their integration intervals. The suprema in \eqref{block_amplitude_constant} can therefore be taken over rational times. Initially $\mathfrak u_n$ is a nonnegative extended integral and finiteness follows by induction below. \\
The forcing estimate in the proof of
Proposition~\ref{prop_weighted_hierarchy}, with the partition formula
\eqref{freq_Bell_poly} and the coefficients in \eqref{hat_N_n_explicit}, gives
$$
\begin{aligned}
 M_n(\bm k_n)
 &\leq\frac{C_Q}{\tau}
 \sum_{l=1}^{n-1}\frac{|(\alpha)_l|(2\pi)^{-dl/2}C_\rho^l}{l!}\sum_{\substack{I_0\dot\cup\cdots\dot\cup I_l=\{1,\ldots,n\}\\I_r\ne\varnothing}}
       \prod_{r=0}^l U_{|I_r|}((k_i)_{i\in I_r}).
\end{aligned}
$$
Here the blocks are ordered and nonempty, with $I_0$ the collision
block. By \eqref{prop_unsteady_cohomology_bound_eq},
$U_n\leq C_{\rm maj}E(\Sigma\bm k_n)M_n$.
Multiply by $\delta(q-\Sigma\bm k_n)\prod_jz(k_j)/n!$ and integrate.
The frequency constraint replaces $E(\Sigma\bm k_n)$ by $E(q)$,
and grouping blocks by their sizes gives, for $n\geq 2$,
\begin{equation}\label{prop82_profile_recursion}
\begin{aligned}
 \mathfrak u_n(q;z)&\leq\frac{C_{\rm maj}C_Q}{\tau}E(q)
 \sum_{l=1}^{n-1}\frac{|(\alpha)_l|(2\pi)^{-dl/2}C_\rho^l}{l!}\sum_{\substack{m_0+\cdots+m_l=n\\m_i\geq 1}}
 (\mathfrak u_{m_0}(\cdot;z)*\cdots*\mathfrak u_{m_l}(\cdot;z))(q).
\end{aligned}
\end{equation}
To see explicitly why the convolution appears, fix an ordered block decomposition and put $m_r=|I_r|$. The definition \eqref{def_integrated_majorant} and the convolution identity \eqref{convolution_rewrite} give
$$
\begin{aligned}
 &\frac1{n!}\int_{{B_c(\varepsilon_B)}^n}\delta(q-\Sigma\bm k_n)
       \prod_{r=0}^l U_{m_r}((k_i)_{i\in I_r})
       \prod_{j=1}^nz(k_j)\,d\bm k_n\\
 &\qquad=\frac{\prod_{r=0}^l m_r!}{n!}
       (\mathfrak u_{m_0}(\cdot;z)*\cdots*\mathfrak u_{m_l}(\cdot;z))(q).
\end{aligned}
$$
Each block integral contributes $m_r!\mathfrak u_{m_r}$, and the
block totals must sum to $q$, which produces the convolution.
For fixed sizes there are $n!/(m_0!\cdots m_l!)$ ordered block
decompositions. Labeling the $l$ density blocks counts each partition
$l!$ times, even when some sizes coincide; the collision block remains
distinguished. Thus the total scalar factor is
$$
 \frac1{n!}\frac1{l!}\frac{n!}{m_0!\cdots m_l!}
       \prod_{r=0}^l m_r!=\frac1{l!}.
$$
Symmetrization preserves this bound, since $\delta(q-\Sigma\bm k_n)\prod_jz(k_j)$ is symmetric in the individual frequencies. Tonelli's theorem justifies these rearrangements even before finiteness is known.\\
Each convolution in \eqref{prop82_profile_recursion} has at least two factors. Applying \eqref{small_divisor_convolution_bound} therefore gives
\begin{equation}\label{integrated_coefficient_majorant}
\begin{aligned}
 \|\mathfrak u_n(\cdot;z)\|_{\mathcal Y}
 &\leq\frac{C_{\rm maj}C_QC_E}{\tau}
 \sum_{l=1}^{n-1}\frac{|(\alpha)_l|(2\pi)^{-dl/2}C_\rho^l}{l!}\sum_{\substack{m_0+\cdots+m_l=n\\m_i\geq 1}}
 \prod_{i=0}^l\|\mathfrak u_{m_i}(\cdot;z)\|_{\mathcal Y},
 \qquad n\geq 2.
\end{aligned}
\end{equation}
Define now the majorization function
$$
 \Psi_\alpha(y):=\sum_{l\geq 1}\frac{|(\alpha)_l|}{l!}y^l,
$$
which is analytic for $|y|<1$  and a polynomial if $\alpha$ is a nonnegative integer, and solve
\begin{equation}\label{integrated_scalar_majorant_equation}
 A=A_1s+\frac{C_{\rm maj}C_QC_E}{\tau}A\Psi_\alpha((2\pi)^{-d/2}C_\rho A),\qquad A(0)=0.
\end{equation}
The derivative with respect to $A$ of the left side minus the right side at $(A,s)=(0,0)$ is one, since $\Psi_\alpha(0)=0$. The analytic implicit-function theorem therefore yields $A$ on a disk of positive radius. At order $n\geq 2$, its coefficient equation is the right-hand side of \eqref{integrated_coefficient_majorant} with each norm $\|\mathfrak u_{m_i}(\cdot;z)\|_{\mathcal Y}$ replaced by $A_{m_i}$. Every $m_i<n$, so this is a recursion in the previously determined coefficients. Starting from $A_1\geq 0$, it gives $A_n\geq 0$. Applying \eqref{integrated_coefficient_majorant} inductively, with $\sum_i m_i=n$, proves \eqref{weighted_profile_majorant_final}, including finiteness of every $\mathfrak u_n$.\\

\noindent
\underline{Convergence of the parametrization}:\\
Choose $R_A>0$ strictly inside the disk of analyticity of $A$, with
$(2\pi)^{-d/2}C_\rho A(R_A)<1$, and then choose $\mathfrak r_*>0$ so that
$$
 C_{\rm comp}^{\mathcal Y}\mathfrak r_*<R_A,\qquad
 (2\pi)^{-d/2}C_\rho A(C_{\rm comp}^{\mathcal Y}\mathfrak r_*)<\tfrac12.
$$
For $1\leq p<2$ we will decrease this radius once more below.
For the differentiability arguments below, we also record the estimate
with one argument distinguished. Positivity and symmetry give the
following mixed-argument form of \eqref{weighted_profile_majorant_final}, for
nonnegative $f,g\in\mathcal Y$ supported in $B_c(\varepsilon_B)$ and $n\geq 2$:
\begin{equation}\label{selected_mixed_bound}
 \left\|\frac1{(n-1)!}\int_{B_c(\varepsilon_B)^n}
 \delta(\,\cdot-\Sigma\bm k_n)U_n(\bm k_n)
 f(k_1)\prod_{j=2}^ng(k_j)\,d\bm k_n\right\|_{\mathcal Y}
 \leq enA_n\|f\|_{\mathcal Y}\|g\|_{\mathcal Y}^{n-1}.
\end{equation}
Indeed, the profile on the left is the coefficient of $\varepsilon$ in
$\mathfrak u_n(\cdot;g+\varepsilon f)$, and hence its norm is at most
$A_n\varepsilon^{-1}(\|g\|_{\mathcal Y}+\varepsilon\|f\|_{\mathcal Y})^n$ for $\varepsilon>0$.
For nonzero norms take $\varepsilon=\|g\|_{\mathcal Y}/((n-1)\|f\|_{\mathcal Y})$
and use $(1+1/(n-1))^{n-1}\leq e$. Note that zero arguments give the assertion directly.\\
By \eqref{def_Un_renormalized}, the normalized order-$n$ integral
contribution to \eqref{defhydro} satisfies
\begin{equation}\label{prop82_coefficient_map}
 \widehat{\mathcal W}^{\mathrm{int}}_n[\xi](k)
 =\frac1{n!}\int_{B_c(\varepsilon_B)^n}
 \delta(k-\Sigma\bm k_n)\hat{\mathfrak U}_n(\bm k_n,\theta[\xi])
 \Pi\zeta[\xi](\bm k_n)\,d\bm k_n.
\end{equation}
Here $\zeta[\xi](-k)=\zeta[\xi](k)^*$, since $\lambda$ is real and even.
The eigenvector satisfies $\hat e(-k)=\hat e(k)^*$. Likewise
$T_{-q}(s)g^*=\bigl(T_q(s)g\bigr)^*$ by \eqref{Duhamel_Tk}, and
$P_{-q}g^*=(P_qg)^*$ on the slow band by
\eqref{Pk_rank_one_corrected}. The forcing in \eqref{hat_N_n_explicit}
has real coefficients, and the cutoff in \eqref{Vn_selected_choice} is
unchanged when all leaf frequencies change sign. Induction in the
Duhamel formula therefore gives
$\hat{\mathfrak U}_n(-\bm k_n,s)=\hat{\mathfrak U}_n(\bm k_n,s)^*$.
Changing $k_j$ to $-k_j$ in \eqref{prop82_coefficient_map} shows that each
integrated term satisfies the same conjugate symmetry in $k$. Since this
condition defines a closed real subspace, it holds for the sum as well.\\
For $n\geq 2$, \eqref{selected_kernel_majorant} gives
\begin{equation}\label{prop82_decay_reconstruction}
 \begin{split}
 \|\widehat{\mathcal W}^{\mathrm{int}}_n[\xi](k)\|_{\mathcal X_v^{(2)}}
 &\leq\frac1{n!}\int_{B_c(\varepsilon_B)^n}
 \delta(k-\Sigma\bm k_n)U_n(\bm k_n)
 e^{-\varkappa\lambda_{\min}(\bm k_n)\theta[\xi]}
 \prod_{j=1}^n|\zeta[\xi](k_j)|\,d\bm k_n.
 \end{split}
\end{equation}
Dropping the decay factor in \eqref{selected_kernel_majorant} and applying
\eqref{weighted_profile_majorant_final} gives, for $z\in\mathcal Y$
supported in $B_c(\varepsilon_B)$,
\begin{equation}\label{prop82_coefficient_norm}
 \left\|\frac1{n!}\int_{B_c(\varepsilon_B)^n}
 \delta(\,\cdot-\Sigma\bm k_n)\hat{\mathfrak U}_n(\bm k_n,\theta)
 \Pi z(\bm k_n)\,d\bm k_n\right\|_{\mathcal Z^{(2)}}
 \leq A_n\|z\|_{\mathcal Y}^n,\qquad \theta\geq 0,
\end{equation}
also for $n=1$, since $\hat{\mathfrak U}_1=e^{\lambda\theta}\hat e$.
Thus \eqref{clock_Y_bound_direct} gives the uniform tail bound
$$
 \sum_{n>N}\|\widehat{\mathcal W}^{\mathrm{int}}_n[\xi]\|_{\mathcal Z^{(2)}}
 \leq\sum_{n>N}A_n(C_{\rm comp}^{\mathcal Y}\mathfrak r)^n
 \longrightarrow0.
$$

\noindent
\underline{Mild invariance}:\\
Proposition~\ref{lem_weighted_semigroup_aux} gives strong continuity
of each fiber semigroup. The bound \eqref{weighted_semigroup_general}
is uniform in $k$ on finite time intervals, so dominated convergence
gives a $C_0$-semigroup on $\mathcal Z^{(2)}$. The block calculation leading to
\eqref{prop82_profile_recursion}, without the factor $C_{\rm maj}E$,
and Young's inequality give, for nonnegative $z\in\mathcal Y$
supported in $B_c(\varepsilon_B)$,
\begin{equation}\label{prop82_forcing_majorant}
 \begin{gathered}
\left\|\frac1{n!}\int_{B_c(\varepsilon_B)^n}
 \delta(\,\cdot-\Sigma\bm k_n)M_n(\bm k_n)
 \Pi z(\bm k_n)\,d\bm k_n\right\|_{\mathcal Y}
 \leq B_n\|z\|_{\mathcal Y}^n,\qquad n\geq 2,\\
 B_1=0,\qquad B(s):=\sum_{n\geq 2}B_ns^n
 =\frac{C_Q}{\tau}A(s)
       \Psi_\alpha((2\pi)^{-d/2}C_\rho A(s)).
 \end{gathered}
\end{equation}
This majorant is analytic for $|s|<R_A$. Comparing with \eqref{integrated_scalar_majorant_equation} gives
$$
 A(s)=A_1s+C_{\rm maj}C_EB(s),\qquad
 A_n=C_{\rm maj}C_EB_n\quad(n\geq 2).
$$
Thus $A$ and $B$ require no separate convergence radii. The same convolution estimates
applied directly to \eqref{def_Nhat} yield
\begin{equation}\label{prop82_nonlinear_bound}
 \begin{gathered}
 \|\widehat{\mathcal N}_\alpha(f)\|_{\mathcal Z^{(2)}}
 \leq\frac{C_Q}{\tau}\|f\|_{\mathcal Z^{(2)}}
       \Psi_\alpha((2\pi)^{-d/2}C_\rho\|f\|_{\mathcal Z^{(2)}}),\qquad
 (2\pi)^{-d/2}C_\rho\|f\|_{\mathcal Z^{(2)}}<1.
 \end{gathered}
\end{equation}
Along the reduced orbit, \eqref{clock_cocycle} gives
$\theta[e^{\lambda t}\xi]=\theta[\xi]+t$ and
$\zeta[e^{\lambda t}\xi]=\zeta[\xi]$.
Apply \eqref{sol_unsteady_cohomol} between $\theta[\xi]$ and
$\theta[\xi]+t$, multiply by
$\delta(k-\Sigma\bm k_n)\Pi\zeta[\xi](\bm k_n)/n!$, and integrate.
By \eqref{weighted_majorant_def} and
\eqref{prop82_forcing_majorant}, the reconstructed forcing at order $n$
has $\mathcal Z^{(2)}$ norm at most
$B_n\|\zeta[\xi]\|_{\mathcal Y}^n$, uniformly in time. These bounds sum to
$B(\|\zeta[\xi]\|_{\mathcal Y})<\infty$.
Together with \eqref{weighted_semigroup_general}, this justifies
interchanging summation and the Bochner integral on every finite time
interval. The partition formula \eqref{hat_N_n_explicit} and absolute
convergence identify the sum of the forcing terms with
$\widehat{\mathcal N}_\alpha(\hat{\mathcal W}[e^{\lambda\sigma}\xi])$.
We therefore obtain
\begin{equation}\label{prop82_mild_identity}
 \begin{split}
 \hat{\mathcal W}[e^{\lambda t}\xi](k)
 ={}&e^{t\mathcal L_k}\hat{\mathcal W}[\xi](k)+\int_0^t e^{(t-\sigma)\mathcal L_k}
       \bigl[\widehat{\mathcal N}_\alpha
       (\hat{\mathcal W}[e^{\lambda\sigma}\xi])\bigr](k)\,d\sigma,
 \qquad t\geq 0.
 \end{split}
\end{equation}
The identity holds in $\mathcal Z^{(2)}$ and hence fiberwise almost
everywhere. Kernel continuity and the same integrable majorants give
continuity of each integral contribution along the orbit; uniform
convergence gives continuity of the sum. Since the reduced flow preserves
the parameter ball, \eqref{prop82_mild_identity} defines a mild kinetic
solution for every $t\geq 0$ and proves forward invariance.

\noindent
\underline{The physical-space conclusions}:\\
In \eqref{prop82_mild_identity} the semigroup acts on the Fourier-space
data before inverse transformation. The bounded map
$\mathcal F^{-1}:L_k^1(\mathcal X_v)\to L_x^\infty(\mathcal X_v)$
therefore carries the convergent series and the Duhamel identity into
$L_x^\infty(\mathcal X_v)$. The Fourier product rule \eqref{def_Nhat}
identifies the nonlinearity with \eqref{defNalpha}. For $p=2$,
Hilbert-space Plancherel gives the same convergence and the
inverse-transformed identity in $L_x^2(\mathcal X_v)$.
For $1\leq p<2$, introduce
$$
 \mathcal H_v=L_v^2(\phi^{p-2}w),\qquad
 \mathcal H_v^{(2)}=\{f:(1+|v|)^2f\in\mathcal H_v\},\qquad
 \|f\|_{\mathcal H_v^{(2)}}=\|(1+|v|)^2f\|_{\mathcal H_v}.
$$
H\"older's inequality applied to
$|f|^pw=(|f|^2\phi^{p-2}w)^{p/2}(\phi^pw)^{(2-p)/2}$ gives
\begin{equation}\label{prop82_hilbert_embedding}
 \|f\|_{\mathcal X_v}\leq c_H\|f\|_{\mathcal H_v},\qquad
 \|f\|_{\mathcal X_v^{(2)}}\leq c_H\|f\|_{\mathcal H_v^{(2)}},\qquad
 c_H=\|\phi\|_{\mathcal X_v}^{1-p/2}.
\end{equation}
The second estimate follows by replacing $f$ with $(1+|v|)^2f$.
The weight $\phi^{p-2}w$ satisfies \eqref{ass_phi_w}: all polynomial
Maxwellian moments are finite by that assumption, and, for $1<p<2$,
$$
 \int\frac{\phi^{2-p}}w\,dv
 \leq\left(\int\phi^pw\,dv\right)^{(2-p)/p}
       \left(\int w^{-1/(p-1)}\,dv\right)^{2(p-1)/p}<\infty.
$$
For $p=1$, use $\int\phi/w\leq\|w^{-1}\|_\infty$.
The preceding construction therefore applies in $\mathcal H_v$, giving
an analytic majorant $A_H(s)=\sum_{n\geq 1}A_{H,n}s^n$ and a positive
radius $\mathfrak r_{*,H}$. Its kernels coincide with the original ones:
\eqref{Duhamel_Tk} and \eqref{Pk_rank_one_corrected} identify the
semigroups and projections under \eqref{prop82_hilbert_embedding},
and this continuous embedding commutes with the convergent integrals
defining the selected homological initial conditions \eqref{Vn_selected_choice}.
Induction in \eqref{sol_unsteady_cohomol} then identifies every order.

Replace $\mathfrak r_*$ by $\min\{\mathfrak r_*,\mathfrak r_{*,H}\}$.
For this common radius, Plancherel and \eqref{prop82_hilbert_embedding}
give the explicit summable estimate
\begin{equation}\label{prop82_physical_hilbert_bound}
 \|\mathcal F^{-1}\widehat{\mathcal W}^{\mathrm{int}}_n[\xi]
       \|_{L_x^2(\mathcal X_v)}
 \leq c_H\|\widehat{\mathcal W}^{\mathrm{int}}_n[\xi]\|_{L_k^2(\mathcal H_v)}
 \leq c_HA_{H,n}(C_{\rm comp}^{\mathcal Y}\mathfrak r)^n.
\end{equation}
This proves uniform absolute convergence in $L_x^2(\mathcal X_v)$.
For $1\leq p<2$ the sum constructed from these same kernels satisfies \eqref{prop82_mild_identity}
in $L_k^1(\mathcal H_v^{(2)})\cap L_k^2(\mathcal H_v^{(2)})$.
Applying Hilbert-space Plancherel and then the continuous embedding
$\mathcal H_v\hookrightarrow\mathcal X_v$ from \eqref{prop82_hilbert_embedding}
gives the inverse-transformed identity in $L_x^2(\mathcal X_v)$ along each
reduced orbit.
\end{proof}

\begin{remark}
The small divisor multiplies a convolution of at least two scalar norm profiles, not a general vector-valued function. The condition $d\geq 3$ enters through $E-1\in L^2$. In $d=2$, $E-1$ is locally integrable but not locally square integrable; in $d=1$ it is not locally integrable. Failure of this particular majorant in lower dimensions is not a nonexistence result.
\end{remark}

\begin{proposition}[Differentiability at equilibrium]\label{prop_selected_tangency}

Assume the hypotheses of Proposition~\ref{prop_majorant_functional_taylor}.
For fixed $0<\mathfrak r<\mathfrak r_*$, set $\hat{\mathcal W}[0]=0$ in \eqref{defhydro}.
The extension of $\hat{\mathcal W}$ is Fr\'echet
differentiable at the origin from the real Banach space
$\Xi(\varepsilon_B)$ with its $C^0$ norm into $\mathcal Z^{(2)}$ with derivative
\begin{equation}\label{selected_tangent_identity}
 D\hat{\mathcal W}[0]h(k,v)
 =\widehat e_h(k,v)
 =h(k)\hat e(k,v)\mathbf1_{B_c(\varepsilon_B)}(k),
 \qquad h\in\Xi(\varepsilon_B),
\end{equation}
where $e_h$ is the linear parametrization \eqref{defexv}, with its
Fourier transform extended by zero outside $B_c(\varepsilon_B)$.
Its range in $\mathcal Z^{(2)}$ is $\{\widehat e_h:h\in\Xi(\varepsilon_B)\}$. After Fourier inversion,
$\mathcal W$ is differentiable at zero into
$L_x^\infty(\mathcal X_v)$, and also into $L_x^2(\mathcal X_v)$ for
$1\leq p\le2$, with $D\mathcal W[0]h=e_h$.
\end{proposition}
\begin{proof}
Retain $\widehat{\mathcal W}^{\mathrm{int}}_n$, $\mathcal Y$, and $A$
from the proof of Proposition~\ref{prop_majorant_functional_taylor}.
By \eqref{W1}, $\widehat{\mathcal W}^{\mathrm{int}}_1[\xi]=\widehat e_\xi$.
The cutoff on the individual frequencies and \eqref{bound_e_hat}, applied in $\mathcal X_v^{(2)}$,
make $h\mapsto\widehat e_h$ bounded from $\Xi(\varepsilon_B)$ to
$\mathcal Z^{(2)}$. It remains to estimate the nonlinear remainder. The strict inequality
$\varkappa>1$ leaves additional decay, which will make the bound for
the normalized remainder tend to zero as $\theta[\xi]\to\infty$,
with $\mathfrak r$ fixed.
Put $t=\theta[\xi]$ and $z=|\zeta[\xi]|$. If
$-\lambda(k_j)=\lambda_{\min}(\bm k_n)$, then
$$
 e^{-\varkappa\lambda_{\min}(\bm k_n)t}\prod_{i=1}^nz(k_i)
 =e^{(\varkappa-1)\lambda(k_j)t}|\xi(k_j)|\prod_{i\ne j}z(k_i).
$$
Since all terms are nonnegative, it follows that
$$
 e^{-\varkappa\lambda_{\min}(\bm k_n)t}\prod_{i=1}^nz(k_i)
 \leq\sum_{j=1}^n e^{(\varkappa-1)\lambda(k_j)t}
                    |\xi(k_j)|\prod_{i\ne j}z(k_i).
$$
Insert this inequality into \eqref{prop82_decay_reconstruction}.
Symmetry makes the $n$ integrated terms equal, so their normalization
is $n/n!=1/(n-1)!$, exactly as in \eqref{selected_mixed_bound}.
Applying that estimate with
$f=e^{(\varkappa-1)\lambda t}|\xi|$ and $g=z$ gives
$$
 \|\widehat{\mathcal W}^{\mathrm{int}}_n[\xi]\|_{\mathcal Z^{(2)}}
 \leq enA_n\|e^{(\varkappa-1)\lambda t}|\xi|\|_{\mathcal Y}
                   \|z\|_{\mathcal Y}^{n-1},\qquad n\geq 2.
$$
Set $s_{\mathfrak r}=C_{\rm comp}^{\mathcal Y}\mathfrak r<R_A$.
Summing and using \eqref{clock_Y_bound_direct} yields
\begin{equation}\label{selected_remainder_bound}
 \frac{\|\hat{\mathcal W}[\xi]-\widehat e_\xi\|_{\mathcal Z^{(2)}}}
      {\|\xi\|_{C^0}}
 \leq e\bigl(A'(s_{\mathfrak r})-A_1\bigr)
       \|e^{(\varkappa-1)\lambda\theta[\xi]}
                  \mathbf1_{B_c(\varepsilon_B)}\|_{\mathcal Y}.
\end{equation}
The first factor is finite and independent of $\xi$.
By \eqref{global_clock_lower_Bound}, $\theta[\xi]\to\infty$ uniformly
as $\|\xi\|_{C^0}\to0$ at fixed $\mathfrak r$.
Since $\varkappa>1$ and $\lambda(k)<0$ for $k\ne0$, dominated convergence
on the ball of individual frequencies makes the last norm tend to zero. This proves
\eqref{selected_tangent_identity} and continuity at zero.\\
Fourier inversion transfers this remainder estimate to
$L_x^\infty(\mathcal X_v)$, and Plancherel does so to
$L_x^2(\mathcal X_v)$ when $p=2$.
For $1\leq p<2$, apply the same estimate in the auxiliary space
$\mathcal H_v^{(2)}$ constructed in the preceding proof, with its
majorant $A_H$ and the common radius fixed there.
The first inequality in \eqref{prop82_physical_hilbert_bound} transfers
the resulting little-$o$ remainder to $L_x^2(\mathcal X_v)$.
\end{proof}

\begin{proposition}[$C^1$ regularity away from the origin]\label{prop_C1_away_origin}
Assume the hypotheses of Proposition~\ref{prop_majorant_functional_taylor}, choose $0<\mathfrak r<\mathfrak r_*$, and suppose
\begin{equation}\label{C1_clock_nu_assumption}
 \nu\geq 2/\tau,\qquad C_{\rm comp}^{\mathcal Y}\mathfrak r<R_A.
\end{equation}
For each $0<\rho<1$, the parametrization is continuously Fr\'echet differentiable, as a map of real Banach spaces,
\begin{equation}\label{C1_punctured_map}
 \hat{\mathcal W}:\mathcal U_\rho^\times\longrightarrow
 \mathcal Z,\qquad
 \mathcal U_\rho^\times=\{\xi\in\Xi(\varepsilon_B):0<\|\xi\|_{C^0}<\rho\mathfrak r\}.
\end{equation}
Its inverse Fourier transform is $C^1$ into $L_x^\infty(\mathcal X_v)$, and also into $L_x^2(\mathcal X_v)$ for $1\leq p\le2$.
\end{proposition}

\begin{proof}
Using the integral contributions \eqref{prop82_coefficient_map}, write
\begin{equation}\label{C1_compensated_expansion}
\hat{\mathcal W}[\xi]
 =\sum_{n\geq 1}\widehat{\mathcal W}^{\mathrm{int}}_n[\xi].
\end{equation}
Differentiating the factors at fixed $\theta$ is controlled by \eqref{selected_mixed_bound},
with order one following from \eqref{W1}.
Equations \eqref{equ_U}, \eqref{C1_generator_weight_loss}, and
\eqref{prop82_forcing_majorant}, together with
$|\Sigma\bm k_n|\leq n(k_c-\varepsilon_B)$, imply
\begin{equation}
 \begin{gathered}
\left\|\frac1{n!}\int_{B_c(\varepsilon_B)^n}
 \delta(\,\cdot-\Sigma\bm k_n)\partial_\theta\hat{\mathfrak U}_n(\bm k_n,\theta)
 \Pi z(\bm k_n)\,d\bm k_n\right\|_{\mathcal Z}
 \leq D_n\|z\|_{\mathcal Y}^n,\\
 D_n=C_{\rm gen}(1+n(k_c-\varepsilon_B))A_n+B_n.
 \end{gathered}
\end{equation}
At fixed order,
Remark~\ref{rem_weighted_majorant_role} gives continuous time derivatives.
Dominated convergence, with the majorants evaluated at
$\mathbf1_{B_c(\varepsilon_B)}$, gives continuity in the multilinear
operator norm from $\Xi(\varepsilon_B)^n$ to $\mathcal Z$.
Thus the integral in \eqref{prop82_coefficient_norm} is $C^1$ in $(\theta,z)$ for $\theta>0$.\\
By Lemma~\ref{lem_theta}, $\theta$ and $\zeta$ are locally $C^1$, with
\begin{equation}\label{C1_zeta_derivative}
 D\zeta[\xi]h=e^{-\lambda\theta[\xi]}h
                  -\lambda\zeta[\xi]D\theta[\xi]h.
\end{equation}
Near each $\xi_0\in\mathcal U_\rho^\times$ their derivatives are bounded, while
\eqref{clock_Y_bound_direct} gives
$\|\zeta[\xi]\|_{\mathcal Y}\leq s_{\mathfrak r}
:=C_{\rm comp}^{\mathcal Y}\mathfrak r<R_A$.
The chain rule consequently bounds the derivative of the $n$th
summand in \eqref{C1_compensated_expansion} by
$$
 C_{\rm clock,loc}(\xi_0)\bigl(D_ns_{\mathfrak r}^n
                   +enA_ns_{\mathfrak r}^{n-1}\bigr).
$$
The constant $C_{\rm clock,loc}(\xi_0)$ depends on the chosen neighborhood of $\xi_0$ and the fixed radius $\mathfrak r$; it is independent of $n$ but need not remain bounded as $\xi_0\to0$.
This sequence is summable by analyticity of $A$ and $B$.
Together with \eqref{prop82_coefficient_norm}, it proves local uniform
convergence of the series and its derivatives, hence $C^1$ regularity.
Fourier inversion and the Hilbert-space argument of
Proposition~\ref{prop_majorant_functional_taylor} give the physical-space
conclusions, using the same estimates in $\mathcal H_v^{(2)}$ for $p<2$.
\end{proof}

\begin{remark}[{The $C^1$ obstruction for unrestricted complex amplitudes}]
\label{rem_selected_not_C1}
For $\alpha>0$, we note that the parametrization need not be $C^1$ at zero
on the larger space of continuous complex amplitudes with the same frequency cutoff and no reflection condition. In this remark we use the same clock and series on that space, regarded as real. The convergence and real derivative estimates above apply there as well, since their majorants use absolute values. Assume $\nu\geq 2/\tau$ and choose a small $\epsilon>0$ as well as frequency vectors
$$
 p=\epsilon(1,-2,0,\ldots,0),\quad
 k=\epsilon(10,10,0,\ldots,0),\quad q_*=p+k.
$$
Take nonnegative smooth bumps $\eta,h$ in balls of radius
$\epsilon/100$ about $p,k$, respectively, positive near their
centres, with
$0<\|\eta\|_{C^0}<\mathfrak r$, and fix $\|h\|_{C^0}=1$.
For small $\epsilon$, both bump supports and their pairwise sums
$p'+k'$ with $p'\in\operatorname{supp}\eta$ and $k'\in\operatorname{supp}h$
lie strictly inside the input and output cutoffs, respectively. Since
$$
 -\lambda(p)=5\tau\epsilon^2+O(\epsilon^4),\quad
 -\lambda(k)=200\tau\epsilon^2+O(\epsilon^4),\quad
 -\lambda(q_*)=185\tau\epsilon^2+O(\epsilon^4),
$$
one has, on the interacting supports,
$$
\begin{gathered}
 -\lambda(p'+k')>-\varkappa_2\lambda(p'),\qquad
 -\lambda(k')+\lambda(p'+k')\geq  c\epsilon^2>0,\\
 -\lambda(p')-\lambda(k')+\lambda(p'+k')\geq  c\epsilon^2>0.
\end{gathered}
$$
In particular,
$\chi\bigl(-\lambda(p'+k')/\min\{-\lambda(p'),-\lambda(k')\}\bigr)=0$
there. Set $\xi_t=e^{\lambda t}\eta$.
The disjoint supports give $D\theta[\xi_t]h=0$.
Using \eqref{Pk_rank_one_corrected}, the eigenvector expansion,
$\rho[2\phi-\hat e(p')-\hat e(k')]=0$, and the Gaussian second moment,
we obtain, uniformly on these supports,
$$
\begin{aligned}
 \frac{\alpha(2\pi)^{-d/2}}{\tau}
 \rho\!\left[P_{p'+k'}\bigl(2\phi-\hat e(p')-\hat e(k')\bigr)\right]
 &=\frac{\alpha(2\pi)^{-d/2}}{\tau}
 \frac{\displaystyle\int
 \frac{2\phi(v)-\hat e(p',v)-\hat e(k',v)}
 {1+\tau\lambda(p'+k')+\ri\tau(p'+k')\cdot v}\,dv}
 {\displaystyle\int
 \frac{\phi(v)}{(1+\tau\lambda(p'+k')+\ri\tau(p'+k')\cdot v)^2}\,dv}\\
 &\quad=\alpha(2\pi)^{-d/2}\tau|p'+k'|^2+O(\epsilon^3)>0.
\end{aligned}
$$
This expression is real by velocity reflection and conjugation.
Differentiating the quadratic term on a small fixed neighborhood
$Q$ of $q_*$ therefore gives the positive integral
$$
\begin{aligned}
 &\rho\!\left[P_q\bigl((D\hat{\mathcal W}[\xi_t]h)(q,\cdot)\bigr)\right]=\frac{\alpha(2\pi)^{-d/2}}{\tau}
 \int\rho\!\left[P_q\bigl(2\phi-\hat e(p')-\hat e(q-p')\bigr)\right]\\
 &\qquad\times
 \frac{e^{[-\lambda(p')-\lambda(q-p')+\lambda(q)](\theta[\eta]+t)}-1}
 {-\lambda(p')-\lambda(q-p')+\lambda(q)}
 e^{\lambda(p')t}\eta(p')h(q-p')\,dp'.
\end{aligned}
$$
There are no other contributions on $Q$. The linear derivative
is supported near $k$. An order-$n$ leaf derivative, $n\geq 3$,
has one individual frequency near $k$ and $n-1$ individual frequencies near $p$, so its first
coordinate is at least
$\epsilon(10-1/100+(n-1)(1-1/100))$, strictly above that of $Q$.
The derivative of the spectral clock vanishes. Local termwise differentiation
is justified by Proposition~\ref{prop_C1_away_origin}.
The convolution $\eta*h$ is positive at $q_*$ and
bounded below on a smaller closed ball in $Q$. On the interacting
supports the minimum of $-\lambda(k')+\lambda(p'+k')$ is positive.
The denominator in the displayed integral is likewise bounded above
and away from zero. For sufficiently large $t$, the numerator is
therefore at least half its exponential term. Since $\theta[\eta]\geq 0$,
multiplication by $e^{\lambda(p')t}$ leaves an exponential with rate
at least that positive minimum. Integration over the smaller ball gives
$$
 \|D\hat{\mathcal W}[\xi_t]h\|_{L_q^2(\mathcal X_v)}
       \geq  C e^{ct},\qquad \|\xi_t\|_{C^0}\longrightarrow0.
$$
Boundedness of $\rho P_q$ on this compact set of summed frequencies justifies the norm lower bound. Thus the derivative is not continuous at zero,
although Proposition~\ref{prop_selected_tangency} proves
differentiability there. Continuity of the derivative at equilibrium on $\Xi(\varepsilon_B)$ remains open for $\alpha>0$.
Note that for $\alpha=0$ the parametrization is linear.
\end{remark}

\begin{remark}[The expansion at auxiliary time zero]
\label{rem_auxiliary_zero_family}
The derivative at equilibrium is the linear lifting
$h\mapsto\widehat e_h$ in \eqref{selected_tangent_identity}.
It is not bounded below from the $C^0$ parameter norm into
$\mathcal Z^{(2)}$, so this tangency does not give a local $C^1$ chart at zero.
Its range is not closed in $\mathcal Z^{(2)}$. Indeed, approximate the
indicator of a smaller ball centered at zero in the interior of the cutoff band by smooth real even
amplitudes in $L^1\cap L^2$. By \eqref{eigenvector} their lifts converge
in $\mathcal Z^{(2)}$, but the density of the limit has no continuous
representative.\\
Use the fixed-clock expansion $\hat{\mathcal W}_\varepsilon[z;\theta]$
of \eqref{defhydro_epsilon} at $\varepsilon=1$. At auxiliary time zero, \eqref{sol_unsteady_cohomol_W} gives
\begin{equation}\label{selected_initial_lifting}
 \hat{\mathcal W}_1[z;0](q)
 =\widehat e_z(q)+\sum_{n\geq 2}\frac1{n!}
 \int_{B_c(\varepsilon_B)^n}
 \delta(q-\Sigma\bm k_n)\hat{\mathfrak V}_n(\bm k_n)
 \Pi z(\bm k_n)\,d\bm k_n.
\end{equation}
For small $\|z\|_{\mathcal Y}$, convergence in $\mathcal Z^{(2)}$
and the integrated mild evolution follow from
\eqref{prop82_coefficient_norm} and the summation of
\eqref{sol_unsteady_cohomol} in Proposition~\ref{prop_majorant_functional_taylor}.
By \eqref{def_Un_renormalized}, this evolution is
$\hat{\mathcal W}_1[e^{\lambda s}z;s]$; at
$z=\zeta[\xi]$ and $s=\theta[\xi]$ it equals
$\hat{\mathcal W}_{\mathfrak r}[\xi]$.
Its initial density is generally nonlinear in $z$.
Section~\ref{sec_slow_spectral_closure} adds the continuous Fourier
estimates and local density coordinates.
\end{remark}

\section{The Kinetic Manifold, Local Density Inversion and the Nonlinear Spectral Closure}\label{sec_slow_spectral_closure}

In this section, we realize the nonlinear spectral closure of Section~\ref{sec_inv_equation_general}
using the invariant parametrization from Section~\ref{sec_proof_main_theorem}.
Theorem~\ref{prop_kinetic_embedding} gives local $C^1$ kinetic charts away from equilibrium, in which the
full kinetic state determines the slow amplitude. On a smaller
forward-invariant domain, Theorem~\ref{prop_density_transversality_s9} proves that the
retracted density derivative is invertible. The mass density alone thus gives local
coordinates and determines the flux in \eqref{dynamcis_rho}. We also prove
positivity and entropy dissipation.\\
We retain $\hat{\mathcal W}$ from \eqref{defhydro} as the parametrization of the hydrodynamic manifold, the spectral clock
\eqref{deftheta}, and the homological initial conditions
\eqref{Vn_selected_choice}. The radius $\mathfrak r$ is fixed throughout, and its implicit occurrence in $\hat{\mathcal W}$ is suppressed as in
Section~\ref{sec_proof_main_theorem}. The argument uses the clock
estimates of Lemma~\ref{lem_theta} and the convergence and regularity
results of Section~\ref{sec_proof_main_thm_proof}.
For the density-inversion argument, we assume the hypotheses of Proposition~\ref{prop_C1_away_origin}, in particular
$\nu\geq 2/\tau$. All derivatives in the following have to be understood as real
Fr\'echet derivatives.\\
To ease notation, throughout this section, $\hat{\mathcal W}[z;s]$ denotes the series
\eqref{defhydro_epsilon} at $\varepsilon=1$ with input $z$ and auxiliary
time $s$, wherever it converges. We write $D_z\hat{\mathcal W}[z;s]$ for the derivative
with respect to the input at fixed $s$. With $\zeta[\xi]=e^{-\lambda\theta[\xi]}\xi$,
Remark~\ref{rem_auxiliary_zero_family} gives
\begin{equation}\label{closure_flow_identity_s9}
\hat{\mathcal W}[\xi]
 =\hat{\mathcal W}[\xi;\theta[\xi]]
 =\left.\hat{\mathcal W}[e^{\lambda s}\zeta[\xi];s]\right|_{s=\theta[\xi]}.
\end{equation}
For the regularity arguments, we again use the spaces $\mathcal Y$ from \eqref{def_majorant_space_Y} and
$\mathcal Z^{(2)}$ from \eqref{def_Xv2_aux}. In addition, density inversion also
requires continuity in the summed frequency, and we therefore introduce the real
spaces
\begin{equation}\label{closure_continuous_spaces_s9}
 \begin{aligned}
 \mathcal Y_\rho&=\{R\in C_0(\mathbb R_k^d;\mathbb C):R(-k)=R(k)^*\},\\
 \mathcal X_{k,v}&=\{f\in C_0(\mathbb R_k^d;\mathcal X_v^{(2)})\cap L_k^1(\mathcal X_v^{(2)}):
                         f(-k,\cdot)=f(k,\cdot)^*\}.
 \end{aligned}
\end{equation}
We equip $\mathcal Y_\rho$ with the supremum norm and $\mathcal X_{k,v}$ with the norm
$$
\|f\|_{\mathcal X_{k,v}}:=\|f\|_{C_0(\mathcal X_v^{(2)})}
                          +\|f\|_{L_k^1(\mathcal X_v^{(2)})}.
$$
Thus $\mathcal X_{k,v}$ is the conjugate-symmetric subspace of
$C_0(\mathbb R_k^d;\mathcal X_v^{(2)})\cap\mathcal Z^{(2)}$, with equivalent sum norms, and
$\|f\|_{\mathcal Z^{(2)}}\leq\tfrac32\|f\|_{\mathcal X_{k,v}}$.
In the notation of Section~\ref{sec_spectral_closure_general}, the kinetic
and observable spaces are here $\mathcal X=\mathcal X_{k,v}$ and
$\mathcal X_{\rm macro}=\mathcal Y_\rho$. By \eqref{selected_tangent_identity}, $\rho[\widehat e_h]$ extends $h$
by zero outside $B_c(\varepsilon_B)$, so $h\mapsto\rho[\widehat e_h]$
is an isometry from $\Xi(\varepsilon_B)$ into $\mathcal Y_\rho$. The following lemma gives the action of the kinetic group $e^{s\mathcal L}$ on $\mathcal X_{k,v}$.

\begin{lemma}[The kinetic group on $\mathcal X_{k,v}$]\label{lem_kinetic_group_s9}
The Fourier-space operator $(\mathcal L f)(k,\cdot)=\mathcal L_k f(k,\cdot)$
of \eqref{Lk} generates a strongly continuous group on $\mathcal X_{k,v}$ with domain
\begin{equation}\label{fourier_generator_domain_s9}
 D(\mathcal L)=\{f\in\mathcal X_{k,v}:-\ri(k\cdot v)f\in\mathcal X_{k,v}\},
\end{equation}
equipped with the graph norm $\|f\|_{\mathcal X_{k,v}}+\|\mathcal L f\|_{\mathcal X_{k,v}}$.
For $s\geq0$, $(e^{s\mathcal L}f)(k)=T_k(s)f(k)$ and
$$
 \|e^{s\mathcal L}f\|_{\mathcal X_{k,v}}
 \leq C_{\mathcal L}^{\rm wgt}(1+s/\tau)^2\|f\|_{\mathcal X_{k,v}}.
$$
\end{lemma}
\begin{proof}
The free multipliers $e^{-\ri s k\cdot v}$ preserve conjugate symmetry and
both norms in $\mathcal X_{k,v}$. On bounded $(k,v)$ sets they converge uniformly
to one as $s\to0$ and profiles with bounded $(k,v)$ support are dense. Continuous even cutoffs in $k$ converge in both norms, and velocity truncation converges in $L_k^1$ by dominated convergence and uniformly
in $k$ because a $C_0$ profile has relatively compact range in $\mathcal X_v^{(2)}$ and $p<\infty$. This proves strong continuity.\\
For $f$ in \eqref{fourier_generator_domain_s9}, the identity
$$
 (e^{-\ri s k\cdot v}-1)f
 =\int_0^s e^{-\ri r k\cdot v}\bigl(-\ri(k\cdot v)f\bigr)\,dr
$$
gives the generator derivative in $\mathcal X_{k,v}$. Conversely, convergence
of the difference quotients in $\mathcal X_{k,v}$ gives convergence in
$\mathcal X_v^{(2)}$ at each $k$, where the pointwise multiplier derivative
identifies the limit as $-\ri(k\cdot v)f$. The collision operator
$\tau^{-1}\mathcal L_Q$ is bounded on $\mathcal X_{k,v}$ by \eqref{def_LQ}
and preserves conjugate symmetry. The bounded perturbation theorem yields
the group on the same domain. Evaluating its Duhamel formula at $k$
gives $T_k(s)$ by uniqueness in \eqref{Duhamel_Tk} and taking the supremum
and the integral in $k$ in \eqref{weighted_semigroup_general} gives the bound.
\end{proof}

\subsection{An invariant domain for density inversion}
Fix $0<\delta<1$. On the punctured ball $\mathcal U_\delta^\times$
of \eqref{C1_punctured_map}, the compensated profile is constant along
the orbit and satisfies \eqref{compensated_clock_normalization_s9}.
The following integral controls the derivative of the clock: 
\begin{equation}\label{closure_clock_data_s9}
 \Gamma_{\mathfrak r}(\zeta[\xi])
 =\int_{B_c(\varepsilon_B)}\frac{\mu(k)}{-\lambda(k)}
 \left(\frac{|\zeta[\xi](k)|}{\mathfrak r}\right)^{-\nu/\lambda(k)-1}\,dk.
\end{equation}
As in \eqref{deftheta}, the value at $k=0$ is irrelevant and can be neglected. For $G>0$,
we define the set 
\begin{equation}\label{invariant_density_domain_s9}
 \mathcal U_{\delta,G}^{(\mathfrak r)}
 =\{\xi\in\mathcal U_\delta^\times:
       \Gamma_{\mathfrak r}(\zeta[\xi])<G\}.
\end{equation}
The following lemma gives an invariant domain on which the clock
derivative can be controlled uniformly.

\begin{lemma}[Invariant domain and derivative of the spectral clock ]
\label{lem_invariant_clock_domain_s9}
The map $\xi\mapsto\Gamma_{\mathfrak r}(\zeta[\xi])$ in
\eqref{closure_clock_data_s9} is finite and continuous on
$\mathcal U_\delta^\times$. For every $G>0$, the domain
\eqref{invariant_density_domain_s9} is open in $\Xi(\varepsilon_B)$
and forward invariant under \eqref{reduced_dyn_linear}.

For every $\xi\in\mathcal U_{\delta,G}^{(\mathfrak r)}$ and
$h\in\Xi(\varepsilon_B)$, the clock derivative satisfies
\begin{equation}\label{compensated_clock_derivative_s9}
 |D\theta[\xi]h|
 \leq \frac{\Gamma_{\mathfrak r}(\zeta[\xi])}{\mathfrak r}
       \|e^{-\lambda\theta[\xi]}h\|_\infty
 \leq \frac{G}{\mathfrak r}\|e^{-\lambda\theta[\xi]}h\|_\infty .
\end{equation}

\end{lemma}
\begin{proof}
Lemma~\ref{lem_theta} and \eqref{C1_zeta_derivative} give local
continuity of $\theta$ and $\zeta$ in the $C^0$ topology, with
$\zeta[\xi](0)=\xi(0)$.
Near any fixed parameter, we choose $0<s<k_c-\varepsilon_B$ and
$0<\beta<1$ such that $|\zeta[\xi](k)|/\mathfrak r\leq\beta$ for
$|k|<s$, uniformly in a parameter neighborhood. Since $y e^{-y}\leq e^{-1}$ for $y\geq 0$, we can estimate 
$$
\begin{aligned}
 \int_{0<|k|<s}\frac{\mu(k)}{-\lambda(k)}
 \left(\frac{|\zeta[\xi](k)|}{\mathfrak r}\right)^{\nu/(-\lambda(k))-1}\,dk &\leq \frac1\beta\int_{0<|k|<s}\frac{\mu(k)}{-\lambda(k)}
   e^{-\nu|\log\beta|/(-\lambda(k))}\,dk\\
& \leq \frac{1}{e\beta\nu|\log\beta|}\int_{|k|<s}\mu(k)\,dk .
\end{aligned}
$$
On $s\leq |k|\leq k_c-\varepsilon_B$, we have $-\lambda(k)\geq -\lambda(s)>0$.
For those values of $k$, the compensated profile  $|\zeta[\xi](k)|/\mathfrak r$ is uniformly bounded in the same
neighborhood of frequencies, possibly after shrinking it further. Thus a fixed multiple
of $\mu$ dominates the integrand on the whole ball. Pointwise
continuity for $k\ne0$ and dominated convergence give finiteness
and continuity of \eqref{closure_clock_data_s9}, implying that the domain $ \mathcal U_{\delta,G}^{(\mathfrak r)}$ is open. It is forward invariant because
$e^{t\lambda}$ contracts the parameter norm, preserves nonzero parameters
at finite times, and leaves $\zeta[\xi]$ unchanged.

In the derivative formula of Lemma~\ref{lem_theta}, write
$\xi=e^{\lambda\theta[\xi]}\zeta[\xi]$ and
$h=e^{\lambda\theta[\xi]}(e^{-\lambda\theta[\xi]}h)$.
The exponential factors in the numerator then cancel
$S_\mu[\xi]=e^{-\nu\theta[\xi]}$ in the denominator, while taking absolute
values gives
$$
 |D\theta[\xi]h|\leq \frac{1}{\mathfrak r}
 \int_{B_c(\varepsilon_B)}\frac{\mu(k)}{-\lambda(k)}
 \left(\frac{|\zeta[\xi](k)|}{\mathfrak r}\right)^{-\nu/\lambda(k)-1}
 \bigl|e^{-\lambda(k)\theta[\xi]}h(k)\bigr|\,dk .
$$
Taking the supremum of the last factor and using
\eqref{closure_clock_data_s9} and \eqref{invariant_density_domain_s9}
proves \eqref{compensated_clock_derivative_s9}.

\end{proof}

Note that, for fixed $\mathfrak r$ and $\delta$, the domains $ \mathcal U_{\delta,G}^{(\mathfrak r)}$ exhaust all of $\mathcal U_\delta^\times$ as $G$ increases, by finiteness of \eqref{closure_clock_data_s9}. We stress that the radius in Theorem~\ref{prop_density_transversality_s9} is chosen after $G$.

\begin{remark}[An explicit forward-invariant class of parameter functions]\label{rem_logarithmic_class_s9}
Let $\eta$  be a non-increasing radial cutoff function with $\eta(0)=1$,
positive for $|k|<k_c-\varepsilon_B$ and vanishing to all orders at
$|k|=k_c-\varepsilon_B$. Fix $m\geq 1$ and $B,c>0$ and consider the following class of logarithmic parameter functions
\begin{equation}
 \mathcal U_{\mathrm{log},m}^{(\mathfrak r)}
 =\left\{\xi\in\Xi(\varepsilon_B):\ \begin{gathered}
 \xi=A\eta e^{\lambda \psi},\quad 0<A<\delta\mathfrak r,\\
 \psi\in C^m_{\rm rad}(\overline{B_c(\varepsilon_B)};\mathbb R),\quad
 \min_{\overline{B_c(\varepsilon_B)}}\psi>c,\quad
 \|\psi-\psi(0)\|_\infty<B
 \end{gathered}\right\}.
\end{equation}
Since $\lambda(0)=0$, $\lambda(k)<0$ for $k\ne0$, and $0\leq\eta\le1$, we have that
\begin{equation}
 \xi(0)=\|\xi\|_{C^0}=A,\qquad
 e^{t\lambda}\xi=A\eta e^{\lambda(\psi+t)}
 \in\mathcal U_{\mathrm{log},m}^{(\mathfrak r)}\quad(t\geq 0),
\end{equation}
and hence, the class $U_{\mathrm{log},m}^{(\mathfrak r)}$ is forward invariant under \eqref{reduced_dyn_linear}.\\
To bound $\Gamma_{\mathfrak r}$ uniformly on the class of logarithmic parameter functions, fix a closed annulus $D$ of
positive measure inside $0<|k|<k_c-\varepsilon_B$. On $D$, both $\eta$
and $-\lambda$ have positive minima and restricting \eqref{deftheta} to $D$ therefore gives the lower bound
$$
 S_\mu[A\eta]\geq 
 \left(\int_D\mu(k)\,dk\right)
 \left(\frac{A\min_D\eta}{\mathfrak r}\right)^{\nu/\min_D(-\lambda)}.
$$
Since the integrand in $S_\mu[\xi]$ acquires the additional factor $e^{-\nu \psi(k)}$, the bound $\|\psi-\psi(0)\|_\infty<B$ therefore gives
\begin{equation}
 |\theta[\xi]-\psi(0)-\theta[A\eta]|\leq B .
\end{equation}
Since $\log(\mathfrak r/A)\geq \log(1/\delta)>0$, combining these estimates gives
$\theta[\xi]-\psi(k)\leq C\log(\mathfrak r/A)$, for some constant $C>0$ depending only on $\delta,B,\eta$ and the fixed clock and spectral data.
For $-\lambda(k)\leq (2C)^{-1}$ we thus have
$$
 \frac{|\zeta[\xi](k)|}{\mathfrak r}
 =\frac{A\eta(k)}{\mathfrak r}
   e^{(-\lambda(k))(\theta[\xi]-\psi(k))}
 \leq \sqrt{A/\mathfrak r}\leq\sqrt\delta .
$$
The integral estimate in the proof of Lemma \ref{lem_invariant_clock_domain_s9} then applies with
$\beta=\sqrt\delta$. On the remaining frequencies, we use that
$x^{q-1}\leq 1+x^q$ for $x\geq 0$ and $q=\nu/(-\lambda(k))>1$.
By \eqref{clock_weight_assumption} and \eqref{compensated_clock_normalization_s9},
\begin{equation}\label{explicit_clock_class_constant_s9}
 \Gamma_{\mathfrak r}(\zeta[\xi])
 \leq \frac{1}{e\sqrt\delta\,\nu|\log\sqrt\delta|}+4C
 =:G_{\rm log}.
\end{equation}
Hence $\mathcal U_{\mathrm{log},m}^{(\mathfrak r)}
\subset\mathcal U_{\delta,G}^{(\mathfrak r)}$ for every $G>G_{\rm log}$,
uniformly in $A$, $\psi(0)$ and $\mathfrak r$.
\end{remark}

\subsection{Preparatory results for local density inversion}

To recover the slow amplitude locally from the mass density in
Theorem~\ref{prop_small_density_chart_section9}, we first establish the
local kinetic embedding in Theorem~\ref{prop_kinetic_embedding} and then
prove that the retracted density map has an invertible derivative in
Theorem~\ref{prop_density_transversality_s9}. The following lemma bounds
the compensated profile $\zeta[\xi]$ and shows that the initial state
$\hat{\mathcal W}[z;0]$ depends $C^1$ on $z$ in the graph norm of
$D(\mathcal L)$ from \eqref{fourier_generator_domain_s9}. These estimates
allow us to differentiate the kinetic evolution.

\begin{lemma}[Compensation and the expansion at auxiliary time zero]
\label{lem_selected_initial_lifting_s9}
Fix the parameters of Proposition~\ref{prop_majorant_functional_taylor},
the clock exponent $\nu\geq 2/\tau$, and the majorant radius $R_A$.
Then one can choose $0<R_{\rm init}<R_A$ and constants
$C_{\rm comp},C_{\rm init}>0$ such that, for every $0<\delta<1$
and every $\xi\in\mathcal U_\delta^\times$, the weighted compensation estimate
\begin{equation}\label{weighted_clock_estimates_s9}
 \|\zeta[\xi]\|_{\mathcal Y}
 +\|(-\lambda)^{-1/2}\zeta[\xi]\|_{\mathcal Y}
 +\|\zeta[\xi]/(-\lambda)\|_1\leq C_{\rm comp}\mathfrak r,
\end{equation}
holds. The fixed-time map $z\mapsto\hat{\mathcal W}[z;0]$ in
\eqref{selected_initial_lifting} is $C^1$ into $D(\mathcal L)$ with its
graph norm, on the open set
$\{z\in\Xi(\varepsilon_B):\|z\|_{\mathcal Y}<R_{\rm init}\}$
with the $C^0$ parameter norm. For sufficiently small $\mathfrak r$, the estimates
\begin{equation}\label{selected_initial_bounds_s9}
 \begin{aligned}
 \|\hat{\mathcal W}[\zeta[\xi];0]-\widehat e_{\zeta[\xi]}\|_{D(\mathcal L)}
 &\leq C_{\rm init}\mathfrak r^2,\\
 \|D_z\hat{\mathcal W}[\zeta[\xi];0]-(h\mapsto\widehat e_h)\|
                  _{\Xi(\varepsilon_B)\to\mathcal X_{k,v}}
 &\leq C_{\rm init}\mathfrak r,
 \end{aligned}
\end{equation}
hold. Here $D_z$ differentiates the argument of the fixed-time map,
before evaluation at $\zeta[\xi]$. The constants are independent of
$\xi$, $\theta[\xi]$ and $\mathfrak r$.
\end{lemma}
\begin{proof}
We first prove the weighted compensation estimate and then apply
the forcing majorants of Proposition~\ref{prop_majorant_functional_taylor}
to the selected initial conditions, followed by a proof of continuity in frequency to obtain the desired graph-norm estimate.\\

\noindent
\underline{The weighted compensation estimate}:\\
By \eqref{C1_clock_nu_assumption} and \eqref{def_lambda_min}, we can bound 
$\nu\geq 2/\tau>2\lambda_*$ and Young's inequality gives, for
$x\geq 0$ and $0<a\leq\lambda_*$,
$$
 \frac{x^2}{a}\leq x^{\nu/a}+a^{-\nu/(\nu-2a)},\qquad
 \sup_{0<a\leq\lambda_*}a^{-2a/(\nu-2a)}<\infty.
$$
Writing the second term as $a^{-1}a^{-2a/(\nu-2a)}$, we observe that the
factor $a^{-2a/(\nu-2a)}$ tends to one as $a\downarrow0$, since its
logarithm is $2a\log(1/a)/(\nu-2a)$ and $a\log(1/a)\to0$.
Because $2\lambda_*<\nu$, this factor extends continuously to
$[0,\lambda_*]$ and has a finite maximum $C$ there, so the second term
is bounded by $C/a$ for every $0<a\leq\lambda_*$.\\

We apply the pointwise inequality with $a=-\lambda(k)$ and
$x=|\zeta[\xi](k)|/\mathfrak r$. Upon integration over the frequency
band, the first term is controlled by the clock normalization
\eqref{compensated_clock_normalization_s9}, since $\mu\geq\mu_0$.
For the second term, the lower bound
$-\lambda(k)\geq(2\tau/\pi)|k|^2$ from \eqref{lambda_quadratic_bounds}
gives
$$
\begin{aligned}
 \frac{\|(-\lambda)^{-1/2}\zeta[\xi]\|_{L_k^2}^2}{\mathfrak r^2}
 &=\int_{B_c(\varepsilon_B)}\frac{1}{-\lambda(k)}
       \left(\frac{|\zeta[\xi](k)|}{\mathfrak r}\right)^2\,dk\\
 &\leq \frac{1}{\mu_0}\int_{B_c(\varepsilon_B)}\mu(k)
       \left(\frac{|\zeta[\xi](k)|}{\mathfrak r}\right)^{\nu/(-\lambda(k))}\,dk
       +\frac{C\pi|\mathbb S^{d-1}|}{2\tau}
          \int_0^{k_c-\varepsilon_B}r^{d-3}\,dr\\
 &=\frac{1}{\mu_0}+\frac{C\pi|\mathbb S^{d-1}|}{2\tau(d-2)}
                  (k_c-\varepsilon_B)^{d-2}.
\end{aligned}
$$

The preceding estimate thus bounds
$\|(-\lambda)^{-1/2}\zeta[\xi]\|_{L_k^2}$ by a fixed multiple of
$\mathfrak r$. Since $\zeta[\xi]$ is supported in $B_c(\varepsilon_B)$,
the Cauchy--Schwarz inequality gives
$$
 \begin{aligned}
 \|(-\lambda)^{-1/2}\zeta[\xi]\|_{L_k^1}
 &\leq |B_c(\varepsilon_B)|^{1/2}
        \|(-\lambda)^{-1/2}\zeta[\xi]\|_{L_k^2},\\
 \|\zeta[\xi]/(-\lambda)\|_{L_k^1}
 &\leq \left(\int_{B_c(\varepsilon_B)}
              \frac{dk}{-\lambda(k)}\right)^{1/2}
        \|(-\lambda)^{-1/2}\zeta[\xi]\|_{L_k^2}.
 \end{aligned}
$$
Since $\|g\|_{\mathcal Y}=\|g\|_{L_k^1}+\|g\|_{L_k^2}$, these bounds
and the estimate for $\|\zeta[\xi]\|_{\mathcal Y}$ in
\eqref{clock_Y_bound_direct} give \eqref{weighted_clock_estimates_s9}.\\

\noindent
\underline{The initial terms and their derivatives}:\\
On the region where we select specific non-zero initial conditions, the estimates from \eqref{Vn_selected_choice},
\eqref{Pk_zero_mass_gain} and \eqref{weighted_majorant_def} give the bounds
\begin{equation}
 \begin{aligned}
 \|\hat{\mathfrak V}_n(\bm k_n)\|_{\mathcal X_v^{(2)}}
 &\leq\frac{C_P|\Sigma\bm k_n|M_n(\bm k_n)}
                {2\lambda_{\min}+\lambda(\Sigma\bm k_n)}
 \leq C_{\mathfrak V}\lambda_{\min}^{-1/2}M_n(\bm k_n),\\
 C_{\mathfrak V}&:=\frac{C_P}{2-\varkappa_2}
                         \sqrt{\frac{\pi\varkappa_2}{2\tau}}.
 \end{aligned}
\end{equation}
Indeed, the denominator is at least $(2-\varkappa_2)\lambda_{\min}$,
and \eqref{lambda_quadratic_bounds} gives
$|\Sigma\bm k_n|^2\leq\pi\varkappa_2\lambda_{\min}/(2\tau)$.
Elsewhere the selected initial condition is zero.\\
To recover $B_n$ in \eqref{prop82_forcing_majorant}, we put
$u_i=\mathfrak u_{m_i}(\cdot;z)$ for $z\geq0$. Each block convolution
has at least two factors. Young's inequality and
\eqref{weighted_profile_majorant_final} give
$$
 \|u_0*\cdots*u_l\|_\infty
 \leq \|u_0\|_2\|u_1\|_2\prod_{i=2}^l\|u_i\|_1
 \leq \prod_{i=0}^l\bigl(A_{m_i}\|z\|_{\mathcal Y}^{m_i}\bigr),
$$
where the empty product is understood as one. Since $\sum_{i=0}^l m_i=n$, summing the
blocks with the coefficients of \eqref{prop82_forcing_majorant} yields
$$
 \left\|\frac1{n!}\int_{B_c(\varepsilon_B)^n}
 \delta(q-\Sigma\bm k_n)M_n(\bm k_n)
 \Pi z(\bm k_n)\,d\bm k_n\right\|_{L^{\infty}_q}
 \leq B_n\|z\|_{\mathcal Y}^n.
$$
Let now $f,g,z\in\mathcal Y$ be nonnegative and supported in
$B_c(\varepsilon_B)$. Introduce two scalar parameters $s,v\geq0$
and replace each factor $z(k_j)$ in the preceding integral by
$z(k_j)+s f(k_j)+v g(k_j)$. For each summed frequency $q$, the
result is a polynomial in $s$ and $v$ with nonnegative coefficients, where the coefficient of $s$ is obtained by choosing one of the $n$ factors to supply $f$, while the coefficient of $sv$ in turn is obtained by choosing one factor to supply $f$ and a different factor to supply $g$. By symmetry
of $M_n$, all $n$ choices in the first case, and all $n(n-1)$ ordered
choices in the second, give the same integral. Thus the coefficient
estimates, in the supremum norm over $q\in\mathbb R^d$, extend to
$$
\begin{aligned}
 &\left\|
 \frac{1}{(n-1)!}\int_{B_c(\varepsilon_B)^n}
 \delta(q-\Sigma\bm k_n)\,M_n(\bm k_n)\,
 f(k_1)\prod_{j=2}^n z(k_j)\,d\bm k_n
 \right\|_{L_q^\infty}
 \leq en B_n\|f\|_{\mathcal Y}\|z\|_{\mathcal Y}^{n-1},\\
 &\left\|
 \frac{1}{(n-2)!}\int_{B_c(\varepsilon_B)^n}
 \delta(q-\Sigma\bm k_n)\,M_n(\bm k_n)\,
 f(k_1)g(k_2)\prod_{j=3}^n z(k_j)\,d\bm k_n
 \right\|_{L_q^\infty}
 \leq e^2n^2 B_n\|f\|_{\mathcal Y}\|g\|_{\mathcal Y}
                   \|z\|_{\mathcal Y}^{n-2}.
\end{aligned}
$$
For $n=2$, the empty product in the second integral equals one.
For the second bound, we divide the diagonal estimate by $sv$ and
choose $s\|f\|_{\mathcal Y}=v\|g\|_{\mathcal Y}
=\|z\|_{\mathcal Y}/(n-2)$ when $n>2$ and then use the elementary estimate
$(1+2/(n-2))^{n-2}\leq e^2$. For $n=2$, let the two equal
quantities tend to infinity, obtaining the smaller factor $4$.\\
Since $\lambda_{\min}^{-1/2}\leq\sum_j(-\lambda(k_j))^{-1/2}$
off the null sets $k_j=0$, symmetry and the one-factor estimate give in turn
$$
 \left\|\frac1{n!}\int_{B_c(\varepsilon_B)^n}
 \delta(q-\Sigma\bm k_n)\hat{\mathfrak V}_n(\bm k_n)
 \Pi z(\bm k_n)\,d\bm k_n\right\|_{L_q^\infty(\mathcal X_v^{(2)})}
 \leq eC_{\mathfrak V}nB_n
       \|(-\lambda)^{-1/2}z\|_{\mathcal Y}\|z\|_{\mathcal Y}^{n-1}.
$$
Differentiating the degree-$n$ term in
\eqref{selected_initial_lifting} gives the normalization $1/(n-1)!$ and the weight $(-\lambda)^{-1/2}$ then falls either on $h$ or on one of
the $n-1$ remaining factors. These alternatives give the
normalizations $1/(n-1)!$ and $(n-1)/(n-1)!=1/(n-2)!$,
respectively, and hence
$$
 \begin{aligned}
 &\left\|\frac1{(n-1)!}\int_{B_c(\varepsilon_B)^n}
 \delta(q-\Sigma\bm k_n)\hat{\mathfrak V}_n(\bm k_n)
 h(k_1)\prod_{j=2}^nz(k_j)\,d\bm k_n
            \right\|_{L_q^\infty(\mathcal X_v^{(2)})}\\
 &\quad\leq C_{\mathfrak V}B_n\bigl(
 en\|(-\lambda)^{-1/2}h\|_{\mathcal Y}\|z\|_{\mathcal Y}^{n-1}
 +e^2n^2\|h\|_{\mathcal Y}\|(-\lambda)^{-1/2}z\|_{\mathcal Y}
                         \|z\|_{\mathcal Y}^{n-2}\bigr).
 \end{aligned}
$$

\noindent
\underline{Summation in the graph norm}:\\
By \eqref{Vn_selected_choice}, every nonlinear term just estimated
lies in $\operatorname{ran}P_k$ and
vanishes outside the fixed compact set
$\{|k|\leq k_0:-\lambda(k)\leq\varkappa_2\lambda_*\}$,
which is contained in $|k|<k_0$ by \eqref{selection_rates} and
the choice of $k_0$. Its supremum norm thus controls its $L^1$ norm, and $\mathcal L_k$ acts simply as multiplication by
$\lambda(k)$. Let $\mathfrak b_n$ denote any of the nonlinear terms just estimated. Once its
continuity is verified below, the support bound and the identity
$\mathcal L_q\mathfrak b_n(q)=\lambda(q)\mathfrak b_n(q)$ on its support give
$$
\begin{aligned}
 \|\mathfrak b_n\|_{D(\mathcal L)}
 &= \|\mathfrak b_n\|_{\mathcal X_{k,v}}+\|\mathcal L\mathfrak b_n\|_{\mathcal X_{k,v}} \\
 &\leq (1+\varkappa_2\lambda_*)
 \bigl(\|\mathfrak b_n\|_{L_q^\infty(\mathcal X_v^{(2)})}
       +\|\mathfrak b_n\|_{L_q^1(\mathcal X_v^{(2)})}\bigr) \\
 &\leq (1+\varkappa_2\lambda_*)(1+|B_{k_0}(0)|)
       \|\mathfrak b_n\|_{L_q^\infty(\mathcal X_v^{(2)})}.
\end{aligned}
$$
Fix $0<R_{\rm init}<R_A$. Since $B$ is analytic and starts at
order two, summation for $\|z\|_{\mathcal Y}<R_{\rm init}$ gives
\begin{equation}\label{selected_initial_mixed_bounds_s9}
 \begin{aligned}
 \|\hat{\mathcal W}[z;0]-\widehat e_z\|_{D(\mathcal L)}
 &\leq C_{\rm init}\|(-\lambda)^{-1/2}z\|_{\mathcal Y}\|z\|_{\mathcal Y},\\
 \|D_z\hat{\mathcal W}[z;0]h-\widehat e_h\|_{D(\mathcal L)}
 &\leq C_{\rm init}\bigl(
 \|(-\lambda)^{-1/2}h\|_{\mathcal Y}\|z\|_{\mathcal Y}
 +\|(-\lambda)^{-1/2}z\|_{\mathcal Y}\|h\|_{\mathcal Y}\bigr).
 \end{aligned}
\end{equation}
The sums $\sum_{n\geq 2}nB_nR_{\rm init}^{n-2}$ and
$\sum_{n\geq 2}n^2B_nR_{\rm init}^{n-2}$ are finite by the majorant argument from the convergence proof.
Moreover, $(-\lambda)^{-1/2}\mathbf1_{B_c(\varepsilon_B)}\in\mathcal Y$
by \eqref{lambda_quadratic_bounds} and $d\geq 3$, so the differentiated
series converges locally uniformly in the $C^0$ parameter norm.\\

\noindent
\underline{Continuity in frequency}.\\
Fix an order $n$ and choose $0<2\epsilon<k_c-\varepsilon_B$. Let
$\eta_\epsilon\colon\mathbb R^d\to[0,1]$ be a smooth radial cutoff function, equal
to zero for $|k|\leq\epsilon$ and to one for $|k|\geq2\epsilon$.
For this continuity argument only, multiply the integrand by
$\prod_{j=1}^n\eta_\epsilon(k_j)$, also in the differentiated integral.
This leaves the summed frequency $q=\Sigma \bm{k}_n$ unrestricted.\\
On the support of the modified integrand, every individual frequency
satisfies $-\lambda(k_j)\geq m_\epsilon:=-\lambda(\epsilon)>0$.
The same lower bound therefore holds for $\lambda_{\min}$ in every
nonempty block. We show inductively that, at each fixed order, $U_n$
is bounded there and the kernel is jointly continuous in the individual
frequencies and auxiliary time on finite time intervals. The assertion
at order one follows from \eqref{eigenvector}. At the next order,
the finite sum in \eqref{hat_N_n_explicit} gives a continuous forcing bounded by $Ce^{-2\lambda_{\min}s}$. For small $q\ne0$ the cutoff
in \eqref{Vn_selected_choice} equals one, since
$\lambda_{\min}\geq m_\epsilon$ and $\lambda(q)\to0$.
By \eqref{Pk_zero_mass_gain}, the selected initial coefficient has norm
at most $C|q|/((2-\varkappa_2)m_\epsilon)$, tending to its prescribed
value zero at $q=0$. Proposition~\ref{prop_weighted_hierarchy} bounds the
complementary part uniformly, and \eqref{prop_unsteady_cohomology_bound_eq}
bounds $U_n$ away from zero. The tail in \eqref{Vn_selected_choice} converges uniformly because $2\lambda_{\min}+\lambda(q)\geq(2-\varkappa_2)m_\epsilon$
where the selected initial condition is nonzero. There the projections are norm-continuous,
and the selected initial condition vanishes near $|\Sigma \bm{k}_n|=k_0$. Strong continuity of $T_{q}(s)$
in $(q,s)$ from \eqref{Duhamel_Tk} now gives joint continuity in
\eqref{sol_unsteady_cohomol}, completing the induction. Substituting
$k_n=q-\sum_{j=1}^{n-1}k_j$, dominated convergence on the compact
integration domain gives continuity in $q$ of the modified integral
and of its derivative applied to any fixed continuous variation $h$.\\
We now remove the auxiliary cutoff. Telescoping
$1-\prod_{j=1}^n\eta_\epsilon(k_j)$ gives terms each containing one
factor $1-\eta_\epsilon(k_j)$, supported in $|k_j|<2\epsilon$.
If the weight $(-\lambda)^{-1/2}$ falls on that factor, use the
one-factor estimate above. If it falls on another factor, use the
two-factor estimate. Bounding the remaining profiles by their
$C^0$ norms on $B_c(\varepsilon_B)$, the errors tend to zero uniformly
in $q$, locally uniformly in $z$, and uniformly for $\|h\|_{C^0}\leq1$,
since
$$
 \|\mathbf1_{\{|k|<2\epsilon\}}\|_{\mathcal Y}
 +\|(-\lambda)^{-1/2}\mathbf1_{\{|k|<2\epsilon\}}\|_{\mathcal Y}
 \longrightarrow 0.
$$
Thus the original integrals are continuous in $q$. The preceding
graph-norm estimate gives convergence also in $D(\mathcal L)$.
Each term is therefore a continuous polynomial into
$D(\mathcal L)$, and local uniform convergence of the series and its
derivatives gives the asserted $C^1$ regularity in the $C^0$ parameter
norm.

\end{proof}

For the following argument, fix once and for all a radius $0<R<R_{\rm init}$ and set
$$
 \mathfrak r_{\rm reg}:=\min\left\{\mathfrak r_*,
                  \frac{R}{2C_{\rm comp}}\right\}.
$$
For the density-inversion argument below we take
$0<\mathfrak r<\mathfrak r_{\rm reg}$, so that
$\|\zeta[\xi]\|_{\mathcal Y}<R$. The fixed constants in this argument
may depend on $R$ and the fixed data of Section~\ref{sec_proof_main_thm_proof},
including the clock exponent $\nu\geq 2/\tau$ and the even weight
$\mu\in L^\infty(B_c)$ satisfying $\int_{B_c}\mu(k)\,dk=1$ and
$\mu(k)\geq\mu_0>0$ almost everywhere, with the clock origin fixed
at $\theta_0=0$, but not on $\xi$, $\theta[\xi]$, $\mathfrak r$,
or $G$.

\begin{theorem}[Local kinetic charts and the differential of the parametrization]
\label{lem_flow_representation_s9}\label{prop_kinetic_embedding}
Under the hypotheses of Proposition~\ref{prop_C1_away_origin}, with the
selection \eqref{Vn_selected_choice} and $d\geq3$, fix a sufficiently small
$0<\mathfrak r<\mathfrak r_{\rm reg}$. Then the parametrization
$$
 \hat{\mathcal W}:
 \{\xi\in\Xi(\varepsilon_B):0<\|\xi\|_{C^0}<\mathfrak r\}
 \longrightarrow\mathcal X_{k,v}
$$
is a local real $C^1$ embedding. For each $\xi_0$ in this domain there
is an open neighborhood $V$ of $\xi_0$ on which $\hat{\mathcal W}$ is a
$C^1$ diffeomorphism onto a submanifold of $\mathcal X_{k,v}$ with the
topology of \eqref{closure_continuous_spaces_s9}. The full parametrized
family is forward invariant, and its inverse spatial Fourier transform consists of real kinetic states. For each such $\xi$, the curve $s\mapsto\hat{\mathcal W}[e^{\lambda s}\zeta[\xi];s]$
is a classical solution of \eqref{maineqfk} in $\mathcal X_{k,v}$ on every
interval $0\leq s\leq T<\infty$.
\end{theorem}

\begin{proof}
Fix $0<\|\xi\|_{C^0}<\mathfrak r$ and write $z=\zeta[\xi]$,
$\theta=\theta[\xi]$, and $\widehat{\mathcal W}_{\rm init}(z):=\hat{\mathcal W}[z;0]$.
By \eqref{closure_flow_identity_s9}, $\hat{\mathcal W}[\xi]=\Phi^\theta(\widehat{\mathcal W}_{\rm init}(z))$
for the local kinetic flow $\Phi$. Varying $\xi$ changes both the initial
profile and the evaluation time.
We first differentiate the kinetic flow in its initial data and in time.
Given the resulting kinetic variation $D\hat{\mathcal W}[\xi]h$, we solve
the linearized equation backwards to time zero. For $\alpha>0$, the
initial-profile term then lies in the slow subspace, while the scalar
$D\theta[\xi]h$ multiplies a vector outside it. A functional vanishing
on the slow subspace recovers this scalar and, after subtracting the clock
term, we recover $h$ from the density of the remaining initial-profile
variation. These operations give the bounded left inverse needed for
the local embedding.

 \noindent\underline{The classical evolution.}
Write $\hat f(s)=\hat{\mathcal W}[e^{\lambda s}\zeta[\xi];s]$ for the
auxiliary trajectory. By Remark~\ref{rem_auxiliary_zero_family}, it
satisfies the Duhamel formula of Section~\ref{sec_proof_main_thm_proof},
$$
 \hat f(s)=e^{s\mathcal L}\hat f(0)
 +\int_0^s e^{(s-\sigma)\mathcal L}
              \widehat{\mathcal N}_\alpha(\hat f(\sigma))\,d\sigma ,
$$
in $\mathcal Z^{(2)}$, where \eqref{prop82_coefficient_norm} and
\eqref{clock_Y_bound_direct} give
$$
 \sup_{s\geq0}\|\hat f(s)\|_{\mathcal Z^{(2)}}
 \leq A(\|\zeta[\xi]\|_{\mathcal Y})\leq C_{\rm evol}\mathfrak r,
 \qquad C_{\rm evol}=A'(R)C_{\rm comp}^{\mathcal Y}.
$$
To read the same formula in $\mathcal X_{k,v}$, it suffices to check the
forcing, since Lemma~\ref{lem_selected_initial_lifting_s9} already places
the initial profile in $D(\mathcal L)$. In each convolution of
\eqref{def_Nhat}, one density factor and the collision factor belong to
$L_k^2$, so their convolution is continuous and vanishes at infinity.
The remaining density factors act by convolution in $L_k^1$, and putting
all factors in $L_k^1$ also bounds the integral norm. Summing these
multilinear estimates on the ball $\|\hat f\|_{\mathcal Z^{(2)}}\leq A(R)$,
where $(2\pi)^{-d/2}C_\rho A(R)<1$, gives a continuous forcing with
$$
 \|\widehat{\mathcal N}_\alpha(\hat f)\|_{\mathcal X_{k,v}}
 \leq C_{\rm nl}\|\hat f\|_{\mathcal Z^{(2)}}^2 .
$$
Lemma~\ref{lem_kinetic_group_s9} now gives continuity of the Duhamel
integral in $\mathcal X_{k,v}$ on every finite interval. Since
\eqref{def_Nhat} is analytic on the open set
$(2\pi)^{-d/2}C_\rho\|\hat f\|_{L_k^1(\mathcal X_v^{(2)})}<1$,
regularity for initial data in $D(\mathcal L)$ gives a classical
solution there \cite[Theorem~1.16]{RevisionSchnaubelt2025}.

 \noindent\underline{Differentiating the parametrization.}\\
For a fixed time, we first vary the initial kinetic state. Given a kinetic direction $h_{\rm kin}\in\mathcal X_{k,v}$, termwise
differentiation of the convolution series \eqref{def_Nhat} gives
\begin{equation}
\begin{aligned}
 D\widehat{\mathcal N}_\alpha(\hat f)h_{\rm kin}
 =\frac1\tau\sum_{\ell\geq1}\binom{\alpha}{\ell}(2\pi)^{-d\ell/2}
 \Bigl((\rho\hat f)^{*\ell}*\mathcal L_Q h_{\rm kin}+\ell(\rho\hat f)^{*(\ell-1)}*(\rho h_{\rm kin})*
                              \mathcal L_Q\hat f\Bigr).
\end{aligned}
\end{equation}
Young's inequality, with the varied factor in $C_0$ and the other
factors in $L_k^1$, gives
\begin{equation}\label{linearized_convolution_bound_s9}
 \|D\widehat{\mathcal N}_\alpha(\hat f)h_{\rm kin}\|_{C_0(\mathcal X_v^{(2)})}
 \leq C_{\rm lin}\|\hat f\|_{L_k^1(\mathcal X_v^{(2)})}
                   \|h_{\rm kin}\|_{C_0(\mathcal X_v^{(2)})}.
\end{equation}
The summation uses the scalar majorant $\Psi_\alpha$ from the proof of
\eqref{prop82_nonlinear_bound} in Proposition~\ref{prop_majorant_functional_taylor}
and its derivative for the terms carrying the factor $\ell$.
For an initial variation $h_{\rm kin}(0)\in\mathcal X_{k,v}$, the
variational equation along $\hat f$ is
\begin{equation}\label{closure_exact_derivative_s9}
\begin{aligned}
 h_{\rm kin}(s,k)
 &=\bigl(D\Phi^s(\hat f(0))h_{\rm kin}(0)\bigr)(k)\\
 &=T_k(s)h_{\rm kin}(0,k)+\int_0^s T_k(s-\sigma)
       \bigl(D\widehat{\mathcal N}_\alpha(\hat f(\sigma))
                                      h_{\rm kin}(\sigma)\bigr)(k)\,d\sigma .
\end{aligned}
\end{equation}

Fix $T<\infty$ and let $\Delta\hat f$ denote the change in the solution
under the initial increment $\Delta\hat f(0)$. The local Lipschitz estimate
\cite[Theorem~1.11(c)]{RevisionSchnaubelt2025} gives
$\sup_{0\leq s\leq T}\|\Delta\hat f(s)\|_{\mathcal X_{k,v}}\leq C_T\|\Delta\hat f(0)\|_{\mathcal X_{k,v}}$.
Subtracting \eqref{closure_exact_derivative_s9} with $h_{\rm kin}(0)=\Delta\hat f(0)$
from the difference of the two Duhamel formulas shows that the Taylor remainder
of $\widehat{\mathcal N}_\alpha$ is $O_T(\|\Delta\hat f(0)\|_{\mathcal X_{k,v}}^2)$
on this interval. Gr\"onwall's inequality then gives
$$
 \sup_{0\leq s\leq T}\|\Delta\hat f(s)-h_{\rm kin}(s)\|_{\mathcal X_{k,v}}
       \leq C_T\|\Delta\hat f(0)\|_{\mathcal X_{k,v}}^2, 
$$
which identifies the solution of \eqref{closure_exact_derivative_s9} as the Fr\'echet derivative applied to the initial increment. Note that the dependence on the initial profile is continuous by the same integral equation. The group in
Lemma~\ref{lem_kinetic_group_s9} lets us solve the linearized equation
backwards on $[0,s]$, so $D\Phi^s(\widehat{\mathcal W}_{\rm init}(z))$ has a bounded inverse.\\
By Lemma~\ref{lem_selected_initial_lifting_s9}, $B=D\widehat{\mathcal W}_{\rm init}(z)$ is continuous
in $z$ in the graph norm, so $s\mapsto e^{s\mathcal L}B$ is
operator-norm continuous. The Volterra iterates in
\eqref{closure_exact_derivative_s9} preserve this continuity (use
$r=s-\sigma$). If $K_T$ bounds
$\|e^{(s-\sigma)\mathcal L}\|\,\|D\widehat{\mathcal N}_\alpha(\hat f(\sigma))\|$
for $0\leq\sigma\leq s\leq T$, integration over the ordered time simplex
bounds the $n$th Volterra power by $(K_T T)^n/n!$, as in
\cite[Lemma~III.1.9]{engel2000one}. The series therefore converges uniformly on $[0,T]$, and the
initial-profile derivative is jointly continuous in $z,s$.
Differentiating $\Phi^s(\hat f(t))=\hat f(t+s)$ at $t=0$ gives
$\dot{\hat f}(s)=D\Phi^s(\widehat{\mathcal W}_{\rm init}(z))\dot{\hat f}(0)$, with $\dot{\hat f}(0)$
continuous in $z$ by graph-norm regularity. The chain rule and
\eqref{C1_zeta_derivative} therefore give, for $h\in\Xi(\varepsilon_B)$,
\begin{equation}\label{selected_variation_initial_s9}
\begin{aligned}
 h_{\rm kin}(\theta)&=D\hat{\mathcal W}[\xi]h,\\
 h_{\rm kin}(0)&=(D\Phi^\theta(\widehat{\mathcal W}_{\rm init}(z)))^{-1}D\hat{\mathcal W}[\xi]h
       =B(e^{-\lambda\theta}h)+\Delta_{\rm clock}(z)\,D\theta[\xi]h,\\
 \Delta_{\rm clock}(z)&=\mathcal L \widehat{\mathcal W}_{\rm init}(z)+\widehat{\mathcal N}_\alpha(\widehat{\mathcal W}_{\rm init}(z))-B(\lambda z).
\end{aligned}
\end{equation}
Thus the scalar clock variation multiplies the difference between the
kinetic time derivative at $\widehat{\mathcal W}_{\rm init}(z)$ and the derivative induced by the
reduced velocity $\lambda z$.

 \noindent\underline{Recovering the parameter variation.}\\
For $\alpha>0$, introduce the
real slow subspace
$$
 E_{{\rm slow},k_0}:=
 \{\widehat e_H:H\in\Xi(k_c-k_0)\}\subset\mathcal X_{k,v}.
$$
It is closed because density recovers $H$ and lifting is bounded on
this band. The selection \eqref{Vn_selected_choice} puts $\widehat{\mathcal W}_{\rm init}(z)$ and
$\operatorname{ran}B$ in $E_{{\rm slow},k_0}$. We claim $\Delta_{\rm clock}(z)\notin E_{{\rm slow},k_0}$.
Otherwise, since $\mathcal L$ preserves $E_{{\rm slow},k_0}$,
\eqref{selected_variation_initial_s9} would place
$\widehat{\mathcal N}_\alpha(\widehat{\mathcal W}_{\rm init}(z))$ in $E_{{\rm slow},k_0}$ as well.
This collision term has zero density by \eqref{def_LQ}, so as a slow
profile it would vanish. In physical space, writing
$f=\mathcal F^{-1}\hat f(0)$, this would give
$$
 \bigl((1+\rho[f])^\alpha-1\bigr)\bigl(\rho[f]\phi-f\bigr)=0 .
$$
In fact, \eqref{compensated_clock_normalization_s9} gives
$\|\zeta[\xi]\|_{C^0}\geq\mathfrak r$, while
\eqref{selected_initial_bounds_s9} bounds the difference between
$\rho[\hat f(0)]$ and $\zeta[\xi]$ by $C_\rho C_{\rm init}\mathfrak r^2$.
Thus the initial density is nonzero for small $\mathfrak r$. Its inverse
Fourier transform is real analytic by compact Fourier
support, and is uniformly small by the $L_k^1$ estimate above. Its
nonzero set is therefore dense, and on that set the first factor is
nonzero for $\alpha>0$. Continuity would then give $f=\rho[f]\phi$
everywhere. Taking Fourier transforms contradicts the nonzero density,
since $\hat f(0)\in E_{{\rm slow},k_0}$ and $\hat e(k,\cdot)\ne\phi$ for $k\ne0$.\\
Put $d_\xi=\operatorname{dist}(\Delta_{\rm clock}(z),E_{{\rm slow},k_0})>0$. The real functional
$\ell(e+c\Delta_{\rm clock}(z))=c$ on $E_{{\rm slow},k_0}+\mathbb R\Delta_{\rm clock}(z)$ has norm $1/d_\xi$.
Hahn--Banach extends it to $\mathcal X_{k,v}$ with this norm
\cite[Corollary~1.2]{BGKproofBrezis}. Thus
$\ell(h_{\rm kin}(0))=D\theta[\xi]h$ and subtracting
$\ell(h_{\rm kin}(0))\Delta_{\rm clock}(z)$ leaves $B(e^{-\lambda\theta}h)$.
To recover its input, define the radial retraction
$\mathscr R:\mathcal Y_\rho\to\Xi(\varepsilon_B)$ for
$k=\kappa\omega$, $0<\kappa\leq k_c-\varepsilon_B$, $|\omega|=1$, by

\begin{equation}\label{closure_retraction_s9}
 (\mathscr R P)(\kappa\omega)=P(\kappa\omega)
   -\frac{\kappa}{k_c-\varepsilon_B}P\bigl((k_c-\varepsilon_B)\omega\bigr),
 \qquad \mathscr R\rho[\widehat e_h]=h .
\end{equation}
At zero take $(\mathscr R P)(0)=P(0)$ and extend by zero outside the
band. The formula preserves conjugate symmetry and gives
$\|\mathscr R\|\leq2$, while the identity on the lifting follows from
\eqref{selected_tangent_identity}. By \eqref{selected_initial_bounds_s9},
$\|\mathscr R\rho B-I\|\leq\|\mathscr R\rho\|C_{\rm init}\mathfrak r<1$
for sufficiently small $\mathfrak r$ and applying
$(\mathscr R\rho B)^{-1}\mathscr R\rho$ to the remaining vector recovers
$e^{-\lambda\theta}h$, while multiplication by $e^{\lambda\theta}$ then recovers $h$.
Each operation is bounded at the fixed $\xi$. For $\alpha=0$ the lifting is linear and $\mathscr R\rho$ recovers $h$ directly. At each nonzero $\xi_0$ the resulting left inverse $L$ splits
$\mathcal X_{k,v}$ into $\operatorname{ran}D\hat{\mathcal W}[\xi_0]$
and $\ker L$. The inverse function theorem for $L\hat{\mathcal W}$
gives local coordinates in which the image is a $C^1$ graph
\cite[Chapter~II, \S2, Proposition~2.2]{BGKproofLang}.
Forward invariance follows from \eqref{prop82_mild_identity}, and
conjugate symmetry gives reality after Fourier inversion.

\end{proof}

Theorem~\ref{prop_kinetic_embedding} recovers $h$ from the full kinetic
variation $D\hat{\mathcal W}[\xi]h$. On
$\mathcal U_{\delta,G}^{(\mathfrak r)}$, the bound on $D\theta[\xi]$
makes the nonlinear correction small enough to preserve invertibility
there. The next theorem follows the variation up to time $\theta[\xi]$,
when its density is $D\mathcal H_\rho(\xi)h$, and proves that the
retracted density still determines $h$. This supplies the derivative
needed for the local density coordinates of
Theorem~\ref{prop_small_density_chart_section9}.

\begin{theorem}[Invertibility of the retracted density derivative]
\label{prop_density_transversality_s9}

Assume the hypotheses of Theorem~\ref{prop_kinetic_embedding}, in particular
$d\geq3$ and $\nu\geq2/\tau$, and let $\mathfrak r_{\rm reg}$ be the
regularity radius used there. Fix
$0<\delta<1$ and $G>0$. Let $\mathcal H_\rho(\xi)=\rho[\hat{\mathcal W}[\xi]]$
be the full Fourier density map of Theorem~\ref{mainthm_formal}, part~3,
with values in $\mathcal Y_\rho$ from \eqref{closure_continuous_spaces_s9},
and let $\mathscr R$ be the radial retraction \eqref{closure_retraction_s9}.
There is a radius $0<\mathfrak r_{\rm inv}(\delta,G)\leq\mathfrak r_{\rm reg}$
such that, for every $0<\mathfrak r<\mathfrak r_{\rm inv}(\delta,G)$ and
every $\xi$ in the forward-invariant domain
$\mathcal U_{\delta,G}^{(\mathfrak r)}$ of \eqref{invariant_density_domain_s9},
the operator
$$
  \mathscr R D\mathcal H_\rho(\xi)\colon
  \Xi(\varepsilon_B)\longrightarrow\Xi(\varepsilon_B)
$$
is a bounded isomorphism of the real Banach space defined in
\eqref{defXi_eps_C0}.
\end{theorem}
\begin{proof}
Fix $\xi\in\mathcal U_{\delta,G}^{(\mathfrak r)}$ and a parameter variation
$h\in\Xi(\varepsilon_B)$ and let $h_{\rm kin}$ solve the linearized kinetic
equation \eqref{closure_exact_derivative_s9} with initial value
\eqref{selected_variation_initial_s9}, whose density at time $\theta[\xi]$ is
$D\mathcal H_\rho(\xi)h$. At time zero the leading term of $h_{\rm kin}(0)$
is just the linear lifting of $e^{-\lambda\theta[\xi]}h$ and we therefore follow
the linear map
$$
 e^{-\lambda\theta[\xi]}h\ \longmapsto\ \mathscr R\rho[h_{\rm kin}(t)],
 \qquad 0\leq t\leq\theta[\xi].
$$
For fixed $\xi$ the
clock is finite, and both $e^{-\lambda\theta[\xi]}$ and its reciprocal
are bounded continuous multipliers on $\Xi(\varepsilon_B)$. Once this map is
invertible at $t=\theta[\xi]$, multiplication by $e^{\lambda\theta[\xi]}$ recovers the original variation $h$.\\

\noindent
\underline{Invertibility at time zero.}\\
The leading contribution to $h_{\rm kin}(0)$ in
\eqref{selected_variation_initial_s9} comes from the initial-profile
term, which, by \eqref{selected_initial_bounds_s9}, differs from the
linear lifting $\widehat e_{e^{-\lambda\theta[\xi]}h}$ by at most
$C_{\rm init}\mathfrak r\|e^{-\lambda\theta[\xi]}h\|_\infty$.
To estimate the clock term, we use the identity
$\mathcal L\widehat e_{\zeta[\xi]}=\widehat e_{\lambda\zeta[\xi]}$,
which cancels the linear contributions to
$\Delta_{\rm clock}(\zeta[\xi])$. The graph-norm bound in
\eqref{selected_initial_bounds_s9} controls the remainder in the
generator term, while the corresponding remainder in the differentiated
lifting is controlled by applying \eqref{selected_initial_mixed_bounds_s9}
in the direction $\lambda\zeta[\xi]$ and using
\eqref{weighted_clock_estimates_s9}. Both remainders are of order
$\mathfrak r^2$, as is the quadratic collision term by the estimate
in Theorem~\ref{prop_kinetic_embedding}. We thus obtain
\begin{equation}
 \|\Delta_{\rm clock}(\zeta[\xi])\|_{\mathcal X_{k,v}}
 \leq C_{\rm defect}\mathfrak r^2 .
\end{equation}
In \eqref{selected_variation_initial_s9}, this defect is multiplied by
$D\theta[\xi]h$, whose absolute value on the stated domain is bounded by
$(G/\mathfrak r)\|e^{-\lambda\theta[\xi]}h\|_\infty$ according to
\eqref{compensated_clock_derivative_s9}. The clock term therefore
contributes at most
$C_{\rm defect}\mathfrak r G\|e^{-\lambda\theta[\xi]}h\|_\infty$.
Together with the estimate for the initial-profile term, this gives
\begin{equation}\label{selected_initial_tangent_cone_s9}
 \|h_{\rm kin}(0)-\widehat e_{e^{-\lambda\theta[\xi]}h}\|_{C_0(\mathcal X_v^{(2)})}
 \leq C_{\rm var,init}\mathfrak r(1+G)
       \|e^{-\lambda\theta[\xi]}h\|_\infty .
\end{equation}
Since $\mathscr R\rho$ recovers $e^{-\lambda\theta[\xi]}h$ exactly
from its linear lifting by \eqref{closure_retraction_s9}, the last
estimate shows that the map
$e^{-\lambda\theta[\xi]}h\mapsto\mathscr R\rho[h_{\rm kin}(0)]$
differs from the identity by at most
$2C_\rho C_{\rm var,init}\mathfrak r(1+G)$ in operator norm.
Choosing $\mathfrak r$ small enough that this bound is at most $1/2$,
we obtain invertibility from the Neumann series, with
$$
 \|e^{-\lambda\theta[\xi]}h\|_\infty
 \leq 2\|\mathscr R\rho[h_{\rm kin}(0)]\|_\infty .
$$

 \noindent\underline{Time integral of the kinetic trajectory.}\par\nobreak
As in the preceding proof, we write
$\hat f(s)=\hat{\mathcal W}[e^{\lambda s}\zeta[\xi];s]$ and recall that the estimate in  \eqref{linearized_convolution_bound_s9} bounds the forcing
in the linearized equation by a constant times
$\|\hat f(s)\|_{L_k^1(\mathcal X_v^{(2)})}
 \|h_{\rm kin}(s)\|_{C_0(\mathcal X_v^{(2)})}$.
Since the clock interval $[0,\theta[\xi]]$ may be arbitrarily long,
we have to control the accumulation of the first factor by proving that
\begin{equation}\label{selected_time_integral_s9}
 \int_0^\infty \|\hat f(s)\|_{L_k^1(\mathcal X_v^{(2)})}\,ds
 \leq C_{\rm time}\mathfrak r .
\end{equation}
To this end, we first establish the linear estimate that will be used for the
zero-density forcing in Duhamel's formula. Indeed, for any initial profile
$\hat g\in\mathcal Z^{(2)}$ with $\rho[\hat g]=0$, we claim that
\begin{equation}\label{zero_density_time_integral_s9}
 \int_0^\infty\|e^{s\mathcal L}\hat g\|_{L_k^1(\mathcal X_v^{(2)})}\,ds
 \leq C\|\hat g\|_{\mathcal Z^{(2)}} .
\end{equation}
The only difficulty occurs at low frequencies, where integration of the
slow component $e^{\lambda(k)s}P_k\hat g(k)$ in time produces the
factor $1/(-\lambda(k))$. For $0<|k|\leq k_c-\varepsilon_B$, the
zero-density estimate \eqref{Pk_zero_mass_gain} supplies a factor
$|k|$ in the projection, so the resulting weight is bounded by
$C/|k|$ by \eqref{lambda_quadratic_bounds}. Since $|k|^{-1}$ is
square integrable on this ball for $d\geq3$, the Cauchy--Schwarz
inequality bounds the slow contribution by
$C\|\hat g\|_{L_k^2(\mathcal X_v^{(2)})}$.\\ 
For the remaining modes, the semigroup decays uniformly at least as
fast as $e^{-\lambda_*s}$. Indeed, Lemma~\ref{lem_selected_semigroup}
gives the faster rate $\gamma>\lambda_*$ for the complementary
component up to $k_0$ and for the full semigroup beyond $k_0$, while
on $k_c-\varepsilon_B\leq |k|\leq k_0$ the slow component also
decays at rate at least $\lambda_*$ by Proposition~\ref{lemma_lambda}
and the projection bound \eqref{Pk_weighted_bound}. Integration in
time therefore bounds the remaining contribution by
$C\|\hat g\|_{L_k^1(\mathcal X_v^{(2)})}$, which together with the
low-frequency estimate proves \eqref{zero_density_time_integral_s9}.\\ 
We now apply Duhamel's formula to $\hat f$, beginning with the
evolution of its initial state $\hat{\mathcal W}[\zeta[\xi];0]$.
By the choice \eqref{Vn_selected_choice}, this state lies in the slow
subspace at each frequency $|k|\leq k_0$ and vanishes outside that
ball, implying that its evolution is given by the slow multiplier
$e^{\lambda(k)s}$. Separating the linear lifting
$\widehat e_{\zeta[\xi]}$ from the remainder, whose supremum norm is
at most $C_{\rm init}\mathfrak r^2$ by
\eqref{selected_initial_bounds_s9}, and integrating the multiplier
in time, we obtain
$$
 \int_0^\infty\|e^{s\mathcal L}\hat f(0)\|_{L_k^1(\mathcal X_v^{(2)})}\,ds
 \leq C\left(\left\|\frac{\zeta[\xi]}{-\lambda}\right\|_1
       +\mathfrak r^2\int_{|k|\leq k_0}\frac{dk}{-\lambda(k)}\right)
 \leq C\mathfrak r .
$$
Here \eqref{weighted_clock_estimates_s9} bounds the first term by
$C\mathfrak r$, while \eqref{lambda_quadratic_bounds} bounds the
integrand in the second by $C|k|^{-2}$, which is integrable on the
ball for $d\geq3$.\\ 
For the nonlinear term in Duhamel's formula, Young's inequality
applied to \eqref{def_Nhat} gives
$$
 \|\widehat{\mathcal N}_\alpha(\hat f(\sigma))\|_{\mathcal Z^{(2)}}
 \leq C\|\hat f(\sigma)\|_{L_k^1(\mathcal X_v^{(2)})}
        \|\hat f(\sigma)\|_{\mathcal Z^{(2)}} .
$$
The proof of Theorem~\ref{prop_kinetic_embedding} bounds the second
factor uniformly by $C_{\rm evol}\mathfrak r$. Since the forcing has
zero density at each time $\sigma$, we can use
\eqref{zero_density_time_integral_s9} to estimate its subsequent
linear evolution. Indeed, after taking norms in Duhamel's formula and
integrating up to time $S$, Tonelli's theorem allows us to integrate
first over $\sigma\leq s\leq S$ with $\sigma$ fixed. Applying
\eqref{zero_density_time_integral_s9} in the elapsed time $s-\sigma$
and combining the result with the bound for the initial state yields
$$
 \int_0^S\|\hat f(s)\|_{L_k^1(\mathcal X_v^{(2)})}\,ds
 \leq C_2\mathfrak r+C_3\mathfrak r
           \int_0^S\|\hat f(\sigma)\|_{L_k^1(\mathcal X_v^{(2)})}\,d\sigma .
$$
Choosing $\mathfrak r$ so that $C_3\mathfrak r\leq1/2$, we absorb
the last term into the left-hand side and obtain a bound independent
of $S$. Letting $S\to\infty$ proves \eqref{selected_time_integral_s9}
and, in particular, gives the required estimate on every clock interval,
\begin{equation}\label{integrated_variational_smallness_s9}
 C_{\rm lin}\int_0^{\theta[\xi]}
       \|\hat f(\sigma)\|_{L_k^1(\mathcal X_v^{(2)})}\,d\sigma
 \leq C_{\rm lin}C_{\rm time}\mathfrak r .
\end{equation}

 \noindent\underline{Recovering the slow amplitude from the density.}\\
Fix now $0\leq t\leq\theta[\xi]$. To recover the initial amplitude
$e^{-\lambda\theta[\xi]}h$ from the retracted density at time $t$, we
first estimate the complementary component forward from time zero.
This allows us to determine, up to a controlled error, how much of
the density at time $t$ comes from the slow component. We then follow
that component back to time zero using its scalar evolution equation.
Throughout the argument, we use the bound
$$
 M=\sup_{0\leq s\leq t} e^{\lambda_*s}
            \|h_{\rm kin}(s)\|_{C_0(\mathcal X_v^{(2)})}.
$$
By the time-zero estimate,
$\|e^{-\lambda\theta[\xi]}h\|_\infty\leq 4C_\rho M$.
The leading lifting in \eqref{selected_initial_tangent_cone_s9} lies
in the range of $P_k$ on the input band and vanishes on its boundary
and outside. Thus only the error in that estimate contributes to the
initial complementary component inside the band or to the initial
variation on and outside its boundary. Both are consequently bounded
by $C\mathfrak r(1+G)M$. For these components, the forward semigroup
from $\sigma$ to $s$ is bounded by $Ce^{-\lambda_*(s-\sigma)}$, so
multiplication of the Duhamel formula \eqref{closure_exact_derivative_s9}
by $e^{\lambda_*s}$ leaves the weight $e^{\lambda_*\sigma}$ in the
forcing integral. By the definition of $M$ and
\eqref{linearized_convolution_bound_s9}, this integral is bounded by
$$
\begin{aligned}
 e^{\lambda_*s}\int_0^s C e^{-\lambda_*(s-\sigma)}
    \bigl\| & D\widehat{\mathcal N}_\alpha(\hat f(\sigma))
               h_{\rm kin}(\sigma)  \bigr\|_{C_0(\mathcal X_v^{(2)})}\,d\sigma\\
 &\leq C C_{\rm lin}\int_0^s
    \|\hat f(\sigma)\|_{L_k^1(\mathcal X_v^{(2)})}
    e^{\lambda_*\sigma}\|h_{\rm kin}(\sigma)\|_{C_0(\mathcal X_v^{(2)})}\,d\sigma\\
 &\leq C C_{\rm lin}M\int_0^s
    \|\hat f(\sigma)\|_{L_k^1(\mathcal X_v^{(2)})}\,d\sigma
 \leq C C_{\rm lin}C_{\rm time}\mathfrak r M .
\end{aligned}
$$
Combining \eqref{integrated_variational_smallness_s9} with the initial
bound therefore gives, for $0\leq s\leq t$,
$$
\begin{aligned}
 e^{\lambda_*s}\|(I-P_k)h_{\rm kin}(s,k)\|_{\mathcal X_v^{(2)}}
     &\leq C\mathfrak r(1+G)M, &&|k|<k_c-\varepsilon_B,\\
 e^{\lambda_*s}\|h_{\rm kin}(s,k)\|_{\mathcal X_v^{(2)}}
     &\leq C\mathfrak r(1+G)M, &&|k|\geq k_c-\varepsilon_B .
\end{aligned}
$$
These two bounds control precisely the terms that distinguish the
retracted density at time $t$ from the slow amplitude
$\rho[P_kh_{\rm kin}(t,k)]$. Indeed, for $k=r\omega$ with
$0<r<k_c-\varepsilon_B$, the retraction formula
\eqref{closure_retraction_s9} gives
$$
\begin{aligned}
 \rho[P_k h_{\rm kin}(t,k)]-(\mathscr R\rho[h_{\rm kin}(t)])(k)
  =-\rho[(I-P_k)h_{\rm kin}(t,k)]+\frac{r}{k_c-\varepsilon_B}\,
       \rho[h_{\rm kin}(t,(k_c-\varepsilon_B)\omega)] .
\end{aligned}
$$
The first bound controls the complementary density in the interior,
while the second controls the boundary value subtracted by the
retraction. Since the boundary term vanishes at $k=0$, the difference
is bounded throughout the input band by
$Ce^{-\lambda_*t}\mathfrak r(1+G)M$.\\
We may now solve the scalar equation for the slow amplitude backward
from $t$ to any $s\in[0,t]$. Since $-\lambda(k)\leq\lambda_*$ on the
input band, its backward multiplier satisfies
$$
 e^{\lambda_*s}e^{\lambda(k)(s-\sigma)}
 \leq e^{\lambda_*\sigma},\qquad s\leq\sigma\leq t .
$$
Thus the same weight that controls the forward evolution also controls
the backward forcing integral. Applying the projection bound
\eqref{Pk_weighted_bound} and the convolution estimate
\eqref{linearized_convolution_bound_s9}, and then using
\eqref{integrated_variational_smallness_s9}, we bound this integral
by $C\mathfrak r M$. Together with the comparison at time $t$, this
yields
\begin{equation}\label{endpoint_comparison_s9}
 e^{\lambda_*s}|\rho[P_kh_{\rm kin}(s,k)]|
 \leq e^{\lambda_*t}\|\mathscr R\rho[h_{\rm kin}(t)]\|_\infty
     +C\mathfrak r(1+G)M,\qquad |k|<k_c-\varepsilon_B .
\end{equation}
To bound the full variation, recall that on the input band
$P_kh_{\rm kin}(s,k)=\rho[P_kh_{\rm kin}(s,k)]\hat e(k)$, with
$\hat e(k)$ uniformly bounded. Adding the forward estimate for the
complementary component inside the band and using the corresponding
estimate for the full variation outside gives, upon taking the
supremum in $s$ and $k$,
$$
 M\leq Ce^{\lambda_*t}\|\mathscr R\rho[h_{\rm kin}(t)]\|_\infty
             +C\mathfrak r(1+G)M .
$$
For fixed $\delta,G$, we choose
$\mathfrak r_{\rm inv}(\delta,G)\leq\mathfrak r_{\rm reg}$ small enough
to absorb the last term into the left-hand side. It remains to compare
the slow amplitude at time zero with the initial parameter variation.
By \eqref{selected_initial_tangent_cone_s9}, their difference on the
input band satisfies
$$
 \bigl|\rho[P_kh_{\rm kin}(0,k)]-(e^{-\lambda\theta[\xi]}h)(k)\bigr|
 \leq C\mathfrak r(1+G)M .
$$
Taking $s=0$ in \eqref{endpoint_comparison_s9} and substituting the
bound just obtained for $M$, we find that
$\|e^{-\lambda\theta[\xi]}h\|_\infty$ is bounded by
$e^{\lambda_*t}\|\mathscr R\rho[h_{\rm kin}(t)]\|_\infty$ times
$1+C\mathfrak r(1+G)$. Decreasing the radius once more therefore gives
$$
 \|e^{-\lambda\theta[\xi]}h\|_\infty
 \leq 2e^{\lambda_*t}\|\mathscr R\rho[h_{\rm kin}(t)]\|_\infty ,
$$
with a choice of radius independent of $t$ and $\xi$.\\

 \noindent\underline{Continuation to the observation time.}\\
It remains to show that the map
$$
 e^{-\lambda\theta[\xi]}h\longmapsto\mathscr R\rho[h_{\rm kin}(t)],
$$
which is invertible at time zero, remains invertible up to the
observation time $t=\theta[\xi]$. For fixed $\xi$, we first verify
that this map depends continuously on $t$ in operator norm.
For the initial-profile term in \eqref{selected_variation_initial_s9},
this continuity follows from Theorem~\ref{prop_kinetic_embedding}.
In the clock term, the functional
$e^{-\lambda\theta[\xi]}h\mapsto D\theta[\xi]h$ is fixed and bounded
by \eqref{compensated_clock_derivative_s9}, while the vector
$D\Phi^t(\hat f(0))\Delta_{\rm clock}(\zeta[\xi])$ depends continuously
on $t$. Their product is therefore continuous in operator norm, and
this property is preserved by the bounded operator $\mathscr R\rho$.\\
We may now apply the method of continuity
\cite[\S5.2]{BGKproofGilbarg}. The set of invertible times contains
zero and is open in $[0,\theta[\xi]]$ by norm continuity. To see that
it is also closed, recall that the preceding estimate bounds the
inverse norm at every invertible time by the same constant
$2e^{\lambda_*\theta[\xi]}$. At any limit point of such times, we can
therefore choose an invertible time so close that the two operators
differ in norm by less than $(2e^{\lambda_*\theta[\xi]})^{-1}$.
The Neumann series then gives invertibility at the limit point.
Since the interval is connected, the map is invertible throughout
$[0,\theta[\xi]]$. At $t=\theta[\xi]$, the retracted density variation is
$\mathscr R D\mathcal H_\rho(\xi)h$ by
\eqref{selected_variation_initial_s9}. Composing the invertible map
at this endpoint with the isomorphism
$h\mapsto e^{-\lambda\theta[\xi]}h$ of $\Xi(\varepsilon_B)$ gives the
asserted invertibility of $\mathscr R D\mathcal H_\rho(\xi)$.

\end{proof}

\subsection{Exact spectral closure and forward continuation}

We now complete the closure argument on the
forward-invariant domain $\mathcal U_{\delta,G}^{(\mathfrak r)}$ from
\eqref{invariant_density_domain_s9}. With the density space $\mathcal Y_\rho$
of \eqref{closure_continuous_spaces_s9} equipped with the supremum norm,
we use Theorem~\ref{prop_density_transversality_s9} to recover the slow
amplitude $\xi$ locally from the full Fourier density. The density then
determines the kinetic state on the nonlinear hydrodynamic manifold and
hence its mass flux, so the continuity equation closes in terms
of the density alone.

\begin{theorem}[Local density coordinates]\label{prop_small_density_chart_section9}
Fix $0<\delta<1$ and $G>0$ and use the radius
$\mathfrak r_{\rm inv}(\delta,G)$ of Theorem~\ref{prop_density_transversality_s9}.
For every $0<\mathfrak r<\mathfrak r_{\rm inv}(\delta,G)$, the full Fourier
density map
\begin{equation}\label{density_chart_section9}
 \mathcal H_\rho:\mathcal U_{\delta,G}^{(\mathfrak r)}
       \longrightarrow\mathcal Y_\rho,\qquad
 \mathcal H_\rho(\xi)=\rho\hat{\mathcal W}[\xi]
\end{equation}
is a local $C^1$ embedding. Every $\xi_0\in\mathcal U_{\delta,G}^{(\mathfrak r)}$
has an open, convex neighborhood $V\subset\mathcal U_{\delta,G}^{(\mathfrak r)}$
in the $C^0$ parameter topology on which $\mathscr R\circ\mathcal H_\rho$ is
a $C^1$ diffeomorphism onto an open set $O\subset\Xi(\varepsilon_B)$.
Here $\mathscr R$ is the retraction \eqref{closure_retraction_s9}.
The full density $R\in\mathcal H_\rho(V)$ determines the parameter by
\begin{equation}\label{closure_density_inverse_s9}
 (\mathcal H_\rho|_V)^{-1}(R)
 =\bigl((\mathscr R\circ\mathcal H_\rho)|_V\bigr)^{-1}(\mathscr R R),
 \qquad R\in\mathcal H_\rho(V).
\end{equation}
The right-hand side gives a $C^1$ extension of the inverse to the open
set $\mathscr R^{-1}(O)\subset\mathcal Y_\rho$.
For the same fixed $\mathfrak r$, every point of the forward orbit
$e^{t\lambda}\xi_0$, $t\geq 0$, has such a chart.
\end{theorem}
\begin{proof}
By Theorem~\ref{prop_kinetic_embedding} and the boundedness of the
density functional \eqref{defCrho}, $\mathcal H_\rho$ is $C^1$.
The derivative of $K_\rho=\mathscr R\circ\mathcal H_\rho$ at $\xi_0$
is an isomorphism by Theorem~\ref{prop_density_transversality_s9}.
The inverse function theorem gives an open neighborhood of $\xi_0$
on which $K_\rho$ is a $C^1$ diffeomorphism, so we may choose an open, convex
neighborhood $V$ inside the guaranteed neighborhood and set $O=K_\rho(V)$. Applying
$(K_\rho|_V)^{-1}$ to $\mathscr R R$ for $R=\mathcal H_\rho(\xi)$,
$\xi\in V$, recovers $\xi$ and proves \eqref{closure_density_inverse_s9}. The same composition
is $C^1$ on $\mathscr R^{-1}(O)$ and, by Theorem~\ref{prop_density_transversality_s9}, the operator
$(\mathscr R D\mathcal H_\rho(\xi_0))^{-1}\mathscr R$ is a bounded left
inverse of $D\mathcal H_\rho(\xi_0)$. The graph argument at the end of
the proof of Theorem~\ref{prop_kinetic_embedding} therefore shows that
$\mathcal H_\rho$ is a local $C^1$ embedding. Finally, Lemma~\ref{lem_invariant_clock_domain_s9} keeps
$e^{t\lambda}\xi_0$ in $\mathcal U_{\delta,G}^{(\mathfrak r)}$ at every
finite forward time.
\end{proof}

We now apply the local inverse to the spectral-closure construction of
Section~\ref{sec_spectral_closure_general}.

\begin{corollary}[The nonlinear density closure]
Let $V$ be a neighborhood from
Theorem~\ref{prop_small_density_chart_section9} and, for $R\in\mathcal H_\rho(V)$, take
$\xi=(\mathcal H_\rho|_V)^{-1}(R)$. On this chart, the constitutive law
\eqref{nonlinear_constitutive_law_density} gives the kinetic perturbation
with Fourier density $R$. Its mass flux and the closed density equation
are given in the same Fourier variables by
\begin{align}
 \widehat{\mathcal C}_{\rm spec}^{\rho}[R]
 &=\hat{\mathcal W}[\xi],\qquad
 \rho\widehat{\mathcal C}_{\rm spec}^{\rho}[R]=R,\quad \widehat j_V[R](k) =\int v\,\widehat{\mathcal C}_{\rm spec}^{\rho}[R](k,v)\,dv,\label{density_spectral_closure_section9}\\
 \partial_{t }R(k)
 &=\bigl(D\mathcal H_\rho(\xi)[\lambda\xi]\bigr)(k)
=-\ri k\cdot\widehat j_V[R](k).
 \label{density_reduced_dynamics_section9}
\end{align}
The constitutive law is $C^1$ into $\mathcal X_{k,v}$ and, for $R_0=\mathcal H_\rho(\xi_0)$ with $\xi_0\in V$, the unique
$C^1$ density trajectory in this chart is
\begin{equation}\label{explicit_global_density_flow_section9}
 R(t )=\mathcal H_\rho(e^{t \lambda}\xi_0)
\end{equation}
for as long as $e^{t \lambda}\xi_0$ stays in $V$.
\end{corollary}
\begin{proof}
The composition in \eqref{density_spectral_closure_section9} is $C^1$
by Theorem~\ref{prop_small_density_chart_section9} and
Theorem~\ref{lem_flow_representation_s9}.
Differentiating \eqref{explicit_global_density_flow_section9}
and using \eqref{reduced_dyn_linear} gives
$\partial_tR=D\mathcal H_\rho(\xi)[\lambda\xi]$ right away. Conversely, for any $C^1$ solution $R(t)\in\mathcal H_\rho(V)$, its inverse
$\xi(t )=(\mathcal H_\rho|_V)^{-1}(R(t ))$ is $C^1$ and satisfies
$$
 D\mathcal H_\rho(\xi(t ))
 \bigl[\partial_{t }\xi(t )-\lambda\xi(t )\bigr]=0.
$$
By Theorem~\ref{prop_density_transversality_s9}, the derivative $D\mathcal H_\rho(\xi(t))$
is injective and thus $\partial_{t }\xi=\lambda\xi$, as in
\eqref{reduced_dyn_linear}, and $\xi(t )=e^{t \lambda}\xi_0$, which proves uniqueness using only $C^1$ regularity.\\
For the flux identity, we take the density of the classical kinetic
equation from Theorem~\ref{lem_flow_representation_s9}.
Since the collision terms have zero density by \eqref{Lk} and \eqref{def_Nhat}, only the transport term remains, which gives
the second equality in \eqref{density_reduced_dynamics_section9}.
The weight in $\mathcal X_v^{(2)}$ controls the first velocity moment,
while membership in $D(\mathcal L)$ gives $-\ri k\cdot\widehat j_V[R]$ in
$\mathcal Y_\rho\cap\mathcal Y$ by
\eqref{fourier_generator_domain_s9}.
Thus the identity also holds after Fourier inversion.
\end{proof}

For a density perturbation $r$ with
$\hat r=R\in\mathcal H_\rho(V)$, the physical-space kinetic perturbation
$\mathcal C_{\rm spec}^{\rho}[r]$ is the inverse spatial Fourier transform
of $\widehat{\mathcal C}_{\rm spec}^{\rho}[R]$ in
\eqref{density_spectral_closure_section9}, with the convention
\eqref{def_Fourier_transform}.
Reintroducing the global Maxwellian background state gives
\begin{equation}\label{physical_spectral_closure_s9}
 \rho_{\rm phys}=\rho_0(1+r),\qquad
 F_{\rm spec}[\rho_{\rm phys}]
 =\rho_0\bigl(\phi+\mathcal C_{\rm spec}^{\rho}
                  [\rho_{\rm phys}/\rho_0-1]\bigr),
\end{equation}
and the exact continuity equation is
$$
 \partial_{t }\rho_{\rm phys}
 +\nabla_x\cdot\int vF_{\rm spec}[\rho_{\rm phys}]\,dv=0.
$$
These formulas use the full convergent family, including all nonlinear
harmonics at the summed frequencies, rather than a truncated constitutive series.\\

\subsection{Positivity and hydrodynamic entropy}

The preceding results show that the family of distribution functions constructed through the successive solution of the homological equation forms a local $C^1$ manifold and that the spectral parameter function $\xi$ can be locally inverted for the mass density $\rho$. We now show that the corresponding distribution functions forming the hydrodynamic manifold remain between two positive multiples of the background
Maxwellian when the parameters are sufficiently small. The
homological initial conditions \eqref{Vn_selected_choice} give these
bounds at the initial time of the auxiliary evolution, and comparison
propagates them along the kinetic solution. We then prove finiteness
and dissipation of the relative entropy, implying that these properties also
hold for the density closure. We note that related comparison results on the torus
are proved in \cite[Theorem~2.1 and Proposition~3.1]{evans2026density}.
We stress again at this point that here we work on $\mathbb R^d$ and use relative entropy because the
homogeneous background has infinite spatial mass.

\begin{theorem}[Positivity and relative entropy on the hydrodynamic manifold]
\label{thm_selected_maxwellian_bounds}
Assume the hypotheses of Proposition~\ref{prop_majorant_functional_taylor},
with the homological initial conditions \eqref{Vn_selected_choice}. There are constants $C_{\rm pos}>0$ and
$0<\mathfrak r_{\rm pos}\leq\mathfrak r_*$, depending only on the fixed
data of Proposition~\ref{prop_majorant_functional_taylor}, with $C_{\rm pos}\mathfrak r_{\rm pos}<1/2$, such that
the following holds. For $0<\mathfrak r<\mathfrak r_{\rm pos}$ and every
$\xi\in\Xi(\varepsilon_B)$ with $\|\xi\|_{C^0}\leq\mathfrak r$,
let $f_{\rm hydro}$ denote the kinetic solution of \eqref{maineqf}
on the hydrodynamic manifold with initial parameter $\xi$. The reduced
flow \eqref{reduced_dyn_linear} gives
$$
 f_{\rm hydro}(0)=\mathcal F^{-1}\hat{\mathcal W}[\xi],\qquad
 f_{\rm hydro}(t)=\mathcal F^{-1}\hat{\mathcal W}[e^{\lambda t}\xi],
 \quad t\geq0.
$$
As in \eqref{deviation_ansatz}--\eqref{maineqf}, the full distribution is given by
$F=\rho_0(\phi+f_{\rm hydro})$, with $\rho_0>0$ and the effective relaxation time
$\tau$. For every $t\geq 0$, almost everywhere,
\begin{equation}\label{selected_uniform_maxwellian_bounds}
 \begin{aligned}
 (1-C_{\rm pos}\mathfrak r)\rho_0\phi
 &\leq F(t)\le(1+C_{\rm pos}\mathfrak r)\rho_0\phi,\\
 (1-C_{\rm pos}\mathfrak r)\rho_0
 &\leq\rho[F(t)]\le(1+C_{\rm pos}\mathfrak r)\rho_0.
 \end{aligned}
\end{equation}
For every admissible original velocity space $\mathcal X_v$,
the same solution also satisfies
\begin{equation*}\tag{\thetheorem a}\label{hydro_uniform_hilbert_s94}
 f_{\rm hydro}\in C([0,\infty);L_x^2L_v^2(\phi^{-1})),\qquad
 \sup_{t\geq 0}\|f_{\rm hydro}(t)\|_{L_x^2L_v^2(\phi^{-1})}\leq C_{\rm pos}\mathfrak r.
\end{equation*}

For $\Phi(y)=y\log y$, the relative entropy \eqref{def_H_Phi_relative}
is finite and satisfies
\begin{equation}\label{selected_entropy_quadratic_bounds}
 \begin{aligned}
 \frac{\rho_0}{2(1+C_{\rm pos}\mathfrak r)}
       \|f_{\rm hydro}(t)\|_{L_x^2L_v^2(\phi^{-1})}^2
 &\leq\mathcal H_{\Phi,\rho_0}^{\rm rel}[F(t)]\\
 &\leq\frac{\rho_0}{2(1-C_{\rm pos}\mathfrak r)}
       \|f_{\rm hydro}(t)\|_{L_x^2L_v^2(\phi^{-1})}^2
 \leq C_{\rm pos}^2\rho_0\mathfrak r^2.
 \end{aligned}
\end{equation}
Along this trajectory the relative entropy is locally absolutely
continuous and nonincreasing. On every finite time interval it satisfies
the integrated form of \eqref{relative_convex_entropy_dissipation}, with
collision rate $\tau^{-1}(\rho[F]/\rho_0)^\alpha$ in the normalization
of \eqref{maineqf}. Since $\mathcal S[F(t)]\geq0$,
the decomposition \eqref{relative_convex_entropy_decomposition} also shows
that the relative entropy of the mass density with respect to $\rho_0$
is bounded above by $\mathcal H_{\Phi,\rho_0}^{\rm rel}[F(t)]$.
On the density charts of
Theorem~\ref{prop_small_density_chart_section9}, these conclusions apply
to \eqref{physical_spectral_closure_s9} and continue along the selected
parameter branch.
\end{theorem}

\begin{proof}
For $\xi=0$ all the assertions follow from $f_{\rm hydro}=0$.
For $\xi\ne0$, let $t_0=-\theta[\xi]$. For the proof we extend
$f_{\rm hydro}$ to $t\geq t_0$ using the mild evolution of
Remark~\ref{rem_auxiliary_zero_family}, with auxiliary time $s=t-t_0$.
\begin{equation}\label{positive_auxiliary_majorant}
 \begin{aligned}
 f_{\rm hydro}(t_0+s)=\mathcal F^{-1}\hat{\mathcal W}[e^{\lambda s}\zeta[\xi];s],
 \qquad s\geq0,\\
 \sup_{s\geq 0}\|\hat{\mathcal W}[e^{\lambda s}\zeta[\xi];s]\|_{\mathcal Z^{(2)}}
 \leq A\bigl(C_{\rm comp}^{\mathcal Y}\mathfrak r\bigr).
 \end{aligned}
\end{equation}
The bound follows from \eqref{prop82_coefficient_norm} and
\eqref{clock_Y_bound_direct}, including $\|\xi\|_{C^0}=\mathfrak r$.
The clock cocycle \eqref{clock_cocycle} identifies this
evolution with the trajectory in the statement for $t\geq0$. Conjugate
reflection commutes with the semigroups, projections and the selection
\eqref{Vn_selected_choice}, so $f_{\rm hydro}$ is real on $[t_0,\infty)$.

\noindent
\underline{The initial Maxwellian bounds}:\\
At $t=t_0$, the selection \eqref{Vn_selected_choice}
and $k_c-\varepsilon_B<k_0$ place every term in the hydrodynamic
eigenspace with summed frequency $|k|<k_0$, thus implying
\begin{equation}
 \hat{\mathcal W}[\zeta[\xi];0](k,v)
 =\rho[\hat{\mathcal W}[\zeta[\xi];0]](k)\hat e(k,v)
 \quad(|k|<k_0),
\end{equation}
and the initial Fourier transform vanishes outside this ball of frequencies.
The density bound \eqref{defCrho}, \eqref{positive_auxiliary_majorant}
and the eigenvector formula \eqref{eigenvector} give
$$
 \|\rho[\hat{\mathcal W}[\zeta[\xi];0]]\|_{\mathcal Y}
 \leq C_\rho A\bigl(C_{\rm comp}^{\mathcal Y}\mathfrak r\bigr),\qquad
 |\hat e(k,v)|\leq\frac{\phi(v)}{1+\tau\lambda(k_0)}.
$$
Here $1+\tau\lambda(k_0)>0$. Fourier inversion and the Hilbert-space
Plancherel identity consequently yield
\begin{equation}\label{positive_initial_bounds}
 \begin{aligned}
 |f_{\rm hydro}(t_0,x,v)|
 &\leq\frac{(2\pi)^{-d/2}C_\rho
                A(C_{\rm comp}^{\mathcal Y}\mathfrak r)}
               {1+\tau\lambda(k_0)}\,\phi(v)
 \leq C_{\rm pos}\mathfrak r\phi(v),\\
 \|f_{\rm hydro}(t_0)\|_{L_x^2L_v^2(\phi^{-1})}
 &\leq\frac{C_\rho A(C_{\rm comp}^{\mathcal Y}\mathfrak r)}
              {1+\tau\lambda(k_0)}
 \leq C_{\rm pos}\mathfrak r.
 \end{aligned}
\end{equation}
To choose $C_{\rm pos}$, fix $0<R_{\rm pos}<R_A$. Since $A(0)=0$ and its coefficients are nonnegative,
$A(r)\leq A'(R_{\rm pos})r$ for $0\leq r\leq R_{\rm pos}$ and both inequalities in \eqref{positive_initial_bounds} therefore hold with
$$
 C_{\rm pos}:=\frac{\max\{1,(2\pi)^{-d/2}\}C_\rho
                 C_{\rm comp}^{\mathcal Y}A'(R_{\rm pos})}
                      {1+\tau\lambda(k_0)}.
$$
We then set
$$
 \mathfrak r_{\rm pos}:=\min\left\{\mathfrak r_*,
       \frac{R_{\rm pos}}{2C_{\rm comp}^{\mathcal Y}},
       \frac1{4C_{\rm pos}}\right\}.
$$
This choice ensures convergence and $C_{\rm pos}\mathfrak r_{\rm pos}<1/2$,
under the compensation hypothesis $\nu\geq2\lambda_*$ of
Section~\ref{sec_proof_main_thm_proof}. The stronger condition
$\nu\geq2/\tau$ in \eqref{C1_clock_nu_assumption} is used for density inversion.

\noindent
\underline{Propagation by comparison}:\\
Set $\varepsilon=C_{\rm pos}\mathfrak r$. Fourier inversion and
\eqref{positive_auxiliary_majorant} give
$$
 \sup_{t\geq t_0}\|\rho[f_{\rm hydro}(t)]\|_{L_x^\infty}
 \leq (2\pi)^{-d/2}C_\rho A(C_{\rm comp}^{\mathcal Y}\mathfrak r)
 \leq\varepsilon<\tfrac12 .
$$
Fourier inversion of the mild equation gives
\eqref{maineqf} in distributions with locally integrable source.
Pullback along $t\mapsto x+(t-t_0)v$ gives an absolutely continuous representative
on almost every characteristic, so variation of constants gives
$$
\begin{aligned}
 &\frac{f_{\rm hydro}(t,x+(t-t_0)v,v)}{\phi(v)}\\
 &\quad=\exp\!\left[-\frac1\tau\int_{t_0}^{t}
      (1+\rho[f_{\rm hydro}(\ell)](x+(\ell-t_0)v))^\alpha\,d\ell\right]
       \frac{f_{\rm hydro}(t_0,x,v)}{\phi(v)}\\
 &\qquad+\int_{t_0}^{t}\exp\!\left[-\frac1\tau\int_{\sigma}^{t}
      (1+\rho[f_{\rm hydro}(\ell)](x+(\ell-t_0)v))^\alpha\,d\ell\right]\\
 &\qquad\quad{}\times\frac{(1+\rho[f_{\rm hydro}(\sigma)](x+(\sigma-t_0)v))^\alpha}{\tau}
                 \,\rho[f_{\rm hydro}(\sigma)](x+(\sigma-t_0)v)\,d\sigma .
\end{aligned}
$$
The integral weight is the $\sigma$-derivative of the exponential and is
nonnegative. Its integral plus the initial weight is one. Thus the
characteristic formula averages values in $[-\varepsilon,\varepsilon]$,
by \eqref{positive_initial_bounds} and the density bound above, giving
$|f_{\rm hydro}(t,x,v)|\leq\varepsilon\phi(v)$ almost everywhere.
For $t\geq0$, this gives the Maxwellian bounds in
\eqref{selected_uniform_maxwellian_bounds} by the definition of $F(t)$, while
integrating them in $v$ gives the corresponding density bounds.

\noindent
\underline{The Hilbert-space estimate}:\\
By \eqref{positive_initial_bounds},
$f_{\rm hydro}(t_0,x,v)/\phi(v)\in L^2(\phi(v)\,dv\,dx)$.
For every $\sigma\geq t_0$,
$$
 \left(\iint |\rho[f_{\rm hydro}(\sigma)](x+(\sigma-t_0)v)|^2\phi(v)\,dv\,dx\right)^{1/2}
 =\|\rho[f_{\rm hydro}(\sigma)]\|_{L_x^2}
 \leq C_\rho A(C_{\rm comp}^{\mathcal Y}\mathfrak r),
$$
where in the first equality we change variables from $x$ to $x+(\sigma-t_0)v$ at
fixed $v$ and use $\int\phi\,dv=1$, while for the inequality we apply the scalar
Plancherel identity to the density and use \eqref{defCrho} as well as
\eqref{positive_auxiliary_majorant}.
Since both exponentials in the characteristic formula are at most one,
Minkowski's inequality yields
$$
 \|f_{\rm hydro}(t)\|_{L_x^2L_v^2(\phi^{-1})}
 \leq\|f_{\rm hydro}(t_0)\|_{L_x^2L_v^2(\phi^{-1})}
   +\frac{(1+\varepsilon)^\alpha}{\tau}
     \int_{t_0}^{t}\|\rho[f_{\rm hydro}(\sigma)]\|_{L_x^2}\,d\sigma .
$$
Along the characteristic, \eqref{maineqf} reads
$$
 \begin{aligned}
 &\frac{d}{dt}\frac{f_{\rm hydro}(t,x+(t-t_0)v,v)}{\phi(v)}
 =\frac{(1+\rho[f_{\rm hydro}(t)](x+(t-t_0)v))^\alpha}{\tau}\\
 &\qquad{}\times\left(\rho[f_{\rm hydro}(t)](x+(t-t_0)v)
        -\frac{f_{\rm hydro}(t,x+(t-t_0)v,v)}{\phi(v)}\right).
 \end{aligned}
$$
The preceding bounds make the right-hand side bounded in
$L^2(\phi(v)\,dv\,dx)$ on each finite time interval. The characteristic
profile is therefore absolutely continuous in that space and strong
continuity of transport translations then gives
$f_{\rm hydro}\in C([t_0,\infty);L_x^2L_v^2(\phi^{-1}))$.
Differentiating the squared norm of the characteristic profile, changing
back to $x$, and using
$\int(f_{\rm hydro}(t)-\rho[f_{\rm hydro}(t)]\phi)\,dv=0$
gives
\begin{equation}\label{positive_hilbert_energy_identity}
 \begin{aligned}
 \|f_{\rm hydro}(t)\|_{L_x^2L_v^2(\phi^{-1})}^2 +\frac2\tau\int_{t_0}^{t}\iint
     (1+\rho[f_{\rm hydro}(\sigma)])^\alpha
     \frac{|\mathcal L_Q f_{\rm hydro}(\sigma)|^2}{\phi}\,dv\,dx\,d\sigma =\|f_{\rm hydro}(t_0)\|_{L_x^2L_v^2(\phi^{-1})}^2.
 \end{aligned}
\end{equation}
Together with \eqref{positive_initial_bounds}, this proves the uniform bound in \eqref{hydro_uniform_hilbert_s94} for $t\geq0$.

\noindent
\underline{The entropy balance}:\\
It remains to justify the entropy identities of Section~\ref{sec_main_equation}. With $F=\rho_0(\phi+f_{\rm hydro})$ and
$\Phi(y)=y\log y$,
Taylor's formula in \eqref{def_H_Phi_relative} and the bounds
\eqref{selected_uniform_maxwellian_bounds} give
\eqref{selected_entropy_quadratic_bounds}. Indeed, on the interval between
$\rho_0$ and $F/\phi$, the second derivative $\Phi''(y)=1/y$ lies between
$[\rho_0(1+\varepsilon)]^{-1}$ and $[\rho_0(1-\varepsilon)]^{-1}$.\\
To prove the bounds on the relative entropy, we abbreviate $\Psi(z)=(1+z)\log(1+z)-z$ and note that
$$
\begin{aligned}
 \mathcal H_{\Phi,\rho_0}^{\rm rel}[F(t)]
 &=\iint\!\left[F(t)\log\frac{F(t)}{\rho_0\phi}-F(t)+\rho_0\phi\right]dv\,dx\\
 &=\rho_0\iint\phi(v)\,\Psi\!\left(\frac{f_{\rm hydro}(t,x+(t-t_0)v,v)}{\phi(v)}\right)dv\,dx ,
\end{aligned}
$$
where the second equality uses translation in $x$. On
$[-\varepsilon,\varepsilon]$, $\Psi'(z)=\log(1+z)$ has Lipschitz
constant $(1-\varepsilon)^{-1}$ and $\Psi'(0)=0$. Thus the integrated
Taylor remainder between two characteristic profiles is bounded by
$\rho_0/[2(1-\varepsilon)]$ times the squared
$L^2(\phi(v)\,dv\,dx)$ distance between them. The linear term is bounded
by Cauchy--Schwarz since $|\Psi'(z)|\leq |z|/(1-\varepsilon)$.
Along the locally absolutely continuous characteristic profile these
bounds give local absolute continuity of the entropy and its derivative
by the chain rule.\\
Returning to the variable $x$, we subtract
$\log(\rho[F]/\rho_0)$ from $\log(F/(\rho_0\phi))$ in the collision
integral, since $\int(F-\rho[F]\phi)\,dv=0$ at each $x$. This is the
cancellation used in Appendix~\ref{app_equilibria_characterization}
to obtain \eqref{eq_stationary_entropy_production}. With $F=F(t)$,
we obtain for almost every $t\geq0$
$$
 \frac{d}{dt}\mathcal H_{\Phi,\rho_0}^{\rm rel}[F(t)]
 =-\frac1\tau\iint\!\left(\frac{\rho[F]}{\rho_0}\right)^\alpha
 (F-\rho[F]\phi)\log\frac{F}{\rho[F]\phi}\,dv\,dx .
$$
The entropy production is nonnegative by monotonicity of the logarithm.
Since the logarithm is Lipschitz on $[1-\varepsilon,1+\varepsilon]$ with
constant $(1-\varepsilon)^{-1}$, the entropy production is integrable on
finite intervals by \eqref{positive_hilbert_energy_identity}.
Integrating the derivative on any finite interval in $[0,\infty)$ gives
the integrated balance \eqref{relative_convex_entropy_dissipation} for
$F$ with effective relaxation time $\tau$.
Restricting these conclusions to the density charts gives the
final assertion.
\end{proof}

\begin{corollary}[Decay along the selected trajectories]\label{cor_selected_decay}
Under the hypotheses of Proposition~\ref{prop_majorant_functional_taylor},
fix $0<\mathfrak r<\mathfrak r_*$ and let
$\xi\in\Xi(\varepsilon_B)$ with $\|\xi\|_{C^0}\leq\mathfrak r$.
There is a constant $C$ independent of $t$ and $\xi$ such that
\begin{equation}\label{hydro_fourier_decay_s94}
 \|\hat{\mathcal W}[e^{\lambda t}\xi]\|_{\mathcal Z^{(2)}}
 \leq C(1+t)^{-d/4}\|\xi\|_{C^0},\qquad t\geq 0 .
\end{equation}
For any admissible $\mathcal X_v$, let $\mathfrak r_{*,H}$ denote the
convergence radius for $H=L_v^2(\phi^{-1})$ with the same clock data,
frequency cutoffs and selection parameters. If $\mathfrak r<\min\{\mathfrak r_{\rm pos},\mathfrak r_{*,H}\}$,
then the corresponding full distribution satisfies
\begin{equation}\label{hydro_entropy_decay_s94}
 \mathcal H_{\Phi,\rho_0}^{\rm rel}[F(t)]
 \leq C\rho_0(1+t)^{-d/2}\|\xi\|_{C^0}^2,\qquad \Phi(y)=y\log y .
\end{equation}
\end{corollary}
\begin{proof}
If $\xi=0$, both assertions follow from $\hat{\mathcal W}[0]=0$ and
$F=\rho_0\phi$. Recall the linear lifting in \eqref{selected_tangent_identity}, 
$\widehat e_\xi(k,\cdot)=\xi(k)\hat e(k,\cdot)$. On the support of $\xi$, the eigenvector formula
\eqref{eigenvector} and the monotonicity in Proposition~\ref{lemma_lambda} give
$|\hat e(k,v)|\leq\phi(v)/(1+\tau\lambda(k_c-\varepsilon_B))$,
where the denominator is strictly positive. Taking first the
$\mathcal X_v^{(2)}$ norm and then the sum of the $L_k^1$ and $L_k^2$
norms, as in \eqref{def_Xv2_aux} and \eqref{def_majorant_space_Y}, thus gives
$$
 \|\widehat e_\xi\|_{\mathcal Z^{(2)}}
 \leq \frac{\|\phi\|_{\mathcal X_v^{(2)}}}
                  {1+\tau\lambda(k_c-\varepsilon_B)}\,\|\xi\|_{\mathcal Y}.
$$
For the nonlinear remainder, the proof of
Proposition~\ref{prop_selected_tangency} applies the mixed estimate
\eqref{selected_mixed_bound} at each order $n\geq2$. One factor is
$e^{(\varkappa-1)\lambda\theta[\xi]}|\xi|$, while the other $n-1$
factors are bounded in $\mathcal Y$ by
$\|\zeta[\xi]\|_{\mathcal Y}\leq C_{\rm comp}^{\mathcal Y}\mathfrak r$
from \eqref{clock_Y_bound_direct}. Summing over $n$, as in the derivation
of \eqref{selected_remainder_bound}, therefore yields
$$
 \|\hat{\mathcal W}[\xi]-\widehat e_\xi\|_{\mathcal Z^{(2)}}
 \leq C\bigl\|e^{(\varkappa-1)\lambda\theta[\xi]}|\xi|\bigr\|_{\mathcal Y}
 \leq C\|\xi\|_{\mathcal Y}.
$$
The majorant series for $A'$ converges at $C_{\rm comp}^{\mathcal Y}\mathfrak r$,
so the constant is uniform on the parameter ball. The last
inequality uses $\lambda\leq0$, $\varkappa>1$ and $\theta[\xi]\geq0$,
which make the exponential at most one. Adding the lifting and remainder
estimates and evaluating at $e^{\lambda t}\xi$ now gives
$\|\hat{\mathcal W}[e^{\lambda t}\xi]\|_{\mathcal Z^{(2)}}
\leq C\|e^{\lambda t}\xi\|_{\mathcal Y}$.
Since $e^{\lambda t}\xi$ is supported in $\overline{B_c(\varepsilon_B)}$,
Cauchy--Schwarz bounds its $L_k^1$ norm by
$|B_c(\varepsilon_B)|^{1/2}\|e^{\lambda t}\xi\|_{L_k^2}$ and, for $t\geq1$, \eqref{lambda_quadratic_bounds} gives the estimate
$$
 \|e^{\lambda t}\xi\|_{L_k^2}^2
 \leq \|\xi\|_{C^0}^2\int_{\mathbb R^d}e^{-4\tau t|k|^2/\pi}\,dk
 \leq C t^{-d/2}\|\xi\|_{C^0}^2 .
$$
The last inequality follows by scaling $k$ by $\sqrt{t}$. On the other hand, for
$0\leq t\leq1$, we have that $\|e^{\lambda t}\xi\|_{L_k^2}
\leq |B_c(\varepsilon_B)|^{1/2}\|\xi\|_{C^0}$. Combining the two estimates gives \eqref{hydro_fourier_decay_s94}.\\
To pass from Fourier decay to entropy decay, we use
\eqref{hydro_fourier_decay_s94} in $H=L_v^2(\phi^{-1})$. Under the
additional radius assumption the series converges both in the original
construction and in the one based on $H$. To identify their sums, compare
the kernels in $L_v^1$, into which both moment spaces embed by
\eqref{defCrho}. On common data the linear semigroups agree by uniqueness
in \eqref{Duhamel_Tk}, and the slow projections agree by
\eqref{Pk_rank_one_corrected}. Starting from the common eigenvector
\eqref{eigenvector}, agreement at lower orders gives the same forcing
in \eqref{eq:renormalized_nonlinearity}, hence the same selected initial
condition in \eqref{Vn_selected_choice} and the same next kernel by
\eqref{sol_unsteady_cohomol}, which implies that the two series give the same
trajectory. By Plancherel, $\|f_{\rm hydro}(t)\|_{L_x^2H}$ is the
$L_k^2(H)$ norm of its Fourier transform, which \eqref{def_Xv2_aux}
bounds by the $\mathcal Z^{(2)}$ norm in \eqref{hydro_fourier_decay_s94}
with $\mathcal X_v=H$. Squaring this decay estimate in the upper bound
of \eqref{selected_entropy_quadratic_bounds} gives
\eqref{hydro_entropy_decay_s94}.
\end{proof}

\begin{remark}
For the three-dimensional vorticity equation, \cite{GallayWayne2002R3} characterize a class of more rapidly decaying solutions by the vanishing of leading asymptotic coefficients, which suggests asking which low-frequency cancellations could sharpen Corollary~\ref{cor_selected_decay}.
\end{remark}

\section{Discussion and Further Perspectives}\label{sec_discussion}

\subsection{Summary and Discussion}

We constructed a nonlinear hydrodynamic manifold for the density-dependent
BGK equation on $\mathbb R^d$. The family is forward invariant, and its
parametrization is differentiable at equilibrium, where its derivative
recovers the slow hydrodynamic modes of the linearized equation. At each
nonzero amplitude in the open parameter ball, the parametrization is
locally a $C^1$ embedding. On a smaller forward-invariant set, the mass
density provides local coordinates on the manifold. In these macroscopic coordinates,
the density determines the kinetic state and hence its flux, yielding an
exact, self-consistent, nonlocal evolution equation for the mass density alone.
For sufficiently small amplitudes, the corresponding distributions remain
between positive multiples of the background Maxwellian and converge to it
at an algebraic rate. Their relative entropy satisfies the kinetic
dissipation balance and also decays algebraically to zero.\\
The principal new analytical and conceptual contribution of this work is the \emph{parametrization of time} for resonances with the kinetic essential spectrum. The classical steady homological equations can already fail at
quadratic order because of essential resonances, see
Proposition~\ref{remark_stationary_quadratic_integrability}, which gives smooth,
arbitrarily small parameters for which the reconstructed stationary quadratic term is not in kinetic $L^2$, although the individual
frequencies remain strictly subcritical. We instead introduce the spectral clock $\theta[\xi]$ and solve the unsteady homological equations by semigroup evolution. The spectral clock estimates compensate for the time dependence of the coefficient kernels and yield a convergent generalized functional expansion by the method of majorants, incorporating the continuous-spectrum resonances in the construction. This mechanism differs from the discrete-resonance corrections in spectral-submanifold theory \cite{cabre2005parameterization,KogelbauerDelaLlave2026Resonant} and complements earlier exact closures for linear BGK models \cite{KogelbauerKarlinClosure2025}.\\
The main result of our work also separates exact kinetic-hydrodynamic reductions from asymptotic hydrodynamics. The leading Chapman--Enskog equation is a nonlinear diffusion equation whose qualitative behavior depends on the diffusion exponent and on the class of initial densities (Section~\ref{sec_CE_wellposedness}), while the exact reduced flow on the hydrodynamic manifold is defined for all forward times. On the invariant domain \eqref{invariant_density_domain_s9}, local density coordinates give an exact autonomous expression, continued along the orbit through successive charts. The leading Chapman--Enskog equation \eqref{prorous_media_eq}
is a truncation in spatial gradients, whereas
\eqref{density_reduced_dynamics_section9} retains the full constitutive
law on the hydrodynamic manifold. At a positive homogeneous background, the Fourier symbol of the linearized diffusion equation agrees to quadratic order in $k$ with the hydrodynamic eigenvalue. This agreement alone gives no estimate of the truncation error for arbitrary times or data. In the context of Hilbert's sixth problem, our contribution is an
invariant-manifold realization of hydrodynamics for a nonlinear kinetic
model at fixed $\tau>0$. Hydrodynamic limits for the Boltzmann equation
address convergence of kinetic solutions to macroscopic equations as
the Knudsen number tends to zero under a prescribed scaling of space
and time \cite{SaintRaymond2009HydrodynamicLimits,GorbanKarlin2014Hilbert}, while
our reduction instead gives exact dynamics on the hydrodynamic manifold without taking this limit, i.e., for any Knudsen number.\\
Several restrictions of the construction remain. The
selected spectral parametrization is differentiable at equilibrium, with continuity of its derivative there still open for physical amplitudes when $\alpha>0$ (Remark~\ref{rem_selected_not_C1}). The local $C^1$ charts
are therefore constructed away from equilibrium, and higher regularity
is not established. Density inversion holds locally on
\eqref{invariant_density_domain_s9}, without a global inverse or a
uniform chart size. The parameters must be small and their individual
frequencies must stay below a fixed subcritical cutoff, although all
nonlinear summed frequencies are retained. The assumption $d\geq 3$
is used in the convolution and clock estimates and the available bounds
do not by themselves establish the corresponding theorem in lower
dimensions. Finally, invariance and tangency describe evolution on
the family but do not prove either exponential or polynomial attraction
of nearby kinetic states. Such a result requires estimates for
deviations from the family and cannot be inferred from its geometric
construction alone. These restrictions raise questions about
regularity, uniqueness, density coordinates and attraction, which we
discuss in the next section.

\subsection{Further Perspectives}

A natural next step is to construct invariant hydrodynamic families for kinetic models that conserve momentum and energy in
addition to mass. The full BGK and Shakhov models provide concrete
starting points, since their hydrodynamic spectral branches have
already been analyzed \cite{KogelbauerKarlinSpectral2024,kogelbauer2024spectral}. Further examples include the ellipsoidal statistical BGK model of
Holway \cite{Holway1966}
and BGK models for gas mixtures \cite{AndriesAokiPerthame2002}.
In these models the scalar density parameter is replaced by several
hydrodynamic variables, and the homological equations couple the
corresponding spectral branches. The full nonlinear Boltzmann equation is, of course,
the principal longer-term goal. While its small-frequency hydrodynamic
branches are classical \cite{EllisPinsky1975}, its collision
operator is no longer a finite-rank perturbation of constant-rate
transport and, for hard spheres, the velocity-dependent collision
frequency also changes the geometry of the essential spectrum.
Extending the present method would therefore require semigroup and
nonlinear estimates that control interactions among several
hydrodynamic branches and their essential resonances.\\
An extension to the Boltzmann equation would also give a way to study how the lifespan of a smooth fluid description depends on the
passage to a hydrodynamic limit. For hard cutoff potentials, the
incompressible Navier--Stokes equations arise rigorously from the Boltzmann
equation in a diffusive scaling~\cite{GolseSaintRaymond2009}. A recent
preprint \cite{AlpogeBuckmaster2026NavierStokes}, however, reports finite-time blow-up from smooth
data for the unforced three-dimensional incompressible equations, while in the compressible case, \cite{MerleEtAl2022ImplosionII} have proved that, for suitable barotropic
pressure laws, smooth finite-energy data can produce a finite-time implosion
at which the density becomes unbounded. Our BGK result thus suggests asking whether the full Boltzmann equation admits a
forward-invariant family near a Maxwellian on which mass, momentum and
energy determine an exact hydrodynamic evolution for all forward times at
fixed Knudsen number, while the first-order Chapman--Enskog equation becomes singular in finite time.\\
Another, more technical question is how the choice of parametrization
affects uniqueness and regularity. The present construction fixes
linear reduced dynamics and selects the homological initial conditions
\eqref{Vn_selected_choice}. Other admissible choices may give different invariant families, or different parametrizations of the same family.
It remains to determine which differences can be removed by a change of parameters and whether an appropriate choice improves regularity
at equilibrium. Allowing nonlinear terms in the reduced dynamics,
as in the parametrization method \cite{cabre2005parameterization},
could absorb internally resonant contributions and might reduce the
regularity loss in the parametrization. The explicit spectral clock in Definition~\ref{def_global_clock},
given by \eqref{deftheta}, uses the linear reduced flow. Nonlinear
reduced dynamics would require a new clock analysis, and essential
resonances would still have to be treated. Extension of density
inversion to the whole punctured parameter ball and global density
injectivity remain open. It is also unresolved whether the failure
of $C^1$ regularity in Remark~\ref{rem_selected_not_C1} persists on
the conjugate-symmetric parameter space. A stable foliation or an
asymptotic-phase map would provide a separate means of relating the
motion on the family to the relaxation of nearby kinetic states, which are further interesting and promising directions to investigate.
Removing the subcritical cutoff $\varepsilon_B$ on the individual frequencies and extending the density chart beyond the small-amplitude regime remain open. Nonlinear summed frequencies at the critical wave number are already retained by the unsteady construction, see Proposition~\ref{remark_stationary_quadratic_integrability}. Addressing these together with attraction would determine how much of the kinetic evolution can be represented by the exact hydrodynamic closure.\\
A different direction could be to change the underlying spatial domain: the role of continuous spatial frequencies could be isolated by studying the model on $\mathbb T^d$. Discrete wave vectors restrict the resonance relations to a lattice, while each nonzero spatial mode of the present model still has essential spectrum arising from the continuous velocity variable.\\
The parametrization of time may extend beyond kinetic systems as it suggests a way to study invariant families for nonlinear PDEs with continuous spectrum on
$\mathbb R^d$. We propose to solve the unsteady homological equations in terms of a 
time variable by Duhamel's formula, as in \eqref{sol_unsteady_cohomol_W},
which keeps the forcing inside a time integral even where the stationary
resonance denominator vanishes. To this end, a suitable representation of time as a functional of the parameter functions has to be defined. Three natural examples are the
incompressible Navier--Stokes equations \cite{Kato1984Strong}, viscous
incompressible magnetohydrodynamics with equal viscosity and magnetic
diffusivity \cite{AgapitoSchonbek2007}, and the parabolic--elliptic
Keller--Segel system \cite{BlanchetDolbeaultPerthame2006}. In each case
the linearization at zero evolves Fourier modes by $e^{-D|k|^2t}$ for
some $D>0$. Orthogonal input frequencies $k$ and $\ell$ already give a
quadratic resonance, since their decay rates add to the decay rate at
$k+\ell$. These rates give a candidate clock of the form \eqref{deftheta}, and
the question is whether the resulting evolution integrals can be bounded
uniformly near the resonant set. Self-similar formulations provide a point of
comparison in the diffusive setting \cite{wayne1997invariant,gallay2002invariant}. We also mention the relation to \cite[Thm.~3]{DuroCarpio2001}, where the authors obtain a logarithm of physical time in the second asymptotic profile of a critical convection-diffusion equation and \cite{LindbladLuhrmannSchlagSoffer2023} where a resonant contribution that decays as $t^{-1/2}\log t$ along
certain rays is isolated for a one-dimensional quadratic Klein--Gordon equation. A diffusive phase mode also governs localized perturbations of homogeneous oscillations in reaction-diffusion systems, see \cite{GallayScheel2011}. We further mention the nonlinear Schr\"odinger equation as a potential extension of the technique, although its free evolution multiplies each Fourier amplitude by a phase factor
of modulus one \cite{Cazenave2003Semilinear}, which implies that a clock depending
only on the moduli of the amplitudes would be constant along each free orbit.
It remains to determine whether the phases can supply a local time
coordinate and whether dispersive estimates can control the resulting
homological integrals, as the analysis of space-time resonances suggests
\cite{GermainMasmoudiShatah2009,Germain2010SpaceTime}. Indeed, \cite[Thm.~1.1]{HaniPausaderTzvetkovVisciglia2015} obtains a resonant profile evolving in $\pi\log t$, an explicit function of physical time, for the cubic Schr\"odinger equation on $\mathbb R\times\mathbb T^d$.

\bibliographystyle{abbrv}
\bibliography{nonlinearkinetic}

\appendix

\section{Proof of Proposition \ref{prop_equilibria}}\label{app_equilibria_characterization}

If $F\equiv0$, the conclusion is immediate.  Assume henceforth that $F$ is not identically zero.  On the set where $\rho[F](x)>0$, multiply \eqref{eq:stationary-bgk} by $\log(F/(\rho_0\phi))$ and integrate in $(x,v)$.  For nonnegative $F$, the calculation is understood by the standard positive regularization of the entropy and the stated finiteness and regularity assumptions permit passage to the limit.  Since $\rho_0\phi$ is independent of $x$, the transport contribution is
\begin{align}
&\int_{\mathbb{R}^d_x\times\mathbb{R}^d_v}
    v\cdot\nabla_x F
    \log\left(\frac{F}{\rho_0\phi}\right)\,dv\,dx \notag =\int_{\mathbb{R}^d_x\times\mathbb{R}^d_v}
    v\cdot\nabla_x
    \left[F\log\left(\frac{F}{\rho_0\phi}\right)-F\right]\,dv\,dx
    =0,
\end{align}
by \eqref{eq:entropy-flux-assumption}. For $F\neq 0$, we can decompose
$$
    \log\left(\frac{F}{\rho_0\phi}\right)
    =
    \log\left(\frac{F}{\rho[F]\phi}\right)
    +
    \log\left(\frac{\rho[F]}{\rho_0}\right),
$$
and the second term makes no contribution because
$$
    \int_{\mathbb R^d_v}(\rho[F]\phi-F)\,dv=0,
$$
implying that
\begin{equation}
0=-\frac1\tau\int_{\mathbb R^d_x}\rho[F](x)^\alpha
   \int_{\mathbb R^d_v}
   (F-\rho[F]\phi)
   \log\left(\frac{F}{\rho[F]\phi}\right)dv\,dx.
\label{eq_stationary_entropy_production}
\end{equation}
The integrand is nonnegative by \eqref{a_b_inequality} and on $\{\rho[F]>0\}$ the prefactor $\rho[F]^\alpha$ is strictly positive.  Hence the inner entropy production vanishes for almost every such $x$, and the equality case in \eqref{a_b_inequality} yields
\begin{equation}
    F(x,v)=\rho[F](x)\phi(v)
    \quad\text{for a.e. }v.
\label{local_maxwellian_equilibrium}
\end{equation}
If $\rho[F](x)=0$, nonnegativity of $F$ and $\rho[F]=\int F\,dv$ imply $F(x,v)=0$ for almost every $v$, so \eqref{local_maxwellian_equilibrium} holds on the vacuum set as well. Substitution into \eqref{eq:stationary-bgk} annihilates the collision term and gives
$$
    \phi(v)\,v\cdot\nabla_x\rho[F](x)=0.
$$
Since $\phi(v)>0$ for every $v$, this identity holds for almost every $v$.  Taking, for example, its $v_i$ moment gives
$$
    \partial_{x_i}\rho[F](x)
       \int_{\mathbb R^d}v_i^2\phi(v)\,dv=0,
    \qquad i=1,\ldots,d,
$$
where we have used that the mixed moments vanish, and therefore $\nabla_x\rho[F]=0$.  Thus $\rho[F]$ is constant on the connected space $\mathbb R^d$, and $F=\rho_*\phi$ for some constant $\rho_*\geq 0$, as claimed.

\section{Proof of Convergence to a Global Maxwellian}\label{proof_global_existence}

\begin{proof}
We use the normalized equation \eqref{maineqf}, the collision operator
$\mathcal L_Q$ from \eqref{def_LQ}, and the Sobolev norms of
Section~\ref{Notation}. The projection $P_\rho$ in \eqref{def_P_rho}
is orthogonal in $L_v^2(\phi^{-1})$, and
the estimate 
$$ \|\rho[f]\|_{H_x^s}\leq\|f\|_{H_x^sL_v^2(\phi^{-1})},
$$
holds. 
For $\|\rho[f]\|_\infty\le1/2$, the Sobolev product and composition
estimates applied to \eqref{defNalpha} give
\begin{equation}
 \begin{aligned}
 \|\mathcal N_\alpha(f)\|_{H_x^sL_v^2(\phi^{-1})} \leq C\|f\|_{H_x^sL_v^2(\phi^{-1})}
          \|\mathcal L_Q f\|_{H_x^sL_v^2(\phi^{-1})}.
 \end{aligned}
\end{equation}
Since $\rho[\mathcal N_\alpha(f)]=0$ and $P_\rho$ commutes with spatial
derivatives, the nonlinear term pairs only with $f-P_\rho f$ in
the Sobolev energy identity. For sufficiently small
$\|f\|_{H_x^sL_v^2(\phi^{-1})}$, it is therefore absorbed by the
linear collision dissipation, giving
\begin{equation}\label{eq_basic_Hs_energy}
 \frac12\frac d{dt}\|f\|_{H_x^sL_v^2(\phi^{-1})}^2
 +\frac1{2\tau}\|\mathcal L_Q f\|_{H_x^sL_v^2(\phi^{-1})}^2\le0.
\end{equation}
For mild solutions, this estimate is justified by translating along
characteristics, which removes the skew-adjoint transport operator.\\
Free transport is an isometric strongly continuous group on
$H_x^sL_v^2(\phi^{-1})$, and the same product and composition estimates
make the collision term locally Lipschitz on $\|\rho[f]\|_\infty<1$.
Choose $\varepsilon_0$ so that absorption in \eqref{eq_basic_Hs_energy}
and the Sobolev bound $\|\rho[f]\|_\infty\le1/2$ hold on the ball of
radius $2\varepsilon_0$. The local mild solution stays in this ball
by \eqref{eq_basic_Hs_energy} and therefore continues uniquely for all
time, with
\begin{equation}\label{eq_global_Hs_bound_final}
 \sup_{t\geq 0}\|f(t)\|_{H_x^sL_v^2(\phi^{-1})}
 \leq\|f_0\|_{H_x^sL_v^2(\phi^{-1})},\qquad
 \sup_{t\geq 0}\|\rho[f(t)]\|_\infty\leq\frac12.
\end{equation}
In particular, the collision rate satisfies
$$
 \frac{2^{-\alpha}}\tau
 \leq\frac{(1+\rho[f])^\alpha}\tau
 \leq\frac{(3/2)^\alpha}\tau.
$$
Sobolev embedding and Cauchy--Schwarz place $f$ in
$L_v^1(\operatorname{BUC}_x)$, where $\operatorname{BUC}_x$ is the space of bounded uniformly
continuous functions on $\mathbb R_x^d$, with the supremum norm, and
$$
 \|f\|_{L_v^1(\operatorname{BUC}_x)}
 =\int_{\mathbb R^d}\sup_{x\in\mathbb R^d}|f(x,v)|\,dv.
$$
Thus the characteristic comparison used
in the proof of Theorem~\ref{thm_selected_maxwellian_bounds} applies
to this rate and preserves the nonnegativity of $F=\rho_0(\phi+f)$.\\
The rate bound above and $L^1_{x,v}$ contractivity of $P_\rho$ give convergence of the transport Duhamel iterates in $L^1_{x,v}$ on finite
intervals. The Sobolev product estimate gives convergence of the same iterates in $H_x^sL_v^2(\phi^{-1})$, so both limits represent $f$. Along characteristics the loss has a positive integrating factor and
the gain proportional to $P_\rho g$ preserves positivity. By linearity and comparison,
the solution with data $|f_0|$ therefore bounds $|f(t)|$ and its integral
equals $\|f_0\|_{L^1_{x,v}}$, since transport preserves the integral
and $\int (1+\rho)^\alpha(P_\rho g-g)\,dv=0$ at each $x$. This gives
\begin{equation}\label{eq_L1_contraction}
 \|f(t)\|_{L^1_{x,v}}\leq\|f_0\|_{L^1_{x,v}}
 \leq\|f_0\|_{L_v^2(\phi^{-1};L_x^1)}, 
\end{equation}
where the last inequality uses $\int\phi\,dv=1$ again. With the Fourier convention
\eqref{def_Fourier_transform}, it follows that
$$
 \sup_k|\rho[\hat f(t,k)]|
 \le(2\pi)^{-d/2}\|f_0\|_{L^1_{x,v}}.
$$
To obtain decay, we need dissipation of the density only at order zero.
Taking the first two velocity moments of \eqref{maineqf}, as in
\eqref{eq_moment_hierarchy}, gives
\begin{equation}\label{eq_macro_system_proof}
 \begin{aligned}
 \partial_t\rho[f]+\nabla_x\cdot\int vf\,dv&=0,\\
 \partial_t\int vf\,dv+\nabla_x\rho[f]
 &=\nabla_x\cdot\int v\otimes v\,\mathcal L_Q f\,dv
   -\frac{(1+\rho[f])^\alpha}\tau\int vf\,dv.
 \end{aligned}
\end{equation}
The Gaussian moments and $\int v\phi\,dv=0$ imply
$$
 \left\|\int vf\,dv\right\|_{L_x^2}
 +\left\|\int v\otimes v\,\mathcal L_Q f\,dv\right\|_{L_x^2}
 \leq C\|\mathcal L_Q f\|_{L_x^2L_v^2(\phi^{-1})}.
$$
For a sufficiently small fixed $\kappa>0$, use the modified energy
\begin{equation}
 E_0(t)=\frac12\|f(t)\|_{L_x^2L_v^2(\phi^{-1})}^2
 +\kappa\left\langle\int vf(t)\,dv,
          \nabla_x(1-\Delta_x)^{-1}\rho[f(t)]\right\rangle_{L_x^2}.
\end{equation}
The preceding moment bound shows that
$E_0\asymp\|f\|_{L_x^2L_v^2(\phi^{-1})}^2$.
Differentiating the inner product with \eqref{eq_macro_system_proof},
the density-gradient term is
$-\|\nabla_x(1-\Delta_x)^{-1/2}\rho[f]\|_{L_x^2}^2$.
The stress and collision terms are bounded by
$$
 C\|\mathcal L_Q f\|_{L_x^2L_v^2(\phi^{-1})}
       \|\nabla_x(1-\Delta_x)^{-1/2}\rho[f]\|_{L_x^2},
$$
while the term from $\partial_t\rho[f]$ is bounded by
$C\|\mathcal L_Q f\|_{L_x^2L_v^2(\phi^{-1})}^2$.
The first two terms use the uniform collision-rate bound and the fact
that the symbols of $\nabla_x(1-\Delta_x)^{-1}$ and
$\nabla_x\nabla_x(1-\Delta_x)^{-1}$ are bounded by
$|k|/(1+|k|^2)^{1/2}$. The last uses boundedness of
$\nabla_x(1-\Delta_x)^{-1}\nabla_x\cdot$ on $L_x^2$.
Young's inequality and the order-zero energy calculation in
\eqref{positive_hilbert_energy_identity} therefore give, for sufficiently
small $\kappa$,
\begin{equation}\label{eq_E0_dissipation}
 \frac d{dt}E_0
 +c\int_{\mathbb R^d}\left(
      \|\mathcal L_Q\hat f(k)\|_{L_v^2(\phi^{-1})}^2
      +\frac{|k|^2}{1+|k|^2}|\rho[\hat f(k)]|^2\right)dk\le0.
\end{equation}
This is the standard hypocoercive correction, see, e.g., \cite{villani2009hypocoercivity}.\\
Apply Fourier splitting \cite{BGKproofSchonbek} to the density
integral in \eqref{eq_E0_dissipation}. For fixed $A>0$ and $t\geq A$,
split at $|k|^2=A/(1+t)$.
On the complementary frequencies,
$|k|^2/(1+|k|^2)\geq  A/(2(1+t))$.
Orthogonality of $P_\rho$ and the Fourier bound following
\eqref{eq_L1_contraction} consequently yield
\begin{equation}
 \begin{aligned}
 E_0'+\frac{cA}{1+t}E_0
 &\leq\frac{CA}{1+t}
       \int_{|k|^2\leq A/(1+t)}|\rho[\hat f(t,k)]|^2\,dk\\
 &\leq C A^{1+d/2}\|f_0\|_{L^1_{x,v}}^2(1+t)^{-1-d/2},
 \end{aligned}
\end{equation}
with $c,C$ independent of $A$. Choose $A$ so that $cA>d/2$
and integrate with the factor $(1+t)^{cA}$.
Since $E_0$ is nonincreasing and
$E_0(0)\leq C\|f_0\|_{L_x^2L_v^2(\phi^{-1})}^2$,
this proves \eqref{decay_estimate} for $t\geq  A$.
For $0\leq t\leq A$, use \eqref{eq_global_Hs_bound_final} and
$(1+t)^{-d/4}\geq (1+A)^{-d/4}$.
The identity $F-\rho_0\phi=\rho_0f$ completes the proof.
\end{proof}

\end{document}